\pdfoutput=1
\documentclass[11pt]{article}
\usepackage[utf8]{inputenc}

\usepackage{amsthm}
\usepackage{amssymb,amsmath}
\usepackage{todonotes} 
\usepackage{hyperref} %
\usepackage{cleveref}
\usepackage{algorithm}
\usepackage{algpseudocode}
\usepackage{enumerate}
\usepackage{accents}
\usepackage{rotating}
\usepackage{pifont} %
\usepackage{multirow} %

\providecommand\leanmapurl{https://mirebeau.github.io/lean-proof-maps/GuaranteedSubsol/}
\providecommand\leancheckmark{\textup{\textsf{\scriptsize Lean}}\ensuremath{^\checkmark}}
\makeatletter
\DeclareRobustCommand\leancheck[1]{%
  \@ifundefined{leancheck@#1}%
    {\PackageWarning{leananchors}{No Lean anchor for the label `#1'}}%
    {\href{\leanmapurl\#label=#1}{\leancheckmark}}}
\DeclareRobustCommand\leancheckdef[1]{%
  \@ifundefined{leancheckdef@#1}%
    {\PackageWarning{leananchors}{No page of definitions `#1'}}%
    {\href{\leanmapurl\#def=#1}{\leancheckmark}}}
\@namedef{leancheck@agsik:until}{1}%
\@namedef{leancheck@algo:Dijkstra}{3}%
\@namedef{leancheck@algo:agsi_parallel}{7}%
\@namedef{leancheck@algo:agsik}{4}%
\@namedef{leancheck@algo:anb}{3}%
\@namedef{leancheck@algo:anb_lazy}{2}%
\@namedef{leancheck@algo:global_iteration}{6}%
\@namedef{leancheck@anb:for}{1}%
\@namedef{leancheck@anb:hatu}{1}%
\@namedef{leancheck@anb:while_sub}{1}%
\@namedef{leancheck@anb:while_sup}{1}%
\@namedef{leancheck@anb_lazy:olu_gets_hatu}{1}%
\@namedef{leancheck@assum:all}{3}%
\@namedef{leancheck@corol:Bn_comp}{2}%
\@namedef{leancheck@corol:G_min_u_t}{2}%
\@namedef{leancheck@corol:NGamma}{6}%
\@namedef{leancheck@corol:max_width}{5}%
\@namedef{leancheck@corol:max_width_withUnknown}{4}%
\@namedef{leancheck@corol:mon_min_max}{2}%
\@namedef{leancheck@corol:subadd_concave}{1}%
\@namedef{leancheck@corol:tol}{5}%
\@namedef{leancheck@enum:Lambda0_alpha}{1}%
\@namedef{leancheck@enum:Lambda0_delta}{2}%
\@namedef{leancheck@eq:Dijkstra_prop}{2}%
\@namedef{leancheck@eq:anb_cost}{1}%
\@namedef{leancheck@eq:causal_old}{2}%
\@namedef{leancheck@eq:comparison_Lambda}{3}%
\@namedef{leancheck@eq:cst_subsol}{2}%
\@namedef{leancheck@eq:lazy_bij_placeholder}{1}%
\@namedef{leancheck@eq:tol_as_Kruzkov}{1}%
\@namedef{leancheck@eqdef:Bpn}{1}%
\@namedef{leancheck@eqdef:Dijkstra_supsol}{1}%
\@namedef{leancheck@eqdef:monotone_extension}{1}%
\@namedef{leancheck@eqdef:stencil}{2}%
\@namedef{leancheck@eqdef:u_less_t}{1}%
\@namedef{leancheck@lem:causal_mix}{12}%
\@namedef{leancheck@lem:cplx_global_iteration}{4}%
\@namedef{leancheck@lem:cplx_global_iteration_withUnknown}{1}%
\@namedef{leancheck@lem:first_iter}{1}%
\@namedef{leancheck@lem:graph_subsol}{1}%
\@namedef{leancheck@lem:max_min}{3}%
\@namedef{leancheck@lem:max_min_mon_causal}{2}%
\@namedef{leancheck@lem:monDom}{3}%
\@namedef{leancheck@lem:pow_cvx}{2}%
\@namedef{leancheck@lem:powr_diff}{2}%
\@namedef{leancheck@lem:u_less_t}{2}%
\@namedef{leancheck@prop:Cpast_formulation}{2}%
\@namedef{leancheck@prop:Dijkstra_band}{4}%
\@namedef{leancheck@prop:Dijkstra_correct}{7}%
\@namedef{leancheck@prop:always_subsol}{4}%
\@namedef{leancheck@prop:causal_formulation}{4}%
\@namedef{leancheck@prop:causal_implies_subadd}{1}%
\@namedef{leancheck@prop:comp_of_causal}{3}%
\@namedef{leancheck@prop:comparison_Lambda}{2}%
\@namedef{leancheck@prop:comparison_Phi}{1}%
\@namedef{leancheck@prop:complexity_lazy}{2}%
\@namedef{leancheck@prop:constr_Gamma}{4}%
\@namedef{leancheck@prop:contractivity}{2}%
\@namedef{leancheck@prop:correct_anb_lazy}{3}%
\@namedef{leancheck@prop:cplx_agsi_parallel}{7}%
\@namedef{leancheck@prop:cst_subsol}{7}%
\@namedef{leancheck@prop:monotonic_Phi}{4}%
\@namedef{leancheck@prop:niter_while}{1}%
\@namedef{leancheck@prop:niter_while_withUnknown}{1}%
\@namedef{leancheck@prop:olu_ge_ustar}{1}%
\@namedef{leancheck@prop:partialMonotone}{2}%
\@namedef{leancheck@prop:pertSubSuper}{4}%
\@namedef{leancheck@prop:strict_causal_implies_subadd}{1}%
\@namedef{leancheck@prop:subadd_alpha_beta}{1}%
\@namedef{leancheck@prop:subadd_local}{4}%
\@namedef{leancheck@prop:t_advances}{1}%
\@namedef{leancheck@prop:t_advances_withUnknown}{1}%
\@namedef{leancheck@rem:anb_npos}{1}%
\@namedef{leancheck@th:complexity}{4}%
\@namedef{leancheck@th:complexity_withUnknown}{2}%
\@namedef{leancheck@th:complexity_with_agsi}{2}%
\@namedef{leancheckdef@AnbAlwaysSubsolWith}{}
\@namedef{leancheckdef@AnbAlwaysSubsolWith/AnbAlwaysSubsol}{}
\@namedef{leancheckdef@AnbAlwaysSubsolWith/AnbAlwaysSubsolWith}{}
\@namedef{leancheckdef@AnbInnerSolver}{}
\@namedef{leancheckdef@AnbInnerSolver/AnbInnerSolver}{}
\@namedef{leancheckdef@AnbInnerSolver/AnbInnerSolver.Valid}{}
\@namedef{leancheckdef@AnbInnerSolver/anbAgsi}{}
\@namedef{leancheckdef@AnbInnerSolver/anbGlobalIter}{}
\@namedef{leancheckdef@AnbParams}{}
\@namedef{leancheckdef@AnbParams.WithUnknown}{}
\@namedef{leancheckdef@AnbParams.WithUnknown/AnbParams.WithUnknown}{}
\@namedef{leancheckdef@AnbParams.WithUnknown/AnbParams.WithUnknown.C}{}
\@namedef{leancheckdef@AnbParams.WithUnknown/AnbParams.WithUnknown.nMax}{}
\@namedef{leancheckdef@AnbParams.WithUnknown/AnbParams.withUnknown}{}
\@namedef{leancheckdef@AnbParams/AnbParams}{}
\@namedef{leancheckdef@AnbParams/AnbParams.nMax}{}
\@namedef{leancheckdef@AnbParams/AnbParams.withUnknown}{}
\@namedef{leancheckdef@Axioms}{}
\@namedef{leancheckdef@Axioms/AdditivelyInvariant}{}
\@namedef{leancheckdef@Axioms/Causal}{}
\@namedef{leancheckdef@Axioms/Comparable}{}
\@namedef{leancheckdef@Axioms/Ext.AdditivelyInvariant}{}
\@namedef{leancheckdef@Axioms/Ext.Causal}{}
\@namedef{leancheckdef@Axioms/Ext.Comparable}{}
\@namedef{leancheckdef@Axioms/Ext.HasComparisonPrinciple}{}
\@namedef{leancheckdef@Axioms/Ext.HasDependency}{}
\@namedef{leancheckdef@Axioms/Ext.HasMemory}{}
\@namedef{leancheckdef@Axioms/Ext.HasStencil}{}
\@namedef{leancheckdef@Axioms/Ext.StrictlyCausal}{}
\@namedef{leancheckdef@Axioms/Ext.SubAdditive}{}
\@namedef{leancheckdef@Axioms/Ext.SubMultiplicative}{}
\@namedef{leancheckdef@Axioms/Global.Causal}{}
\@namedef{leancheckdef@Axioms/Global.Ext.Causal}{}
\@namedef{leancheckdef@Axioms/Global.Ext.HasMemory}{}
\@namedef{leancheckdef@Axioms/Global.Ext.HasStencil}{}
\@namedef{leancheckdef@Axioms/Global.Ext.StrictlyCausal}{}
\@namedef{leancheckdef@Axioms/Global.Ext.SubAdditive}{}
\@namedef{leancheckdef@Axioms/Global.Ext.revStencil}{}
\@namedef{leancheckdef@Axioms/Global.HasMemory}{}
\@namedef{leancheckdef@Axioms/Global.StrictlyCausal}{}
\@namedef{leancheckdef@Axioms/Global.SubAdditive}{}
\@namedef{leancheckdef@Axioms/HasComparisonPrinciple}{}
\@namedef{leancheckdef@Axioms/HasDependency}{}
\@namedef{leancheckdef@Axioms/HasMemory}{}
\@namedef{leancheckdef@Axioms/HasStencil}{}
\@namedef{leancheckdef@Axioms/StrictlyCausal}{}
\@namedef{leancheckdef@Axioms/SubAdditive}{}
\@namedef{leancheckdef@Axioms/SubMultiplicative}{}
\@namedef{leancheckdef@Definitions}{}
\@namedef{leancheckdef@Definitions/Global.kruzkovOp}{}
\@namedef{leancheckdef@Definitions/Global.monotoneExtension}{}
\@namedef{leancheckdef@Definitions/bandExtend}{}
\@namedef{leancheckdef@Definitions/bandOp}{}
\@namedef{leancheckdef@Definitions/constrGamma}{}
\@namedef{leancheckdef@Definitions/cutoffLe}{}
\@namedef{leancheckdef@Definitions/cutoffLt}{}
\@namedef{leancheckdef@Definitions/kruzkov}{}
\@namedef{leancheckdef@Definitions/kruzkovInv}{}
\@namedef{leancheckdef@Definitions/kruzkovOp}{}
\@namedef{leancheckdef@Definitions/monotoneExtension}{}
\@namedef{leancheckdef@Definitions/reflectOp}{}
\@namedef{leancheckdef@Definitions/rescaleOp}{}
\@namedef{leancheckdef@Definitions/restrictOp}{}
\@namedef{leancheckdef@Definitions/stencilDist}{}
\@namedef{leancheckdef@Definitions/stencilReach}{}
\@namedef{leancheckdef@OrderedPartialMap}{}
\@namedef{leancheckdef@OrderedPartialMap/OrderedPartialMap}{}
\@namedef{leancheckdef@OrderedPartialMap/bot}{}
\@namedef{leancheckdef@OrderedPartialMap/botVal}{}
\@namedef{leancheckdef@OrderedPartialMap/empty}{}
\@namedef{leancheckdef@OrderedPartialMap/entries}{}
\@namedef{leancheckdef@OrderedPartialMap/entriesIco}{}
\@namedef{leancheckdef@OrderedPartialMap/erase}{}
\@namedef{leancheckdef@OrderedPartialMap/erase1}{}
\@namedef{leancheckdef@OrderedPartialMap/pop}{}
\@namedef{leancheckdef@OrderedPartialMap/set}{}
\@namedef{leancheckdef@OrderedPartialMap/set1}{}
\makeatother

\def\diff{\mathrm{d}}
\usepackage{bbm} \def\bbone{\mathbbm{1}}
\def\ve{\varepsilon}
\def\<{\langle} \def\>{\rangle}
\def\interp{\mathrm I}
\DeclareMathOperator\Lip{Lip}
\def\sm{\setminus}

\def\bN{{\mathbb N}}
\def\bR{{\mathbb R}}
\def\bZ{{\mathbb Z}}

\def\cD{{\mathcal D}}
\def\cF{{\mathcal F}}
\def\cO{{\mathcal O}}
\def\cS{{\mathcal S}}
\def\cT{{\mathcal T}}
\def\cV{{\mathcal V}}
\def\bV{{\mathbb V}}

\newtheorem{theorem}{Theorem}[section]
\newtheorem{remark}[theorem]{Remark}
\newtheorem{lemma}[theorem]{Lemma}
\newtheorem{proposition}[theorem]{Proposition}
\newtheorem{corollary}[theorem]{Corollary}

\newtheorem{assumption}[theorem]{Assumption}

\usepackage{url}
\def\Kruzkov{Kruzk\"ov}
\def\ol{\overline}
\def\ul{\underline}
\def\Int{\mathring}
\def\past{\mathrm{past}}
\def\rE{\mathrm{E}}

\def\lin{\mathrm{lin}}
\def\loc{\mathrm{loc}}
\def\LF{\mathrm{LF}}
\def\SL{\mathrm{SL}}
\def\UD{\mathrm{UD}}
\def\bd{\mathrm{bd}}
\def\sign{\mathrm{sign}}

\def\init{\mathrm{init}}

\DeclareMathOperator\Cost{Cost}

\title{A narrow band method for non-causal discretizations of anisotropic eikonal equations, with adaptive width, sub- and super-solution construction, and quasi-linear complexity}

\author{
Jean-Marie Mirebeau\footnote{
Universit\'e Paris-Saclay, ENS Paris-Saclay, CNRS, Centre Borelli, F-91190 Gif-sur-Yvette, France.
},
Yanis Mouheb\footnotemark[1]\ \footnote{CEA, DAM, DIF, F-91297 Arpajon, France}
}

\begin{document}
\maketitle

\begin{abstract}
We describe a narrow band method for numerically solving anisotropic eikonal equations, and more generally static first order Hamilton-Jacobi equations with possibly non-convex Hamiltonians, discretized using standard schemes (semi-Lagrangian, upwind-differences, Lax-Friedrichs).
The causality property is \emph{not} needed, in contrast with the fast marching method.
The proposed approach simultaneously builds a numerical sub-solution $\ol u$ and super-solution $\ul u$ to the problem, and adapts the width of the narrow band based on their comparison. 
We prove that $(1+\ve) \ol u \geq \ul u$ at termination, where $\ve>0$ is any given tolerance parameter, and that the worst case complexity is $\cO(N |\ln \ve| \ln (N/\ve))$, where $N$ denotes the number of discretization points, and the underlying constant is explicit in terms of those associated with the scheme properties. 
Numerical experiments illustrate our results.
\end{abstract}

\section{Introduction}
\label{sec:intro}

The eikonal partial differential equation (PDE) is a simple model of front propagation, in which the front speed at any point of the domain is given in advance as a fixed positive value. 
\emph{Anisotropic} variants of the eikonal PDE take the form of static first order Hamilton-Jacobi-Bellman (HJB) PDEs, whose Hamiltonian is convex %
and dual to a metric $\cF$ encoding a direction dependent speed profile at each point. A further generalization, Hamilton-Jacobi-Isaacs (HJI) PDEs, allows a non-convex Hamiltonian $\cF^*$, often referred to as the slowness profile. 
Numerical solutions to the eikonal, HJB and HJI PDEs have countless applications in motion planning, image processing, traveltime tomography, dynamic games theory, etc, which motivates the design of efficient and general numerical methods.

Discretization schemes for such eikonal type PDEs can be split into two categories: \emph{Causal schemes} are solvable exactly, with $\cO(N\ln N)$ complexity using the fast-marching method and variants \cite{tsitsiklis1995EfficientTrajectories,sethian1996fast, kimmel1998Manifold,sethian2003OUM,vladimirsky2008LabelSetting,alton2012MAOUM,mirebeau2014marchingLattice,mirebeau2014EfficientFinsler,mirebeau2018curvature,mirebeau2019riemannian,dahiya2018ordered}. 
\emph{Non-causal schemes}, in contrast, are solvable up to a user prescribed tolerance $\ve > 0$, using a variety of iterative methods - fast sweeping (FS) \cite{tsai2003fastSweeping,zhao2005FastSweeping,kao2004lax,qian2007fastSweeping,luo2016convergenceFSM}, adaptive Gauss-Seidel iteration (AGSI) \cite{bornemann2006FiniteElement}, the fast iterative method (FIM) \cite{jeong2008fim} - or narrow band methods \cite{cristiani2009fast,cacace2011local,cacace2014can,mansour2024solving}.
Depending on the geometry of the speed profile and on the problem dimension $d$, designing a causal scheme can be complicated, or costly, or come with drawbacks such as the use of very large stencils \cite{mirebeau2016minimalStencils}, or simply be impossible, which motivates the study of iterative or narrow band methods applicable to non-causal schemes.
The numerical complexity $\Cost(N,\ve)$ of these methods is empirically between $\cO(N)$ and $\cO(N^{1+1/d})$ \cite{zhao2005FastSweeping,bornemann2006FiniteElement}, depending on the method and the test case, 
without complexity guarantee in terms of $N$ and $\ve$ for general non-causal schemes, except for the narrow band (NB) method \cite{mansour2024solving} by one of the authors, with complexity $\cO(N\ln^2(N/\ve))$.
Let us acknowledge that an important research effort has been devoted to developing eikonal solvers with linear $\cO(N)$ complexity, but theoretical guarantees remain partial, see \cref{rem:linear_complexity}.

We develop in this paper an \emph{adaptive} narrow band (ANB) method, which applies three standard non-causal schemes (semi-Lagrangian without the acuteness condition, Lax-Friedrichs relaxation, upwind-differences with drift), and has a quasi-linear complexity:
\begin{equation}
	\label{eq:N_lneps_lnNeps}
	\cO(N|\ln \ve|\ln (N/\ve)),
\end{equation}
see \cref{th:complexity,th:complexity_withUnknown} and \eqref{eq:ln_N_eps}.
The complexity is similar to the NB method \cite{mansour2024solving}, but we address a key shortcoming of this earlier work.
One limitation of the NB method is indeed that the narrow band width (among other parameters) must be \emph{fixed} a-priori. The theoretically established parameter values ensuring the method's convergence and complexity guarantees are unfortunately quite pessimistic and sub-optimal, leading to poor numerical efficiency in practice, whereas the use of ad-hoc empirical parameters improves practical efficiency but voids the mathematical guarantees.
While the present paper is clearly a spiritual continuation of \cite{mansour2024solving}, \emph{none of the results presented here depends on that earlier work}, which is still under review at the time of writing.
Some of its material is nevertheless recalled in \cref{sec:effective}.

\textbf{Lean formalization and AI usage.}
The original version of this manuscript and the C++ codes were written by hand by the authors, over the course of 2026. On October of that year, after being stupefied like many in the mathematical community by news of the efficiency of AI systems in mathematics, the first author decided to use them to validate the work before submission. 
All the mathematical results presented in this paper (lemmas, propositions, theorems) have thus been formalized in the Lean language, using Claude$^{\small\textregistered}$ AI together with Scott Armstrong's Lean skills\footnote{\url{https://github.com/scottnarmstrong/LeanAutoformalizationSkills}}. The formalization process was surprisingly smooth and fast (around three days, reviewing the formal statements as they were produced), but beyond correctness we make no claims about the quality or reusability of the Lean code. Some minor typos and mistakes were corrected in the manuscript in the process, and the most significant one is acknowledged in the proof of \cref{prop:complexity_lazy}. 
The C++ codes were audited, and some additional references were suggested. Beyond this, no AI was used in writing the manuscript and C++ codes.\\

In order to present the adaptive narrow band (ANB) method, whose narrow band width is adjusted at runtime based on a-posteriori error estimates, rather than fixed a-priori as in \cite{mansour2024solving}, we need to introduce some notation. Similarly to the previously mentioned methods \cite{zhao2005FastSweeping,bornemann2006FiniteElement,jeong2008fim,cristiani2009fast,mansour2024solving}, the PDE domain is approximated with a finite set $X$.
We denote by $\bR^X$ the collection of all mappings $u : X \to \bR$, equipped with the following partial order relation: 
for all $u,v \in \bR^X$
\begin{equation*}
	u \leq v \quad \Leftrightarrow \quad (\,\forall x \in X,\ u(x) \leq v(x)\,).
\end{equation*}
Our objective is to compute or approximate the fixed point $u_* \in \bR^X$ of a given operator $\Lambda : \bR^X \to \bR^X$,
which is a discretization scheme for a PDE of interest:
\begin{equation}
\label{eq:fixedPoint}
	\Lambda u_* = u_*
\end{equation}
The proposed ANB method simultaneously constructs a \emph{sub-solution} $\ol u$ %
and a \emph{super-solution} $\ul u\in \bR^X$, whose relative difference is at most a given $\ve>0$ at termination: 
\begin{align}
\label{eq:cvtol_intro}
	\ol u &\leq \Lambda \ol u, &
	\Lambda \ul u &\leq \ul u, &
	(1+\ve) \ol u \geq \ul u.
\end{align}
Here, $\ol u, \ul u, \Lambda \ol u, \Lambda \ul u \in \bR^X$ are mappings from the finite set $X$ to the set of real numbers $\bR$, and all inequalities hold pointwise on $X$. %
Note that classical methods \cite{zhao2005FastSweeping,jeong2008fim,bornemann2006FiniteElement,mansour2024solving} only generate a super-solution $\ul u$, see \cref{rem:subsol_construction} for the specific difficulties associated with the construction of a sub-solution.
The discrete comparison principle, see \cref{prop:comparison_Lambda}, yields
\begin{align}
\label{eq:tol_fixed_point}
	\ol u &\leq u_* \leq \ul u, &
	\text{hence }
	(1+\ve)^{-1} \ul u \leq u_* \leq (1+\ve) \ol u,
\end{align}
and therefore $\max \{\|u_*-\ol u\|_\infty, \|\ul u-u_*\|_\infty\} \leq \ve \|u_*\|_\infty$.
See \cref{corol:tol} for further discussion of convergence tolerance.

A consequence of our convergence criterion \eqref{eq:cvtol_intro} is the non-negativity of the sub-solution: indeed $(1+\ve) \ol u \geq \ul u \geq \ol u$ by the convergence tolerance criterion and comparison principle, hence $\ve \ol u \geq 0$, and therefore $0 \leq \ol u \leq u_*$; this is why we assume positive boundary conditions below, which is fortunately not restrictive thanks to the additive invariance of the problem. 
Let us add that the $L^\infty$ error between the fixed point $u_*$ of the discretization scheme \eqref{eq:fixedPoint}, and the solution $U$ of the eikonal, HJB or HJI PDE of interest, can usually also be estimated: $\|u_*-U\|_\infty = \cO(\sqrt h)$ where $h>0$ denotes the discretization grid scale, under standard assumptions using the doubling of variables technique.
See \cite[Chapter 10]{evans2010pde} for a general overview, \cite{shum2016rateOUM,mirebeau2019riemannian} and references therein for the application to various discretizations of eikonal and HJB PDEs, and \cite[Appendix C]{mansour2024solving} for the specific schemes considered in this paper.

We describe and study the proposed ANB method in an abstract setting, by assuming that a hypothetical operator $\Lambda$ satisfies a series of axioms.
In practice, these axioms are satisfied by a variety of discretizations of the eikonal PDE, including the Semi-Lagrangian (SL), Upwind-Differences (UD), and Lax-Friedrichs (LF) schemes described in \cref{sec:effective}.
For concreteness, we also present in \cref{rem:FE} below the finite element discretization \cite{bornemann2006FiniteElement} of HJB PDEs, which is a special case of a Semi-Lagrangian discretization.

\begin{remark}[Finite element discretization of eikonal and HJB PDEs]
\label{rem:FE}
Consider an open boun\-ded domain $\Omega \subset \bR^d$, and a slowness profile $\cF^*\in C^0(\overline \Omega \times \bR^d , [0,\infty[)$, $(x,\zeta) \mapsto \cF^*_x(\zeta)$, often referred to as the Hamiltonian, subject to assumptions specified later. 
Our objective is to compute the viscosity solution $U : \overline \Omega \to \bR$ \cite{bardi2008viscosity} to
\begin{align}
\label{eq:eikonal}
	\forall x \in \Omega,\, \cF_x^*(\nabla U(x)) &= 1,  &
	\forall x \in \partial \Omega,\, U(x) &= U_\bd(x).
\end{align}
An upper bound on the Lipschitz constant of the boundary data $U_\bd : \partial \Omega \to \bR$ is required to ensure the continuity of the solution, see \cite{bardi2008viscosity,bornemann2006FiniteElement} for a discussion, which is out of the scope of this paper.
An important special case\footnote{This condition implies the convexity of $\zeta \mapsto \cF^*_x(\zeta)$, thus excludes two-player differential games, 
see \cref{subsec:pursuit}.} is when $\cF_x^*(\zeta) = \max \{\<\zeta,\xi\> \mid \cF_x(\xi) \leq 1\}$ is dual to a \emph{Finsler metric}%
\footnote{I.e.\ $\cF_x$ is a norm on $\bR^d$, possibly asymmetric, for all $x \in \overline\Omega$.}
$\cF \in C^0(\overline \Omega \times \bR^d, [0,\infty[)$, $(x,\xi)\mapsto \cF_x(\xi)$,
on the domain. Bellman's optimality principle then states that, for any point $x$ and neighborhood $V$ such that $x\in V \subset \overline \Omega$, one has
\begin{align}
\label{eq:BellmanOpt}
	U(x) &= \min_{y \in \partial V} \Big(d_\cF(y, x) + U(y)\Big),\\ 
	\nonumber
	\text{where} \quad
	d_\cF(x,y) &:= \min \Big\{\int_0^1 \cF_{\gamma(t)}(\gamma'(t)) \diff t \mid
	\gamma \in \Lip([0,1], \overline \Omega),\ \gamma(0) = x,\ \gamma(1) = y\Big\}.
\end{align}
Taking $V = \overline \Omega$ in \eqref{eq:BellmanOpt}, one sees in particular that the viscosity solution to \eqref{eq:eikonal} is the path-length distance from $\partial \Omega$, w.r.t.\ the metric $\cF$, and with initial penalty $U_\bd$.

We now turn to the so-called finite element discretization \cite{bornemann2006FiniteElement} of \eqref{eq:eikonal}, and for that purpose we consider a triangulation $\cT$ of $\overline \Omega$,
whose interior (resp.\ boundary) vertices are denoted $\Int X$ (resp.\ $\partial X$). Let $X := \Int X \sqcup \partial X$, and denote by $V(x)$ the union of the simplices of $\cT$ containing a given point $x\in \Int X$. Define for all $u \in \bR^X$ and all $x \in \Int X$ %
\begin{align}
\label{eqdef:FE}
	\Lambda u(x) &:= \min_{y \in \partial V(x)} \Big(\cF_x(x-y) + \interp_\cT u(y)\Big), &
	\Gamma u(x) &:= c(x) + \min_{y \in X \cap \partial V(x)} u(y),
\end{align}
where $\interp_\cT u(y)$ denotes the piecewise linear interpolation on the triangulation $\cT$ of the mapping $u : X \to \bR$, evaluated at the point $y\in \overline \Omega$, and where $c(x) := \min\{\cF_x(x-y) \mid y \in \partial V(x)\}$.
We also let $\Lambda u (x) := \Gamma u(x) := u_\bd(x)$ for all $x \in \partial X$, where $u_\bd : \partial X \to \bR$ approximates the boundary data $U_\bd$. The fixed point equation \eqref{eq:fixedPoint}, characterizing the discrete scheme solution $u_*$, thus mimics Bellman's optimality principle \eqref{eq:BellmanOpt}, using the polygonal neighborhood $V(x)$ and the approximations $d_\cF(y,x) \approx \cF_x(x-y)$ and $U(y) \approx \interp_\cT u(y)$.
The operator (\ref{eqdef:FE}, right) is intended as a cheap approximation of $\Lambda$, solvable via Dijkstra's algorithm, which is used in a subroutine of the proposed ANB method.
\end{remark}

We fix throughout this paper a finite set $X$, which is the disjoint union $X = \Int X \sqcup \partial X$ of two subsets $\Int X$ and $\partial X$. 
As announced, we consider an abstract operator $\Lambda : \bR^X \to \bR^X$ for which we \emph{always implicitly assume Dirichlet boundary conditions}: $\Lambda u := u_\bd$ on $\partial X$, where the boundary data $u_\bd \in \bR^{\partial X}$ is fixed and independent of $u\in \bR^X$.
We often assume the following axioms, %
where $C_\past\geq 0$ is a constant
\begin{itemize}
	\item (Monotonicity) $u \leq v \Rightarrow \Lambda u \leq \Lambda v$, for all $u,v\in \bR^X$.
	\item (Additive invariance \leancheckdef{Axioms/AdditivelyInvariant}) $\Lambda (u + t) = \Lambda u + t$ on $\Int X$, for all $t \in \bR$, $u \in \bR^X$.
	\item (Positive boundary conditions) $\min(u_\bd) > 0$. %
	\item (Super-solution) There exists $u \in \bR^X$ such that $\Lambda u \leq u$.
	\item ($C_\past$-memory \leancheckdef{Axioms/HasMemory}\footnote{We denote by $\max\{u,t\}\in \bR^X$ the elementwise maximum $x \in X \mapsto \max \{u(x), t\}$, and likewise $\min \{u , t\}$.}) For all $u,v\in  \bR^X$, and all $t \in \bR$, 
	\begin{equation}
	\label{eqdef:Cpast_memory}
	\max \{u,t\} = \max\{v,t\}
	\ \Rightarrow \
	\max \{\Lambda u, t+C_\past\} = \max \{\Lambda v, t+C_\past\}.
	\end{equation}
\end{itemize}
Some brief preliminary comments: Monotonicity is a classical assumption in this context, which does however exclude discretizations of accuracy order $>1$ \cite{bornemann2006FiniteElement}. Additive invariance is expected since it holds at the continuous level \eqref{eq:eikonal}.
The positivity of the boundary condition is needed for our convergence tolerance criterion \eqref{eq:cvtol_intro} to make sense; as already mentioned, it is not restrictive, in view of the additive invariance of the problem.
The existence of a super-solution needs to be assumed, otherwise the problem may have no finite solution. 
The $C_\past$-memory axiom means that $\Lambda u(x)$ does not depend on excessively old neighbor values $u(y) \leq \Lambda u(x)-C_\past$, see \cref{prop:Cpast_formulation}, which is convenient in the design of narrow band methods \cite{mansour2024solving}.

Denote by $\bR_\infty := ]-\infty,\infty]$, the real line extended \emph{to the right} only.
Any monotone operator $\Lambda : \bR^X \to \bR^X$ admits an extension
$\Lambda_\infty : \bR_\infty^X \to \bR_\infty^X$ defined as follows: for all $u \in \bR^X_\infty$
\begin{equation}
\label{eqdef:monotone_extension}
	\Lambda_\infty u := \lim_{t \to \infty} \Lambda \min\{u,t\}.
\end{equation}
Note that the r.h.s.\ of \eqref{eqdef:monotone_extension} is pointwise non-decreasing w.r.t.\ the variable $t \in \bR$, by monotonicity, hence the limit exists.
Moreover $\min\{u,t\} = u$ for all $t \geq \max (u)$, hence $\Lambda = \Lambda_\infty$ on $\bR^X$. 
For this reason we often conflate $\Lambda$ with $\Lambda_\infty$ in the paper, dropping the $\infty$ subscript. The Lean formalization is necessarily stricter in this regard, and in particular it features counterparts of the above axioms for operators $\bR_\infty^X \to \bR_\infty^X$, as well as straightforward transfer theorems between them by monotone extension \eqref{eqdef:monotone_extension} and by restriction to $\bR^X$.
The \Kruzkov\ coordinate transform $\rE_\alpha : \bR_\infty \to \bR$ and its left inverse $\rE_\alpha^{-1} : \bR \to \bR_\infty$ are defined\footnote{%
Some authors define instead the \Kruzkov\ transform as $1-\exp(-u/\alpha)$, in such way that it takes non-negative values for all $u \geq 0$, but this was not convenient for our purposes.} as
\begin{align}
\label{eqdef:Kruzkov}
	\rE_\alpha u &:= \begin{cases}
		-\exp(-u/\alpha) &\text{if } u<\infty\\
		0 &\text{else}
	\end{cases} &
	\rE_\alpha^{-1} v &:= \begin{cases}
		-\alpha\ln(-v) &\text{if } v<0\\
		\infty &\text{else},
	\end{cases}
\end{align}
where $\alpha>0$ is a constant. Note that $\rE_\alpha^{-1} \rE_\alpha u = u$ and $\rE_\alpha \rE_\alpha^{-1} v = \min\{v,0\}$. The mapping $\rE_\alpha$ (resp.\ $\rE_\alpha^{-1}$) applies elementwise to $u \in \bR_\infty^X$ (resp.\  $v\in \bR^X$).
We define %
\begin{align}
\label{eqdef:Phi}
	\Phi_\alpha &: \bR^X \to \bR^X, &
	\Phi_\alpha := \rE_\alpha \Lambda_\infty \rE_\alpha^{-1}.
\end{align}
We often assume the following axiom, where $0 \leq \rho < 1$ is a constant:
\begin{itemize}
	\item ($\rho$-sub-additivity \leancheckdef{Axioms/SubAdditive}) $\Phi_\alpha(v+t) \leq \Phi_\alpha(v) + \rho t$, for all $t\geq 0$, $v \in \bR^X$.
\end{itemize}
The above condition is essentially equivalent to the contractivity of the operator $\Phi_\alpha$ in the $\ell^\infty$ norm, with contraction factor $\rho$, 
see \cref{prop:contractivity}. 
Similarly to \cite{mansour2024solving}, it is a key ingredient of the convergence analysis.

Finally, we consider an additional operator $\Gamma$, which like (\ref{eqdef:FE}, right) is typically associated to a positively weighted graph. It is regarded as an approximation of $\Lambda$, whose evaluation is cheap and whose solutions can be computed using Dijkstra's method since it satisfies the \emph{causality property}, see \cref{prop:Dijkstra_correct}.
The \emph{same boundary conditions as $\Lambda$} are always implicitly assumed: $\Gamma u := u_\bd$ on $\partial X$, for all $u \in \bR^X$.
In addition, we often assume the following axioms, %
where $C_\Gamma\geq 1$ and $C'_\Gamma\geq 0$ are constants:
\begin{itemize}
	\item (Monotonicity) $u \leq v \Rightarrow \Gamma u \leq \Gamma v$, for all $u,v \in \bR^X$.
	\item (Additive invariance \leancheckdef{Axioms/AdditivelyInvariant}) $\Gamma (u + t) = \Gamma u + t$ on $\Int X$, for all $t \in \bR$, $u \in \bR^X$.
	\item (Causality \leancheckdef{Axioms/Causal}) For all $u,v \in \bR^X$, $t \in \bR$,
	\begin{equation}
	\label{eqdef:causality}
		\min \{u,t\} = \min\{v,t\} \ \Rightarrow \ \exists \delta>0,\, \min\{\Gamma u,t+\delta\} = \min \{\Gamma v,t+\delta\}.
	\end{equation}
	\item (Under-estimation\footnote{I.e.\ $\forall u \in \bR^X,\, \forall x \in X,\, \Gamma u(x) \leq \Lambda u(x)$.}) $\Gamma \leq \Lambda$.
	\item ($C_\Gamma,C'_\Gamma$-Comparison \leancheckdef{Axioms/Comparable}%
) 
$C_\Gamma \Lambda u - C'_\Gamma \leq \Lambda [C_\Gamma u] \leq C_\Gamma \Gamma u$ on $\Int X$, for all $u \in \bR^X$.
\end{itemize}
A similar graph-based approximation of $\Lambda$ is considered in \cite{mansour2024solving}, but in that case it is assumed to be \emph{over-estimating} rather than \emph{under-estimating}.

Our main result, \cref{th:complexity} on the complexity of the proposed Adaptive Narrow Band (ANB) method, regards the above properties as assumptions. Of course, one also needs to establish that these properties hold for the concrete numerical schemes of interest, and how the related constants depend on the grid scale and other geometrical properties, see \cref{sec:effective} and \cite{mansour2024solving}.
See also \cref{rem:homog_csts,rem:validityFE,rem:causality} for an early discussion, focused on the case of the finite element scheme \eqref{eqdef:FE}. 
The formulation given above of the causality (resp.\ $C_\past$-memory) property differs from \cite{mirebeau2019riemannian,mirebeau2018curvature} (resp.\ \cite{mansour2024solving}), but it is actually equivalent as established in \cref{prop:causal_formulation} (resp.\ \cref{prop:Cpast_formulation}).
This reformulation simplifies and clarifies the statements and proofs, and is also more convenient for the formalization in the Lean language.

\begin{remark}[Other properties of interest]
\label{rem:submult_strictcaus}
The following properties, of an operator $\Lambda : \bR^X \to \bR^X$, are not used in our main results, but are often discussed: with $\delta,\nu>0$
\begin{itemize}
	\item ($\nu$-sub-multiplicativity \leancheckdef{Axioms/SubMultiplicative}\footnote{%
	The sub-multiplicativity parameter is denoted $\sigma$ in \cite{mansour2024solving}, but here we reserve this symbol for another use.
	})
	$\Lambda[(1+t) u] \leq (1+t) \Lambda u - \nu t$ on $\Int X$, for all $t \geq 0$, $u \in \bR^X$.
	\item ($\delta$-strict-causality \leancheckdef{Axioms/StrictlyCausal}) Like causality, but $\delta>0$ is independent of $u,v\in \bR^X$ and $t\in \bR$.
\end{itemize}
\end{remark}

A specificity of the proposed ANB method is the construction of a \emph{sub-solution} to the discretized problem, see \cref{rem:subsol_construction}. 
The next result shows that the null map is a sub-solution to $\Lambda$ and $\Gamma$, under various combinations of the above axioms, hence can be used to initialize the ANB.

\begin{proposition}[Constant sub-solutions \leancheck{prop:cst_subsol}]
\label{prop:cst_subsol}
Consider two additively invariant operators $\Lambda, \Gamma : \bR^X \to \bR^X$, and parameters $\nu,\alpha,\rho,C_\Gamma>0$, with $\rho<1$. 
\begin{enumerate}
	\item If $\Lambda$ is $\nu$-sub-multiplicative, then $\Lambda[0] \geq \nu$ on $\Int X$.
	\label{enum:Lambda0_nu}
	\item If $\Phi_\alpha := \rE_\alpha \Lambda_\infty \rE_\alpha^{-1}$ is $\rho$-sub-additive, then $\Lambda[0] \geq -\alpha \ln \rho$ on $\Int X$.
	\label{enum:Lambda0_alpha}
	\item If $\Gamma$ is causal, then $\Gamma[0] > 0$ on $\Int X$. (If $\Gamma$ is $\delta$-strictly causal, then $\Gamma[0] \geq \delta$ on $\Int X$.) 
		\label{enum:Lambda0_delta}
\end{enumerate}
If $\Lambda[0] \geq 0$ and the upper estimate $\Lambda[C_\Gamma u] \leq C_\Gamma \Gamma u$ of $C_\Gamma,C'_\Gamma$-comparison holds (resp.\ $\Gamma[0] \geq 0$ and $\Gamma\leq \Lambda$), then
\begin{align}
\label{eq:cst_subsol}
	t &\leq \Lambda[t],&
	t &\leq \Gamma[t], 
	&\text{for all } t \leq \min(u_\bd).
\end{align}
\end{proposition}

\begin{proof}
	Case \ref{enum:Lambda0_nu}: apply the sub-multiplicativity axiom to $u=0$ and $t=1$. 
	Case \ref{enum:Lambda0_alpha}: we have $\Phi_\alpha(t) = t \exp(-\Lambda[0]/\alpha)$ for all $t<0$, using the additive invariance of $\Lambda$, hence $\exp(-\Lambda[0]/\alpha) \leq \rho$, as announced.
	Case \ref{enum:Lambda0_delta}: causality and additive invariance, applied to $u=0$, $v=1$, and $t=0$, yield $\min \{\Gamma[0],\delta\} = \min\{1+\Gamma[0],\delta\}$, for some $\delta>0$ depending on $u,v,t$ (resp.\ the fixed $\delta$ in the strictly-causal case), hence $\Gamma[0] \geq \delta$ as announced.
	
	Additional claim \eqref{eq:cst_subsol}: note that $0 \leq \Lambda[0] \leq C_\Gamma \Gamma[0]$ (resp.\ $0 \leq \Gamma[0] \leq \Lambda[0]$) on $\Int X$, hence both $\Lambda[0], \Gamma[0] \geq 0$ on $\Int X$, and therefore $\Lambda[t], \Gamma[t] \geq t$ on $\Int X$ by additive invariance, for all $t\in \bR$. We conclude recalling that $\Lambda[t]=\Gamma[t]=u_\bd$ on $\partial X$.
\end{proof}
Conversely, $C_\past$-memory and additive invariance imply the upper bound $\Lambda[0] \leq C_\past$ on $\Int X$ (indeed, choosing $u=0$, $v=-1$ and $t=0$ yields $\max\{\Lambda[0],C_\past\} = \max \{\Lambda[0]-1,C_\past\}$).

The underestimating operator (\ref{eqdef:FE}, right) associated to the finite element scheme (\ref{eqdef:FE}, left) is a special case of a generic construction, see the next result, which also applies to the Semi-Lagrangian, Lax-Friedrichs, and Upwind-Differences schemes considered in \cref{sec:effective}.
On the positive side, this construction is straightforward, satisfies all the required axioms, and is used successfully in our numerical experiments. On the negative side, an ad-hoc construction depending on the numerical scheme of interest may yield sharper constants $C_\Gamma$ and $C_\Gamma'$, hence possibly improved complexity both in theory and practice. In addition, the following construction relies on an assumption \eqref{eqdef:Cpast_dependency} which is a strengthened version of $C_\past$-memory \eqref{eqdef:Cpast_memory}, and which excludes for instance the operator $\Lambda u(x) := C_\past + \max_{\cV(x)} u$, where $\emptyset \neq \cV(x)\subset X$ for all $x \in \Int X$ is an arbitrary stencil, as well as similar schemes defined as \emph{maxima} of simpler operators and encountered in e.g.\ some discretization of Isaacs equations associated to differential games.

\begin{proposition}[Constructing $\Gamma$ from $\Lambda$ \leancheck{prop:constr_Gamma}]
\label{prop:constr_Gamma} 
	Let $\Lambda : \bR^X \to \bR^X$ be monotone, additively invariant, and such that $\Lambda[0] \geq \delta > 0$ on $\Int X$.
	Assume that $\Lambda u(x)$ only depends on the values of $u$ on a stencil $\emptyset\neq\cV(x) \subset X$, for all $x \in \Int X$, and is bounded in terms of the smallest of these values:
	\begin{align}
	\label{eqdef:Cpast_dependency}
		u_{|\cV(x)} &= v_{|\cV(x)} \ \Rightarrow \ \Lambda u(x) = \Lambda v(x),&
		\Lambda u(x) &\leq C_\past + \min_{\cV(x)} u,
	\end{align}
	for all $u,v\in \bR^X$ and all $x\in \Int X$, where $C_\past\geq 0$ is a constant.
	Then $\Lambda$ has $C_\past$-memory, and admits a super-solution if the stencils are connected to the boundary\footnote{In the sense that: $\forall x_0 \in \Int X$, $\exists n \geq 0$, $\exists x_1,\cdots,x_n \in \Int X$, $\exists x_{n+1} \in \partial X$, $\forall 0 \leq i \leq n$, $x_{i+1}\in \cV(x_i)$}. In addition, the following operator $\Gamma$ is monotone, additively invariant, $\delta$-strictly-causal, underestimates $\Lambda$, and is $C_\Gamma, C'_\Gamma$-comparable. Given $u \in \bR^X$, set $\Gamma u = u_\bd$ on $\partial X$, and for all $x \in \Int X$
	\begin{align}
	\label{eqdef:constr_Gamma}
		\Gamma u(x) &:= c(x) + \min_{\cV(x)} u, &
		c(x) &:= \Lambda[0](x), &
		C_\Gamma &:= \frac {C_\past} \delta, &
		C'_\Gamma &:= C_\past C_\Gamma-\delta.
	\end{align}
\end{proposition}

The proof is simple but postponed to \cref{sec:effective}, to streamline the introduction.

\begin{remark}[Homogeneity properties of the constants]
\label{rem:homog_csts}
The following intrinsic grid scale $h>0$ can be associated to the finite element scheme \eqref{eqdef:FE}:
\begin{equation}
\label{eqdef:hFE}
h = \max \{ \cF_x(x-y)	\mid x \in \Int X,\ y \in X \cap \partial V(x)\}.
\end{equation}
A grid scale $h>0$ is also associated to the Semi-Lagrangian, Upwind-Differences, and Lax-Friedrichs schemes, see \cref{sec:effective}.
For these schemes, and under mild geometric non-degeneracy assumptions, $\alpha/h$ is positively bounded above and below, $1-\rho$ is positively bounded below, and $(C_\past/h, C_\Gamma, C_\Gamma'/h)$ are bounded above, \emph{independently} of $h>0$, see \cref{sec:effective}. The final complexity estimate \eqref{eq:N_lneps_lnNeps} assumes these homogeneity properties, see \eqref{eq:homog_scheme_csts} and \eqref{eq:ln_N_eps} below.
\end{remark}

\begin{remark}[Validity of the scheme axioms]
\label{rem:validityFE}
Consider the finite element scheme \eqref{eqdef:FE}, from \cite{bornemann2006FiniteElement}.
Monotonicity and additive invariance hold by construction. 
Property \eqref{eqdef:Cpast_dependency} is satisfied with $C_\past := h$ from \eqref{eqdef:hFE}, hence \cref{prop:constr_Gamma} yields the remaining desired properties of $\Lambda$ and $\Gamma$, except the positivity of boundary conditions, which must be assumed, and the $\rho$-sub-additivity of $\Phi_\alpha$, which is further discussed in \cref{rem:causality} and is established in \cref{sec:effective}.
\end{remark}

\begin{remark}[Causality, and \Kruzkov-contractivity]
\label{rem:causality}
A monotone and causal \eqref{eqdef:causality} operator admits a unique fixed point, which can be computed with complexity $\cO(N\ln N)$ using the so-called dynamic programming principle, also referred to as Dijkstra's algorithm or the fast marching method depending on the context, see \cref{prop:Dijkstra_correct,prop:comp_of_causal}.
Causality holds for the operator associated to a graph distance with positive edge lengths, such as $\Gamma$ as defined in (\ref{eqdef:FE}, right) with $\delta:=\min\{c(x)\mid x \in \Int X\}$. It also holds for the operator $\Lambda$ defined in (\ref{eqdef:FE}, left), provided a  geometrical \emph{acuteness condition} is satisfied \cite{kimmel1998Manifold,sethian2003OUM,vladimirsky2008LabelSetting}, and this observation is at the foundation of numerous variants \cite{alton2012MAOUM,mirebeau2014marchingLattice,mirebeau2014EfficientFinsler,dahiya2018ordered} of the \emph{semi-Lagrangian} fast marching algorithm \cite{tsitsiklis1995EfficientTrajectories}. Finally, causality holds for numerical schemes defined using suitable \emph{upwind-differences}, see 
the foundational work
\cite{sethian1996fast} and its anisotropic variants \cite{mirebeau2018curvature,mirebeau2019riemannian} also discussed in \cref{sec:effective}.

On the other hand, causality does \emph{not} hold for the scheme \eqref{eqdef:FE} if the mesh has obtuse angles (measured w.r.t.\ the metric), or for the upwind-differences scheme augmented with a drift term, or the Lax-Friedrichs scheme. As a substitute, we use the $\rho$-sub-additivity of the operator $\Phi_\alpha$ expressed in \Kruzkov\ coordinates \eqref{eqdef:Phi}, which is established in \cite{mansour2024solving} for the schemes of interest in this paper, see \cref{sec:effective}. Here, crucially, $\alpha$ is proportional to (resp.\ $\rho$ is independent of) the grid scale $h$, see \cref{rem:KruzkovFixed}.
This property can be reformulated in terms of the following \emph{local contraction ratio} if the scheme $\Lambda$ is everywhere differentiable: for all $x \in \Int X$ and all $u \in \bR^X$, 
\begin{align}
\label{eq:subadd_diff}
\rho^\loc(\alpha; u,x) &\leq \rho < 1& \text{where } 
\rho^\loc(\alpha; u,x) &:= 
	\sum_{y \in X} \frac{\partial \Lambda u(x)} {\partial u(y)} \exp\big(\frac {u(y) - \Lambda u(x)} \alpha \big).
\end{align}
See \cref{prop:subadd_local} for a similar statement without the differentiability assumption. Noting that $0 \leq \frac{\partial \Lambda u(x)} {\partial u(y)}$ by monotony we obtain $\frac{\partial \Lambda u(x)} {\partial u(y)} \leq \rho \exp(-\frac {u(y) - \Lambda u(x)} \alpha)$. In other words, $\Lambda u(x)$ may depend on later, non-causal values $u(y) > \Lambda u(x)$, but this dependency decays exponentially fast with the difference $u(y) - \Lambda u(x)$, with the timescale $\alpha>0$.
We regard this property as a weakened form of the causality property, which itself implies $\frac{\partial \Lambda u(x)} {\partial u(y)} = 0$ whenever $u(y) > \Lambda u(x)$.

Say that an operator $\Lambda$ is $\alpha$-\Kruzkov-contractant if there exists $\rho<1$ such that $\Phi_\alpha := \rE_\alpha \Lambda_\infty \rE_\alpha^{-1}$ is $\rho$-sub-additive, see \cref{prop:contractivity}. Then for any $\delta \geq \mu > 0$ and $0 < \alpha \leq \beta$:
\begin{align*}
	\!
	\delta\text{-strictly-causal}
	\Rightarrow \mu\text{-strictly-causal}
	\Rightarrow\text{causal}
	\Rightarrow \alpha\text{-\Kruzkov-contractant}
	\Rightarrow \beta\text{-\Kruzkov-contractant}
\end{align*}
The first two implications are trivial, and the last two are established in \cref{prop:causal_implies_subadd,prop:subadd_alpha_beta} respectively, under mild assumptions.
The classical fast marching method relies on causality, and Dial's algorithm uses the even stronger strict-causality property \cite{tsitsiklis1995EfficientTrajectories}. In contrast, the ANB method introduced in this paper uses the weaker $\alpha$-\Kruzkov-contractivity, or the even weaker $\beta$-\Kruzkov-contractivity for some \emph{unknown} parameter $\beta$ see \cref{th:complexity_withUnknown}.
\end{remark}

\begin{remark}[\Kruzkov\ transformation with $\alpha$ independent of $h$]
\label{rem:KruzkovFixed}
Consider a viscosity solution $u$ to some HJB or HJI degenerate elliptic PDE $\cF^*_x(\nabla u(x)) = 1$, where the Hamiltonian $\cF^*$ is positively $1$-homogeneous. Then $v := -\exp(-u/\alpha)$ is a viscosity solution to the elliptic PDE $v(x)/\alpha + \cF^*_x(\nabla v(x)) = 0$, where $\alpha>0$ is a fixed timescale. The presence of a zeroth order term with a positive discount coefficient $1/\alpha$ makes the latter PDE much easier to study \cite{bardi2008viscosity}. 

In the discrete setting, a \Kruzkov\ transformation with a fixed timescale $\alpha>0$ independent of the grid scale $h$ is considered in \cite{akian2018minimax,kushner1990numerical}. This yields a contracting scheme in the $\ell^\infty$ norm, with contraction factor $\rho \approx 1-c h/\alpha$ as $h\to 0^+$, where $c>0$ is a constant.
While this property is sufficient to prove the existence and uniqueness of the scheme solution, the contraction ratio is too close to $1$ to be helpful in the complexity analysis.
In contrast, we apply in this paper the \Kruzkov\ change of coordinates with a timescale $\alpha = \Theta( h )$ proportional to the grid scale, so as to obtain a contraction ratio $\rho<1$ independent of $h$.
\end{remark}

\begin{algorithm}
\caption{Adaptive Narrow Band (ANB) \leancheck{algo:anb}}
\label{algo:anb}
\begin{algorithmic}[1]
\State \label{anb:input}%
\textbf{Input:} Tolerance $\ve>0$, constants $\tau>0$, $\alpha>0$, $\rho<1$.
\State \label{anb:ensures}%
\textbf{Ensures at termination:} $\ol u\leq \Lambda \ol u$, $\Lambda \ul u \leq \ul u$, and $(1+\ve) \ol u \geq \ul u$.
\State \label{anb:init}%
\textbf{Initialization:} Set $\ol u=0$ and $\ul u = \infty$, everywhere on $X$.
\State \label{anb:for}%
\textbf{For} $n=1,2,\cdots$ \textbf{until} $(1+\ve) \ol u \geq \ul u$ \textbf{do}:
\State \phantom{bl} \label{anb:hatu}%
Compute the solution of $\hat u = \max \{\ol u, \min\{ n\tau, \Gamma \hat u\} \}$ on $X$ \hfill (See \cref{prop:Dijkstra_band})
\State \phantom{bl} \label{anb:hatu2}%
$\ol u \gets \hat u$
\State \phantom{bl} \label{anb:tn}%
$t_n \gets \min \{n\tau, \min \{\ol u(x) \mid x\in X,\, (1+\ve) \ol u(x) < \ul u(x)\}\}$ %
\State \phantom{bl} \label{anb:Deltan}%
$\Delta_n \gets n\tau-t_n$
\State \phantom{bl} \label{anb:Bn}%
$B_n \gets \{x \in X \mid n\tau-2\Delta_n \leq \ol u(x) < n \tau\}$.
\hfill (Note: 
$n\tau-2 \Delta_n = t_n-\Delta_n$,\ $n\tau = t_n+\Delta_n$)
\State \label{algoline:etan}%
\phantom{bl} $\eta_n \gets (1-\rho) \exp(-n\tau/\alpha)$
\State \phantom{bl} \label{anb:while_sub}
\textbf{Until} $\exp(-\frac{\ol u} \alpha) \leq \exp(-\frac{\Lambda \ol u} \alpha) + \eta_n$ on $B_n$ \textbf{do}:
\textbf{For all} $x \in B_n$ \textbf{in parallel}: $\ol u(x) \gets \Lambda \ol u(x)$
\State \phantom{bl} \label{anb:while_sup}
\textbf{Until} $\exp(-\frac{\Lambda\ul u} \alpha) \leq \exp(-\frac{\ul u} \alpha) + \eta_n$ on $B_n$ \textbf{do}:
\textbf{For all} $x \in B_n$ \textbf{in parallel}: $\ul u(x) \gets \Lambda \ul u(x)$
\end{algorithmic}
\end{algorithm}

\begin{figure}[!t]
\centering
	\def\h{9cm}
	\includegraphics[height=\h]{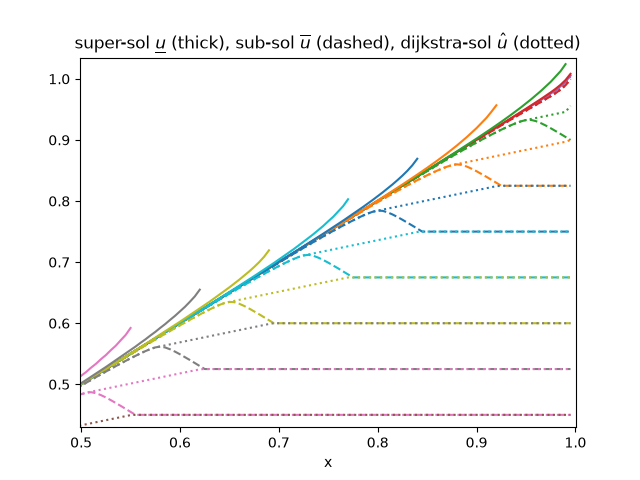}
	\caption{
	Super-solution $\underline u_n \in \bR_\infty^X$ (thick), sub-solution $\overline u_n\in \bR^X$ (dashed), and Dijkstra extension $\hat u_n\in \bR^X$ (dotted), in some successive steps of the ANB method (colors change for each $n$). The ANB method is here applied to the toy one-dimensional problem $|u'|=1$ on the domain $]0,1[$, with source $u(0)=0$, and outflow condition $u(1)=\infty$, discretized using a diffusive Lax-Friedrichs scheme, see \cref{subsec:1D} for more discussion and \cref{fig:1d} for the complete picture. 
	}
	\label{fig:sub_sup_dij_1D}
\end{figure}

\begin{remark}[Challenges of the sub-solution construction]
\label{rem:subsol_construction}
Narrow band methods turn the global problem \eqref{eq:fixedPoint} on the finite set $X = \Int X \sqcup \partial X$ into a sequence of sub-problems posed on successive and possibly overlapping subdomains $B_n \subset X$, $n \geq 0$. 
On the $n$-th iteration one typically has $X = A_n \sqcup B_n \sqcup C_n$, where the numerical solution is essentially frozen in $A_n$, and is completely unknown in $C_n$. 
This raises the issue of choosing the boundary condition on $C_n$ for the sub-problem posed on $B_n$. When building a super-solution, the outflow boundary condition $+\infty$ on $C_n$ yields a small and acceptable boundary layer. When building a sub-solution, the simplest counterpart is to choose the boundary value $0$ on $C_n$, but this creates a huge, domain-wide boundary layer, which is not acceptable, see the discussion in \cite{cristiani2009fast}. 
For this reason, our sub-solution construction involves an additional ingredient: an extension by solving an obstacle-type problem via Dijkstra's method, see line \ref{anb:hatu} of \cref{algo:anb} and \cref{fig:sub_sup_dij_1D}, which (i) does \emph{not} require boundary conditions on $C_n$, (ii) preserves the sub-solution property, and (iii) provides boundary conditions for the next sub-problem.
Iterative methods, such as the FS, AGSI and FIM \cite{zhao2005FastSweeping,bornemann2006FiniteElement,jeong2008fim} face a similar issue: the super-solution $+\infty$ is a much better initialization than the sub-solution $0$, which is why they work with super-solutions in practice.
\end{remark}

\paragraph{Adaptive Narrow Band algorithm.}
The main contribution of this paper is the  numerical solver of the fixed point equation \eqref{eq:fixedPoint} described in \cref{algo:anb}, and its numerical analysis summarized in \cref{th:complexity,th:complexity_withUnknown}.
\Cref{algo:anb} is inspired by the Banded Fast Marching (BFM) method of Cristiani et al.\ \cite{cristiani2009fast,cacace2011local,cacace2014can}, but it is (heavily) modified in such way that a quasi-linear complexity guarantee can be established. On the other hand, compared to \cite{mansour2024solving}, the key improvement is that the narrow band width $\Delta_n$ is adaptive, rather than prescribed a-priori via the scheme parameters - even though we eventually show that $\Delta_n \leq \Delta_{\max}$, where the upper bound $\Delta_{\max}$ depends on the constants associated to the scheme properties. %

For simplicity, \cref{algo:anb} is presented using high level pseudocode, which is sufficient to establish the \emph{update count} complexity result \cref{th:complexity}.
\Cref{algo:anb_lazy}, which is mathematically equivalent by \cref{prop:correct_anb_lazy}, contains a detailed implementation of lines \ref{anb:hatu}, \ref{anb:hatu2}, \ref{anb:tn} and \ref{anb:Bn} involving sorted lists implemented as balanced binary search trees, which ensures quasi-linear complexity in terms of usual \emph{elementary operations}.
Let us comment the successive lines of \cref{algo:anb}.
We denote by $u_*$ the fixed point \eqref{eq:fixedPoint} of $\Lambda$, see  \cref{prop:comparison_Lambda} for existence and uniqueness.

\begin{enumerate}
	\item The user defined tolerance $\ve>0$ expresses the desired solution accuracy. %
	Our first complexity analysis  \cref{th:complexity} assumes that $\Phi_\alpha$ is $\rho$-sub-additive, which puts strong constraints on the method parameters $\alpha$ and $\rho$, similarly to \cite{mansour2024solving}. Strikingly, however, we can \emph{lift this assumption} in \cref{th:complexity_withUnknown}, requiring instead that $\Phi_\beta$ is $\sigma$-sub-additive for some \emph{unknown} constants $\beta>0$ and $\sigma<1$, and still establish quasi-linear complexity \eqref{eq:ln_N_eps}.

	\item The inequalities $\ol u\leq \Lambda \ol u$ and $\Lambda \ul u \leq \ul u$ mean that $\ol u$ and $\ul u$ are respectively a sub- and super-solution to the problem. These properties are maintained throughout the iterations, see \cref{prop:always_subsol}, and imply that $\ol u \leq u_* \leq \ul u$ by the comparison principle, see \cref{prop:comparison_Lambda}. The additional property $(1+\ve) \ol u \geq \ul u$ at termination yields  \eqref{eq:tol_fixed_point}.

	\item
	Let us mention that $\ol u$ (resp.\ $\ul u$) monotonically increases (resp.\ decreases) throughout the iterations, see \cref{prop:always_subsol}.
	The evaluation of the monotone operator $\Lambda : \bR^X \to \bR^X$ on the super-solution $\ul u \in \bR_\infty^X$, which possibly features infinite values, is defined by \eqref{eqdef:monotone_extension}.

	\item The continuation criterion of the for loop ensures that $(1+\ve)\ol u \geq \ul u$ holds at termination.
	We provide in \cref{th:complexity,th:complexity_withUnknown} an upper bound $n_{\max}$ on the number of iterations.

	\item The computation of $\hat u$ is based on Dijkstra's method, see \cref{algo:Dijkstra}. Thanks to a well-chosen initial condition, only points within a small band $B'_n$ are updated by Dijkstra's method in the $n$-th iteration, see \cref{prop:Dijkstra_band}. %

	\item This step increases $\ol u$, and preserves the sub-solution property, see \cref{prop:always_subsol}.
	In the detailed implementation \cref{algo:anb_lazy}, which uses a special symbol $\diamond$ to represent the special cutoff value $n\tau$, 
	the assignment $\ol u(x) \gets \hat u(x)$ actually only needs to be performed at the points $x\in B'_n$ of a small band, see the previous point.

	\item The center $t_n$ of the narrow band time range is defined as the first time when the accuracy tolerance $(1+\ve)\ol u \geq \ul u$ fails (or $n\tau$ if this condition fails later, in which case $\Delta_n = 0$, $B_n = \emptyset$, and the method essentially jumps to the next iteration).
	The detailed \cref{algo:anb_lazy} maintains \eqref{eqdef:M_heap} the list $\{(\ol u(x),x) \mid  \ol u(x) < n\tau \text{ and } (1+\ve) \ol u(x) < \ul u(x)\}$, sorted w.r.t.\ the first component $\ol u(x)$ and implemented using a balanced binary search tree, whose first element yields $t_n$.

	\item As already mentioned, a key result of the complexity analysis is that the (half-)width $\Delta_n$ of the narrow band time range is uniformly bounded by a constant $\Delta_{\max}$, itself of the order $h \ln (1/\ve)$.
	See \cref{corol:max_width,corol:max_width_withUnknown} for the proof of this property.

	\item The narrow band $B_n$ is defined by a range of values $[t_n-\Delta_n,t_n+\Delta_n[$ of the current sub-solution $\ol u$. We show in  \cref{corol:Bn_comp} that $B_n$ is contained in a similar range of values of the exact solution $u_*$. 
	A similar band was defined through values of a super-solution in \cite{mansour2024solving}, but using the sub-solution $\ol u$ is more convenient for the purposes of this paper.	
	
	\item The tolerance parameter $\eta_n$ works best if $\Phi_\alpha$ is $\rho$-sub-additive. 

	\item The termination criterion of the inner while loop is expressed in \Kruzkov\ coordinates \eqref{eqdef:Kruzkov}.
	The contraction property of the operator $\Phi_\alpha$, or $\Phi_\beta$ depending on the assumption, allows bounding the number of iterations, see \cref{prop:niter_while,prop:niter_while_withUnknown}.

	In addition, this \emph{global iteration} within the narrow band $B_n$, can be replaced with \emph{adaptive Gauss-Seidel iteration} \cite{bornemann2006FiniteElement} within $B_n$ again, see \cref{sec:AGSI}. This modification does not improve the theoretical worst case complexity, see \cref{th:complexity_with_agsi}, but the empirical numerical cost is often strongly reduced, see \cref{sec:num}.

	\item Similar to the previous line, but addresses the super-solution, rather than the sub-solution.
\end{enumerate}

We now gather the assumptions scattered in this introduction, and state our main results.
\begin{assumption}
\label{assum:all}
	The operator $\Lambda$ is monotone, additively invariant, has positive boundary conditions, admits a finite super-solution, and has $C_\past$-memory. 
	The operator $\Gamma$ is monotone, additively invariant, causal, underestimates $\Lambda$, and is $C_\Gamma,C'_\Gamma$-comparable to $\Lambda$. We denoted by $C_\past\geq 0$, $C_\Gamma\geq 1$, $C'_\Gamma\geq 0$ some arbitrary constants.
\end{assumption}

\begin{theorem}[\leancheck{th:complexity}]
\label{th:complexity}
Let \cref{assum:all} hold, and further assume that $\Phi_\alpha := \rE_\alpha \Lambda_\infty \rE_\alpha^{-1}$ is $\rho$-sub-additive for some $\alpha>0$, $0 < \rho < 1$, and that $0 < \ve \leq \ve_{\max}$ and $2\Delta_{\max} \leq \min (u_\bd)$.

The operator $\Lambda$ admits a unique fixed point $u_* \in \bR^X$, and satisfies the comparison principle.
\Cref{algo:anb}, with parameters $0 < \ve \leq \ve_{\max}$ and $\tau>0$, terminates using at most $n_{\max}$ iterations.
In the $n$-th such iteration, one has $\Delta_n < \Delta_{\max}$, and the number of inner loops in lines \ref{anb:while_sub} and \ref{anb:while_sup} is at most $N_\Lambda''$.
In addition, any given $x\in X$ appears at most $N'_\Lambda$ times in the narrow band: $\# \{n = 1,2,\cdots \mid x \in B_n\} \leq N'_\Lambda$.
A sub- and a super-solution $\ol u,\ul u \in \bR^X$ are produced, satisfying $(1+\ve) \ol u \geq \ul u$ at termination.

Finally, \cref{algo:anb} uses, at any given point $x \in X$, at most $N_\Lambda$  updates $\ol u(x) \gets \Lambda \ol u(x)$, at most $N_\Lambda$ updates $\ul u(x) \gets \Lambda \ul u(x)$, and at most $N_\Gamma$ updates $\hat u(x) \gets \Gamma \hat u(x)$ assuming that line \ref{anb:hatu} of \cref{algo:anb} is implemented as in \cref{prop:Dijkstra_band}. 
We denoted
\begin{align*}
	\Delta_{\max} &:= 2 \alpha  \ln \big(\frac{2 e} \ve\big) + 2\tau, &
	n_{\max} &:= \Big\lceil \frac{\|u_*\|_\infty + \Delta_{\max}} \tau \Big\rceil, &
	\ve_{\max} &:= \min \big\{\frac \alpha{\|u_*\|_\infty+\Delta_{\max}}, 2\big\},
\end{align*}
\vspace{-0.5cm}
\begin{align*}
	\Delta^+_{\max} &:= (C_\Gamma-1) \Delta_{\max} + \ve \|u_*\|_\infty + C_\past + C'_\Gamma,&
	N'_\Lambda &:= \Big\lceil \frac{2\Delta_{\max} + \Delta^+_{\max}} \tau \Big\rceil,&
\end{align*}
\vspace{-0.5cm}
\begin{align*}
N_\Gamma &:= N'_\Lambda+1, &
N''_\Lambda &:= \Big\lceil \frac{2 \Delta_{\max}/\alpha - \ln (1-\rho)}{-\ln \rho}	\Big\rceil, &
N_\Lambda &:= N_\Lambda' N_\Lambda''.
\end{align*}
\end{theorem}

We now state a variant of \cref{th:complexity} where $\Phi_\beta$ is $\sigma$-sub-additive, for some \emph{arbitrary and unknown constants} $\beta > 0$ and $\sigma < 1$, a property which is always satisfied for the schemes of interest.
The complexity guarantee applies to \cref{algo:anb} whose input parameters $\alpha>0$ and $\rho < 1$ \emph{do not need to satisfy any property}. 
This is convenient for applications, since evaluating a valid \Kruzkov\ timescale $\alpha$ and contraction factor $\rho$ for the considered scheme does have a cost, see \cref{sec:effective}, and since the guaranteed theoretical values may be too pessimistic for practical use.
It improves on \cref{th:complexity} and \cite{mansour2024solving} which require that $\Phi_\alpha$ is $\rho$-sub-additive, for some known constants used as input parameters to the numerical methods. 
 
By \cref{prop:subadd_alpha_beta}, the unknown constant $\beta$ (resp.\ $\sigma$) can be chosen arbitrarily large (resp.\ arbitrarily close to $1$), hence the assumption $\beta \geq\alpha$ (resp.\ $(1-\sigma)^\beta \leq (1-\rho)^\alpha$) does not constrain the choice of the input parameter $\alpha$ (resp.\ $\rho$).

\begin{theorem}[\leancheck{th:complexity_withUnknown}]
\label{th:complexity_withUnknown}
Let \cref{assum:all} hold, and further assume that $\Phi_\beta$ is $\sigma$-sub-additive with $0<\alpha \leq \beta$, $\rho,\sigma \in ]0,1[$ and $(1-\rho)^\alpha \geq (1-\sigma)^\beta$, and that $0 < \ve \leq \ve_{\max}$ and $2\Delta_{\max} \leq \min (u_\bd)$. The conclusions of \cref{th:complexity} hold, using the modified constants:
\begin{align*}
	\Delta_{\max} &:= 2 \beta \ln \big(\frac{2 e C} \ve\big) + 2\tau  \quad 
	\text{where } C := \frac{(1-\rho)^{\frac \alpha \beta}}{1-\sigma}, &
	N_\Lambda'' &:= \Big\lceil \frac{2 \Delta_{\max}/\alpha + \ln (\beta/\alpha) - \ln (1-\rho)}{-\ln \sigma} \Big\rceil.
\end{align*}
\end{theorem}
The constants associated to the scheme properties, and involved in \cref{th:complexity,th:complexity_withUnknown}, are under mild geometric non-degeneracy assumptions related to the grid scale $h>0$ by
\begin{align}
\label{eq:homog_scheme_csts}
	\rho\leq \sigma &\leq  \sigma_{\max} < 1, &
	\tau,\alpha,\beta &= \Theta(h), &
	C'_\Gamma,C_\past &= \cO(h), &
	C_\Gamma &= \cO(1), &
	\|u_*\|_\infty &= \Theta(1),
\end{align}
where the upper bound $\sigma_{\max}$ is independent of $h$, and where $a = \Theta(b)$ means: $a=\cO(b)$ and $b=\cO(a)$, as $h\to 0^+$.
These homogeneity properties were already mentioned in \cref{rem:homog_csts}, and are established in \cref{sec:effective} which is based on  \cite{mansour2024solving}. Thus, by \cref{th:complexity,th:complexity_withUnknown}: %
\begin{align} \label{eq:bigO_h_ln_eps}
	\Delta_{\max},\Delta^+_{\max} &= \cO(-h \ln \ve), &
	n_{\max} &= \cO(h^{-1} - \ln \ve), &
	\ve_{\max} = \Theta(h), \\ \nonumber
	N'_\Lambda, N_\Gamma, N''_\Lambda &= \cO(-\ln \ve), &
	N_\Lambda &= \cO(\ln^2 \ve).
\end{align}

\begin{remark}[Valid range of the tolerance parameter $\ve>0$]
\label{rem:tol_range}
The assumptions of \cref{th:complexity,th:complexity_withUnknown} impose
$\ve \leq \ve_{\max} = \Theta(h)$ and $\min(u_\bd) \geq 2\Delta_{\max} = \cO(-h\ln \ve)$.
For sufficiently small $h>0$, there is therefore a non-empty range of acceptable tolerances $\ve \in [\exp(-c/h), c' h]$, where $c, c' >0$ are positive constants.
\end{remark}

Denote by $N$ the number of discretization points, and by $\widetilde N$ the number of connections between the scheme stencils:
\begin{align}
\label{eqdef:stencil_card}
	N &:= \#(X), & %
	\widetilde N &:= \#\{(x,y) \in \Int X \times X \mid \Gamma u(x) \text{ depends on } u(y)\}. &
	\text{Thus } N, \widetilde N &= \Theta(h^{-d}),
\end{align}
where the last estimate holds under mild geometric assumptions.
The complexity of \cref{algo:anb} in terms of \emph{elementary operations}, when implemented as specified in \cref{algo:anb_lazy}, is as announced, by \cref{prop:complexity_lazy}.
\begin{align}
\nonumber
&\cO(\widetilde N \ln \widetilde N) 
+ 2 N N_\Lambda U_\Lambda 
+ (\widetilde N N_\Gamma + (N+\widetilde N) N'_\Lambda) U_\Gamma 
+ (N+\widetilde N)(N'_\Lambda + N_\Gamma) \cO( \ln N) +  \cO(n_{\max}\ln N)
	\\
	&= \cO( N \ln^2 \ve - N \ln \ve \ln N).
	\label{eq:ln_N_eps}
\end{align}
We denoted by $U_\Lambda$, $U_\Gamma$ and $\cO(\ln N)$ the unit cost of, respectively, the evaluation of $\Lambda u(x)$, $\Gamma u(x)$, and the insertion or deletion operations in a balanced binary search tree with at most $N$ elements. We assumed (\ref{eqdef:stencil_card}, right) in \eqref{eq:ln_N_eps}, and introduced the estimates \eqref{eq:bigO_h_ln_eps}. We also made in \eqref{eq:ln_N_eps} the assumption  $U_\Lambda, U_\Gamma = \cO(1)$, which is justified if $\Lambda u(x)$ admits a closed form expression, e.g.\ for the Lax-Friedrichs and upwind differences schemes, as well as the semi-Lagrangian scheme with an isotropic, Riemannian or Randers metric \cite{sethian2003OUM}. On the other hand, when the semi-Lagrangian scheme is applied to complex anelliptic norms arising in seismology as in \cite{desquilbet2021single}, the computation of $\Lambda u(x)$ involves numerically solving an inner convex optimization problem via a Newton-like method, which has a significant cost.

\paragraph{Comparison with the state of the art.}
The computation of minimal paths is a ubiquitous problem in applied mathematics, for which numerous methods have been considered; a popular approach, followed in this paper, is to numerically solve a Hamilton-Jacobi-Bellman (HJB) PDE, from which the minimal paths can be extracted by backtracking. 
\emph{Non-HJB} methods include linearization techniques based on the Hopf-Cole transform \cite{crane2013GeodesicsHeat}, exact methods on triangulated surfaces \cite{surazhsky2005exactDistance}, a variant of Dijkstra's method with short term memory \cite{campen2013PracticalGeodesy}, and many others see the survey \cite{crane2020geodesicSurvey} and references therein. The advantages of HJB methods include their strong mathematical guarantees based on the theory of viscosity solutions to degenerate elliptic PDEs \cite{bardi2008viscosity}, and their extension to 
degenerate Hamiltonians, associated to non-holonomic vehicle models \cite{mirebeau2018curvature}, and to non-convex Hamiltonians associated to differential games \cite{cacace2011local} (Hamilton-Jacobi-Isaacs PDEs).

The proposed ANB method estimates the numerical solution of a HJB or HJI PDE of interest up to a tolerance $\ve>0$ given by the user \eqref{eq:tol_fixed_point}, and has quasi-linear complexity \eqref{eq:ln_N_eps} w.r.t.\ the number $N$ of discretization points.
These guarantees hold even if the constants associated to the scheme properties are unknown, as established in \cref{th:complexity_withUnknown}.
Let us compare this result with alternative methods in the literature:
\begin{itemize}
	\item The fast marching method (FMM) \cite{tsitsiklis1995EfficientTrajectories,sethian1996fast} computes the fixed point $u_*$ of the operator $\Lambda$ \emph{exactly}, and has quasi-linear complexity $\cO(N\ln N)$. However, this algorithm is intrinsically sequential, and is only applicable if $\Lambda$ obeys the causality property \eqref{eqdef:causality}. The latter constraint typically  imposes the use of wide and complex stencils in anisotropic settings, and excludes non-convex Hamiltonians, see \cite{kimmel1998Manifold,sethian2003OUM,vladimirsky2008LabelSetting,alton2012MAOUM,mirebeau2014marchingLattice,mirebeau2014EfficientFinsler,mirebeau2018curvature,mirebeau2019riemannian,dahiya2018ordered}.
	As already discussed, some variants of the FMM achieve linear $\cO(N)$ complexity, at the price of stronger assumptions or of small additional numerical errors: see Dial's method \cite{tsitsiklis1995EfficientTrajectories}, the untidy FMM \cite{yatziv2006untidyFMM}, and group marching \cite{kim2001groupMarching}. %

	\item Iterative solvers such as the adaptive Gauss-Seidel iteration (AGSI), fast iterative method (FIM), and fast sweeping (FS) \cite{bornemann2006FiniteElement,jeong2008fim,zhao2005FastSweeping,tsai2003fastSweeping,qian2007fastSweeping} produce a super-solution $\ul u \in \bR^X$ such that $0 \leq \ul u-\Lambda \ul u \leq \ve$.
	This tolerance criterion, on the scheme residual, slightly differs from ours, see \cite[Remark 3.5]{mansour2024solving} and \cref{corol:tol}.
	The empirically observed complexity is between linear $\cO(N)$ in favorable cases \cite{zhao2005FastSweeping,luo2016convergenceFSM}, and super-linear $\cO(N^{1+\frac 1 d})$ if the scheme is diffusive, if it uses a non-standard stencil, or if the minimal paths are complex, %
	see \cite[Remark 4.4]{mansour2024solving} and \cref{sec:num}.
	
	\item
	Narrow band methods turn a global problem into a sequence of sub-problems posed on possibly overlapping subdomains	\cite{adalsteinsson1995fast,cristiani2009fast,cacace2011local,cacace2014can,mansour2024solving}.
	Only the latter \cite{mansour2024solving} has a complexity guarantee: it produces a super-solution $\ul u\in \bR^X$ such that $u_* \leq \ul u \leq u_* + \ve$, with quasi-linear complexity $\cO(N \ln^2 (N/\ve))$. However, for this guarantee to hold, the algorithm parameters need to be set to specific theoretical values depending on the constants $(\rho,\alpha,C_\past,%
	\nu)$ associated to the scheme properties.
\end{itemize}
In practice, the constants associated with the scheme properties are rarely evaluated, since their computation is not free, and since for various reasons the associated parameters and estimates are often excessively pessimistic: indeed, the theoretical analysis is often sub-optimal, the estimates correspond to the worst case scenario, and the constants may reflect a localized difficulty in a large inhomogeneous model.
If these constants are regarded as unknown, then the ANB and FMM are the only methods with guarantees of uniform accuracy and quasi-linear complexity.

The narrow band methods closest to the proposed ANB are banded fast marching (BFM) \cite{cristiani2009fast}, and the non-adaptive narrow band (NB) method introduced in \cite{mansour2024solving}.
BFM contains like the ANB the idea of advancing the narrow band based on an \emph{a-posteriori} error estimate obtained by comparing a super-solution \emph{with a sub-solution};
however it does not come with a complexity guarantee.
Two additional works are close in spirit to the ANB:
the method \cite{clawson2014causalDomain}, related to the $A^*$ algorithm, uses both a lower and a upper estimates of the solution to restrict the computational domain of an eikonal solver; the two-scale method \cite{chacon2012twoScale} combines marching on a coarse scale with sweeping on a fine scale, i.e.\ an ordered Dijkstra-like method with local iterative solvers, similarly to the NB and ANB.

The NB method \cite{mansour2024solving} relies like the ANB on the $\rho$-sub-additivity of the operator $\Phi_\alpha$ in \Kruzkov\ coordinates, to establish a quasi-linear complexity guarantee $\cO(N \ln^2 (N/\ve))$, the first one to the authors' knowledge for a method applicable to non-causal schemes. 
A key difference is that the width of the narrow band is fixed a-priori in the NB method, whereas it is adaptive in the ANB. 
For this reason, the theoretical convergence and complexity guarantees of the NB require specific input parameters depending on the scheme properties, raising the issues discussed above, while the ANB can in principle accommodate arbitrary parameters, see \cref{th:complexity_withUnknown} and the discussion in \cref{sec:num}. 

Finally, let us acknowledge that GPU accelerated solvers of eikonal equations \cite{kimmel2008ParallelDistance,jeong2008fim,mirebeau2023massively} would likely vastly outperform the proposed ANB method in terms of computation time when applicable, thanks to parallelization over thousands of threads, and despite their lack of complexity guarantee. Efficiently parallelizing the ANB on such architectures is non-trivial, and will be the object of future work.

\begin{remark}[Eikonal solvers with $\cO(N)$ complexity]
\label{rem:linear_complexity}
Dial's method \cite{tsitsiklis1995EfficientTrajectories} achieves this if the scheme is \emph{strictly causal}, by replacing the priority queue of the fast marching method with a list of buckets. If the scheme is \emph{causal} but not strictly causal, and one is willing to tolerate small numerical errors, then the implementation \cite{yatziv2006untidyFMM} of the FMM using an untidy priority queue similarly has linear complexity; see also the group marching variant \cite{kim2001groupMarching}. 
The fast sweeping method has $\cO(k N)$ complexity, where $k = k(N,\ve)$ denotes the number of sweeps, which is often presented as being bounded $k=\cO(1)$ independently of the problem size and target tolerance, resulting in linear complexity.
For the eikonal equation with unit speed, which characterizes the Euclidean distance function, discretized using the standard $2d$-point scheme, \cite{zhao2005FastSweeping} proves indeed that $k = 2^d = \cO(1)$ sweeps suffice to converge exactly ($\ve=0$) if the source is a single point, and for a general source set to converge approximately ($\ve>0$, unspecified) and reach an accuracy $\cO(h |\ln h|)$ comparable to the discretization error, where $h \approx N^{-\frac 1 d}$ denotes the grid scale. For a general speed function, \cite{zhao2005FastSweeping} relates the number $k$ of sweeps to the number of times the tangent to a characteristic curve (i.e.\ a minimal geodesic to the source set) changes orthant, which is itself related to its curvature\footnote{This is a legitimate heuristic, but not a complete argument, since for instance the curve $t \mapsto (t, \delta^2 \sin(t/\delta))$ has curvature at most $1$, for any $\delta > 0$, while its tangent changes quadrant every $\pi\delta$.}. 
A strength of the FS method is that it can also solve \emph{non-causal}\footnote{``Causality'' in \cite{luo2016convergenceFSM} is related to our monotonicity assumption, and \emph{differs} from the causality property defined below and used by FMMs \cite{tsitsiklis1995EfficientTrajectories,sethian1996fast, kimmel1998Manifold,sethian2003OUM,vladimirsky2008LabelSetting,alton2012MAOUM,mirebeau2014marchingLattice,mirebeau2014EfficientFinsler,mirebeau2018curvature,mirebeau2019riemannian,dahiya2018ordered}.} monotone schemes \cite{tsai2003fastSweeping, kao2004lax,qian2007fastSweeping,luo2016convergenceFSM}. 
Convergence is established in that context \cite{luo2016convergenceFSM} as $k \to \infty$, and empirically a fixed number of sweeps $k = \cO(1)$ still suffices in many cases. However, a number of sweeps $k = \cO(h^{-1})$ growing as the grid scale $h$ is refined seems hard to avoid if one uses the diffusive Lax-Friedrichs scheme, or a scheme with a non-standard stencil, or if a characteristic curve exhibits a complex space filling behavior, see \cref{sec:num}; similarly, $k = \cO(h^{-1})$ in the fully periodic linear example considered in \cite{luo2016convergenceFSM}.
We refer to \cite{gomez2019fastMethods} for a survey and an experimental comparison of these various methods.
\end{remark}

\paragraph{Outline.}
Our main complexity results \cref{th:complexity,th:complexity_withUnknown} are established in \cref{sec:complexity}.
Numerical experiments illustrate our findings in \cref{sec:num}. A small modification of \cref{algo:anb}, using adaptive Gauss-Seidel iterations within the narrow band subdomains, is discussed in \cref{sec:AGSI}. Semi-Lagrangian, Lax-Friedrichs, and upwind-differences discretization schemes are described in \cref{sec:effective}, and it is shown that they satisfy \cref{assum:all} with explicit constants. \Cref{sec:subadd} establishes the results tying the causality property with the $\rho$-sub-additivity of $\Phi_\alpha$ announced in \cref{rem:causality}.

\paragraph{Notations.}
We fix a finite set $X = \Int X \sqcup \partial X$, which is the union of two disjoint subsets $\Int X$ and $\partial X$. For all $u \in \bR^X$, $t\in \bR$, we denote by $\min \{u,t\} \in \bR^X$ the mapping $x \in X \mapsto \min \{u(x),t\}$, and likewise $\max\{u,t\}$.
We denote $\{u<t\} := \{x \in X \mid u(x) < t\}$, when the domain $X$ is clear from context. $\bR_\infty:=]-\infty,\infty]$.

Operators denoted by the symbols $\Lambda,\Gamma : \bR^X \to \bR^X$ are subject to the boundary condition $\Lambda u = \Gamma u = u_\bd$ on $\partial X$, for all $u\in \bR^X$, where $u_\bd\in \bR^{\partial X}$ is some fixed boundary data. Sub-solutions $\ol u \leq \Lambda \ol u$ and super-solutions $\Lambda \ul u \leq \ul u$ of $\Lambda$, are usually distinguished with an over- and under-line respectively, whereas solutions $\hat u = \Gamma \hat u$ of the auxiliary operator $\Gamma$ are usually distinguished with a hat.

\section{Complexity estimate}
\label{sec:complexity}
This section is devoted to the proof of our main complexity estimates \cref{th:complexity,th:complexity_withUnknown}.
We recall in \cref{subsec:comp_principle} two classical approaches for computing the fixed point of a monotone operator: (i) fixed point iterations, if the operator is in addition $\rho$-sub-additive, or (ii) Dijkstra's algorithm, if the operator is in addition causal.
We study in \cref{subsec:dijkstra_under} the improvement $\hat u$ of the subsolution $\ol u$ via Dijkstra's method, see line \ref{anb:hatu} in \cref{algo:anb}. We establish \cref{th:complexity} in \cref{subsec:cplx_subadd}, and \cref{th:complexity_withUnknown} in \cref{subsec:unknown_params}. Finally, we describe and study in \cref{subsec:heaps} a lower level implementation \cref{algo:anb_lazy} of the proposed method \cref{algo:anb}, which is proved to be equivalent and to achieve quasi-linear complexity \eqref{eq:ln_N_eps} in terms of elementary operations, rather than updates.

\subsection{General properties of monotone, $\rho$-sub-additive, and causal operators}
\label{subsec:comp_principle}

This section recalls some basic and foundational properties of discrete monotone operators, see e.g.\ \cite{oberman2006convergent} for more discussion of this rich topic.
We first discuss their extension to positively infinite values and their restriction to a subdomain. We then turn to the comparison principle, the existence and uniqueness of a fixed point, and its computation. For that purpose, monotonicity alone is not enough, but needs to be combined with either $\rho$-sub-additivity, with $\rho<1$, or causality.
Some elementary and classical statements are reproduced from \cite{mansour2024solving} for completeness, together with their (straightforward) proofs.

\paragraph{Value extension, and domain restriction, of monotone operators.}
We first gather elementary properties of the extension \eqref{eqdef:monotone_extension} of a monotone operator to $\bR_\infty := ]-\infty,\infty]$, i.e.\ to \emph{positively infinite} values. This is a prerequisite for the definition of the \Kruzkov\ transform \eqref{eqdef:Phi}, which is used in the analysis of the proposed ANB method. 

\begin{proposition}[\leancheck{prop:monotonic_Phi}]
\label{prop:monotonic_Phi}
	Let $\Lambda : \bR^X \to \bR^X$ be a monotone operator. Define $\Lambda_\infty : \bR_\infty^X \to \bR_\infty^X$  as in \eqref{eqdef:monotone_extension}, and $\Phi_\alpha := \rE_\alpha \Lambda_\infty \rE_\alpha^{-1} : \bR^X \to \bR^X$ as in \eqref{eqdef:Phi}. Then $\Lambda_\infty u = \Lambda u$ for all $u \in \bR^X$, and $\Lambda_\infty$ and $\Phi_\alpha$ are monotone.
	If in addition $\Lambda$ is additively invariant, then $\Phi_\alpha$ is multiplicatively invariant in the sense that $\Phi_\alpha(\gamma v) (x) = \gamma \Phi_\alpha v (x)$ for all $\gamma>0$, $v \in \bR^X$, and $x \in \Int X$. If $\Lambda$ is monotone and causal, then so is $\Lambda_\infty$.
\end{proposition}

\begin{proof}
For any $u \in \bR_\infty^X$, $x \in X$, the mapping $t \in \bR \mapsto \Lambda \min \{u,t\} (x) \in \bR$ is non-decreasing by monotonicity of $\Lambda$. The following limit is therefore well-defined, for any $u \in \bR^X_\infty$ and $x \in X$
\begin{equation*}
	\Lambda_\infty u(x) := \lim_{t \to \infty} \Lambda \min \{u,t\} (x) = \sup_{t \in \bR} \Lambda \min \{u,t\} (x) \in ]-\infty,\infty].
\end{equation*}
Furthermore $\Lambda_\infty$ is monotone, as a supremum of monotone functions. Thus $\Phi_\alpha$ is monotone by composition with \Kruzkov's monotone changes of variables $\rE_\alpha$ and $\rE_\alpha^{-1}$.
If $u \in \bR^X$, then $\min\{u,t\} = u$ for all $t \geq \max (u)$, hence $\Lambda_\infty u = \Lambda u$.
The multiplicative invariance of $\Phi_\alpha$ immediately follows from the properties of \Kruzkov's transform \eqref{eqdef:Kruzkov}, namely the exponential $\rE_\alpha$ and logarithm $\rE_\alpha^{-1}$, and the additive invariance of $\Lambda$. Indeed for all $\gamma>0$, $v \in \bR^X$, $x \in \Int X$:
\begin{equation*}
	\Phi_\alpha (\gamma v) (x) = \rE_\alpha \Lambda_\infty (\rE_\alpha^{-1} v - \alpha \ln \gamma)(x) = \rE_\alpha( \Lambda_\infty \rE_\alpha^{-1} v - \alpha \ln \gamma) (x)= \gamma \rE_\alpha \Lambda_\infty \rE_\alpha^{-1} v(x).
\end{equation*}
Finally, assume that $\Lambda$ is monotone and causal, and let $u,v\in \bR_\infty^X$ and $T \in \bR$ such that $\min \{u,T\} = \min \{v,T\}$. 
Then $\min \{\min \{u,t\},T\} = \min \{\min \{v,t\},T\}$ for any $t \in \bR$, hence by causality one has (\ref{eqdef:causality}, right), which can be reformulated as follows: for any $x \in X$, either
$\Lambda \min \{u,t\} (x) = \Lambda \min \{v,t\}(x)$ or both terms are $> T$. Since $\Lambda \min \{u,t\}(x)$ and $\Lambda \min \{v,t\}(x)$ are \emph{non-decreasing} w.r.t.\ $t\in \bR$, by monotony, we can pass this property to the limit as $t \to \infty$, and obtain that $\Lambda_\infty u(x) = \Lambda_\infty v(x)$ or both terms are $>T$. Since $x$ is arbitrary and $X$ is finite, this yields the causality property for $\Lambda_\infty$, as announced.
\end{proof}

Given an operator $\Lambda : \bR^X \to \bR^X$, and a sub-domain $B\subset X$ (resp.\ and some boundary data $v \in \bR^{X\sm B}$), we define the restriction $\Lambda_B : \bR^X \to \bR^X$ (resp. $\Lambda_B^v : \bR^B \to \bR^B$) as follows 
\begin{align}
\label{eqdef:Lambda_B}
	\Lambda_B u(x) &:= \begin{cases}
		\Lambda u(x) & \text{if } x \in B,\\
		u(x) & \text{else,}
	\end{cases}
	&
	\text{resp. }
	\Lambda_B^v w(x) &:= 
	\Lambda [\bbone_B w + \bbone_{X\sm B} v](x)
\end{align}
for all $u \in \bR^X$ and $x \in X$ (resp.\ $w\in \bR^B$ and $x \in B$).
Note that lines \ref{anb:while_sub} and \ref{anb:while_sup} of the proposed \cref{algo:anb} implicitly involve $\Lambda_{B_n}$, where $B_n$ is the computational narrow band. 
The following is clear. 

\begin{lemma}[\leancheck{lem:monDom}]
\label{lem:monDom}
Using the notations of \eqref{eqdef:Lambda_B}.
If $\Lambda : \bR^X \to \bR^X$ is monotone, then so are $\Lambda_B$ and $\Lambda_B^v$. If in addition $\Lambda$ is $\rho$-sub-additive, then so is $\Lambda_B^v$.
\end{lemma}

\begin{proposition}[Partial application preserves sub- and super-solutions, {\cite[Prop 2.9]{mansour2024solving}} \leancheck{prop:partialMonotone}]
\label{prop:partialMonotone}
Let $\Lambda : \bR^X \to \bR^X$ be a monotone operator, and let $B \subset X$ be a subset.
Let $u \in \bR^X$ and $u' := \Lambda_B u$.  If $\Lambda u \leq u$ then $\Lambda u' \leq u'$. Likewise if $\Lambda u \geq u$ then $\Lambda u' \geq u'$.
\end{proposition}
\begin{proof} %
Consider the case $\Lambda u \leq u$ of a super-solution. Then $\Lambda u \leq u' \leq u$ by assumption, thus $\Lambda u' \leq \Lambda u$ by monotonicity, and therefore $\Lambda u' \leq u'$ as announced. The case $\Lambda u \geq u$ is similar. 
\end{proof}

\paragraph{Monotonicity and $\rho$-sub-additivity.}
These two properties together imply a contraction property, hence the existence, uniqueness, and computability of a fixed point of the operator by the Banach-Picard theorem. We also recall a discrete version of the comparison principle.

\begin{proposition}[Contraction, see {\cite[Corollary 2.3]{mansour2024solving}} \leancheck{prop:contractivity}]
\label{prop:contractivity}
	Let $\Phi : \bR^X \to \bR^X$ be monotone and $\rho$-sub-additive. Then $\|\Phi u - \Phi v\|_\infty \leq \rho \|u-v\|_\infty$ for all $u,v \in \bR^X$. If $\rho<1$, then $\Phi$ admits a unique fixed point $v_*\in \bR^X$, and $\|\Phi^n v - v_*\|_\infty \leq \rho^n \|\Phi v - v\|_\infty/(1-\rho)$ for all $v \in \bR^X$.
\end{proposition}

\begin{proof}
	One has $v-t \leq u \leq v+t$ with $t := \|u-v\|_\infty$, hence $\Phi v- \rho t \leq \Phi u \leq \Phi v + \rho t$ by $\rho$-sub-additivity, i.e.\ $\|\Phi u-\Phi v\|_\infty \leq \rho t$ as announced. 
	Then apply Picard-Banach's theorem.
\end{proof}

\begin{proposition}[{\cite[Proposition 2.6]{mansour2024solving}} \leancheck{prop:pertSubSuper}]
\label{prop:pertSubSuper}
	Let $\Phi : \bR^X\to \bR^X$ be a $\rho$-sub-additive operator, for some $\rho<1$. Let $u \in \bR^X$ and $\ol c,\ul c \geq 0$ be such that $u-\ol c \leq \Phi u \leq u + \ul c$.
	Then $\ol u:= u-\ol c/(1-\rho)$ and $\ul u :=u+\ul c/(1-\rho)$ are a sub-solution and a super-solution respectively: $\ol u \leq \Phi \ol u$ and $\Phi \ul u \leq \ul u$.
\end{proposition}
\begin{proof} 
	For any $t \geq 0$ one has
	$\Phi (u+t)
	\leq \Phi u+\rho t
	\leq u + \ul c +\rho t
	= (u+t) + \ul c-(1-\rho) t
	$, hence $\ul u := u+t$ is a super-solution provided $t \geq \ul c/(1-\rho)$, as announced, and likewise for $\ol u$.
\end{proof}

\begin{proposition}[Comparison principle, see {\cite[Proposition 2.7]{mansour2024solving}} \leancheck{prop:comparison_Phi}]
\label{prop:comparison_Phi}
Let $\Phi : \bR^X \to \bR^X$ be monotone and $\rho$-sub-additive, with $\rho<1$.
Then for any partition $X = \Int Y \sqcup \partial Y$ and $u, v \in \bR^X$:
\begin{equation*}
	(u - \Phi u \leq v - \Phi v \text{ on } \Int Y, \text{ and } u \leq v \text{ on } \partial Y)\ \Rightarrow \ u \leq v.
\end{equation*}
\end{proposition}

\begin{proof}
	Let $x \in X$ be such that $t:=u(x)-v(x)$ is maximal, and assume for contradiction that $t>0$, thus $\Phi u(x) \leq \Phi (v+t)(x) \leq \Phi v(x) + \rho t$ by monotony and $\rho$-sub-additivity.
	On the other hand one has $x \in \Int Y$, since $u \leq v$ on the complement $\partial Y$, and therefore $t=u(x)-v(x) \leq \Phi u(x)-\Phi v(x)$ by assumption, which contradicts our previous estimate.
\end{proof}

The comparison principle for an operator $\Lambda$, and the existence of a fixed point, can be deduced from the same properties for the operator $\Phi_\alpha$ written in \Kruzkov\ coordinates. The next result also shows that $\Lambda$ admits a unique fixed point under the assumptions of \cref{th:complexity,th:complexity_withUnknown}.
\begin{proposition}[\leancheck{prop:comparison_Lambda}]
\label{prop:comparison_Lambda}
	Assume that $\Lambda : \bR^X \to \bR^X$ is monotone, and that $\Phi_\alpha := \rE_\alpha \Lambda_\infty \rE_\alpha^{-1}$ is $\rho$-sub-additive for some $\alpha>0$ and $\rho<1$.
	Then for any partition $X = \Int Y \sqcup \partial Y$ and any $u, v \in \bR^X$:
\begin{equation}
\label{eq:comparison_Lambda}
	(u \leq \Lambda u \text{ and } \Lambda v \leq v \text{ on } \Int Y, \text{ and } u \leq v \text{ on } \partial Y)\ \Rightarrow \ u \leq v.
\end{equation}
If $\Lambda \ul u \leq \ul u$ for some $\ul u \in \bR^X$, 
then $\Lambda$ admits a unique fixed point $u_* \in \bR^X$, and $u_* \leq \ul u$.
\end{proposition}

\begin{proof}
	By \cref{prop:monotonic_Phi}, the operator $\Phi_\alpha$ is monotone. Since by assumption $\Phi_\alpha$ is $\rho$-sub-additive, it admits a unique fixed point $v_*\in \bR^X$ by \cref{prop:contractivity}, and satisfies the comparison principle by \cref{prop:comparison_Phi}.
	Denoting $\ul v := \rE_\alpha \ul u$, and observing that $\ul v \geq \Phi_\alpha \ul v = \rE_\alpha \Lambda \ul u$ by monotonicity of the \Kruzkov\ transform, we obtain	
	\begin{align*}
		v_* - \Phi_\alpha v_* &= 0 \leq \ul v - \Phi_\alpha \ul v \text{ on } X. &
		\text{Therefore } v_* &\leq \ul v \leq \rE_\alpha(\max (\ul u)) < 0,
	\end{align*}
	by the comparison principle for $\Phi_\alpha$ with  $\Int Y := X$  and $\partial Y := \emptyset$.
	It follows that $u_* := \rE_\alpha^{-1} v_*$ takes finite values, hence it is a fixed point of $\Lambda$.

	Considering $u,u' \in \bR^X$, and denoting $v := \rE_\alpha u$ and $v' := \rE_\alpha u'$, we have the following sequence of equivalences and implications, which establishes the comparison principle for $\Lambda$:
	\begin{align*}
		(u \leq \Lambda u  \text{ and } & \Lambda u' \leq u' \text{ on } \Int Y, \text{ and } u \leq u' \text{ on } \partial Y)\\
		\Rightarrow\ & (v \leq \Phi_\alpha v \text{ and } \Phi_\alpha v' \leq v' \text{ on } \Int Y, \text{ and } v \leq v' \text{ on } \partial Y)\ \Rightarrow\ v \leq v' \ \Leftrightarrow\ u \leq u'.
	\end{align*}
	It immediately follows that $u_*$ is the unique fixed point of $\Lambda$, and that $u_* \leq \ul u$.
\end{proof}

The next result compares the convergence tolerance of the ANB, on the relative error between a sub- and a super-solution, with the usual tolerance on the scheme residual underlying the Fast-Sweeping, AGSI, and FIM methods \cite{zhao2005FastSweeping,bornemann2006FiniteElement,jeong2008fim}. 
All schemes considered in \cref{sec:effective} are $\nu$-sub-multiplicative, with $\nu=\Theta(h)$ proportional to the grid scale $h$.
\begin{corollary}[Convergence tolerance \leancheck{corol:tol}]
\label{corol:tol}
Consider an operator $\Lambda : \bR^X \to \bR^X$ which obeys the comparison principle \eqref{eq:comparison_Lambda}, and admits a fixed point $u_*\in \bR^X$. Let also $\ve\geq 0$.
\begin{itemize}
	\item Assume that $\Lambda$ is monotone. If $\ol u, \ul u \in \bR^X$ obey $\ol u \leq \Lambda \ol u$, $\Lambda \ul u \leq \ul u$, and $(1+\ve) \ol u \geq \ul u$, then 
	\begin{align*}
		(1+\ve)^{-1} u_* & \leq \ol u \leq u_*, &
		u_* & \leq \ul u \leq (1+\ve) u_*, &
		0 \leq \Lambda \ol u - \ol u & \leq \tfrac \ve {1+\ve} u_*, &
		0 \leq \ul u-\Lambda \ul u &\leq \ve u_*.
	\end{align*}
	\item 
	Assume in addition that $\Lambda$ is $\nu$-sub-multiplicative, for some $\nu >0$, and has non-negative boundary conditions $u_\bd \geq 0$.
	If $\ul u \in \bR^X$ obeys $0 \leq \ul u - \Lambda \ul u \leq \ve$ (resp.\ 
	$\ol u \in \bR^X$ obeys $0 \leq \Lambda \ol u - \ol u\leq \ve $, where $0 \leq \ve < \nu$), and $\ul u = u_\bd$ (resp.\ $\ol u=u_\bd$) on $\partial X$ then 
	\begin{align*}
		u_* &\leq \ul u \leq \big(1+\frac \ve \nu\big) u_*,  &
		\Big(\text{resp.}\ \big(1-\frac \ve \nu\big) u_* \leq \ol u \leq u_*\Big).
	\end{align*}
\end{itemize}
\end{corollary}

\begin{proof}
	\emph{Proof of the first point.} By the comparison principle, $\ol u \leq u_* \leq \ul u$, and by the additional assumption $(1+\ve)^{-1} \ul u \leq u_* \leq (1+\ve) \ol u$, which implies the first two inequalities.
	Since $\Lambda$ is monotone, $\Lambda \ol u$ and $\Lambda \ul u$ are respectively a sub- and a super-solution, 
	hence $\Lambda \ol u \leq u_* \leq \Lambda \ul u$, and therefore $\ul u - \Lambda \ul u \leq (1+\ve) u_* - u_* = \ve u_*$ and $\Lambda \ol u - \ol u \leq u_* - (1+\ve)^{-1} u_* = \frac \ve{1+\ve} u_*$ as announced. 
	
	\emph{Proof of the second point.}
	One has $\Lambda [(1-\eta) \ul u] \geq (1-\eta) \Lambda \ul u +\nu \eta\geq (1-\eta) (\ul u-\ve) +\nu \eta$ on $\Int X$, for any $\eta \in ]0,1[$ by $\nu$-sub-multiplicativity with $1+t := (1-\eta)^{-1}$, as also observed in \cite[Lemma 3.2]{mansour2024solving}; 
	choosing $\eta := \ve/(\nu + \ve)$  we find that $(1-\eta)\ul u$ is a sub-solution. Therefore $(1-\eta) \ul u \leq u_*$, which is the announced inequality, by the comparison principle since $(1-\eta) \ul u = (1-\eta) u_\bd \leq u_\bd$ on $\partial X$. 
	Likewise, one has $\Lambda [(1+\eta) \ol u] \leq (1+\eta) \Lambda \ol u - \nu \eta \leq (1+\eta) (\ol u+\ve) - \nu \eta$ on $\Int X$, for any $\eta \geq 0$, and choosing $\eta := \ve/(\nu-\ve)$ we find that $(1+\eta) \ol u$ is a super-solution, hence $(1+\eta) \ol u \geq u_*$ as announced.	
\end{proof}

\paragraph{Monotony and causality.}
These two properties together imply the existence and computability of a fixed point, via Dijkstra's algorithm which is an instantiation of the dynamic programming principle, see \cref{algo:Dijkstra,prop:Dijkstra_correct}. We also prove a version of the comparison principle.
To each $x \in \Int X$ we associate a stencil $\cV(x)\subset X$, such that $\Gamma u(x)$ only depends on the values $(u(y))_{y \in \cV(x)}$. The reversed stencil $\widetilde \cV(y)$, for all $y \in X$, is obtained by reversing this neighborhood relation:
\begin{align}
\label{eqdef:stencil}
	\forall u, v \in \bR^X, \big(u_{|\cV(x)} &= v_{|\cV(x)} \Rightarrow \Gamma u(x) = \Gamma v(x)\big), &
	\widetilde \cV(y) &:= \{x \in \Int X \mid y\in \cV(x)\}.
\end{align}
In order to break ties in Dijkstra's method, and make it deterministic, we assume in \cref{algo:Dijkstra} that $X$ is linearly ordered. This order is arbitrary and has no influence on the method's output.
Given $T \in \bR_\infty$, $u\in \bR_\infty^X$, we define the cutoffs:
\begin{align}
\label{eqdef:u_less_t}
	u^{<T} (x) &:=
	\begin{cases}
		u(x) \!\!& \text{if } u(x) < T\\
		\infty  &\text{else}
	\end{cases} &
	u^{\leq T} (x) &:=
	\begin{cases}
		u(x) \!\!& \text{if } u(x) \leq T\\
		\infty & \text{else}
	\end{cases}
\end{align}
We also consider the operator $\Gamma^{<T} : u \mapsto (\Gamma u)^{<T}$. Note that $u^{<\infty} = u^{\leq \infty} = u$, and likewise $\Gamma$.
The sorted list in \cref{algo:Dijkstra} should be implemented as a heap, or as a balanced binary search tree if range queries and deletions at arbitrary positions are needed as in \cref{algo:anb_lazy}.

\begin{algorithm}
\caption{Generalized Dijkstra algorithm, with stopping time $T \in \bR_\infty$. \leancheck{algo:Dijkstra}}
\label{algo:Dijkstra}
\begin{algorithmic}[1]
\State \textbf{Input}: monotone and causal  $\Gamma$, reversed stencils $\widetilde \cV$, initial guess $u \in \bR_\infty^X$ s.t. $(\Gamma u)^{<T} \leq u$.
\State\textbf{Ensures at termination}: $u = (\Gamma u)^{<T}$.
\label{Dijkstra:ensure}
\State \textbf{Variables}:
 $L$ a list of pairs $(t,x) \in \bR \times X$ sorted lexicographically.
\State \textbf{Initialization}:
$L \gets \{(\Gamma u(x),x) \mid x \in X,\ \Gamma u(x) < u(x)\}$. \hfill (Removed in \cref{algo:anb_lazy}.)
\State \textbf{While} $L$ contains $(t,x)$ such that $t<T$ \textbf{do}: 
\State \phantom{bla} \label{Dijkstra:accept}
Remove from $L$ the pair $(t,x) \in \bR \times X$ which is lexicographically minimal.
\State \phantom{bla} $u(x) \gets t$ \hfill (Dijkstra update at $x$)
\State \phantom{bla} \textbf{For all} $y \in \widetilde \cV(x)$ \textbf{do}:
\State \phantom{blabla} Compute $s := \Gamma u(y)$.
\State \phantom{blabla} \textbf{If} $s < u(y)$ \textbf{then}:
\State \phantom{blablabla} Insert $(s,y)$ in $L$. If $L$ contains an outdated key $(s',y) \in L$, then remove it.
\end{algorithmic}
\end{algorithm}

\begin{proposition}[Correctness of generalized Dijkstra \leancheck{prop:Dijkstra_correct}]
\label{prop:Dijkstra_correct}
Let $\Gamma : \bR^X \to \bR^X$ be monotone and causal, extended to $\bR_\infty^X$ as in \eqref{eqdef:monotone_extension}. Consider a stopping time $T \in \bR_\infty$, and let $u \in \bR^X_\infty$ satisfy $(\Gamma u)^{<T} \leq u$.
Denote by $L_n, u_n$ the state of $L, u$ before executing line \ref{Dijkstra:accept} of the $n$-th iteration, and by $(t_n,x_n)$ the pair removed from $L$ in that line.
Then %
\begin{align}
\label{eq:Dijkstra_prop}
	\Gamma u_n &\leq u_n, &
	L_n &= \{(\Gamma u_n(x),x) \mid x  \in X,\ \Gamma u_n(x) < u_n(x)\},
\end{align}
\begin{align}
\label{eq:Dijkstra_prop2}
	u_{n+1} &\leq u_n, &
	\min\{t_n, u_{n+1}\} &= \min \{t_n, u_n\}, &
	t_{n+1} &\geq t_n.
\end{align}
The points $(x_n)_{n \geq 0}$ are pairwise distinct, thus \cref{algo:Dijkstra} terminates, at some iteration $N$.
The map $u_* := u_N$ is a fixed point of $\Gamma^{<T}$, and $\{x_0,\cdots,x_{N-1}\} = \{x\in X \mid u_*(x) < u_0(x)\}$.
\end{proposition}

\begin{proof}
Dijkstra's algorithm applied to $\Gamma$ with stopping time $T$, is clearly equivalent to the same algorithm applied to $\Gamma^{<T}$, except that the elements $\{ (\Gamma u(x), x) \mid T \leq \Gamma u(x) < u(x)\}$ are removed from $L$. We thus assume $T=\infty$, i.e.\ no stopping time, w.l.o.g.\ up to considering $\Gamma^{<T}$.

Obviously, \eqref{eq:Dijkstra_prop2} is only considered if the $n$-th iteration is not the last one.
Let us observe that (\ref{eq:Dijkstra_prop}, $n=0$) holds at initialization. Here and below, by (\ref{eq:Dijkstra_prop}, $n=0$) (resp.\ (\ref{eq:Dijkstra_prop}, i, $n=0$)) we refer to the statements (resp.\ the first statement) of \eqref{eq:Dijkstra_prop} specialized with $n=0$.

By \cref{prop:partialMonotone} with $B_0 := \{x_0\}$, one has $\Gamma u_1 \leq u_1$, proving (\ref{eq:Dijkstra_prop}, i, $n=1$).
We have $u_1 \leq u_0$ by construction, with equality except at $x_0$, thus $\Gamma u_1 \leq \Gamma u_0$ by monotony, with equality except possibly on $\widetilde \cV(x_0)$, which implies (\ref{eq:Dijkstra_prop}, ii, $n=1$).
One has $\min \{t_0, u_0\} = \min \{t_0, u_1\}$ by construction since $t_0 = u_1(x_0) = \Gamma u_0(x_0) < u_0(x_0)$,
hence $\min \{t_0, \Gamma u_0\} = \min \{t_0, \Gamma u_1\}$
by the causality property, which implies $t_1 \geq t_0$, and concludes the proof of (\ref{eq:Dijkstra_prop2}, $n=0$). By an immediate induction argument, \eqref{eq:Dijkstra_prop} and \eqref{eq:Dijkstra_prop2} hold in all iterations.

If by contradiction $x_n = x_m$ for some $n<m$, then denoting $x := x_m = x_n$ we obtain
\begin{equation*}
	t_n = \Gamma u_n(x) = u_{n+1}(x) \geq u_m(x) > \Gamma u_m(x) = t_m,
\end{equation*}
contradicting the non-decreasing property of $(t_n)_{n \geq 0}$. Thus the points $(x_n)$ are pairwise distinct, and therefore \cref{algo:Dijkstra} terminates in a finite number $N \leq \# X$ of iterations. It follows that $L_N = \emptyset$, hence $u_N$ is a fixed point of $\Gamma$ by (\ref{eq:Dijkstra_prop}, $n=N$). We conclude noting that $m-n = \# \{x \in X \mid u_m (x) < u_n(x)\}$.
\end{proof}

We next state a version of the comparison principle \eqref{eq:comparison_Lambda} for causal operators. 
\begin{proposition}[Comparison principle, via causality \leancheck{prop:comp_of_causal}]
\label{prop:comp_of_causal}
Let $\Gamma : \bR^X \to \bR^X$ be monotone and causal. Then for any partition $X = \Int Y \sqcup \partial Y$ and any $u, v \in \bR_\infty^X$
\begin{equation*}
	(u \leq \Gamma u \text{ and } \Gamma v \leq v \text{ on } \Int Y, \text{ and } u \leq v \text{ on } \partial Y)\ \Rightarrow \ u \leq v,
\end{equation*}
abusively denoting by $\Gamma$ the monotone extension $\Gamma_\infty$ for readability, since they coincide on finite maps by \cref{prop:monotonic_Phi}.
Given a super-solution $\ul u \in \bR_\infty^X$ as initial guess (for instance $\ul u = \infty$), Dijkstra's algorithm thus computes the \emph{unique} fixed point $u_*$ of $\Gamma$, using exactly $\#\{u_* < \ul u\} \leq \# X$ iterations. In addition, $u_* \leq \ul u$, hence $u_*$ takes finite values if $\ul u$ does.
\end{proposition}

\begin{proof}
	By contradiction, consider $x \in X$ such that $u(x) > v(x)$ and $t := v(x)$ is minimal.
	Then $\min\{t,u\} \leq \min \{t,v\}$, hence $\Gamma \min \{t,u\} \leq \Gamma \min \{t,v\}$ by monotony, thus $\min \{t+\delta, \Gamma u\} \leq \min \{t+\delta, \Gamma v\}$ for some $\delta>0$ by causality (using that $\min \{t+\delta, \Gamma u\} = \min\{t+\delta, \Gamma \min\{t,u\}\}$ and likewise $v$), and therefore $\min \{t+\delta, u\} \leq \min\{t+\delta, v\}$, contradicting $u(x) > v(x) = t$. 
	
The statements regarding Dijkstra's method immediately follow from \cref{prop:Dijkstra_correct}.
\end{proof}

The maximum of two operators is defined as $\max\{\Gamma,\Lambda\} u(x)  := \max \{\Gamma u(x), \Lambda u(x)\}$ for all $u \in \bR^X$, $x \in X$, and likewise the minimum.

\begin{lemma}[\leancheck{lem:causal_mix}]
\label{lem:causal_mix}
	If two operators $\Gamma, \Lambda : \bR^X \to \bR^X$ are monotone (resp.\ causal, resp.\ $\rho$-sub-additive, resp.\ $\nu$-sub-multiplicative), then so are $\max \{\Gamma, \Lambda\}$ and $\min\{\Gamma, \Lambda\}$.
\end{lemma}
\begin{proof}
The announced monotonicity follows from the monotonicity of the $\min$ and $\max$ functions.
Causality follows from the distributivity $\min \{t,\max \{a,b\}\} = \max \{\min \{t,a\}, \min \{t,b\}\}$ and  $\min \{t,\min \{a,b\}\} = \min \{\min \{t,a\}, \min \{t,b\}\}$ of $\min$ over $\max$ and over itself, where $t, a, b \in \bR$.
The proof of $\rho$-sub-additivity (resp.\ $\nu$-sub-multiplicativity) is trivial noting that $\min \{a+t,b+t\} = \min\{a,b\}+t$, and $\min\{t a, t b\} = t \min \{a,b\}$ for all $t \geq 0$, and likewise for the $\max$ function.
\end{proof}

\paragraph{Reformulations of causality and $C_\past$-memory.}
The axiomatizations of causality \eqref{eqdef:causality} and $C_\past$-memory \eqref{eqdef:Cpast_memory} considered in this paper are chosen for their simplicity and symmetry. The following results show that they are equivalent with previous formulations from \cite{mirebeau2019riemannian,mansour2024solving}, and illustrate some immediate consequences.
\begin{proposition}[Reformulation of causality \leancheck{prop:causal_formulation}]
\label{prop:causal_formulation}
For any  $u,v \in \bR^X$ and any $t \in \bR$, one has
\begin{align*}
	\min\{u,t\} &= \min\{v,t\} \Leftrightarrow u^{<t} = v^{<t}, &
	(\exists \delta>0,\, \min\{u,t+\delta\} = \min\{v,t+\delta\}) \Leftrightarrow u^{\leq t} = v^{\leq t}
\end{align*}
Therefore, the causality property \eqref{eqdef:causality} is equivalent to: 
$\forall u,v \in \bR^X$, $\forall t\in \bR$,
\begin{equation}
\label{eq:causal_old}
	u^{<t} = v^{<t} \Rightarrow (\Gamma u)^{\leq t} = (\Gamma v)^{\leq t}.
\end{equation}
If $\Gamma$ is monotone (resp.\ causal), then so is $\Gamma^{<T} : u \mapsto (\Gamma u)^{<T}$, for any $T \in \bR$. 
\end{proposition}

\begin{proof}
	The causality of $\Gamma^{<T}$ follows from the commutation property $(a^{\leq t})^{<T} = (a^{<T})^{\leq t}$, where $a,t,T \in \bR$. The rest of the proof is omitted since obvious.
\end{proof}

The last point was implicitly used in \cref{prop:Dijkstra_correct}.
Letting $v := \min \{u,t\}$ in \eqref{eq:causal_old} we get: %
\begin{corollary}[\leancheck{corol:G_min_u_t}]
\label{corol:G_min_u_t}
	Let $\Gamma : \bR^X \to \bR^X$ be causal. Then for all $u\in \bR^X$ and all $x \in X$,
	\begin{align*}
		\Gamma u(x) &= \Gamma \min \{u,t\}(x), &
		\text{for all } t &\geq \Gamma u(x).
	\end{align*}
\end{corollary}

The $C_\past$ memory property is presented in a similar way \eqref{eq:Lambda_max_u_t} in \cite{mansour2024solving}.
\begin{proposition}[\leancheck{prop:Cpast_formulation}]
\label{prop:Cpast_formulation}
Let $\Lambda : \bR^X \to \bR^X$ be monotone and additively invariant.
The $C_\past$-memory property is equivalent to: for all $u \in \bR^X$ and all $x \in X$
\begin{align}
\label{eq:Lambda_max_u_t}
	\Lambda u(x) &= \Lambda \max \{u, t\}(x), &
	\text{for all } t &\leq \Lambda u(x) - C_\past.
\end{align}
\end{proposition}

\begin{proof}
\emph{Proof of \eqref{eq:Lambda_max_u_t} assuming $C_\past$-memory}. 
We first consider the case $t<\Lambda u(x) - C_\past$, and define $v := \max \{u,t\}$ which satisfies $\max \{u,t\} = \max \{v,t\}$, hence $\max \{\Lambda u,t+C_\past\} = \max \{\Lambda v,t+C_\past\}$, and therefore $\Lambda u (x) = \Lambda v(x)$, as announced. We now turn to the case $t = \Lambda u(x)-C_\past$, and note that (\ref{eq:Lambda_max_u_t}, left) holds for all $x\in \partial X$ in view of the Dirichlet boundary conditions, whereas for all $x \in \Int X$ and all $\ve>0$ one has 
\begin{equation*}
	\Lambda \max \{u,t\}(x) \geq \Lambda u(x) = \Lambda \max \{u,t-\ve\}(x) \geq \Lambda [\max \{u,t\}-\ve](x) = \Lambda \max \{u,t\}(x) - \ve
\end{equation*}
using successively (i) monotony, (ii) the previously established case since $t-\ve < \Lambda u(x) - C_\past$, (iii) monotony again, and (iv) additive invariance. Letting $\ve\to 0^+$ we obtain \eqref{eq:Lambda_max_u_t} as announced.

\emph{Proof of $C_\past$-memory assuming  \eqref{eq:Lambda_max_u_t}.} Consider $u,v\in \bR^X$ and $t \in \bR$ such that $\max \{u,t\} = \max \{v,t\}$. Let also $x \in \Int X$, since the announced result is obvious on $\partial X$ in view of the Dirichlet boundary conditions. Our objective is to establish that $\max \{\Lambda u(x),t+C_\past\} = \max \{\Lambda v(x),t+C_\past\}$. We assume that $\Lambda u(x) \geq \Lambda v(x)$, w.l.o.g.\ up to exchanging the roles of $u$ and $v$, and we also assume that $\Lambda u(x) > t+C_\past$, otherwise there is nothing to prove.
Denoting $\ve := \Lambda u(x) - (t+C_\past) > 0$, and $Y := \{u>t\}$, we obtain
\begin{equation}
\label{eq:Lambda_eps_Y_geq}
	\Lambda [u+\ve\bbone_Y] (x)= \Lambda \max \{u+\ve \bbone_Y,t+\ve\}(x) \geq \Lambda [u+\ve](x) = \Lambda u(x) + \ve,
\end{equation}
using successively (i) \eqref{eq:Lambda_max_u_t} along with the observation that $t+\ve+C_\past = \Lambda u(x) \leq \Lambda [u+\ve\bbone_Y](x)$ by monotony, (ii) monotony along with the observation that $u(y) \leq t$ for all $y \notin Y$, (iii) additive invariance. 
Note that equality actually holds in \eqref{eq:Lambda_eps_Y_geq}, by monotonicity and since $u+\ve \bbone_Y \leq u+\ve$.
Therefore, denoting by $Y^c  := \{u\leq t\}$ the complement of $Y$, one obtains
\begin{equation*}
	\Lambda [u-\ve \bbone_{Y^c}](x) = \Lambda [u+\ve\bbone_Y - \ve] (x) = \Lambda [u+\ve\bbone_Y] (x) - \ve  = \Lambda u(x)
\end{equation*}
using (i) $\bbone_Y + \bbone_{Y^c} = 1$, (ii) additive invariance, and (iii) our previous estimate.
Thus 
\begin{align*}
	\forall n \in \bN,\, \Lambda u (x) &= \Lambda [u- n \ve \bbone_{Y^c}](x).&
	\text{Hence } 
	\Lambda u (x) &  \leq \Lambda v(x).
\end{align*}
We used (left) an immediate induction argument, and (right) monotonicity and an integer $n\geq 0$ large enough such that $u- n \ve \bbone_{Y^c} \leq v$, which is possible since $\max \{u,t\} = \max \{v,t\}$. Recalling that $\Lambda u(x) \geq \Lambda v(x)$, we obtain $\Lambda u(x) = \Lambda v(x)$, as desired. 
\end{proof}

\begin{corollary}[\leancheck{corol:mon_min_max}]
\label{corol:mon_min_max}
Let $\Lambda: \bR^X \to \bR^X$ be monotone, causal, additively invariant, with $C_\past$-memory.
Let $u \in \bR^X$, $x \in X$, $t_0 \leq \Lambda u(x)-C_\past$ and $t_1 \geq \Lambda u(x)$.
Let $v\in \bR^X$ be such that $v \geq u$ on $\{t_0 \leq u < t_1\}$, and $v \geq t_1$ on $\{u\geq t_1\}$. Then $\Lambda u(x) \leq \Lambda v(x)$.
\end{corollary}

\begin{proof}
Denote $\tilde u := \min \{t_1,u\}$, $\tilde v := \min \{t_1,v\}$, and choose $t'_0< t_0$ such that $\{t'_0 \leq u\} = \{t_0\leq u\}$, i.e.\ $t'_0$ is close enough to $t_0$ such that no value of  $u$ belongs to $[t'_0,t_0[$. Then, as announced
\begin{align*}
	t'_0 & +C_\past < \Lambda u(x) = \Lambda \tilde u(x) = \Lambda \max \{t'_0, \tilde u\} (x)  \leq \Lambda \max \{t'_0, \tilde v\} (x) = \Lambda \tilde v(x) \leq \Lambda v(x),
\end{align*}
using successively 
(i) the assumptions $t'_0 < t_0 \leq \Lambda u(x)-C_\past$,
(ii) causality as in \cref{corol:G_min_u_t}, 
(iii) $C_\past$-memory as in \cref{prop:Cpast_formulation}, 
(iv) monotony, 
(v) observing that $C_\past$-memory \eqref{eqdef:Cpast_memory} implies 
$\max \{t'_0+C_\past, \Lambda \max \{t'_0,\tilde v\}(x)\} = 
\max \{t'_0+C_\past, \Lambda \tilde v(x)\}$, and noting that $\Lambda \max \{t'_0,\tilde v\}(x) \geq \Lambda u(x) > t'_0+C_\past$ we obtain this equality, 
(vi) monotony.
\end{proof}

\subsection{Improvement of a sub-solution using Dijkstra's algorithm}
\label{subsec:dijkstra_under}

We study the extension $\hat u$ obtained by Dijkstra's algorithm of the subsolution $\ol u$, and the resulting narrow band $B_n$, see lines \ref{anb:hatu} and \ref{anb:Bn} of \cref{algo:anb}.
Our first objective is to establish, in \cref{prop:always_subsol}, that this construction respects the monotony properties underlying \cref{algo:anb}. Given $\Lambda, \Gamma : \bR^X \to \bR^X$, we recall that $\Gamma \leq \Lambda \Leftrightarrow (\forall u \in \bR^X,\, \forall x \in X,\, \Gamma u(x) \leq \Lambda u(x))$.

\begin{lemma}[\leancheck{lem:graph_subsol}]
\label{lem:graph_subsol}
Consider operators $\Lambda, \Gamma: \bR^X \to \bR^X$. Assume that $\Lambda$ is monotone, and that $\Gamma \leq \Lambda$.
If $\ol u,\hat u \in \bR^X$ satisfy $\ol u \leq \Lambda \ol u$ and $\hat u = \max \{\ol u, \Gamma\hat u\}$, then $\hat u \leq \Lambda \hat u$.
\end{lemma}

\begin{proof}
		One has $\ol u \leq \hat u$ by construction, hence $\Lambda \ol u \leq \Lambda \hat u$ by monotony of $\Lambda$.
		We conclude that $\hat u = \max \{\ol u, \Gamma \hat u\} \leq \max \{\Lambda \ol u, \Lambda \hat u\} = \Lambda \hat u$, using $\Gamma \leq \Lambda$ and $\ol u \leq \Lambda \ol u$ for the inequality. 
\end{proof}

\begin{remark}[\leancheck{rem:anb_npos}]
\label{rem:anb_npos}
The iteration numbers in \cref{algo:anb} are positive: $n = 1,2,...$ However, it would make no difference to start at $n=0$, since the $0$-th iteration would be vacuous: one would have $\ol u_0 = \hat u_0 = 0$, $\ul u_0 = \infty$, hence $t_0 = \Delta_0 = 0$, and therefore an empty  narrow band $B_0 = \emptyset$. (The Lean implementation of \cref{algo:anb} contains this additional $0$-th iteration.)
\end{remark}

\begin{proposition}[\leancheck{prop:always_subsol}]
\label{prop:always_subsol}
Let $\Lambda, \Gamma : \bR^X \to \bR^X$. Assume that $\Lambda$ is monotone, $\Lambda[0] \geq 0$, and $\Gamma \leq \Lambda$.
Denote by $\ol u_n, \ul u_n, \hat u_n$ the state of the variables $\ol u, \ul u, \hat u$ in line \ref{anb:hatu} of the $n$-th iteration of \cref{algo:anb}. Then
\begin{align}
\label{eq:always_subsol}
	\ol u_n &\leq \Lambda \ol u_n, &
	\Lambda \ul u_n &\leq \ul u_n, &
	\hat u_n &\leq \Lambda \hat u_n, &
	\ol u_n &\leq \hat u_n \leq \ol u_{n+1}, &
	\ul u_{n+1} &\leq \ul u_n.
\end{align}
\end{proposition}

\begin{proof}
	Note that $n\geq 1$, see \cref{rem:anb_npos}, and that the properties regarding $\ol u_{n+1}$ and $\ul u_{n+1}$ are of course only considered if the $n$-th iteration is not the last one.
	At initialization, one has $\ol u_1 = 0$ and $\ul u_1 := \infty$ by definition, hence $\ol u_1 \leq \Lambda \ol u_1$ by assumption and $\Lambda \ul u_1 \leq \ul u_1$ by construction, which yields (\ref{eq:always_subsol}, i and ii, $n=1$).
	Applying \cref{lem:graph_subsol} to the operator $\Gamma' := \min \{n \tau, \Gamma\}$, and recalling that $\hat u_1 = \max \{\ol u_1, \Gamma' \hat u_1\}$ by definition, we find that $\hat u_1 \leq \Lambda \hat u_1$, and obviously $\ol u_1 \leq \hat u_1$, which yields  (\ref{eq:always_subsol}, iii, $n=1$) and (\ref{eq:always_subsol}, iv, l.h.s., $n=1$).
	By \cref{prop:partialMonotone}, the partial updates of lines \ref{anb:while_sub} and \ref{anb:while_sup} in \cref{algo:anb} preserve the sub- and the super-solution property, hence $\ol u_2 \leq \Lambda \ol u_2$, $\Lambda \ul u_2 \leq \ul u_2$ and $\hat u_1 \leq \ol u_2$, $\ul u_2 \leq \ul u_1$, which yields (\ref{eq:always_subsol}, iv, r.h.s., $n=1$), (\ref{eq:always_subsol}, v, $n=1$) and (\ref{eq:always_subsol}, i and ii, $n=2$). The result immediately follows by induction.
\end{proof}

Note that the assumption $\Lambda[0] \geq 0$ of \cref{prop:always_subsol} follows from \cref{assum:all}, by \cref{prop:cst_subsol} (case \ref{enum:Lambda0_delta}), underestimation $\Gamma \leq \Lambda$, and the positivity of the boundary conditions.

We next describe the effective computation of $\hat u$, which is based on Dijkstra's algorithm, with a well-chosen initial guess \eqref{eqdef:Dijkstra_supsol}. The latter needs to be large enough to be a super-solution of the adequate operator $\Gamma' := \max \{\ol u_n, \min \{n\tau,\Gamma\}\}$, but not too large in such way that Dijkstra's algorithm only visits a thin band of points \eqref{eqdef:Bpn}. 
	Given $\ol u\in \bR^X$ and $T \in \bR$, we denote by $\Gamma' := \max \{\ol u, \min \{T, \Gamma\}\} : \bR^X \to \bR^X$, the operator defined by $\Gamma' u(x) := \max \{\ol u(x), \min \{T, \Gamma u(x)\}\}$ for all $u \in \bR^X$, $x \in X$.

\begin{lemma}[\leancheck{lem:max_min_mon_causal}]
\label{lem:max_min_mon_causal}
	If $\Gamma : \bR^X \to \bR^X$ is monotone and causal, and $\ol u \in \bR^X$, $T \in \bR$, then $\Gamma' := \max \{\ol u, \min \{T, \Gamma\}\}$ is monotone and causal.
\end{lemma}

\begin{proof}
	Monotonicity and causality follow from \cref{lem:causal_mix}, and from the fact that the \emph{constant} operators $u \mapsto \ol u$ and $u \mapsto T$ are monotone and causal.
\end{proof}

\begin{proposition}[\leancheck{prop:Dijkstra_band}]
\label{prop:Dijkstra_band}
Assume that $\Gamma:\bR^X \to \bR^X$ is monotone and causal, and use the notations $\ol u_n, \hat u_n$ and the assumptions of \cref{prop:always_subsol}.
Given $n\geq 1$, define the mapping $\hat u^\init_n \in \bR^X$ by
\begin{align}
\label{eqdef:Dijkstra_supsol}
	\hat u_n^\init(x) &:= \begin{cases}
		\hat u_{n-1}(x) & \text{if } \hat u_{n-1}(x) < t_{n-1}^-, \\ 
		\max\{ \ol u_{n} (x), n\tau\} & \text{otherwise},
	\end{cases} &
	\text{where } t_n^- &:= t_n - \Delta_n = n\tau-2 \Delta_n,
\end{align}
and by convention $\hat u_0 := 0$ and $t_0^- := 0$.
Then $\Gamma' \hat u_n^\init \leq \hat u_n^\init$, where $\Gamma' := \max \{\ol u_n, \min \{n\tau,\Gamma\}\}$.
Applying Dijkstra's \cref{algo:Dijkstra} to $\Gamma'$ and to the initial guess $\hat u_n^\init$, we obtain $\hat u_n$.
All the points updated in the course of \cref{algo:Dijkstra} belong to
\begin{equation}
\label{eqdef:Bpn}
	B'_n := \{x \in X \mid t^-_{n-1} \leq \hat u_{n-1}(x) \text{ and } \hat u_n(x) < n\tau\}.
\end{equation}
\end{proposition}

\begin{proof}
We denote $u := \hat u_n^\init$ for readability.
Fix $n\geq 1$, and first observe that
\begin{align}
\label{eq:u_nm1_eq_max}
\hat u_{n-1} &= \max \{\ol u_{n-1}, \min \{(n-1) \tau, \Gamma \hat u_{n-1}\}\}, \\
\label{eq:min_tnm1_hat_sub_u}
\min \{t_{n-1}^-, \hat u_{n-1}\} &= \min \{t_{n-1}^-, \ol u_n\} = \min \{t_{n-1}^-, u\}, \\
\label{eq:min_t_nm1_Gamma}
\min \{t_{n-1}^-, \Gamma \hat u_{n-1}\} &= \min\{t_{n-1}^-, \Gamma u\}
\end{align}
using successively (i) the definition of $\hat u_{n-1}$, (ii) the fact that the values updated in line \ref{anb:while_sub} of \cref{algo:anb} belong to $B_{n-1}$ hence are larger than $t_{n-1}^-$, and (iii) the causality property.

\emph{Proof that $\Gamma' u  \leq u$}. Let $x\in X$ be arbitrary.
In the \emph{otherwise} case of \eqref{eqdef:Dijkstra_supsol}, one has by construction $\Gamma' u(x) = \max \{\ol u_n(x), \min \{n\tau, \Gamma u(x)\}\} \leq \max \{\ol u_n(x), n\tau\} = u(x)$ as desired.

We may thus assume that $\hat u_{n-1}(x) < t_{n-1}^-$. This implies $\hat u_{n-1}(x) = \ol u_n(x) = u(x)$ by \eqref{eq:min_tnm1_hat_sub_u}, and 
$\hat u_{n-1}(x) = \max \{\ol u_{n-1}(x), \Gamma \hat u_{n-1}(x)\}$ by \eqref{eq:u_nm1_eq_max} since $t^-_{n-1} \leq (n-1) \tau$.
It then follows that $\Gamma \hat u_{n-1}(x) \leq \hat u_{n-1}(x) < t_{n-1}^-\leq (n-1) \tau$, thus $\Gamma \hat u_{n-1}(x) = \Gamma u(x)$ in view of \eqref{eq:min_t_nm1_Gamma}.
We conclude as desired that 
\begin{equation*}
u(x) 
= \hat u_{n-1}(x) 
= \max \{\ol u_{n-1}(x), \Gamma \hat u_{n-1}(x)\}
= \max \{\ol u_n(x), \Gamma u(x)\} 
= \Gamma' u(x),
\end{equation*}
using successively (i) the definition \eqref{eqdef:Dijkstra_supsol} of $u$, (ii, iii) previously observed identities, and (iv) $\Gamma u(x) = \Gamma \hat u_{n-1}(x) \leq (n-1) \tau \leq n\tau$. 
The super-solution property follows.

\emph{Proof that the updated points belong to $B'_n$.}
\Cref{prop:Dijkstra_correct} shows that the fixed point $\hat u_n = \Gamma' \hat u_n$ is indeed computed by Dijkstra's algorithm, since $\Gamma'$ is monotone and causal by \cref{lem:max_min_mon_causal}, and since the initial guess satisfies $\Gamma' u \leq u$ by the above.
A point $x\in X$ is updated in the course of Dijkstra's algorithm iff $\hat u_n(x) < u(x)$, see \cref{prop:Dijkstra_correct}.
This implies $\hat u_n(x)<n\tau$, since otherwise $\hat u_n(x) = u(x) = \max\{\ol u_n(x), n\tau\}$.
This also implies $t_{n-1}^- \leq \hat u_{n-1}(x)$, since otherwise one would have $u(x) = \hat u_{n-1}(x)$ by \eqref{eqdef:Dijkstra_supsol}, and therefore $\hat u_n(x) < u(x) = \hat u_{n-1}(x)$ which contradicts (\ref{eq:always_subsol}, iv). The result follows.
\end{proof}

The objective of the rest of this subsection is to bound above and below the values of the exact fixed point $u_*$ of $\Lambda$ on the narrow bands $B_n$ and $B_n'$, see \cref{corol:Bn_comp}. For that purpose we first describe the effect of the saturation threshold $n\tau$ in the defining equation of $\hat u_n$, see \cref{lem:max_min}, and we compare the growth of $u_*$ and of a sub-solution of $\Gamma$, see \cref{prop:olu_ge_ustar}.

\begin{lemma}[\leancheck{lem:max_min}]
\label{lem:max_min}
	Let $\Gamma : \bR^X \to \bR^X$ be monotone and causal, and $\ol u \in \bR^X$, $T \in \bR$.
If $u \in \bR^X$ satisfies $u = \max \{\ol u, \Gamma u\}$, then $\hat u := \max \{\ol u, \min \{T, u\}\}$ is the unique fixed point of $\Gamma' := \max \{\ol u, \min \{T, \Gamma\}\}$.
In addition, $\min \{u,T\} = \min \{\hat u, T\}$, hence $\{u < t\} = \{\hat u < t\}$, $\forall t\leq T$.
\end{lemma}

\begin{proof}
\emph{Proof of the additional property.}
One has $\min \{T,\hat u\} = \max \{\min \{T, \ol u\},\min \{T,u\}\} = \min \{T,u\}$ by distributivity of $\min$ and since $u \geq \ol u$. 
\emph{Proof of the first statement.}
One has 
\begin{equation}
	\label{eq:min_T_Gammap_hatu_eq_min_T_u}
	\min\{T, \Gamma' \hat u\} = \min \{T, \max \{\ol u, \Gamma \hat u\}\} = \min \{T, \max \{\ol u,  \Gamma u\}\} = \min \{T,u\}
\end{equation}
using successively (i) the expression of $\Gamma'$, (ii) the equality $\min \{T, \Gamma \hat u\} = \min \{T,\Gamma u\}$ which holds by causality since $\min\{T,\hat u\} = \min \{T,u\}$, (iii) the assumption on $u$. 
	On the set $Y := \{u<T\}$, one thus has $\Gamma'\hat u = u$ by \eqref{eq:min_T_Gammap_hatu_eq_min_T_u}, and in addition $\hat u = \max \{\ol u, \min \{T,u\}\} = \max \{\ol u, u\} = u$ since $\ol u \leq u < T$, thus $\Gamma'\hat u = \hat u$ as desired.
	On the other hand, on the complement $X \sm Y = \{u \geq T\}$, 
\begin{equation*}
	\hat u 
	:= \max \{\ol u, \min \{T,u\}\} 
	= \max \{\ol u, T\}
	= \Gamma' \hat u
\end{equation*}
using successively (i) the definition of $\hat u$, (ii) this case's assumption $u\geq T$, (iii) the definition of $\Gamma'$ and the observation that $\Gamma' \hat u \geq T$ by \eqref{eq:min_T_Gammap_hatu_eq_min_T_u}. This shows that $\hat u$ is a fixed point of $\Gamma'$, and the uniqueness follows by the comparison principle \cref{prop:comp_of_causal}, since this operator is monotone and causal by \cref{lem:max_min_mon_causal}.
\end{proof}

\begin{proposition}[\leancheck{prop:olu_ge_ustar}]
\label{prop:olu_ge_ustar}
Let $\Lambda : \bR^X \to \bR^X$ be monotone and additively invariant, and let $\Gamma : \bR^X \to \bR^X$ obey the $(C_\Gamma, C'_\Gamma)$-comparison property for some $C_\Gamma \geq 1$, $C'_\Gamma \geq 0$.
Assume furthermore that $\Lambda$ has $C_\past$-memory, with $C_\past \geq 0$, satisfies the comparison principle \eqref{eq:comparison_Lambda},
and that $u_* \in \bR^X$ satisfies $u_* \leq \Lambda u_*$.
Let $u\in \bR^X$, $t\geq 0$ and $\ve>0$ obey
\begin{itemize}
	\item $(1+\ve) u \geq u_*$ on $\partial Y := \{u < t\}$.
	\item $\Gamma u \leq u \leq \Lambda u$.
\end{itemize}
Then $u \geq t + (u_*-(1+\ve) t-C_\past -C'_\Gamma)/C_\Gamma$ on $\Int Y := X \sm \partial Y = \{u \geq t\}$.
\end{proposition}

\begin{proof}
	Let us recall the implicit additional assumption that $\Lambda$ and $\Gamma$ share the same Dirichlet boundary data $u_\bd$ on $\partial X$.
	Define $v,w \in \bR^X$ by 
	\begin{align}
	\label{eqdef:vw_Cpast}
		v &:= C_\past + C'_\Gamma + (1+\ve)t+ C_\Gamma (u-t) = T + C_\Gamma u, &
		w &:= \max \{v,(1+\ve) t\},
	\end{align}
	where $T := C_\past + C_\Gamma' + (1+\ve)t - C_\Gamma t$.
	On the subdomain $\Int Y$, one has $u \geq t$, hence $v\geq C_\past + C_\Gamma'+(1+\ve) t$, and therefore $v=w$ since $C_\past, C'_\Gamma\geq 0$.
	In addition, on the subdomain $\Int Y \cap \Int X$
	\begin{equation}
	\label{eq:get_Cpast}
	\Lambda v
	= T + \Lambda [C_\Gamma u]
	\geq T + C_\Gamma \Lambda u - C'_\Gamma
	\geq T + C_\Gamma u- C'_\Gamma = v - C'_\Gamma \geq (1+\ve) t+C_\past,
	\end{equation}
	using successively (i) the additive invariance of $\Lambda$ on $\Int X$, (ii) the first part of the $(C_\Gamma,C'_\Gamma)$-comparison assumption: $C_\Gamma\Lambda u - C'_\Gamma\leq \Lambda [C_\Gamma u]$ on $\Int X$, (iii) the assumption $u \leq \Lambda u$, (iv) the defining expression \eqref{eqdef:vw_Cpast} of $v$, and (v) the fact that $u\geq t$ on $\Int Y$.
	
	Therefore, again on the subdomain $\Int Y \cap \Int X$
	\begin{equation}
	\label{eq:supsol_Cpast}
		\Lambda \max \{v,(1+\ve) t\} =
		\Lambda v
		 = \Lambda [T + C_\Gamma u]
		 = T + \Lambda [C_\Gamma u]
		 \leq T+ C_\Gamma \Gamma u
		 \leq T + C_\Gamma u = v
		 = w, %
	\end{equation}
	using successively (i) the $C_\past$-memory property with cutoff time $t' := (1+\ve) t \leq \Lambda v(x) - C_\past$ by \eqref{eq:get_Cpast}, (ii) the defining expression \eqref{eqdef:vw_Cpast} of $v$, (iii) the additive invariance of $\Lambda$ on $\Int X$, (iv) the second part of the $(C_\Gamma,C'_\Gamma)$-comparison assumption: $\Lambda [C_\Gamma u] \leq C_\Gamma \Gamma u$ on $\Int X$, (v) the assumption $\Gamma u \leq u$, (vi) the definition \eqref{eqdef:vw_Cpast} of $v$ again, and (vii) $v=w$ on $\Int Y$ as already observed.

	We have shown that $\Lambda w \leq w$ on $\Int Y \cap \Int X$ by \eqref{eq:supsol_Cpast}, and on the other hand 
	$\Lambda w = u_\bd = \Gamma u \leq u \leq v \leq w$ 
	on $\Int Y \cap \partial X$, using successively 
	(i and ii) the Dirichlet boundary conditions on $\partial X$, always implicitly assumed, 
	(iii) the assumption, 
	(iv) $v \geq t + C_\Gamma (u-t)$ by (\ref{eqdef:vw_Cpast}, left), and $u \geq t$ on $\Int Y$, 
	and (v) the defining expression (\ref{eqdef:vw_Cpast}, right) of $w$. 
		
	On $\partial Y$ one has $w \geq (1+\ve) t \geq (1+\ve) u \geq u_*$, using  (i) the defining expression of $w$,  (ii) the defining expression of $\partial Y$, and (iii) the assumption on $\partial Y$.
	By the comparison principle, we obtain $w \geq u_*$. Thus $v = w \geq u_*$ on $\Int Y$, which is the announced inequality.
\end{proof}

\begin{corollary}[\leancheck{corol:Bn_comp}]
\label{corol:Bn_comp}
	Let \cref{assum:all} hold, and further assume that $\Lambda$ admits a fixed point $u_* \in \bR^X$, and satisfies the comparison principle, and let $\ve,\tau>0$. Then
	\begin{align*}
		B_n &\subset \{t_n^- \leq u_* < t_n^+\},& 
		B'_n &\subset \{t_{n-1}^- \leq u_* < t_n^+\}, &
		\text{with } t_n^+ := (1+\ve) t_n + C_\past + C'_\Gamma + C_\Gamma \Delta_n,
	\end{align*}
	in the $n$-th iteration of \cref{algo:anb}. The time $t_n^-$ and band $B'_n$ are from \cref{prop:Dijkstra_band}.
\end{corollary}

\begin{proof}
Let us recall that $B_n := \{t_n^- \leq \hat u_n < n\tau\}$ and $B'_n := \{t_{n-1}^- \leq \hat u_{n-1} \text{ and } \hat u_n < n\tau\}$.

\noindent
\emph{Proof of the lower bound $u_* \geq t_n^-$ on $B_n$ (resp.\ $u_* \geq t_{n-1}^-$ on $B'_n$).}
By \cref{prop:always_subsol}, one has $\hat u_n \leq \Lambda \hat u_n$, hence $\hat u_n \leq u_*$ by the comparison principle, and likewise $\hat u_{n-1} \leq u_*$, as desired. 
	
\emph{Proof of the upper bound $u_* < t_n^+$ on $B_n$ and $B'_n$.} Denote by $u \in \bR^X$ the solution of $u = \max \{\ol u_n, \Gamma u\}$, which exists by \cref{prop:Dijkstra_correct} (indeed, the monotone and causal operator $\max \{\ol u_n, \Gamma \}$ admits the finite super-solution $u_*+t$, for sufficiently large $t$, since $\max \{\ol u_n, \Gamma (u_*+t)\} = \max \{\ol u_n, \Gamma u_* + t\} = \Gamma u_*+t \leq \Lambda u_* + t = u_*+t$ on $\Int X$, and $\max \{\ol u_n, u_\bd\} \leq u_*+t$ on $\partial X$).
One has $\Gamma u \leq u$ by construction, and $u \leq \Lambda u$ by \cref{lem:graph_subsol}. 
In addition $\hat u_n = \max \{\ol u_n, \min \{n\tau, u\}\} \leq \max \{\ol u_n, u\} = u$ and $\{\hat u_n < t\} = \{u < t\}$ for any $t \leq n\tau$, by \cref{lem:max_min}.
In particular, $(1+\ve) u \geq (1+\ve) \hat u_n \geq \ul u_n \geq u_*$ on $\{\hat u_n < t_n\} = \{u < t_n\}$ where $t_n$ is defined in line \ref{anb:tn} of \cref{algo:anb}.
Thus $u$ satisfies the assumptions of \cref{prop:olu_ge_ustar} with $t:=t_n$, which is non-negative by construction since it is either $n\tau$ or a value of $\hat u_n \geq \hat u_0 = 0$. Hence denoting $\partial Y := \{u<t_n\}$ and $\Int Y := X\sm \partial Y$
	\begin{equation*}
		\{\hat u_n < n \tau\} \cap \Int Y
		= \{u<n\tau\} \cap \Int Y
		\subset \{t_n + (u_* - (1+\ve) t_n - C_\past - C_\Gamma')/C_\Gamma < n\tau \}
		= \{u_* < t_n^+\},
	\end{equation*}
	using successively (i) the equality $\{u<n\tau\} = \{\hat u_n < n \tau\}$ by \cref{lem:max_min}, (ii) the conclusion of \cref{prop:olu_ge_ustar}, and (iii) recalling that $t_n+\Delta_n = n\tau$.
	On $\partial Y$ one has on the other hand $u_* \leq (1+\ve) u < (1+\ve) t_n \leq t_n^+$, thus $\{\hat u_n < n \tau\} \cap \partial Y \subset \partial Y \subset \{u_* < t_n^+\}$. The result follows.
\end{proof}

\subsection{Complexity assuming that $\Phi_\alpha$ is $\rho$-sub-additive}
\label{subsec:cplx_subadd}

We conclude in this subsection the proof of our complexity result \cref{th:complexity}, which assumes that $\Phi_\alpha$ is $\rho$-sub-additive, in other words that the constants $(\alpha,\rho)$ associated to this key property of the scheme are \emph{known} and are used as the parameters of \cref{algo:anb}.
The key ingredient for that purpose is a uniform bound on the (half-)width $\Delta_n$ of the narrow band, established in \cref{corol:max_width} below, which is preceded with a technical result.

For convenience, we reformulate the target tolerance criterion using \Kruzkov\ coordinates:
\begin{align}
\label{eq:tol_as_Kruzkov}
	(1+\ve) \ol u &\geq \ul u \ \Leftrightarrow \ \ol v \exp(-\ve \ol u/\alpha) \geq \ul v, &
	\text{where } \ol v &:= \rE_\alpha \ol u,\ \ul v := \rE_\alpha \ul u,
\end{align}
and $\ol u \in \bR^X$, $\ul u \in \bR^X_\infty$. 
The assumptions of \cref{prop:t_advances} informally state that (a) the scheme residue of the sub- and super-solutions $\ol v$ and $\ul v$ is sufficiently small on a band of points denoted $\Int Y$, (b) $\ol v$ is sufficiently large on the upper boundary of the band, (c) the tolerance criterion \eqref{eq:tol_as_Kruzkov} is satisfied on the lower boundary of the band, and \eqref{eq:constr_tadvances} the band is sufficiently large. We show under these conditions, that the tolerance criterion \eqref{eq:tol_as_Kruzkov} is satisfied on a large subset of the band.
Rephrasing this result with the notations of the ANB method \cref{algo:anb}, we obtain in \cref{corol:max_width} that the time $t_n$, associated to the first violation of the constraint \eqref{eq:tol_as_Kruzkov} in the $n$-th iteration, satisfies $t_{n+1} \geq t_n + \Delta_n/2$, if the narrow band width $\Delta_n$ is sufficiently large. Since 
$t_{n+1}+\Delta_{n+1} = (n+1) \tau =  \tau + t_n+\Delta_n$ by construction, 
this implies $\Delta_{n+1} \leq \Delta_n/2+\tau$, and eventually yields the boundedness of the narrow band width $\forall n,\ \Delta_n \leq \Delta_{\max}$ as desired.

\begin{proposition}[\leancheck{prop:t_advances}]%
\label{prop:t_advances}
Let $\Lambda : \bR^X\to \bR^X$ be a monotone and additively invariant operator,   and assume that $\Phi := \rE_\alpha \Lambda_\infty \rE_\alpha^{-1}$ is $\rho$-sub-additive with $0 \leq \rho < 1$.
Let $\ol u \in \bR^X$, $\ol v := \rE_\alpha \ol u$, and $\ul v \in ]-\infty,0]^X$. Let $0 \leq T^- < T^+$, $\xi \geq 0$, $\ve>0$.
Assume that
\begin{enumerate}[(a)]
	\item $-(1-\rho) \xi \leq \ol v-\Phi \ol v\leq 0$ and $0 \leq \ul v-\Phi \ul v \leq (1-\rho)\xi$ on $\Int Y := \{T^- \leq \ol u < T^+\}$.
	\item $\ol v \geq -\xi$ on $\partial Y_1 := \{T^+ \leq \ol u\}$.
	\item $\ol v \exp(-\ve \ol u/\alpha) \geq \ul v$ on $\partial Y_2 := \{ \ol u < T^-\}$.
\end{enumerate}
Denote $t := (T^++T^-) / 2$ and $\Delta := (T^+-T^-) / 2$, and assume furthermore that for some $C \geq 1$
\begin{align}
\label{eq:constr_tadvances}
	\xi &\leq C \exp(-T^+/\alpha), &
	\ve T^+ &\leq \alpha, &
	\Delta \geq 2 \alpha \max \{1, \ln (2 e C/\ve)\}.
\end{align}
Then $\ol v \exp(-\ve \ol u/\alpha) \geq \ul v$ on
$\{t-\Delta/2 \leq \ol u \leq t+\Delta/2\}$.
\end{proposition}

\begin{proof}
Define $v := \ol v \exp(-\ve T^-/\alpha) + 2\xi$, and observe that $\ol v < 0$ and $\ul v \leq 0$. One has on $\Int Y\sm \partial X$
\begin{equation*}
	v - \Phi v
	\geq \exp(-\ve T^-/\alpha) (\ol v - \Phi \ol v) + 2 \xi (1-\rho)
	\geq (\ol v - \Phi \ol v) + 2 \xi (1-\rho)
	\geq  \xi (1-\rho)
	\geq \ul v-\Phi \ul v,
\end{equation*}
using successively (i) the $\rho$-sub-additivity and the multiplicative invariance of $\Phi$, see \cref{prop:monotonic_Phi}, and (ii, iii, iv) assumption (a) and the fact that $\exp(-\ve T^-/\alpha) \leq 1$.
On the other hand, one has on $\Int Y \cap \partial X$
\begin{equation*}
	v-\Phi v = -(1-\exp(-\ve T^-/\alpha)) \ol v + 2 \xi + \ol v-\Phi \ol v \geq 2 \xi -(1-\rho) \xi
	\geq  \xi (1-\rho)
	\geq \ul v-\Phi \ul v,
\end{equation*}
using successively (i) the expression of $v$ and the fact that $\Phi v = \Phi \ol v$ on $\partial X$, in view of the Dirichlet boundary conditions on $\partial X$, always implicitly assumed, (ii) inequalities $\ol v\leq 0$, $\exp(-\ve T^-/\alpha) \leq 1$, and assumption (a), (iii) the assumption $\rho,\xi \geq 0$, and (iv) assumption (a) again.

By (b), one has $v \geq 0 \geq \ul v$ on $\partial Y_1$.
By (c), one has
	$\ul v \leq \ol v\exp(-\ve T^-/\alpha) \leq v$ on $\partial Y_2$.
Therefore $\ul v \leq v$ on $X$, by the comparison principle, see \cref{prop:comparison_Phi}.

	We now remark that on the subdomain $\Int B := \{t-\Delta/2 \leq \ol u \leq t+\Delta/2\}$ one has
	\begin{align}
		\frac{2\xi} {-\ol v} \nonumber
		\leq 2C\exp(\frac{\ol u-T^+} \alpha)
		\leq 2 C \exp( -\frac{\Delta}{2 \alpha})
		\leq  \frac \ve e
		\leq  \frac {\ve\Delta}{2 \alpha e}
		\leq \exp(&-\frac{\ve T^-} \alpha) - \exp(-\frac{\ve (T^- +  \Delta/ 2)}\alpha) \\
		&\leq \exp(-\frac{\ve T^-} \alpha) - \exp(-\frac {\ve \ol u} \alpha).
		\label{ineq:t_advances}
	\end{align}

	We used successively (i) the upper bound (\ref{eq:constr_tadvances}, left) on $\xi$, and the assumption $\ol v = \rE_\alpha \ol u = -\exp(-\ol u/\alpha)$, (ii) $\ol u \leq t+\Delta/2$ on $\Int B$, and the observation that $T^\pm = t\pm \Delta$, (iii, iv) the lower bound (\ref{eq:constr_tadvances}, right), (v) the observation that $\exp(-a) - \exp(-b) \geq (b-a)/e$ for all $a \leq b\leq 1$, which is obtained by integration, with $a:=\ve T^-/\alpha \leq b:=\ve(T^-+\Delta/2)/\alpha \leq 1$ by (\ref{eq:constr_tadvances}, center), and (vi) $\ol u \geq t-\Delta/2 = T^-+\Delta/2$ on $\Int B$.

	The obtained inequality \eqref{ineq:t_advances} can be reformulated as $\ol v \exp(-\ve\ol u/\alpha) \geq v$.
	Combining it with the estimate $\ul v \leq v$ previously obtained by the comparison principle, we conclude the proof.
\end{proof}

In the rest of this subsection, we use the notations $t_n,\Delta_n,B_n,\eta_n$ of \cref{algo:anb}, and we denote by $\ol u_n$ and $\ul u_n$ the state of the variables $\ol u$ and $\ul u$ in the beginning of the $n$-th iteration, 
and by $\hat u_n$ the map $\hat u$ computed in line \ref{anb:hatu} of that iteration,
similarly to \cref{prop:always_subsol}. The next lemma shows that, because of the positive boundary conditions, the first iterations of \cref{algo:anb} are essentially vacuous.

\begin{lemma}[First iterations \leancheck{lem:first_iter}]
\label{lem:first_iter}
	Under \cref{assum:all}:
	in any iteration $1 \leq n \leq \min(u_\bd)/\tau$ of \cref{algo:anb}, one has 
	$\ol u_n = (n-1) \tau$, $\ul u_n = \infty$, $\hat u_n = n\tau$, $t_n = n\tau$, $\Delta_n = 0$, 
	$B_n = \emptyset$.
\end{lemma}

\begin{proof}
	One has $\ol u_1 = 0$ and $\ul u_1 = \infty$ by initialization. Let us assume, for induction, that $1 \leq n \leq \min(u_\bd)/\tau$, $\ol u_n = (n-1) \tau$ and $\ul u_n = \infty$. 
	By \cref{prop:cst_subsol} (case \ref{enum:Lambda0_delta}) one has $n\tau \leq \Gamma[n\tau]$, hence $n\tau = \max\{(n-1) \tau, \min\{n\tau, \Gamma[n\tau]\}\}$. The operator $\Gamma' := \max\{(n-1) \tau,\min \{n\tau,\Gamma\}\}$ is monotone and causal by \cref{lem:max_min_mon_causal}, hence admits a \emph{unique} fixed point by \cref{prop:comp_of_causal}, and therefore $\hat u_n = n\tau$. 	
It immediately follows that $t_n = n\tau$, hence $\Delta_n = 0$ and $B_n = \emptyset$, which means that the method jumps to the next iteration with $\ol u_{n+1} = n\tau$ and $\ul u_{n+1} = \infty$, as announced.
\end{proof}

Recall that $\lceil a \rceil \in \bZ$ denotes the smallest integer above a real number $a \in \bR$.

\begin{corollary}[Maximum narrow band width \leancheck{corol:max_width}]
\label{corol:max_width}
Let \cref{assum:all} hold, and further assume that $\Phi_\alpha := \rE_\alpha \Lambda_\infty \rE_\alpha^{-1}$ is $\rho$-sub-additive with $\alpha>0$ and $0 \leq \rho < 1$, and that
\begin{align*}
	0 &< \ve \leq \ve_{\max}, & 2\Delta_{\max} \leq \min(u_\bd).
\end{align*}
The number of iterations of \cref{algo:anb} is bounded by $n_{\max}$, and in each iteration one has $\Delta_n < \Delta_{\max}$.
We denoted:
\begin{align*}
	\Delta_{\max} &:= 2 \alpha  \ln (2 e/\ve) + 2\tau, &
	n_{\max} &:= \Big\lceil \frac{\|u_*\|_\infty + \Delta_{\max}} \tau \Big\rceil, &
	\ve_{\max} &:= \min \big\{\frac \alpha{\|u_*\|_\infty+\Delta_{\max}}, 2\big\}.
\end{align*}
\end{corollary}

\begin{proof}
	We proceed by contradiction, and denote by $n+1$ the least integer such that  $\Delta_{n+1} \geq \Delta_{\max}$. Note that $\Delta_n < \Delta_{\max}$, and that $n \geq 1$ since $\Delta_1 = 0$.
	Recall that $\hat u_n \leq \ol u_{n+1} \leq \hat u_{n+1}$ and $\ul u_{n+1} \leq \ul u_n$, see \cref{prop:always_subsol}, hence
	\begin{equation}
	\label{eq:Enp1_subset_En}
	E_{n+1} := \{(1+\ve) \hat u_{n+1} < \ul u_{n+1}\} \subset \{(1+\ve) \ol u_{n+1} < \ul u_{n+1}\} \subset E_n := \{(1+\ve) \hat u_n < \ul u_n\} \subset \{\hat u_n \geq t_n\},
	\end{equation}
	where the last inclusion holds by definition of $t_n := \min\{n \tau, \min_{E_n} \hat u_n\}$, see line \ref{anb:tn} of \cref{algo:anb}. 
	Therefore $t_{n+1} := \min\{(n+1) \tau, \min_{E_{n+1}} \hat u_{n+1}\} \geq t_n$, since $E_{n+1} \subset E_n$ and $\hat u_{n+1} \geq \hat u_n$.
	Thus
	\begin{equation*}
		\Delta_n := n \tau-t_n
		\geq n\tau-t_{n+1}
		= \Delta_{n+1} - \tau
		\geq \Delta_{\max}-\tau
		= 2 \alpha \ln (2e/\ve) + \tau
		> 0.
	\end{equation*}

	We now apply \cref{prop:t_advances}, with $C=1$,  $\ol u := \ol u_{n+1}$, $\ul v := \rE_\alpha \ul u_{n+1}$, $T^\pm := t_n \pm \Delta_n$, and $\xi := \exp(-n\tau/\alpha)$. 
	One has $n\tau > \min (u_\bd)$ by \cref{lem:first_iter} and since $\Delta_n>0$, hence $T^- = n\tau-2 \Delta_n > \min(u_\bd)-2\Delta_{\max} \geq 0$ as required. 
	Assumption (a) holds in view of the stopping criterion of lines \ref{anb:while_sub} and \ref{anb:while_sup}, and since
	$B_n := \{T^- \leq \hat u_n < T^+\} \supset \{T^- \leq \ol u_{n+1} < T^+\}$, recalling that $\hat u_n \leq \ol u_{n+1}$ and observing that the values of $\hat u_n$ smaller than $T^-$ are not modified in line \ref{anb:while_sub} of \cref{algo:anb}. Assumption (b) holds by choice of $\xi := \exp(-T^+/\alpha)$ and since $\ol v = \rE_\alpha \ol u = -\exp(-\ol u /\alpha)$.
	Assumption (c) holds noting that
	\begin{equation}
	\label{eq:olv_exp_olu_iff}
		\ol v \exp(-\ve \ol u/\alpha) \geq \ul v \ \Leftrightarrow \ (1+\ve) \ol u_{n+1} \geq \ul u_{n+1},
	\end{equation}
	and since $\{(1+\ve) \ol u_{n+1} \geq \ul u_{n+1}\} \supset \{\hat u_n < t_n\} \supset \{\ol u < T^-\}$, using successively \eqref{eq:Enp1_subset_En}, $\ol u := \ol u_{n+1} \geq \hat u_n$ and $T^- < t_n$.
	Assumption (\ref{eq:constr_tadvances}, left) holds by choice of $\xi := \exp(-T^+/\alpha)$. Assumption (\ref{eq:constr_tadvances}, center) holds by choice of $\ve_{\max}$ and since
	\begin{equation}
	\label{eq:n_le_nmax}
		T^+ = n\tau = t_n + \Delta_n < \|\hat u_n\|_\infty + \Delta_{\max} \leq \|u_*\|_\infty + \Delta_{\max},
	\end{equation}
	noting that $\hat u_n\leq u_*$ by the comparison principle. Finally, assumption (\ref{eq:constr_tadvances}, right) holds since $\ln (2 e/\ve) \geq \ln e = 1$.

	Consider $x \in X$. We have the implications and equivalences
	\begin{align}
		\label{eq:tnp1_ge_tnpD2}
		(1 &+ \ve)  \hat u_{n+1}(x) < \ul u_{n+1}(x)
		\ \Rightarrow \
		(1+\ve) \ol u_{n+1}(x) < \ul u_{n+1}(x)
		\ \Leftrightarrow \
		\ol v (x) \exp(-\ve \ol u(x)/\alpha) < \ul v(x) \\
		\label{eq:u_nm1_le_tn_or}
		& \Rightarrow \
		(\ol u_{n+1}(x) < t_n -\Delta_n/2 \text{ or } \ol u_{n+1}(x) > t_n + \Delta_n/2)
	\end{align}
	using successively (i) $\ol u_{n+1} \leq \hat u_{n+1}$, (ii) the observation \eqref{eq:olv_exp_olu_iff}, and (iii) the conclusion of \cref{prop:t_advances}.
	Recalling the inclusion \eqref{eq:Enp1_subset_En}, which implies $\ol u_{n+1}(x) \geq \hat u_n(x) \geq t_n > t_n-\Delta_n/2$, we exclude the first clause of the disjunction \eqref{eq:u_nm1_le_tn_or} and find that 
	\begin{align*}
		x\in E_{n+1}\ \Rightarrow \ t_n + \Delta_n/2 < \ol u_{n+1}(x) \leq \hat u_{n+1}(x).
	\end{align*}
	Thus $t_{n+1} := \min\{(n+1) \tau,\, \min_{E_{n+1}} \hat u_{n+1}\} \geq \min \{(n+1)\tau,\, t_n+\Delta_n/2\} = t_n+\Delta_n/2$.
	It follows that $\Delta_{\max} \leq \Delta_{n+1} = \tau + \Delta_n+t_n-t_{n+1} \leq \tau + \Delta_n - \Delta_n/2 = \Delta_n/2 + \tau < \Delta_{\max}/2 + \tau$, noting that $\Delta_n < \Delta_{\max}$ by choice of $n$ for the last inequality, hence $\Delta_{\max} < 2\tau$, which is the desired contradiction.
	Incidentally, we proved \eqref{eq:n_le_nmax} for each iteration where $E_n\neq \emptyset$; taking into account the final iteration, where $E_n = \emptyset$, yields the announced $n_{\max}$.
\end{proof}

\begin{corollary}[\leancheck{corol:NGamma}]
\label{corol:NGamma}
	Under the assumptions of \cref{corol:max_width}.
	In any iteration $n\geq 1$
\begin{align}
\label{eq:Bn_subset_Bpn_subset}
	B_n &\subset \{ n\tau -2 \Delta_{\max} \leq u_* < n \tau + \Delta^+_{\max} \}, &
	B'_n &\subset \{ (n-1) \tau - 2 \Delta_{\max} \leq u_* < n\tau + \Delta^+_{\max}\},
\end{align}
with $\Delta^+_{\max} := (C_\Gamma-1) \Delta_{\max} + \ve \|u_*\|_\infty + C_\past + C'_\Gamma$.
Therefore, for any $x \in X$,
\begin{align}
\label{eq:count_x_in_Bn}
	\# \{n\geq 1 \mid x \in B_n\} &\leq N'_\Lambda := \Big\lceil \frac{2\Delta_{\max} + \Delta^+_{\max}} \tau \Big\rceil, &
	\# \{n\geq 1 \mid x \in B'_n\} &\leq N_\Gamma := N'_\Lambda+1.
\end{align}
\end{corollary}

\begin{proof} Recalling the notations $t^\pm_n$ of \cref{corol:Bn_comp} one has
\begin{align*}
	t^-_n &:= t_n - \Delta_n = n\tau - 2\Delta_n \geq n\tau-2 \Delta_{\max}, \\
	t^+_n &:= (1+\ve) t_n + C_\past + C'_\Gamma + C_\Gamma \Delta_n =
	n \tau + \ve t_n + C_\past + C'_\Gamma + (C_\Gamma-1) \Delta_n \leq n \tau +\Delta^+_{\max},
\end{align*}
using that
$t_n = n\tau-\Delta_n$, and $\Delta_n \leq \Delta_{\max}$, and $t_n \leq \|\hat u_n\|_\infty \leq \|u_*\|_\infty$.
The inequality $t_n \leq \|\hat u_n\|_\infty$ uses that $t_n := \min \{n\tau,\min_{E_n} \hat u_n\}$ is bounded by any value of $\hat u_n$ on the set $E_n := \{(1+\ve) \hat u_n < \ul u_n\}$, provided it is non-empty. 
From this point, the announced inclusions immediately follow those established in \cref{corol:Bn_comp}.
If $x\in B_n$, then $-\Delta^+_{\max}/\tau < n - u_*(x)/\tau \leq 2 \Delta_{\max}/\tau$, and the number of such integers $n$ is at most $N'_\Lambda$, as announced. (Indeed, $\#(]a,b] \cap \bZ) \leq \lceil b-a \rceil$ for all $a \leq b$.)
Proceeding likewise for $N_\Gamma$, we conclude the proof when $E_n\neq \emptyset$.

Finally, we briefly discuss the case $E_n = \emptyset$, which in view of the stopping criterion can only happen in the very last iteration. This implies $B_n = \emptyset$, and on $B'_n$ one has $u_* \leq \ul u_n(x) \leq (1+\ve) \hat u_n(x) < n\tau + \ve \|u_*\|_\infty \leq n\tau + \Delta^+_{\max}$ as before, which concludes.
\end{proof}

The last ingredient of the proof of \cref{th:complexity} is an upper bound on the number of inner iterations, using the contraction property of the update operator $\Phi_\alpha$ in \Kruzkov\ coordinates.
The next lemma also ensures that the inner iterations always terminate, which was used implicitly.
We denote by $\lceil t\rceil_+ := \min \{ n \in \bN_{\geq 0} \mid n \geq t\}$ the smallest non-negative integer above $t$.

\begin{algorithm}
\caption{Global iteration, for computing a super-solution. \leancheck{algo:global_iteration}\newline
Variant for computing a sub-solution: replace all instances of $v - \Phi v$ below with $\Phi v-v$.
}
\label{algo:global_iteration}
\begin{algorithmic}[1]
\State 
\textbf{Input:} a target tolerance $\eta>0$, a mapping $v\in \bR^X$. 
\State
\textbf{Ensures at termination:} $0 \leq v-\Phi v\leq \eta$.
\State \textbf{While} $\max(v-\Phi v) > \eta$ \textbf{do}: \textbf{For all} $x\in X$  \textbf{do in parallel}: 
$v(x) \gets \Phi v(x)$.
\end{algorithmic}
\end{algorithm}

\begin{lemma}[\leancheck{lem:cplx_global_iteration}]
\label{lem:cplx_global_iteration}
	Let $\Phi : \bR^X \to \bR^X$ be monotone and $\rho$-sub-additive, with $0< \rho < 1$, and let $v \in \bR^X$ satisfy $0 \leq v-\Phi v \leq C$ (resp.\ $0 \leq \Phi v-v \leq C$). Then \cref{algo:global_iteration} (resp.\ the sub-solution variant) terminates using at most $\lceil \frac{\ln (C/\eta)} {-\ln \rho} \rceil_+$ iterations.
\end{lemma}

\begin{proof}
	One has $0 \leq \Phi^n v - \Phi^{n+1} v \leq C \rho^n$ (resp.\ $0 \leq \Phi^{n+1} v - \Phi^{n} v \leq C \rho^n$ for sub-solutions) for all $n\geq 0$, by an immediate induction argument, using (i) monotonicity and (ii) $\rho$-sub-additivity. The result follows. %
\end{proof}

\begin{proposition}[\leancheck{prop:niter_while}]
\label{prop:niter_while}
Under the assumptions of \cref{corol:max_width}, and $\rho>0$.
The number of iterations in the inner while loops lines \ref{anb:while_sub} and \ref{anb:while_sup} of \cref{algo:anb}, is at most
\begin{equation}
\label{eq:NLambdapp}
	N''_\Lambda := \Big \lceil \frac{2 \Delta_{\max}/\alpha - \ln (1-\rho)}{-\ln \rho} \Big \rceil
\end{equation}
\end{proposition}

\begin{proof}
	Denote by $\hat u_n$ and $\ul u_n$ the sub-and super solution at the beginning of these inner loops, in the $n$-th iteration of \cref{algo:anb}, and recall that $\hat u_n \leq \Lambda \hat u_n \leq u_* \leq \Lambda \ul u_n \leq \ul u_n$, see \cref{prop:always_subsol}, and that $B_n := \{n\tau-2\Delta_n \leq \hat u_n < n\tau\}$. Let also $\hat v_n := \rE_\alpha \hat u_n$ and $\ul v_n := \rE_\alpha \ul u_n$.
	One has $\hat u_n \geq n\tau-2 \Delta_{\max}$ on $B_n$, and therefore on this set
	\begin{align*}
		0 &\leq \Phi \hat v_n - \hat v_n \leq -\hat v_n \leq \exp(\frac {2\Delta_{\max} - n\tau}\alpha), &
		0 \leq \ul v_n - \Phi \ul v_n \leq - \Phi \ul v_n \leq - \hat v_n \leq \exp(\frac{2 \Delta_{\max}-n\tau}{\alpha}).
	\end{align*}
	We conclude by applying \cref{lem:cplx_global_iteration} to the restricted map $[\ul v_n]_{|B_n} \in \bR^{B_n}$ and the restricted \eqref{eqdef:Lambda_B} operator $[\Phi_\alpha]_{B_n}^{\ul v_n} : \bR^{B_n} \to \bR^{B_n}$, which is monotone and $\rho$-sub-additive by \cref{lem:monDom} 
	(resp.\ apply the sub-solution variant of \cref{lem:cplx_global_iteration} to $[\hat v_n]_{|B_n}$ and $[\Phi_\alpha]_{B_n}^{\hat v_n}$).
\end{proof}

\begin{proof}[Proof of \cref{th:complexity}]
The existence and uniqueness of a fixed point, and the comparison principle for $\Lambda$, are established in \cref{prop:comparison_Lambda}. The termination of \cref{algo:anb} follows from the bound $N_\Lambda''$ of \cref{prop:niter_while} on the number of inner loops, and the bound $n_{\max}$ of \cref{corol:max_width} on the number of outer loops. The number of updates $\ol u(x) \gets \Lambda \ol u(x)$ (resp.\ $\ul u(x) \gets \Lambda \ul u(x)$), at any given point $x\in X$, is bounded by the product $N_\Lambda := N'_\Lambda N_\Lambda''$ of the number $N'_\Lambda$ of iterations such that $x\in B_n$, see \eqref{eq:count_x_in_Bn}, with the number $N''_\Lambda$ of updates \eqref{eq:NLambdapp} in the inner loops. 

Finally, the update $\hat u(x)\gets\Gamma\hat u(x)$ is applied at most once to any given point $x\in X$ in the $n$-th iteration of \cref{algo:anb}, and in that case $x \in B_n'$. Therefore, at most $N_\Gamma$ such updates are performed overall by \eqref{eq:count_x_in_Bn}, which concludes. 
\end{proof}

\subsection{Complexity with unknown property constants}
\label{subsec:unknown_params}

We establish \cref{th:complexity_withUnknown}, which assumes that $\Phi_\beta$ is $\sigma$-sub-additive for some \emph{unknown} parameters $0 < \beta$, $\sigma < 1$. The proof is similar to \cref{th:complexity} in \cref{subsec:cplx_subadd}, up to some rough estimates for the difference of power functions, see \cref{lem:powr_diff,lem:pow_cvx}. 
Indeed, \Kruzkov's exponential change of coordinates $\rE_\alpha := -\exp(-\cdot/\alpha)$, where $\alpha>0$, satisfies for all $u\in \bR^X$
\begin{align*}
	w &= -(-v)^r, &
	\text{if } 
	v &:= \rE_\alpha u, \quad
	w := \rE_\beta u, \quad 
	r := \alpha/\beta.
\end{align*}

\begin{lemma}[\leancheck{lem:powr_diff}]
\label{lem:powr_diff}
	For all $a,b\geq 0$, $r\in ]0,1]$, one has $|a^r-b^r| \leq |a-b|^r$.
\end{lemma}
\begin{proof}
	We assume $a \geq b$, w.l.o.g.\ by symmetry. The announced inequality is equivalent to $(b+c)^r \leq b^r + c^r$, where $c := a-b\geq 0$. We assume $b+c = 1$, w.l.o.g.\ by homogeneity. Thus $b^r+c^r \geq b+c = 1 = (b+c)^r$, as announced, which concludes the proof.
\end{proof}

\begin{proposition}[\leancheck{prop:t_advances_withUnknown}]
\label{prop:t_advances_withUnknown}
Let $\Lambda : \bR^X \to \bR^X$ be a monotone and additively invariant operator. Assume that $\Phi_\beta := \rE_\beta \Lambda_\infty \rE_\beta^{-1}$ is $\sigma$-sub-additive, where the parameters $\alpha,\beta,\rho,\sigma$ satisfy 
\begin{align}
\label{eq:t_advances_withUnknown_params}
	\alpha, \beta & >0 , &
	\rho,\sigma &\in ]0,1[, &
	1 \geq r &:= \alpha/\beta, &
	1 \leq C &:= (1-\rho)^r /(1-\sigma).
\end{align}
Let $\ol u \in \bR^X$, $\ol v := \rE_\alpha \ol u$, and let $\ul v \in ]-\infty,0]^X$. Let $0 \leq T^- < T^+$, $\xi \geq 0$, $\ve>0$.
Assume that
\begin{enumerate}[(a)]
	\item $-(1-\rho) \xi \leq \ol v-\Phi_\alpha \ol v \leq 0$ and $0 \leq \ul v-\Phi_\alpha \ul v \leq (1-\rho)\xi$ on
	$\Int Y := \{T^- \leq \ol u < T^+\}$.
	\item $\ol v \geq -\xi$ on $\partial Y_1 := \{T^+ \leq \ol u\}$.
	\item $\ol v \exp(-\ve \ol u/\alpha) \geq \ul v$ on $\partial Y_2 := \{ \ol u < T^-\}$.
\end{enumerate}
Denote $t := \frac {T^++T^-} 2$ and $\Delta := \frac{T^+-T^-} 2$, and assume furthermore that
\begin{align}
\label{eq:constr_tadvances_withUnknown}
	\xi &\leq \exp(-T^+/\alpha), &
	\ve T^+ &\leq \alpha, & 
	\Delta \geq 2 \beta \max \{1, \ln (2 e C/\ve)\}.
\end{align}
Then $\ol v \exp(-\ve \ol u/\alpha) \geq \ul v$ on
$\{t-\Delta/2 \leq \ol u \leq t+\Delta/2\}$.
\end{proposition}

\begin{proof}
We denote $\ol w := \rE_\beta \ol u = -(-\ol v)^r$, $\ul w := -(-\ul v)^r$, and observe that pointwise on $X$
\begin{align}
\label{eq:olv_olw}
	\ol v \exp(-\ve \ol u/\alpha) &\geq \ul v \Leftrightarrow \ol w \exp(-\ve \ol u/\beta) \geq \ul w, &
	\Phi_\beta \ol w &= -(-\Phi_\alpha \ol v)^r, &
	\Phi_\beta \ul w &= -(-\Phi_\alpha \ul v)^r.
\end{align}
Thus $\ol w \leq \Phi_\beta \ol w$ and $\Phi_\beta \ul w \leq \ul w$ on $\Int Y$.
	In addition, one has by \cref{lem:powr_diff}
	\begin{align*}
		\ul w - \Phi_\beta \ul w &= (-\Phi_\alpha \ul v)^r -(-\ul v)^r \leq (\ul v-\Phi_\alpha \ul v)^r \leq  ((1-\rho) \xi)^r = (1-\sigma) \zeta,
	\end{align*}
	with $\zeta := C \xi^r \leq C \exp(-T^+/\beta)$, where $C$ is defined in (\ref{eq:t_advances_withUnknown_params}, right).
	Likewise, $\Phi_\beta \ol w - \ol w \leq (1-\sigma) \zeta$. Note also that $\ol w = -(-\ol v)^r \geq -\xi^r \geq -\zeta$ on $\partial Y_1$, since $C \geq 1$ by assumption, and that $\ve T^+ \leq \alpha \leq \beta$, since $r \leq 1$ by assumption.
	We have shown that $(\beta, \sigma,\ol u, \ol w, \ul w,T^+,T^-,\zeta,\ve,C)$ obey the assumptions of \cref{prop:t_advances}, and the result follows using (\ref{eq:olv_olw}, left).
\end{proof}

\begin{corollary}[\leancheck{corol:max_width_withUnknown}]
\label{corol:max_width_withUnknown}
Let \cref{assum:all} hold, and further assume that $\Phi_\beta$ is $\sigma$-sub-additive where $0 < \alpha \leq \beta$, and $\rho,\sigma \in ]0,1[$ satisfy $0 < (1-\sigma)^\beta \leq (1-\rho)^\alpha$. Assume furthermore that $0 < \ve \leq \ve_{\max}$ and $2 \Delta_{\max} \leq \min(u_\bd)$. The number of iterations of \cref{algo:anb} is bounded by $n_{\max}$, and in each iteration $\Delta_n < \Delta_{\max}$. We denoted
\begin{align*}
	\Delta_{\max} &:= 2 \beta  \ln (2 e C/\ve) + 2\tau,&
	n_{\max} &:= \Big\lceil \frac{\|u_*\|_\infty + \Delta_{\max}} \tau \Big\rceil,&
	\ve_{\max} &:= \min \big\{\frac \alpha {\|u_*\|_\infty + \Delta_{\max}}, 2\big\},
\end{align*}
where $C$ is from \cref{prop:t_advances_withUnknown}.
Likewise, the statement of \cref{corol:NGamma} holds with this modified value of $\Delta_{\max}$.
\end{corollary}

\begin{proof}
	The proof is identical to \cref{corol:max_width}, except that we need to use \cref{prop:t_advances_withUnknown} with $C = (1-\rho)^r / (1-\sigma)$, rather than \cref{prop:t_advances} with $C=1$.
	The statement of \cref{corol:NGamma} follows, since its proof only uses \cref{assum:all}. Note that the existence of a solution $u_*$, and the comparison principle, now come from \cref{prop:comparison_Lambda} applied to the unknown constants $\beta,\sigma$. 
\end{proof}

\begin{lemma}[\leancheck{lem:pow_cvx}] \label{lem:pow_cvx}
	For any $0 \leq a \leq b$, $s \geq 1$, one has $b^s \leq a^s + s b^{s-1} (b-a)$.
\end{lemma}
\begin{proof}
Convexity of $f(t) := t^s$, where $t \geq 0$, yields $f(a) \geq f(b) + f'(b) (a-b)$, as announced.
\end{proof}

\begin{lemma}[\leancheck{lem:cplx_global_iteration_withUnknown}] \label{lem:cplx_global_iteration_withUnknown}
Let $\Lambda : \bR^X \to \bR^X$ be monotone and such that $\Phi_\beta := \rE_\beta \Lambda_\infty \rE_\beta^{-1}$ is $\sigma$-sub-additive, where $0 < \alpha \leq \beta$ and $0<\sigma<1$.
Let also $\ol u, \ul u \in \bR_\infty^X$ be such that $\ol u \leq \Lambda \ol u$, $\Lambda \ul u \leq \ul u$ and $\rE_\alpha \ol u \geq -C$ where $C \geq 0$. Then \cref{algo:global_iteration} (global iteration) applied to $\Phi_\alpha$ and $\rE_\alpha \ol u$ (sub-solution version) or $\rE_\alpha \ul u$ (super-solution version), with tolerance $\eta>0$, uses at most 
$\lceil \frac{\ln (C/\eta) + \ln (\beta/\alpha)} {-\ln \sigma} \rceil_+$
iterations.
\end{lemma}

\begin{proof}
	Denote $\ol v := \rE_\alpha \ol u$, $\ul v := \rE_\alpha \ul u$, $\ol w := \rE_\beta \ol u$, $\ul w := \rE_\beta \ul u$, and $r := \alpha/\beta \in ]0,1]$. 
	One has 
	\begin{align}
	\label{eq:cplx_global_iteration_withUnknown}
		\ol w &= - (-\ol v)^r \geq -C^r,&
		0 &\leq \Phi_\beta \ol w - \ol w \leq -\ol w \leq C^r,&
		0 &\leq \ul w - \Phi_\beta \ul w\leq -\Phi_\beta \ul w \leq -\ol w \leq C^r.
	\end{align}
	We used the growth of $c \mapsto -(-c)^r$ in the first inequality, and to show that $\ol w$ and $\ul w$ are respectively a sub- and super-solution of $\Phi_\beta$. It follows that $\Phi_\beta \ul w$ is also a super-solution, hence $\ol w \leq \Phi_\beta \ul w$ since $\Phi_\beta$ satisfies the comparison principle by \cref{prop:comparison_Phi}, which proves \eqref{eq:cplx_global_iteration_withUnknown}.
	By an immediate induction argument, we obtain for all $n \in \bN$
	\begin{align*}
		0 &\leq \Phi^{n+1}_\beta \ol w - \Phi^n_\beta \ol w \leq C^r \sigma^n, &
		0 &\leq \Phi^{n}_\beta \ul w - \Phi^{n+1}_\beta \ul w \leq C^r \sigma^n, &
		-C^r &\leq \Phi^n_\beta \ol w \leq \Phi^n_\beta \ul w.
	\end{align*}
	Finally, denoting $s := \beta/\alpha = 1/r$, we obtain using \cref{lem:pow_cvx}
	\begin{equation*}
		0 
		\leq \Phi^{n+1}_\alpha \ol v - \Phi^n_\alpha \ol v 
		= (-\Phi^{n}_\beta \ol w)^s - (-\Phi^{n+1}_\beta \ol w)^s
		\leq s  (-\Phi^{n}_\beta \ol w)^{s-1} (\Phi^{n+1}_\beta \ol w -\Phi^{n}_\beta \ol w)
		\leq s (C^r)^{s-1} C^r \sigma^n = s C \sigma^n.
	\end{equation*}
	The iteration bound is the least integer such that $sC\sigma^n \leq \eta$. Likewise for the super-solution.
\end{proof}

\begin{proposition}[\leancheck{prop:niter_while_withUnknown}]
\label{prop:niter_while_withUnknown}
Under the assumptions of \cref{corol:max_width_withUnknown}.
The number of iterations in the inner while loops lines \ref{anb:while_sub} and \ref{anb:while_sup} of \cref{algo:anb} is at most
\begin{equation*}
	N_\Lambda'' :=
	\Big\lceil
	\frac{2 \Delta_{\max}/\alpha + \ln (\beta/\alpha)  - \ln (1-\rho)}{-\ln \sigma}
	\Big\rceil
\end{equation*}
\end{proposition}

\begin{proof}
Denote by $\hat u_n$ and $\ul u_n$ the sub- and super-solution at the beginning of these inner loops, in the $n$-th iteration of \cref{algo:anb}, and note that $\hat u_n \geq n\tau - 2 \Delta_{\max}$ on $B_n$.
We now apply \cref{lem:cplx_global_iteration_withUnknown} to the restricted maps $[\hat u_n]_{|B_n}$ and $[\ul u_n]_{|B_n}$, and the restricted operator
$\Lambda_{B_n}^{\hat u_n}$ (resp.\ $\Lambda_{B_n}^{\ul u_n}$) on $\bR^{B_n}$, whose \Kruzkov\ transform is $\sigma$-sub-additive by \cref{lem:monDom}.
The lower bound is 
$C := \exp( -(n\tau - 2 \Delta_{\max})/\alpha)$ and the tolerance criterion $\eta_n = (1-\rho) \exp(-n\tau/\alpha)$. The result follows.
\end{proof}

\begin{proof}[Proof of \cref{th:complexity_withUnknown}]
Same as the proof of \cref{th:complexity}, presented at the end of \cref{subsec:cplx_subadd}, replacing
\cref{corol:max_width,prop:niter_while} with
\cref{corol:max_width_withUnknown,prop:niter_while_withUnknown}.
\end{proof}

\subsection{Lazy implementation of the ANB}
\label{subsec:heaps}

The ANB method, as described in \cref{algo:anb}, includes several operations which incur a super-linear cost if naïvely implemented: the array copy $\ol u \gets \hat u$, the computation of the threshold value $t_n$, the enumeration of the narrow band elements $B_n$, and evaluation of the stopping criterion, on lines \ref{anb:hatu2}, \ref{anb:tn}, \ref{anb:Bn} and \ref{anb:for} of \cref{algo:anb}.
In this subsection, we describe a lazy implementation \cref{algo:anb_lazy} of the ANB, whose complexity is quasi-linear in terms of \emph{elementary operations} (rather than updates), while remaining fully equivalent to \cref{algo:anb} via an explicit bijection between their state variables, see \cref{prop:correct_anb_lazy} and \cref{prop:complexity_lazy}.
This reformulation requires two main ingredients: lists that remain sorted throughout the execution of \cref{algo:anb_lazy} (implemented as balanced binary search trees), and a placeholder value $\diamond$ representing the time $n\tau$ in the sub-solution.

\paragraph{A placeholder value ``$\diamond$'' representing $n\tau$.} 
The ANB assigns the value $n\tau$ to the sub-solution $\ol u$, and to Dijkstra's solution $\hat u$, in a large part of the domain, see \cref{algo:anb} lines \ref{anb:hatu} and \ref{anb:hatu2}.
As a result, simply replacing $(n-1) \tau$ with $n\tau$ in all relevant entries would have a super-linear cost $\cO(N n_{\max})$, where $n_{\max}$ bounds the number of iterations of the main loop (typically $n_{\max} = \cO(N^{1/d})$, see  \eqref{eq:bigO_h_ln_eps}). We use two different tricks for $\ol u$ and $\hat u$. 
Denote $\bR_\diamond := \bR \sqcup \{\diamond\}$, where $\diamond$ is an arbitrary additional symbol, and for any $u \in \bR_\diamond^X$, $T \in \bR$
\begin{align}
\label{eqdef:u_diamond_to}
	u^{\diamond \to T}(x) &:=
	\begin{cases}
		u(x) \!\!&\text{if } u(x) \neq \diamond\\
		T & \text{else}
	\end{cases} 
\end{align}
The sub-solution denoted $\ol u \in \bR^X$ and $\tilde u \in \bR_\diamond^X$, in the original and lazy ANB respectively, are related by $\ol u = \tilde u^{\diamond \to (n-1)\tau}$ at the start of the $n$-th iteration, see \cref{prop:correct_anb_lazy}.
On the other hand, Dijkstra's solution denoted $\hat u \in \bR^X$ and $\check u \in \bR_\infty^X$, in the original and the lazy ANB respectively, solve the slightly different problems below, taking $T=n\tau$ in the $n$-th step. The cutoffs $u^{< T}$ and $u^{\leq T}$ are defined in \eqref{eqdef:u_less_t}.

\begin{lemma}[\leancheck{lem:u_less_t}]
\label{lem:u_less_t}
		Let $\Gamma : \bR^X_\infty \to \bR^X_\infty$ be monotone and causal, let $\ol u \in \bR^X$, and let $T \in \bR$.
	Then so are the operators $\Gamma' := \max \{\ol u, \Gamma\}$, $\check \Gamma := \max\{\ol u,\Gamma\}^{<T}$ and $\hat \Gamma := \max \{\ol u, \min \{T,\Gamma\}\}$, and as a result each has a unique fixed point in $\bR_\infty^X$, respectively denoted $u, \check u, \hat u$. They are related as follows: 
	\begin{align*}
		\check u &= u^{<T} = \hat u^{<T}, &
		\hat u &= \max \{\ol u, \min \{T, u\}\} = \max \{\ol u, \min \{T, \check u\}\}.
	\end{align*}	 
\end{lemma}
\begin{proof}
The three operators are monotone and causal, since these properties transfer to the minimum or maximum of two such operators by \cref{lem:causal_mix}, and to a thresholded operator by \cref{prop:causal_formulation}. The existence and uniqueness of a fixed point follows by \cref{prop:comp_of_causal}.

\emph{Proof that $\check \Gamma u^{<T} = u^{<T}$.} %
Letting $v:= u^{<T}$ one has $(\Gamma' v)^{\leq T} = (\Gamma' u)^{\leq T} = u^{\leq T}$, 
using successively (i) causality formulated as in \eqref{eq:causal_old}, and (ii) the assumption $u = \Gamma' u$. 
Therefore, as announced, $(\Gamma' v)^{< T} = ((\Gamma' v)^{\leq T})^{<T} = (u^{\leq T})^{<T} = u^{<T} = v$,
using successively (i) cutoff definitions \eqref{eqdef:u_less_t}, 	 (ii) the previously established identity $(\Gamma' v)^{\leq T} = u^{\leq T}$, (iii) again \eqref{eqdef:u_less_t}, and (iv) the definition of $v$. 

\emph{Proof of the other relations.}
The unique fixed point of the operator $\hat \Gamma$ is $\hat u = \max \{\ol u, \min \{T, u\}\}$, by \cref{lem:max_min}. The result follows noting that $\min \{T,u\} = \min\{T,\check u\}$, and that $\hat u^{<T} =  \max \{\ol u, \min \{T, u\}\}^{<T} =  \max \{\ol u, u\}^{<T} = u^{<T}$ since $u \geq \ol u$.
\end{proof}

\paragraph{Sorted lists maintained throughout \cref{algo:anb_lazy}.} 
The lazy ANB \cref{algo:anb_lazy}, maintains four sorted lists $L$, $M$, $S$ and $S^c$ of pairs $(t,x) \in \bR_\infty \times X$, sorted lexicographically, using an arbitrary linear ordering of $X$. Similar to Dijkstra's method \cref{algo:Dijkstra}, the order on $X$ is only used to break ties w.r.t.\ the first component $t$ and make the algorithm deterministic; it has no impact on the method's output. 

The list $L$ is used by Dijkstra's \cref{algo:Dijkstra}, called as a sub-routine. 
The list $M$ allows computing threshold value $t_n$, on which the narrow band is centered.
The list $S$ allows enumerating the narrow band elements $B_n$. Finally, $S^c$ allows checking the termination criterion.
See lines \ref{anb_lazy:Dijkstra}, \ref{anb_lazy:tn}, \ref{anb_lazy:Bn} and \ref{anb_lazy:main_for} of \cref{algo:anb_lazy} respectively. 
Let us emphasize that a naïve brute force computation of $t_n$, $B_n$, or verification of the stopping criterion, would yield a super-linear complexity. %

Throughout the execution of \cref{algo:anb_lazy}, regular updates to $L$, $M$, $S$ and $S^c$ ensure that the following properties are maintained at appropriate keypoints:
\begin{align} 
\label{eqdef:L_heap}
	L &= \{ (t,x)\in \bR \times X \mid x \in X,\, t := \max \{\tilde u^{\diamond \to (n-1)\tau}(x),\, \Gamma \check u(x)\} < \check u(x) \,\}, \\
\label{eqdef:M_heap}
	M &= \{ (\tilde u(x), x) \in \bR \times X \mid x \in X,\, \tilde u(x) \neq \diamond,\ (1+\ve) \tilde u(x) < \undertilde u(x) \},\\
\label{eqdef:S_heap}
	S &= \{ (\tilde u(x), x) \in \bR \times X \mid x \in X,\, \tilde u(x) \neq \diamond\}, \\
\label{eqdef:Sc_heap}
	S^c & = \{(\undertilde u(x),x) \in \bR_\infty \times X \mid x \in X,\ \tilde u(x) = \diamond\}.
\end{align}

The most natural data structure for $L$, $M$, $S$ and $S^c$ is a balanced binary search tree, which is a sorted data structure with an $\cO(\ln N)$ cost for insertion and deletion, where $N$ is an upper bound for the number of elements. Those are available in e.g.\ the C++ standard library as \verb!std::set<std::pair<float,int>>!, where the floating point value represents $t\in \bR$, and the integer represents the index of the discretization point $x\in X$.

\paragraph{Correctness and complexity.} We show that the implementation \cref{algo:anb_lazy} of the ANB is equivalent to \cref{algo:anb}, and we estimate its complexity in terms of elementary operations, pointwise evaluations of the operators $\Lambda$ and $\Gamma$, and insertions and deletions in binary search trees.
The next result fully reconstructs the state of the ANB from the lazy ANB, and vice versa.
\begin{proposition}[Equivalence of \cref{algo:anb,algo:anb_lazy} \leancheck{prop:correct_anb_lazy}]
\label{prop:correct_anb_lazy}
Under \cref{assum:all},
denote by $\tilde u_n$, $\undertilde u_n$, $\check u_n^\init$ (resp.\ and $\check u_n$) the state of the variable $\tilde u$, $\undertilde u$, $\check u$ (resp.\ also $\check u$) in line \ref{anb_lazy:main_for} (resp.\ line \ref{anb_lazy:tn}) of the $n$-th iteration of \cref{algo:anb_lazy}. Then 
\begin{align}
\label{eq:lazy_bij}
	\ol u_n &= \tilde u_n^{\diamond \to (n-1)\tau}, &
	\undertilde u_n &= \ul u_n, &
	\check u_n^\init &= (\hat u_n^\init)^{< n\tau}, & 
	\hat u_n^\init &= \min\{\check u_n^\init, \max \{\overline u_n, n\tau\}\}, 
\end{align}
\begin{align}
\label{eq:lazy_bij_hat}
	\check u_n &= \hat u_n^{<n\tau}, &
	\hat u_n = \max \{\ol u_n, \min \{n\tau,\check u_n\}\},
\end{align}
\begin{equation}
\label{eq:lazy_bij_placeholder}
	\tilde u_n = 
	\begin{cases}
	\diamond \text{ where } \ol u_{n-1} <(n-1) \tau \text{ and } \hat u_{n-1} = (n-1) \tau, \\
	\ol u_n \text{ otherwise}.  	
	\end{cases}
\end{equation}
where $\ol u_n$, $\ul u_n$, $\hat u_n^\init$ and $\hat u_n$ are from \cref{algo:anb}, as described in \cref{prop:always_subsol,prop:Dijkstra_band}.
In addition, \eqref{eqdef:L_heap}, \eqref{eqdef:M_heap}, \eqref{eqdef:S_heap} and \eqref{eqdef:Sc_heap} hold on line \ref{anb_lazy:main_for} and after each update of $M,S,S^c$ in \cref{algo:anb_lazy}.
\Cref{algo:anb,algo:anb_lazy} have equivalent stopping criteria, and thus stop at exactly the same iteration.
\end{proposition}

\begin{proof}
The sorted lists $M$, $S$, $S^c$ are initialized consistently with $\tilde u$ and  $\undertilde u$, and are updated each time these mappings are modified, maintaining \eqref{eqdef:M_heap}, \eqref{eqdef:S_heap} and \eqref{eqdef:Sc_heap} throughout execution as announced.
The stopping criterion is equivalent to $(1+\ve) \tilde  u^{\diamond \to (n-1)\tau} \geq \undertilde u$, thanks to $M$ at the points where $\tilde u \neq \diamond$, and to $S^c$ at the points where $\tilde u = \diamond$. This is consistent with the stopping criterion $(1+\ve)\ol u \geq \ul u$ of the ANB \cref{algo:anb}, provided (\ref{eq:lazy_bij}, i and ii) hold which we prove below. 
We establish the other properties by induction on the iteration number $n$. 

\emph{Initialization $n=1$}. One has
$\undertilde u_1 = \infty = \ul u_1$,
$\tilde u_1^{\diamond \to 0} = 0 = \ol u_1$. Evaluating \eqref{eqdef:Dijkstra_supsol} yields $\hat u_1^\init = \tau$, given that $t_0^- = \hat u_0 = \ol u_1 = 0$, hence $(\hat u_1^\init)^{<n\tau} = \tau^{<\tau} = \infty = \check u^\init_1$ as desired. 
In addition \eqref{eqdef:L_heap} holds since $\Gamma[\infty] = \infty$ on $\Int X$ by additive invariance, and since $\max \{\tilde u^{\diamond \to 0}, \Gamma[\infty]\} = \max \{0,u_\bd\} = u_\bd$ on $\partial X$ recalling that the boundary conditions are non-negative.

\emph{Induction on $n$}.
We consider the $n$-th iteration of \cref{algo:anb_lazy}, and assume as an induction hypothesis that \eqref{eq:lazy_bij} and \eqref{eqdef:L_heap} hold on line \ref{anb_lazy:main_for}. 

We first prove that $\check u_n^\init$ is a super-solution for the operator $\max \{\ol u_n, \Gamma\}^{<n\tau}$: 
\begin{equation*}
	(\hat u_n^\init)^{< n\tau} 
	\geq \max \{\ol u_n, \min \{n\tau, \Gamma \hat u_n^\init\}\}^{< n\tau}
	= \max \{\ol u_n, (\Gamma \hat u_n^\init)^{< n\tau}\}^{<n\tau}
	= \max \{\ol u_n, \Gamma (\hat u_n^\init)^{< n\tau} \}^{< n\tau}.
\end{equation*}	
using (i) that $\hat u_n^\init$ is a super-solution for the operator $\max\{\ol u_n, \min \{n\tau, \Gamma\}\}$, see \cref{prop:Dijkstra_band}, (ii) $\min \{t, a\}^{<t}=a^{<t}$ for any $t,a \in \bR$, (iii) the causality of $\Gamma$.

The assumptions of Dijkstra's method are thus satisfied, for the operator $\max \{\tilde u^{\diamond \to (n-1)\tau}, \Gamma\} = \max \{\ol u_n, \Gamma\}$, and the stopping time $n\tau$, see \cref{prop:Dijkstra_correct}. Therefore it returns $\check u_n$ satisfying \eqref{eq:lazy_bij_hat} by \cref{lem:u_less_t}. After Dijkstra's method, $L$ satisfies
\begin{equation}
\label{eq:Lmid_anb_lazy}
	L = \{ (t,x) \mid x \in X,\, t := \max \{\ol u_n(x),\, \Gamma \check u_n(x)\} <  \check u_n(x) \,\}, 
\end{equation}
and $t\geq n\tau$ for all $(t,x)\in L$
by the stopping criterion of Dijkstra's method. 
Next we assign %
\begin{equation}
\label{eq:tilde_u_assign}
	\tilde u_{n+1}^\init(x) :=
	\begin{cases}
		\diamond \quad \text{if } \tilde u_n^{\diamond \to (n-1)\tau}(x) < n\tau \text{ and } \check u_n(x) = \infty \\
		\max \{\tilde u_n^{\diamond \to (n-1)\tau}(x), \min \{n\tau, \check u_n(x)\}\} \quad \text{else}
	\end{cases} \!\!\!
	=
	\begin{cases}
		\diamond \quad \text{if } \ol u_n(x) < n\tau \text{ and } \hat u_n(x) = n\tau \\
		\hat u_n(x) \quad \text{else}
	\end{cases}
\end{equation}
to the variable $\tilde u$
at a selection of points $x\in X$, and we claim that $\tilde u = \tilde u_{n+1}^\init$ everywhere. The equality (\ref{eq:tilde_u_assign}, right) follows from $\hat u_n = \max \{\ol u_n, \min \{n\tau, \check u_n\}\}$, and the equivalence ($\ol u_n(x) < n\tau$ and $\hat u_n(x) = n\tau$) $\Leftrightarrow$ ($\tilde u_n^{\diamond \to (n-1)\tau}(x) < n\tau$ and $\check u_n(x) = \infty$), given that $\check u_n(x)\in ]-\infty,n\tau[ \cup \{\infty\}$.
The points $x\in X$ where \eqref{eq:tilde_u_assign} is performed cover all those where the r.h.s.\ may have changed since the previous iteration as shown by the following case analysis: 
(a) if $\tilde u_n(x) = \diamond$, then either $\check u_n(x) = \infty$ and no update is needed since $\tilde u_{n+1}^\init = \diamond = \tilde u_n(x)$, or $\check u_n(x)<\infty$ and $x$ has just been selected by Dijkstra's method. We now assume that $\tilde u_n(x) \in \bR$, and distinguish 
(b) if $\tilde u_n(x) \geq n\tau$, then $\tilde u^\init_{n+1}(x) = \tilde u_n(x)$ and no update is needed, 
(c) if $\tilde u_n(x) < (n-1) \tau$ then either $x\in B_{n-1}$ is covered, or $\check u_n(x) = \tilde u_n(x) = \tilde u_{n+1}^\init(x)$ and no update is needed, and finally  
(d) the values $\tilde u_n(x) \in [(n-1)\tau,n\tau[$ are caught by the range query on $S$.

This modification of $\tilde u$ also does not invalidate $L$, which still satisfies \eqref{eq:Lmid_anb_lazy}, since
\begin{equation*}
	\max \{(\tilde u_{n+1}^\init)^{\diamond \to n\tau}, \Gamma \check u_n\} = \max \{\hat u_n,\Gamma \check u_n\} = \max \{ \ol u_n, \min \{n\tau, \Gamma \hat u_n\}, \Gamma \check u_n\} = \max \{\ol u_n, \Gamma \check u_n\},
\end{equation*}
using (i) $(\tilde u_{n+1}^\init)^{\diamond \to n\tau} = \hat u_n$, (ii) the fixed point property of $\hat u_n$, and (iii) $\Gamma \hat u_n \leq \Gamma (\hat u_n)^{<n\tau} = \Gamma \check u_n$.
The identity $(\tilde u_{n+1}^\init)^{\diamond \to n\tau} = \hat u_n$ also implies that $t_n$, $\Delta_n$ and $B_n$ are identical to those of the ANB, from which we deduce that $\undertilde u_{n+1} = \ul u_{n+1}$ and $\tilde u_{n+1}^{\diamond \to n\tau} = \ol u_{n+1}$ after the inner loops, establishing (\ref{eq:lazy_bij}, i and ii, iteration $n+1$). 
The band $B_n$ contains no placeholder points, since they correspond to the value $n\tau$, hence $\{\tilde u_{n+1} = \diamond\} = \{\tilde u_{n+1}^\init = \diamond\}$, which establishes (\ref{eq:lazy_bij_placeholder}, iteration $n+1$). Finally, $\check u_{n+1}^\init$ is obtained by resetting $\check u_n$ to $\infty$ on $B_n = \{t_n^- \leq \check u_n < n\tau\}$, but $\check u_n$ takes no value in $[n\tau,\infty[$, hence $\check u_{n+1}^\init = \check u_n^{<t_n^-} = \hat u_n^{<t_n^-}$, and (\ref{eq:lazy_bij}, iii and iv, iteration $n+1$) follow by comparison with \eqref{eqdef:Dijkstra_supsol}. The last update of $L$ ensures that \eqref{eqdef:L_heap} still holds after these updates to $\tilde u$ and $\check u$, and the proof is complete.
\end{proof}

\begin{proposition}[\leancheck{prop:complexity_lazy}]
\label{prop:complexity_lazy}
Under the assumptions of \cref{th:complexity} or \ref{th:complexity_withUnknown},
\cref{algo:anb_lazy} has complexity
\begin{equation}
\label{eq:anb_cost}
\cO(\widetilde N \ln \widetilde N) 
+ 2 N N_\Lambda U_\Lambda 
+ (\widetilde N N_\Gamma + (N+\widetilde N) N'_\Lambda) U_\Gamma 
+ (N+\widetilde N)(N'_\Lambda + N_\Gamma) \cO( \ln N) +  \cO(n_{\max} \ln N),
\end{equation}
where $N := \#(X)$ and $\widetilde N := \sum_{x \in X} \#(\cV(x)) = \sum_{x \in X} \#(\widetilde \cV(x))$, is the total cardinality of the stencils of $\Gamma$, see \eqref{eqdef:stencil_card} and \eqref{eqdef:stencil}.
We denoted by $U_\Lambda$, $U_\Gamma$ and $\cO(\ln N)$ the unit cost of, respectively, the evaluation of $\Lambda u(x)$, $\Gamma u(x)$, and the insertion or deletion into a balanced binary search tree with at most $N$ elements.
$N_\Lambda$, $N'_\Lambda$ and $N_\Gamma$ are from \cref{th:complexity} or \ref{th:complexity_withUnknown}.
\end{proposition}

\begin{proof}
	The computation of the reversed stencils $\widetilde \cV$, in the pre-processing step, amounts to reversing the $\widetilde N$ edges of the oriented graph defined by $\cV$. This is achieved by sorting these edges w.r.t.\ their endpoint rather than their source, which has cost $\cO(\widetilde N \ln \widetilde N)$. For the initialization of $L$, the Lean implementation recovers $u_\bd$ and $\partial X$ by computing $\Gamma[\infty]$, which adds $N$ evaluations of $\Gamma$ to the estimate \eqref{eq:anb_cost}.

A given point $x \in X$ is selected at most $N_\Gamma$ times in Dijkstra's subroutine line \ref{anb_lazy:Dijkstra}, overall in the course of \cref{algo:anb_lazy}, by \cref{th:complexity,th:complexity_withUnknown}. Each time, this triggers an update of $L$, an update of $M$, $S$ and $S^c$, and $\# \widetilde \cV(x)$ evaluations of the operator $\Gamma$ and insertions in the list $L$, with total cost $\#\widetilde \cV(x)\, U_\Gamma + (\#\widetilde \cV(x)+4) \cO(\ln N)$.
Summing over all points $x\in X$, each selected at most $N_\Gamma$ times by Dijkstra's method, the overall cost of this part of the algorithm is thus
$N_\Gamma (\widetilde N U_\Gamma + (\widetilde N + N) \cO(\ln N))$.

A given point $x \in X$ is contained within the narrow band $B_n$ in at most $N'_\Lambda$ iterations in the course of \cref{algo:anb_lazy}, by \cref{th:complexity,th:complexity_withUnknown}. Each time, there is a cost $\cO(\ln N)$ for extracting this element in line \ref{anb_lazy:Bn}, at most $N''_\Lambda$ evaluations of $\Lambda \ol u(x)$ and $\Lambda \ul u(x)$ in lines \ref{anb_lazy:while_sub} and \ref{anb_lazy:while_sup}, one update of $M$ and of $S$ in line \ref{anb_lazy:update_MS}, an assignment in line \ref{anb_lazy:parallel_for2}, and $1+\#\widetilde \cV(x)$ evaluations of $\Gamma$ and updates of $L$ in line \ref{anb_lazy:L}, and finally one update of $M$ and $S$ in line \ref{anb_lazy:olu_gets_hatu}, with cost $5 \cO(\ln N) + 2 N''_\Lambda U_\Lambda + (\#\widetilde\cV(x)+1) (U_\Gamma + \cO(\ln N))$. Summing over all $x\in X$, and over all iterations in which they belong to $B_n$, the overall cost of this part of the algorithm is at most $N N'_\Lambda (2 N''_\Lambda U_\Lambda + \cO(\ln N)) + (N+\widetilde N) N'_\Lambda (U_\Gamma + \cO(\ln N))$. Recall that $N_\Lambda := N'_\Lambda N''_\Lambda$.

Line \ref{anb_lazy:olu_gets_hatu} introduces the set $W_n := \{ (t,x) \in S \mid t \in [(n-1)\tau,n\tau[\}$. Since $S$ is a balanced binary search tree, constructing $W_n$ has cost $\#(W_n) \cO(\ln N)$. 
One has $x \in W_n$ iff $\tilde u_n(x) \in [(n-1)\tau,n\tau[$ (in particular $\tilde u_n(x) \neq \diamond$), say in iterations $n_1,\cdots,n_K$.
Since Dijkstra's method in iteration $m$ only produces values 
$<m\tau$,
a real value $\tilde u_n(x) \geq (n-1)\tau$ at the start of iteration $n$ comes from an inner loop of an earlier iteration $k$ with $x\in B_k$, and is consumed when $x\in W_n$; hence $K \leq N'_\Lambda$.
In each occurrence, an update of $M$, $S$, $S^c$ is triggered, with $\cO(\ln N)$ cost. The overall associated cost, over all points and iterations, is thus $4 N  N'_\Lambda \cO(\ln N)$. Acknowledgement: our original manuscript did not include the set $W_n$, and this mistake was spotted by Claude$^{\small\textregistered}$ during formalization, which proposed on the fly the fixes to \cref{algo:anb_lazy,prop:correct_anb_lazy,prop:complexity_lazy}; this is the main gap discovered during formalization, and the main AI contribution to this manuscript.

The remaining non-trivial computations are the extraction of the smallest element of $M$ to evaluate $t_n$, the largest element of $S^c$ to check the stopping criterion, and we still need to count two range queries on $S$ even if the set $B_n$ is empty, which amounts to $4 \cO(\ln N)$ cost per iteration. 
\end{proof}

\begin{algorithm}[!t]
\caption{Adaptive Narrow Band (ANB), complete implementation \leancheck{algo:anb_lazy}}
\label{algo:anb_lazy}
\begin{algorithmic}[1]
\State
\textbf{Input:} Tolerance $\ve>0$, constants $\tau>0$, $\alpha>0$, $\rho<1$, positive boundary conditions $u_\bd$
\State
\textbf{Ensures at termination:} 
$\ol u\leq \Lambda \ol u$, 
$\Lambda \undertilde u \leq \undertilde u$, 
and $(1+\ve) \ol u \geq \undertilde u$, 
where $\ol u := \tilde u^{\diamond \to (n-1) \tau}$
and $n$ is the index of the last iteration.
\State 
\textbf{Initialization:} Set $\tilde u\gets \diamond$ and $\check u, \undertilde u \gets \infty$, everywhere on $X$. Compute the reversed stencils $\widetilde \cV$, see \eqref{eqdef:stencil}. Set $M, S, B_0 \gets \emptyset$, and 
$L\gets \{(u_\bd(x),x) \mid x \in \partial X\}$, and $S^c \gets \{(\infty,x) \mid x \in X\}$.
\State \label{anb_lazy:main_for}
\textbf{For} $n=1,2,\cdots$ \textbf{until} 
$M = \emptyset$ 
\State \phantom{blablablablablabla}
and ($S^c = \emptyset$ or the largest $(t,x)\in S^c$ satisfies $t \leq (1+\ve)(n-1)\tau$) \textbf{do}:
\State \phantom{bla} \label{anb_lazy:Dijkstra}
Call Dijkstra's method, see \cref{algo:Dijkstra}, with the following parameters: 
\State
\phantom{bla} initial guess $\check u$, list $L$ (already initialized), stopping time $n\tau$, operator $\max \{\tilde u^{\diamond\to (n-1) \tau}, \Gamma\}$.
\State \phantom{bla} \label{anb_lazy:olu_gets_hatu}
\textbf{For all}  $x$ in $B_{n-1}$, or updated by Dijkstra, or s.t. $(t,x) \in S$ with $(n-1)\tau \leq t < n\tau$
 \textbf{do}:
\State \phantom{blabl}
$\tilde u(x) \gets \diamond$ if $\tilde u^{\diamond \to (n-1)\tau}(x) < n\tau$ and $\check u(x)=\infty$, else $\max \{\tilde u^{\diamond \to (n-1)\tau}(x), \min\{n\tau, \check u(x)\}\}$
\State \phantom{bla}
Update $M$, $S$, $S^c$, accordingly.
\State \phantom{bla} \label{anb_lazy:tn}
$t_n \gets \min \{n\tau, \min \{t \mid (t,x) \in M\}\}$
\State \phantom{bla}
$\Delta_n \gets n\tau-t_n$
\State \phantom{bla} \label{anb_lazy:Bn}
$B_n \gets \{x \in X \mid (t,x) \in S \text{ with } n\tau-2\Delta_n \leq t < n \tau\}$.
\State \phantom{bla} \label{anb_lazy:etan}
$\eta_n \gets (1-\rho) \exp(-n\tau/\alpha)$
\State \phantom{bla} \label{anb_lazy:while_sub}
\textbf{Until} $\exp(-\tilde u/ \alpha) \leq \exp(-\Lambda \tilde u^{\diamond \to n\tau}/ \alpha) + \eta_n$ on $B_n$ \textbf{do}:
\State \phantom{blabla} \label{anb_lazy:while_sub_updt}
\textbf{For all} $x \in B_n$ \textbf{in parallel}: $\tilde u(x) \gets \Lambda \tilde u^{\diamond \to n\tau}(x)$
\State \phantom{bla} \label{anb_lazy:while_sup}
\textbf{Until} $\exp(-\frac{\Lambda\undertilde u} \alpha) \leq \exp(-\frac{\undertilde u} \alpha) + \eta_n$ on $B_n$ \textbf{do}:
\State \phantom{blabla} \label{anb_lazy:while_sup_updt}
\textbf{For all} $x \in B_n$ \textbf{in parallel}: $\undertilde u(x) \gets \Lambda \undertilde u(x)$
\State \phantom{bla} \label{anb_lazy:update_MS}
Update the lists $M$, $S$, $S^c$, using the new values of $\tilde u(x)$ and $\undertilde u(x)$ for all $x\in B_n$.
\State \phantom{bla} \label{anb_lazy:parallel_for2}
\textbf{For all} $x \in B_n$ \textbf{do in parallel}: $\check u(x) \gets \infty$.
\State \phantom{bla} \label{anb_lazy:L}
Update $L$ using the new value $\max \{\tilde u^{\diamond \to n\tau}(y), \Gamma \check u(y)\}$ for all $y \in \cup_{x \in B_n} (\{x\} \cup \widetilde\cV(x)) $
\end{algorithmic}
\end{algorithm}

\paragraph{Additional implementation notes for \cref{algo:anb_lazy}.}
In order to avoid floating point underflow, the tests of lines \ref{anb_lazy:while_sub} and \ref{anb_lazy:while_sup} are reformulated as follows: for all $a,b,n\tau\in \bR$ and $\alpha>0$
$$
\exp(-\frac a \alpha) \leq \exp(-\frac b \alpha) + (1-\rho) \exp(-\frac{n\tau} \alpha)
\ \Leftrightarrow \
1 \leq \exp( \frac {a-b} \alpha) + (1-\rho) \exp(\frac {a-n\tau} \alpha).
$$
One may implement the modified update $\undertilde u(x) \gets \min\{\Lambda \undertilde u(x), \undertilde u(x)\}$ in line \ref{anb_lazy:while_sup_updt}; this is mathematically equivalent since $\Lambda \undertilde u(x) \leq \undertilde u(x)$ by the super-solution property, but it limits some issues related with floating-point round-off errors.
Similarly $\tilde u(x) \gets \max \{\Lambda \tilde u^{\diamond \to n\tau}(x), \tilde u(x)\}$ in line \ref{anb_lazy:while_sub_updt}.

Many optimizations can improve practical performance, even though the quasi-linear complexity guarantee established above may be lost.
For instance one may replace the balanced binary search trees $L$ and $M$ with priority queues, which are computationally cheaper but only allow the removal of the smallest element (in which case one needs to reject outdated keys). One may also implement $S$ as a simple array, and sort appropriate tail subsets of that array. One can get rid of $S^c$, and use the termination test ``$M=\emptyset$ and $L = \emptyset$'', %
in which case it may continue for a few more iterations than necessary, with empty bands $B_n$, producing the same super-solution and a slightly larger (still valid) sub-solution.
See also \cref{sec:AGSI}.

\section{Numerical results}
\label{sec:num}
We compare the methods introduced in this paper (ANB, ANB$^*$, AGSIK), with reference iterative methods from the literature (FS, FSd, AGSI, FIM, GI, NB), 
see the next paragraph for the description of these acronyms. 
The test cases are similar to \cite{mansour2024solving}, which already gathered classical reference test cases from the literature \cite{cacace2011local,mirebeau2014EfficientFinsler,mirebeau2014marchingLattice}. 
They include a one-dimensional toy experiment, a pursuit-evasion game involving a non-convex Hamiltonian, an anelliptic problem based on geological data, a maze with slightly permeable walls,
a strongly anisotropic test case related with tubular structure segmentation, and a strongly asymmetric instance of Zermelo's navigation problem.
Expectedly, some numerical results for the legacy methods are identical to those reported in \cite{mansour2024solving}. C++ codes are available on demand.
All test cases are discretized on regular Cartesian grids of shape $n_1 \times \cdots \times n_d$, limiting ourselves to dimension  $d=2$ for simplicity and visualization purposes, even though nothing prevents running the methods in higher dimension.

\paragraph{Considered fixed point solvers.}
\emph{Iterative methods.}
We consider the fast-sweeping method (FS) \cite{zhao2005FastSweeping}, its popular variant using diagonal sweeps (FSd) and often referred to as alternating diagonal sweeping, the adaptive Gauss-Seidel iteration method (AGSI) \cite{bornemann2006FiniteElement}, the fast iterative method (FIM) \cite{jeong2008fim}, and global iteration (GI) see \cref{algo:global_iteration}.

\emph{Narrow band methods.}
We implement the adaptive narrow band (ANB) method described in \cref{algo:anb,algo:anb_lazy}, as well as the variant (ANB$^*$) described in \cref{sec:AGSI} which uses the AGSI for the inner iterations within the narrow band, thus avoiding many redundant updates. 
\emph{For practical applications, in particular with semi-Lagrangian schemes, we highly recommend the ANB$^*$ over the ANB, even though its guarantees are slightly weaker, see \cref{table:guarantees,table:niter}.}
We also consider the (NB) method introduced by one of the authors \cite{mansour2024solving}, again using the AGSI for the inner iterations within the corresponding narrow bands. Finally, we introduce the (AGSIK) method, see \cref{algo:agsik} of \cref{sec:AGSI}, which is a  progressive variant of the AGSI, featuring intermediate objectives whose stopping criterion is expressed in \Kruzkov\ coordinates.

\begin{figure}
\def\h{3cm}
\includegraphics[height=\h]{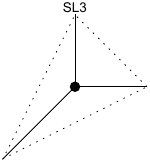}
\includegraphics[height=\h]{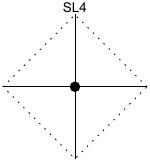}
\includegraphics[height=\h]{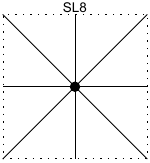}
\hspace{1cm}
\includegraphics[height=3.5cm]{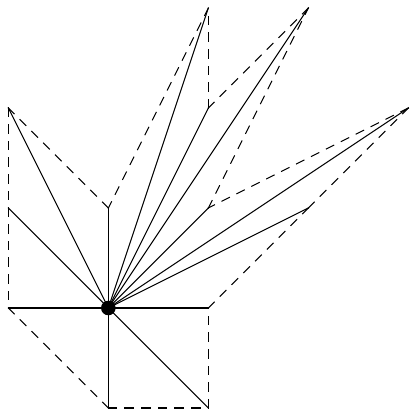}
\hspace{-0.3cm}
\includegraphics[height=\h]{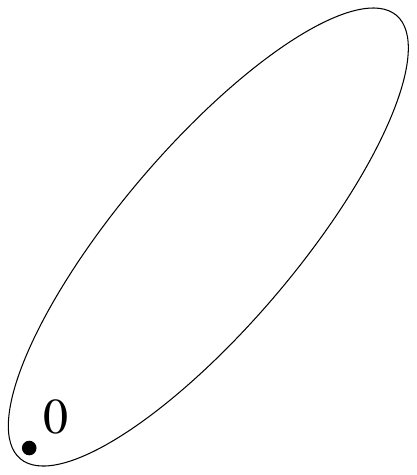}
\caption{
(i, ii, iii) Stencils of the three-, four- and eight-point Semi-Lagrangian method (SL3, SL4, SL8). 
SL4 and SL8 are common in applications, whereas SL3 is introduced to mimic unstructed triangulations with ill-shaped elements. (iv, v) Example of an adaptively refined anisotropic stencil, used to ensure the causality property in the FM-ASR method \cite{mirebeau2014EfficientFinsler}.
}
\label{fig:stencils}
\end{figure}

\emph{Fast-Marching type methods.} 
Several of the considered problems could be addressed using variants of the fast marching method, such as  \cite{kimmel1998Manifold,sethian2003OUM,alton2012MAOUM,vladimirsky2008LabelSetting,dahiya2018ordered}.
Our experience with the variants introduced by one of the authors, namely fast marching using anisotropic stencil refinement (FM-ASR) \cite{mirebeau2014EfficientFinsler,mirebeau2014marchingLattice}, fast marching using Voronoi's first reduction (FM-VR1) \cite{mirebeau2019riemannian}, and a FMM tailored for the pressure wave traveltime in geological media \cite{desquilbet2021single}, is that they often outperform the alternatives considered above \emph{both in terms of computation time and accuracy}. 
However, these approaches typically solve a \emph{heavily modified discretization scheme} based on \emph{ad-hoc adaptive stencils ensuring the causality property}, as illustrated on \cref{fig:stencils}, hence we are not really comparing apples to apples, as discussed in the next paragraph.
 We refer to \cite{mansour2024solving} for numerical comparisons between the NB method and such anisotropic variants of the fast marching method, on the same test cases as in this paper, and we do not consider FMMs further in this section. 

\paragraph{Choice of discretization scheme.}
Our objective is to illustrate the fact that the ANB, ANB$^*$ and AGSIK, maintain a quasi-linear complexity even when solving very ``bad'' numerical schemes, in the sense that they are \emph{highly non-causal}. Said otherwise, these methods can efficiently solve schemes which are very \emph{low on the causality ladder} described in \cref{rem:causality}: they are not causal, but only $\alpha$-\Kruzkov-contractant, and the timescale $\alpha>0$ is large and often unknown. 
In applications we would rather recommend using a causal scheme, or a moderately non-causal scheme, but that is not always a possibility as discussed below.
The numerical schemes are described in \cref{sec:effective}, and the semi-Lagrangian stencils are illustrated on \cref{fig:stencils}.

The only \emph{causal} scheme used in these experiments is the standard upwind differences (UD) discretization of the \emph{isotropic} eikonal PDE \cite{sethian1996fast}, in the maze test case \cref{subsec:maze}. 
There also exists causal schemes for the PDEs considered in \cref{subsec:seismic,subsec:snail,subsec:zermelo}, see \cite{desquilbet2021single,mirebeau2019riemannian,mirebeau2014EfficientFinsler}, but they are not considered here, see   \cite{mansour2024solving}.
Such causal schemes should usually be preferred when available, but it is not always the case since their design and implementation may require significant effort depending on the metric structure and domain representation.

The \emph{moderately non-causal} schemes include the upwind differences scheme (UD) when a drift term is present, and the four-point semi-Lagrangian scheme SL4 (resp.\ eight point SL8), see \cref{fig:stencils}, when the metric is anisotropic (resp.\ strongly anisotropic).
They are used in the test cases inspired by tubular segmentation \cref{subsec:snail}, and Zermelo's navigation problem \cref{subsec:zermelo}. 
Let us acknowledge that the diagonal Fast-Sweeping scheme FSd (resp.\ the original fast sweeping scheme FS) is particularly efficient at solving the SL4 and SL8 schemes (resp.\ the SL8 scheme only), since the sweeping directions match well with the grid-aligned stencil.
 
The \emph{highly non-causal} schemes used in our experiments are the Lax-Friedrichs (LF) scheme, due to the strong diffusivity of the elliptic regularization term, and the semi-Lagrangian scheme when it is combined with the small and ill-shaped $3$-point stencil (SL3). These two schemes also suffer from poor overall accuracy, and the LF scheme additionally creates a spurious boundary layer near the domain boundary. Needless to say, we recommend \emph{against} using these schemes in applications. 
Nevertheless, the Lax-Friedrichs scheme is hard to avoid when the Hamiltonian is non-convex, see the differential game test case \cref{subsec:pursuit}, or when it has a complex algebraic structure, see the seismic travel time test case \cref{subsec:seismic}.
The semi-Lagrangian scheme, on the other hand, is often used to discretize a PDE on a \emph{given unstructured mesh} (rather than a Cartesian grid), which dictates the stencil and often features some ill shaped elements; it is this unfavorable configuration that we attempt to emulate with the SL3 stencil, see \cref{fig:stencils}.

\paragraph{Choice of convergence tolerance.}
The convergence tolerance for the proposed ANB method, $(1+\ve) \ol u \geq \ul u$ where $\ol u$ and $\ul u$ are respectively a sub- and a super-solution, differs from the classical one $0 \leq \ul u - \Lambda \ul u \leq \tilde \ve$ used in the FS, AGSI, FIM \cite{zhao2005FastSweeping,bornemann2006FiniteElement,jeong2008fim}, or the AGSIK.
If one wants to achieve a similar accuracy in the $\ell^\infty(X)$ norm, then one should use $\tilde \ve \approx \nu \ve $ in view of \cref{corol:tol}, where the parameter $\nu>0$ is associated with the $\nu$-sub-multiplicativity property. 
For the numerical schemes of interest, $\nu$ is a small constant times the grid scale $h>0$, see \cref{sec:effective}.
Based on these considerations we use the tolerance parameters 
\begin{align}
\label{eq:default_tol}
	\ve &= 10^{-3}, &
	\tilde \ve &= 10^{-4} h,
\end{align}
in the following experiments, unless otherwise mentioned. Note that the numerical cost of the ANB only depends logarithmically on $\ve$, as established in \cref{th:complexity}, and illustrated on \cref{fig:iter_wrt_eps}.

\paragraph{Computational cost.}
We report the number of pointwise evaluations of the update operator $\Lambda$ for the various methods, see \cref{table:niter,fig:iter_wrt_alpha,fig:iter_wrt_eps,fig:iter_wrt_rho,fig:iter_wrt_tau}.
The number of evaluations of the graph based approximation $\Gamma$ is reported in \cref{fig:iter_wrt_tau}, and their numerical cost is in practice negligible except in the maze test case \cref{subsec:maze}.
The reported CPU times, see \cref{table:elapsed}, were obtained using a sequential (single thread) C++ implementation on an Apple$^{\footnotesize\textregistered}$
laptop equipped with an M1-pro$^{\footnotesize\textregistered}$ processor. 
Note that the ANB and ANB$^*$ construct both a sub-solution and a super-solution, whereas the other methods only produce a super-solution; for fairness, we thus count separately the evaluations of $\Lambda$ on the sub-solution, and the corresponding CPU time, in \cref{table:niter,table:elapsed}. 

Let us acknowledge that achieving state of the art computation time would require a parallel implementation, typically using thousands of threads on a GPU architecture following \cite{jeong2008fim,kimmel2008ParallelDistance,mirebeau2023massively}. Significant algorithmic modifications would be required, such as grouping the discretization points into $8 \times 8$ or $4 \times 4 \times 4$ blocks depending on the dimension following \cite{jeong2008fim}, and finding parallelizable counterparts for the ANB subroutines which include Dijkstra's method as well as the maintenance operations of the various heap structures used in \cref{algo:anb_lazy}. We regard this extension as an opportunity for future work.

\begin{table}
	\begin{tabular}{|c|l|cccc|} \hline
		& & ANB & ANB$^*$ & NB & AGSIK,\,FS,\,FIM,\,AGSI \\ \hline
		\multirow{2}{*}{Tolerance guarantee} &$\alpha, \rho$ known & \checkmark & \checkmark & \checkmark & \checkmark \\
		 & $\alpha, \rho$ unknown & \checkmark & \checkmark  & \ding{55} & \checkmark \\ \hline
		\multirow{2}{*}{Complexity guarantee} &$\alpha, \rho$ known & \checkmark & \checkmark & \checkmark & \ding{55} \\
		 &$\alpha, \rho$ unknown & \checkmark & \ding{55} & \ding{55} & \ding{55} \\ \hline
	\end{tabular}
	\caption{Summary of the guarantees of the various methods, see \cref{th:complexity,th:complexity_withUnknown,th:complexity_with_agsi} for a complete statement.
	Note that the FS, FIM and AGSI methods ignore the parameters $\alpha,\rho$.
	}
	\label{table:guarantees}
\end{table}

\paragraph{Summary of tolerance and complexity guarantees.} When running a computer program, two questions are fundamental: Can the output be trusted? Can the running time be reasonably bounded? 
We summarize below, and in \cref{table:guarantees}, the theoretical guarantees of the different methods, focusing on a particular issue: does it require as input a valid timescale $\alpha>0$, and contraction ratio $0<\rho<1$ for the discretization scheme in \Kruzkov\ coordinates?
\begin{itemize}
	\item \emph{Tolerance guarantee}, i.e.\ can the method output be trusted (if it terminates)? 
	The termination criterion of the ANB method, and of its variant the ANB$^*$, ensures that at termination the sub-solution $\ol u$ and super-solution $\ul u$ satisfy $(1+\ve) \ol u \geq \ul u$, see \cref{algo:anb}. 
	Likewise, the termination criterion of the FS, FIM, AGSI, GI, and the proposed AGSIK, ensures that at termination the numerical solution satisfies $0 \leq u-\Lambda u \leq \tilde \ve$. 
	This is in contrast with the NB method, for which an inadequate choice of parameters may yield an incorrect numerical solution at termination; when this failure mode is diagnosed, i.e.\ $u-\Lambda u \not \leq \tilde \ve$, the remedy is to increase the parameter $\alpha$ \cite[Remark 4.3]{mansour2024solving}, as in \cref{subsec:snail,subsec:zermelo}.
	
	\item \emph{Quasilinear complexity guarantee}, i.e.\ can the running time be reasonably bounded? 
	This property is guaranteed for the ANB with any choice of parameters, see \cref{th:complexity_withUnknown} (up to mild constraints on the tolerance $\ve>0$, see \cref{rem:tol_range}), even though the underlying multiplicative constant increases if e.g.\ $\alpha$ is chosen excessively small or large.
	On the other hand, the quasi-linear complexity of the ANB$^*$ is only guaranteed if $\Phi_\alpha$ is $\rho$-sub-additive, see \cref{th:complexity_with_agsi}, and likewise for the NB method.
	Finally, the classical FS, FIM, AGSI, and the proposed AGSIK, do not offer any complexity guarantee, even though they often work well in practice. Note that Fast-Marching type methods have guaranteed quasi-linear complexity, but they are not considered here, see the discussion above. 
\end{itemize}

\paragraph{Practical choice of parameters, and linearized contraction factor.}
The convergence analysis of the ANB method is based on a contraction property $\|\Phi_\alpha u-\Phi_\alpha v\|_\infty \leq \rho \|u-v\|_\infty$, for the non-linear discretization scheme $\Phi_\alpha := \rE_\alpha \Lambda_\infty \rE_\alpha^{-1}$ expressed in \Kruzkov\ coordinates, where $u,v \in \bR^X$ are arbitrary, see \cref{rem:causality}. Our theoretical estimates of $\alpha>0$ and $0 < \rho < 1$ are rather pessimistic, because they are based on the rough estimation of composed nonlinear effects, see \cref{sec:effective}.

In our test cases, these theoretical estimates typically suggest $\alpha \geq 100h$, either due to the strong diffusivity of the Lax-Friedrichs scheme, or to the strong anisotropy of the metric.
However, the 
minimal running time of the ANB and ANB$^*$ is
typically obtained with $h \leq \alpha \leq 5h$, and choosing $\alpha$ outside this sweet spot can significantly increase computational cost, see \cref{fig:iter_wrt_alpha}.

\begin{figure}
\def\h{4cm}
	\includegraphics[height=\h]{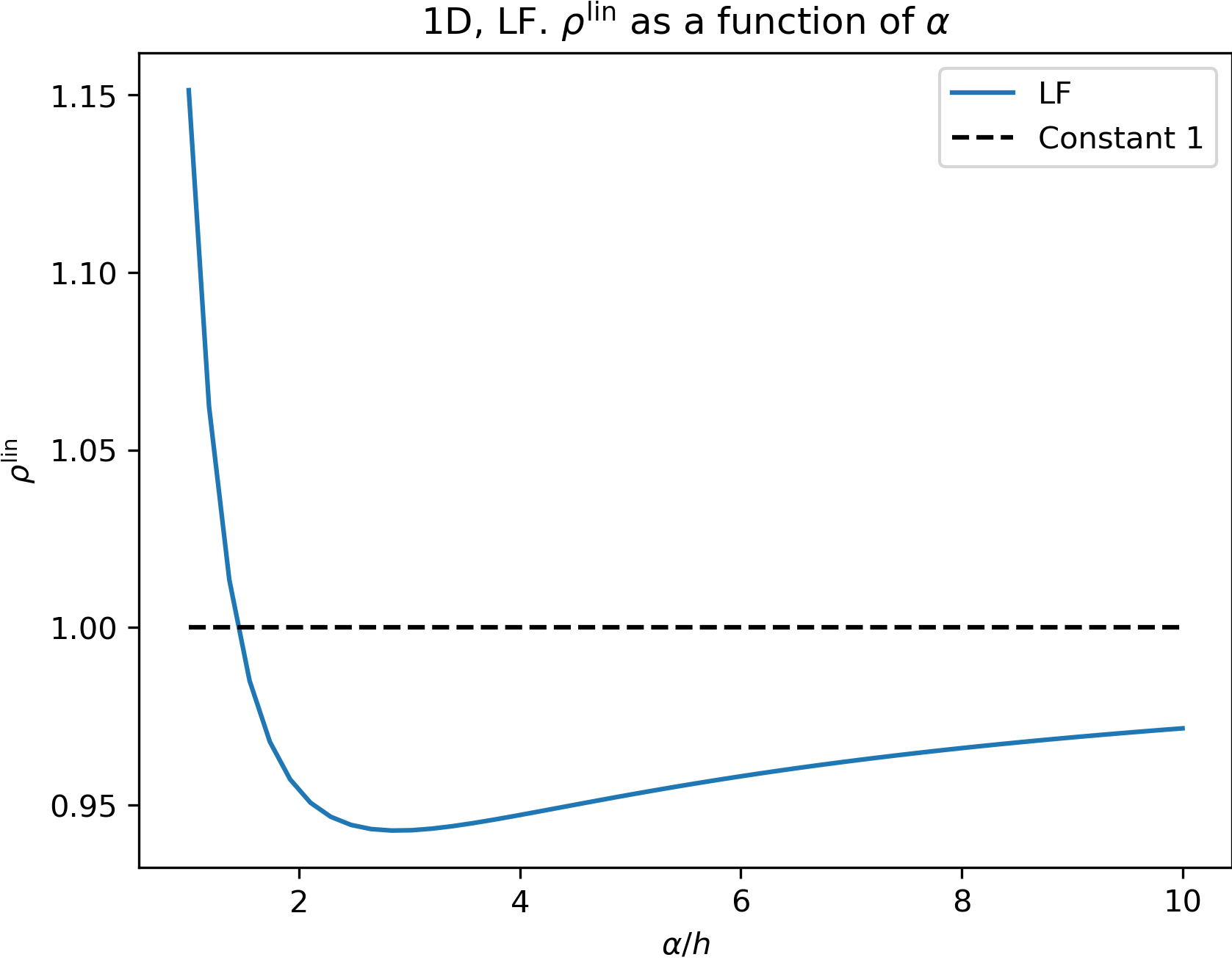}
	\includegraphics[height=\h]{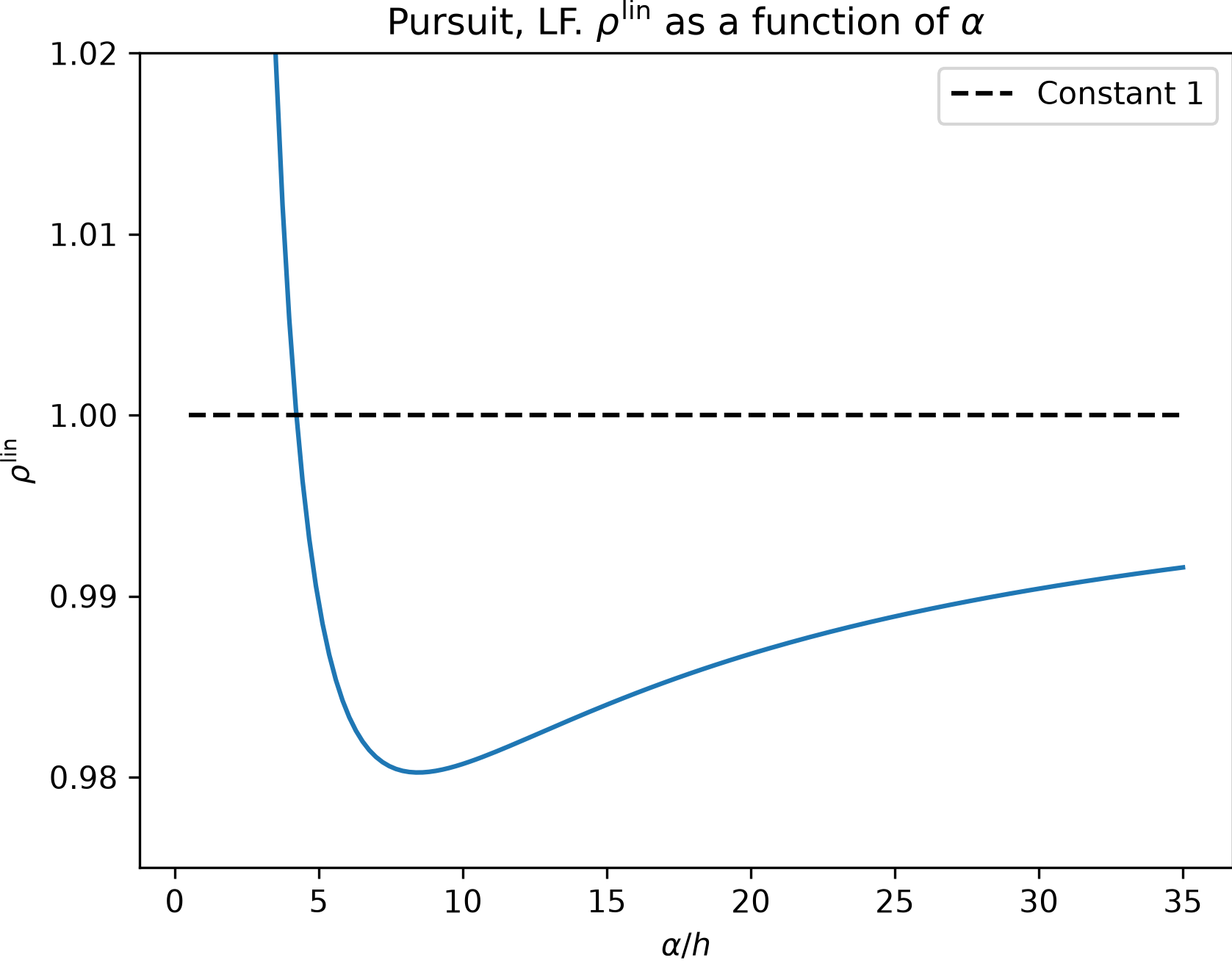} 
	\includegraphics[height=\h]{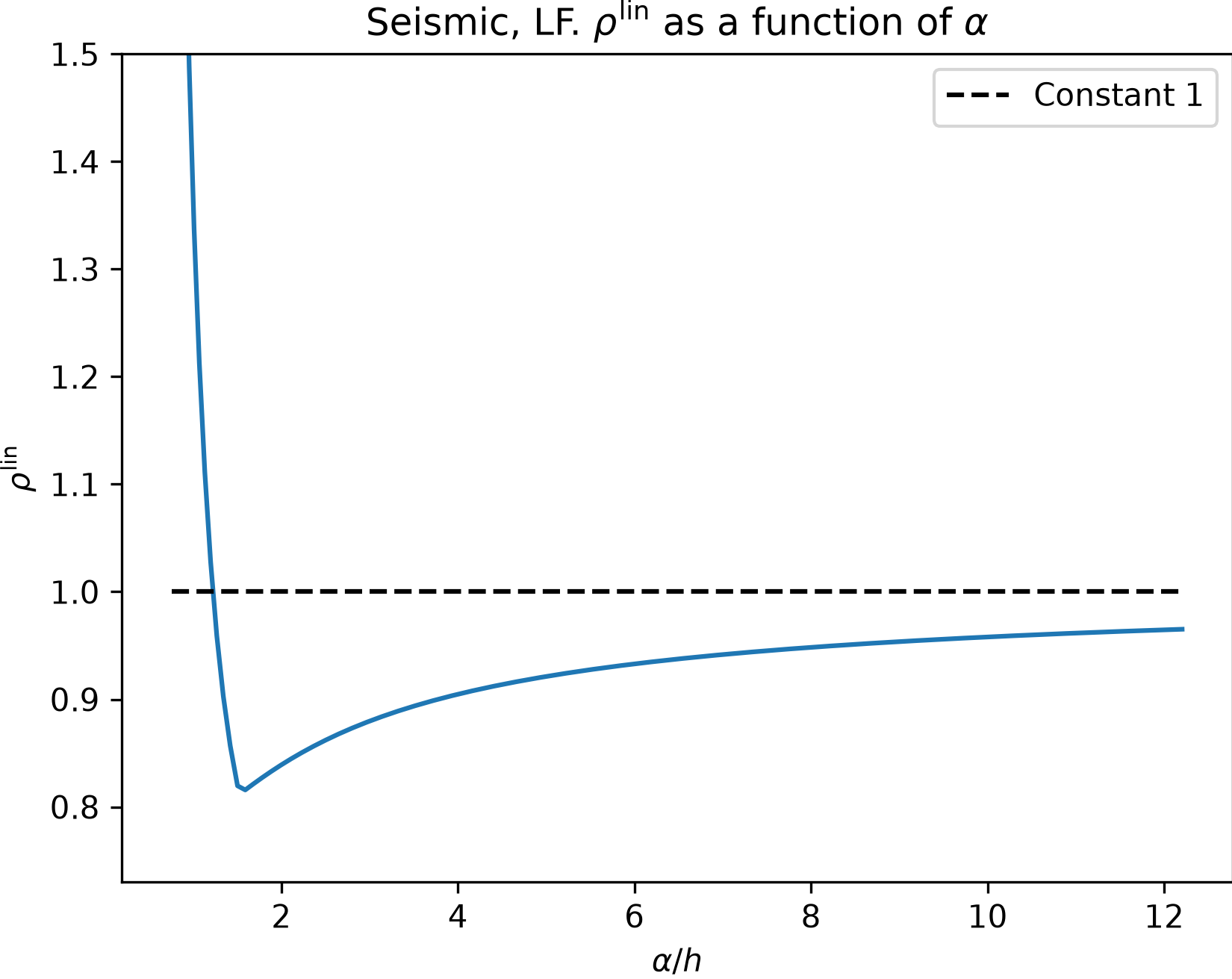} \\
	\includegraphics[height=\h]{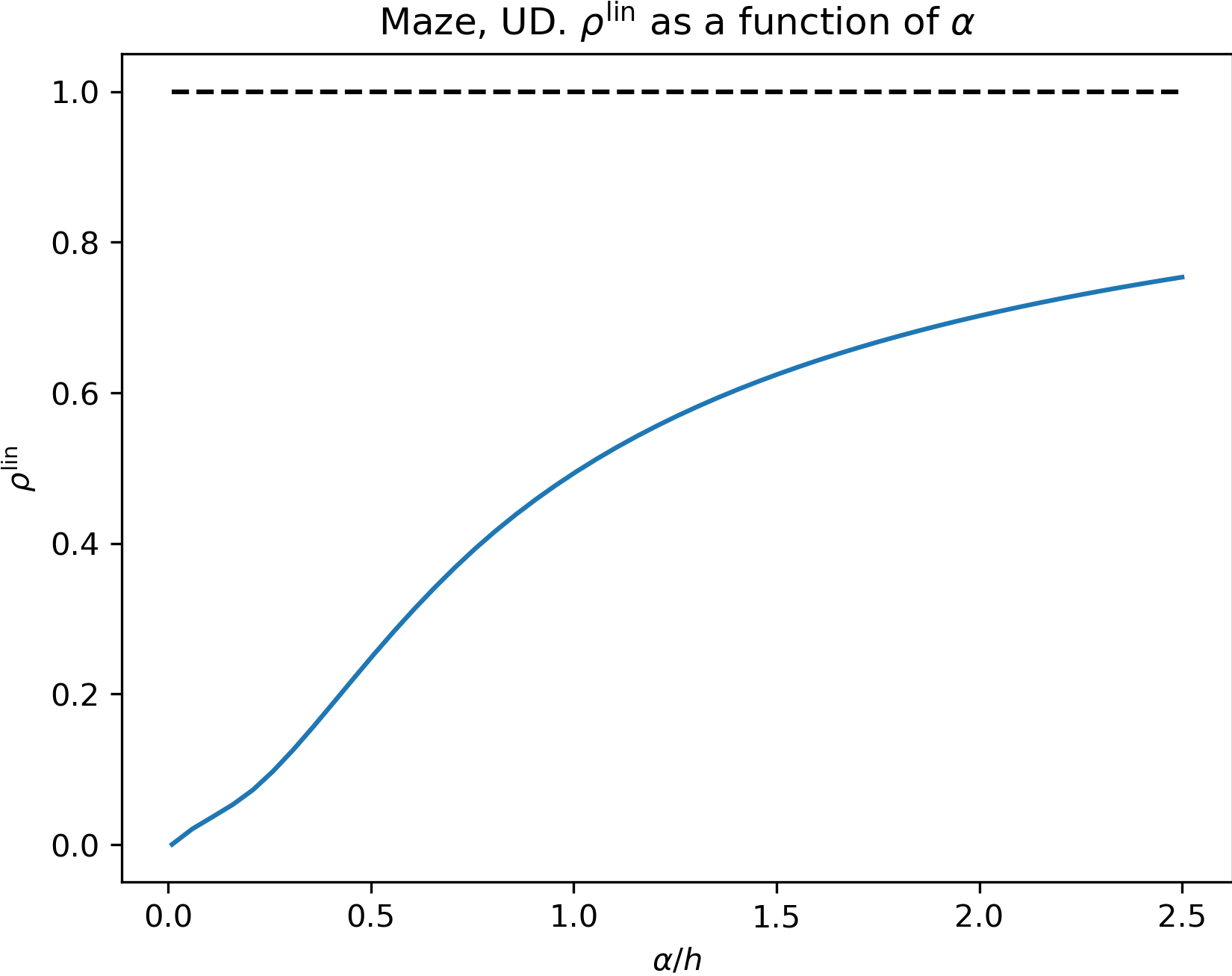} 
	\includegraphics[height=\h]{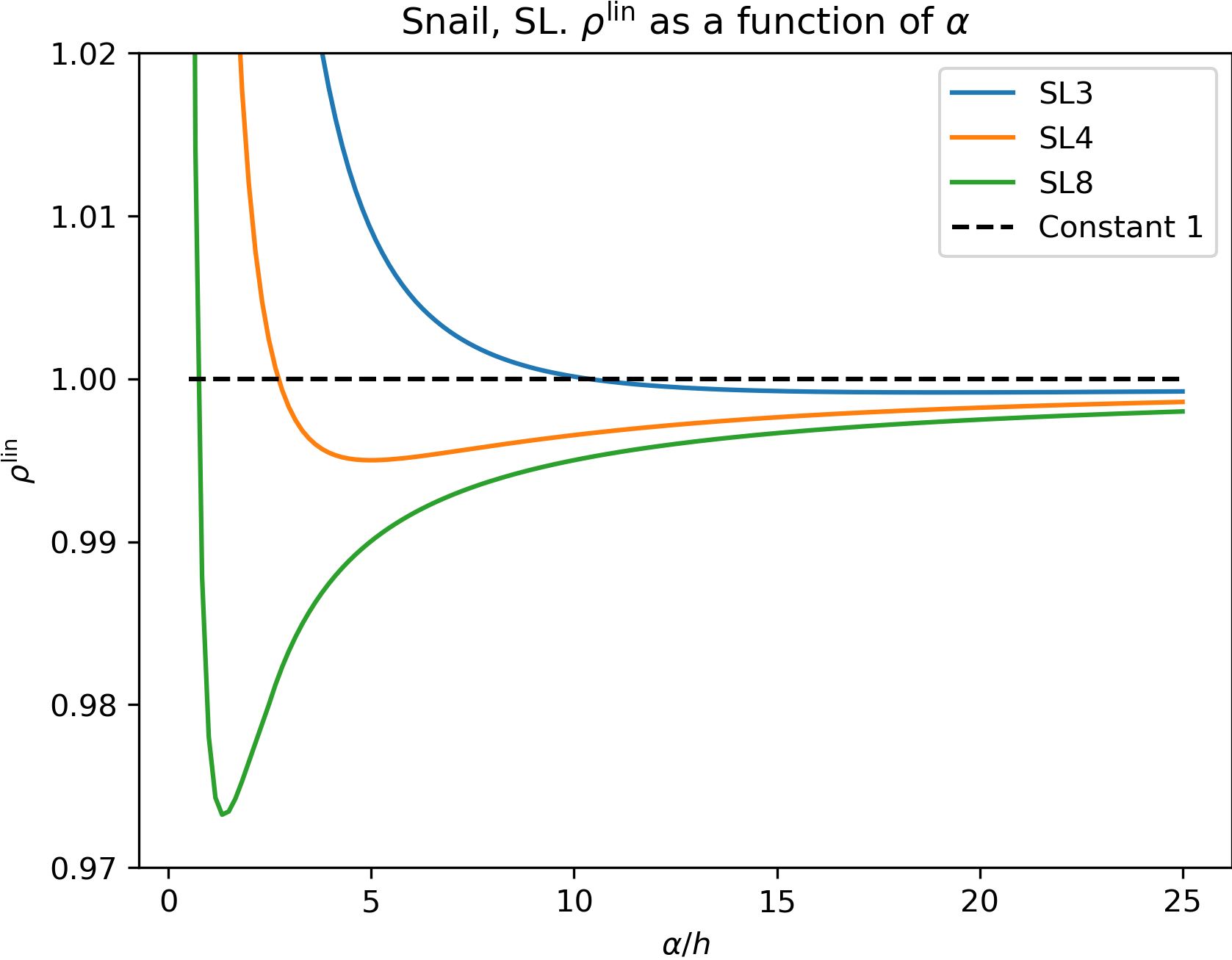}	
	\includegraphics[height=\h]{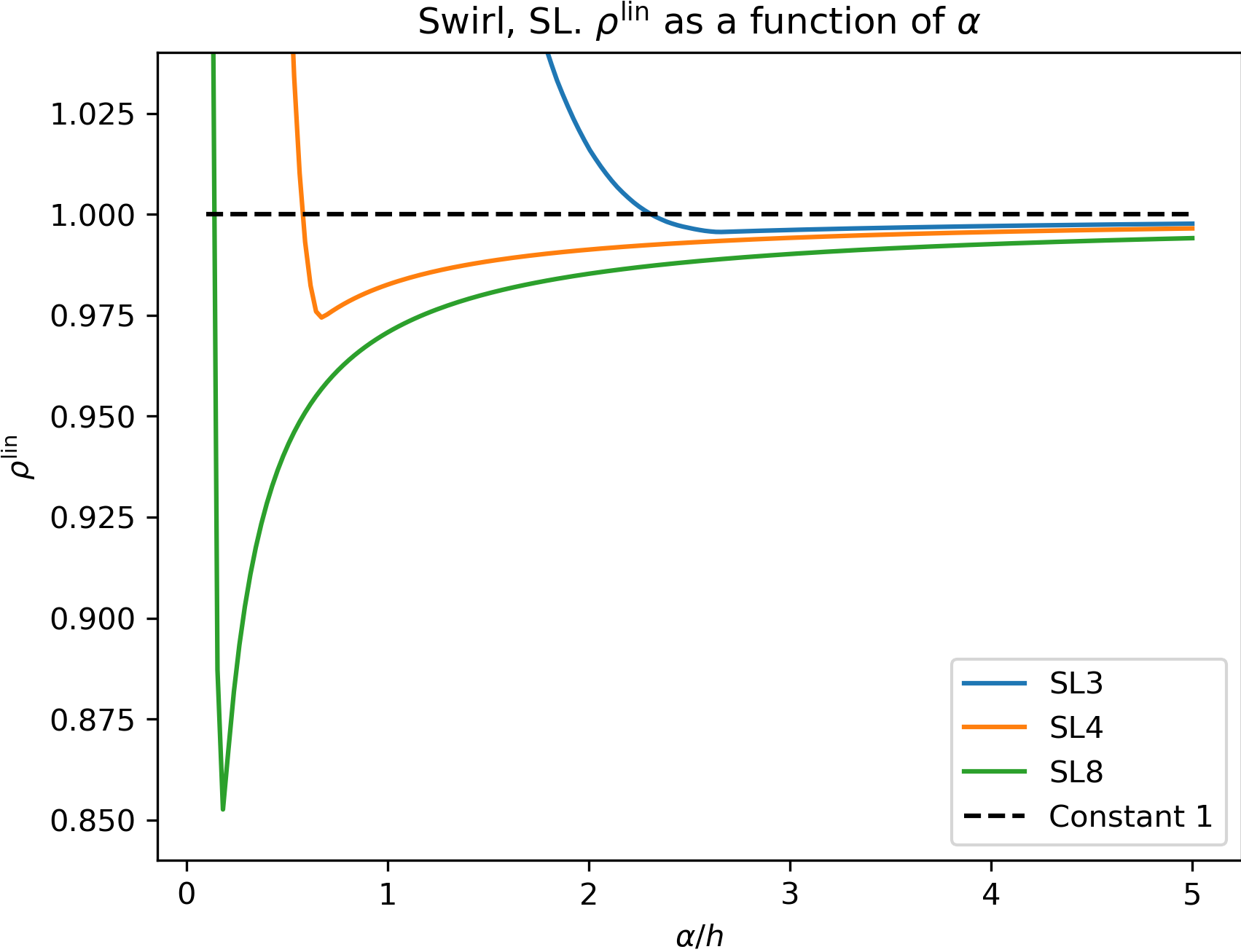} 
	\caption{Local contraction ratio near linear maps $\rho^\lin$ as a function of the \Kruzkov\ timescale $\alpha$, normalized by the gridscale $h$, in the different test cases. 
	In the maze case (bottom left), one has $\rho^\lin(\alpha)<1$ for all $\alpha>0$, since the upwind difference scheme is causal, see \cref{prop:causal_implies_subadd}. In contrast,  $\rho^\lin(\alpha) > 1$ when  $\alpha/h$ is too small in the other cases. In all cases, $\rho^\lin(\alpha) < 1$ when $\alpha/h$ is sufficiently large, and $\rho^\lin(\alpha) \to 1$ as $\alpha\to \infty$. 
	The best timescale $\alpha>0$ for the numerical efficiency of the ANB is often the one which minimizes $\rho^\lin(\alpha)$, see \cref{fig:iter_wrt_alpha}. %
	}
	\label{fig:rho_lin}
\end{figure}

In order to resolve this apparent contradiction, and to have usable estimates of the adequate \Kruzkov\ timescale and corresponding contraction ratio, 
we introduce the 
\emph{local contraction ratio near linear maps}, defined as follows:
\begin{align}
\label{eq:rho_lin}
	\rho^\lin(\alpha; x) &:= 
	\sup_{\cF_x^*(\eta) = 1} \rho^\loc\big(\alpha; \<\eta,\cdot\>, x\big), &
	\rho^\lin(\alpha) &:= 
	\max_{x \in \Int X} \rho^\lin(\alpha; x).
\end{align}
This is a specialization of the local contraction ratio $\rho^\loc(\alpha; u,x)$, see \eqref{eq:subadd_diff}, near linear maps $u = \<\eta,\cdot\> \in \bR^X$, $\eta\in \bR^d$, which satisfy the PDE of interest $1 = \cF^*_x(\nabla u(x)) = \cF_x^*(\eta)$ at the given point $x\in X$.
This quantity is easy to compute (or at least to approximate well) numerically in a preprocessing step, using the explicit expression of $\rho^\loc\big(\alpha; \<\eta,\cdot\>, x\big)$ derived in \cref{eq:rho_loc_LF,eq:rho_loc_SL,eq:rho_loc_UD} of \cref{sec:effective} for the schemes of interest.
The parameter $\alpha>0$ which minimizes $\rho^{\lin}(\alpha)$ often also approximately minimizes numerical cost (measured in evaluations of $\Lambda$) in our numerical experiments, compare \cref{fig:rho_lin,fig:iter_wrt_alpha}.

The impact of the remaining ANB parameters $\rho$ and $\tau$ on numerical cost was on the other hand quite modest in our experiments, see \cref{sec:num_secondary,fig:iter_wrt_rho,fig:iter_wrt_tau} for a discussion and illustration. Therefore, unless otherwise specified, we choose for simplicity and consistency
\begin{align}
	\label{eq:default_rho_tau}
	\rho &= 0.95, &
	\tau &= \alpha /2. 
\end{align}

\paragraph{Boundary conditions.}
In all our numerical examples, we numerically solve HJB or HJI PDEs $\cF_x^*(\nabla u(x)) = 1$ over a domain $x \in \Omega \sm \{x_0\}$, where $x_0 \in \Omega$ is referred to as the source point. We use outflow boundary conditions $u=\infty$ on $\partial \Omega$, and following standard practice, we set $u(x_0) = 0$ at the source point, which is always a discretization grid point. 

Alternatively, as suggested in the introduction, we could have set $u(x_0)=1$ so as to satisfy the \cref{assum:all} of positive boundary conditions, and then taken advantage of the additive invariance of the problem to recover the null-source solution. This trick was omitted in our experiments since it led to no observable difference on solution accuracy or computational cost. 

\begin{figure}
\def\w{15cm}
\includegraphics[width=\w]{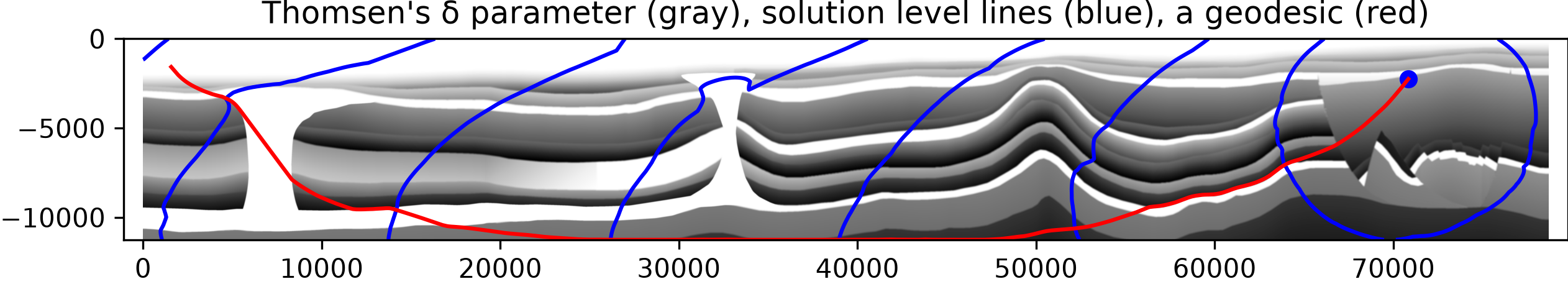}\\
\includegraphics[width=\w]{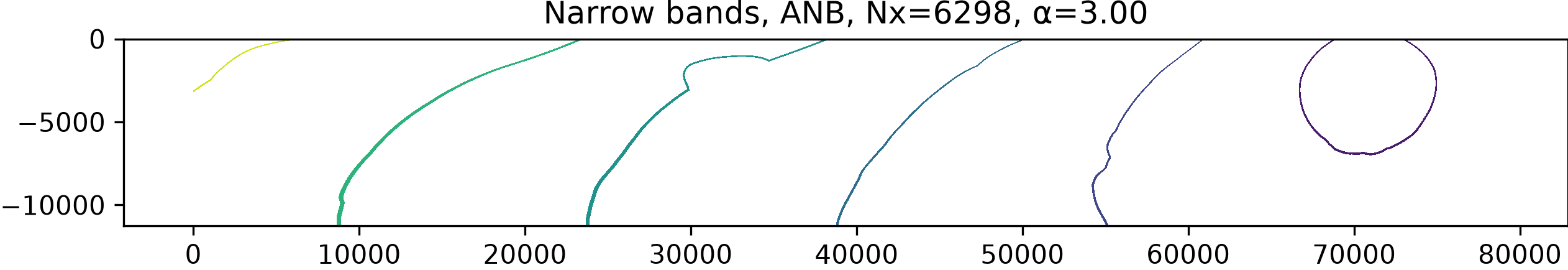}\\
\includegraphics[width=\w]{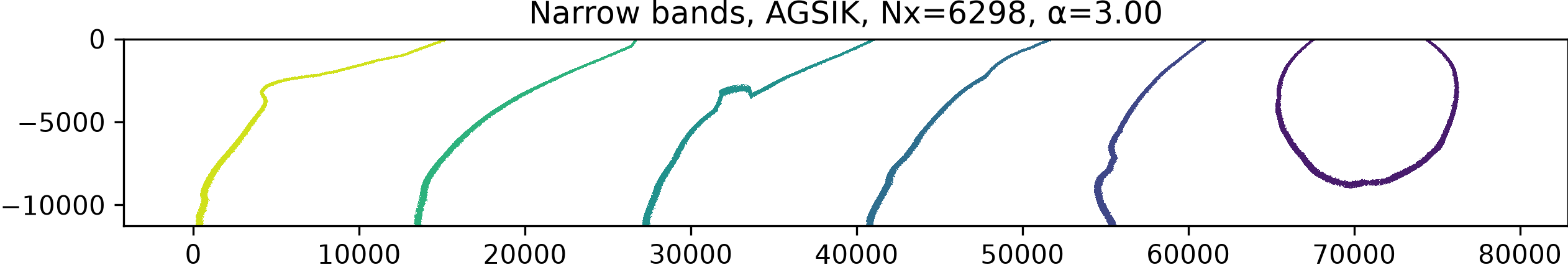}
\caption{Seismic test case, Lax-Friedrichs scheme, see \cref{subsec:seismic}. 
(Top) Visualization of Thomsen's geophysical parameter $\delta$ (gray), some level lines of the solution (blue), and a geodesic (red).
(Middle) Narrow bands produced in some iterations of the ANB$^*$. (Bottom) Narrow bands produced in some iterations of the AGSIK. Problem size $6298 \times 900$ (i.e.\ dataset subsampled by a factor two), 
\Kruzkov\ timescale $\alpha=3h$, other parameters defaulted \eqref{eq:default_tol}, \eqref{eq:default_rho_tau}.
}
\label{fig:seismic}
\end{figure}

\begin{figure}
	\def\h{3.9cm}
	\includegraphics[height=\h]{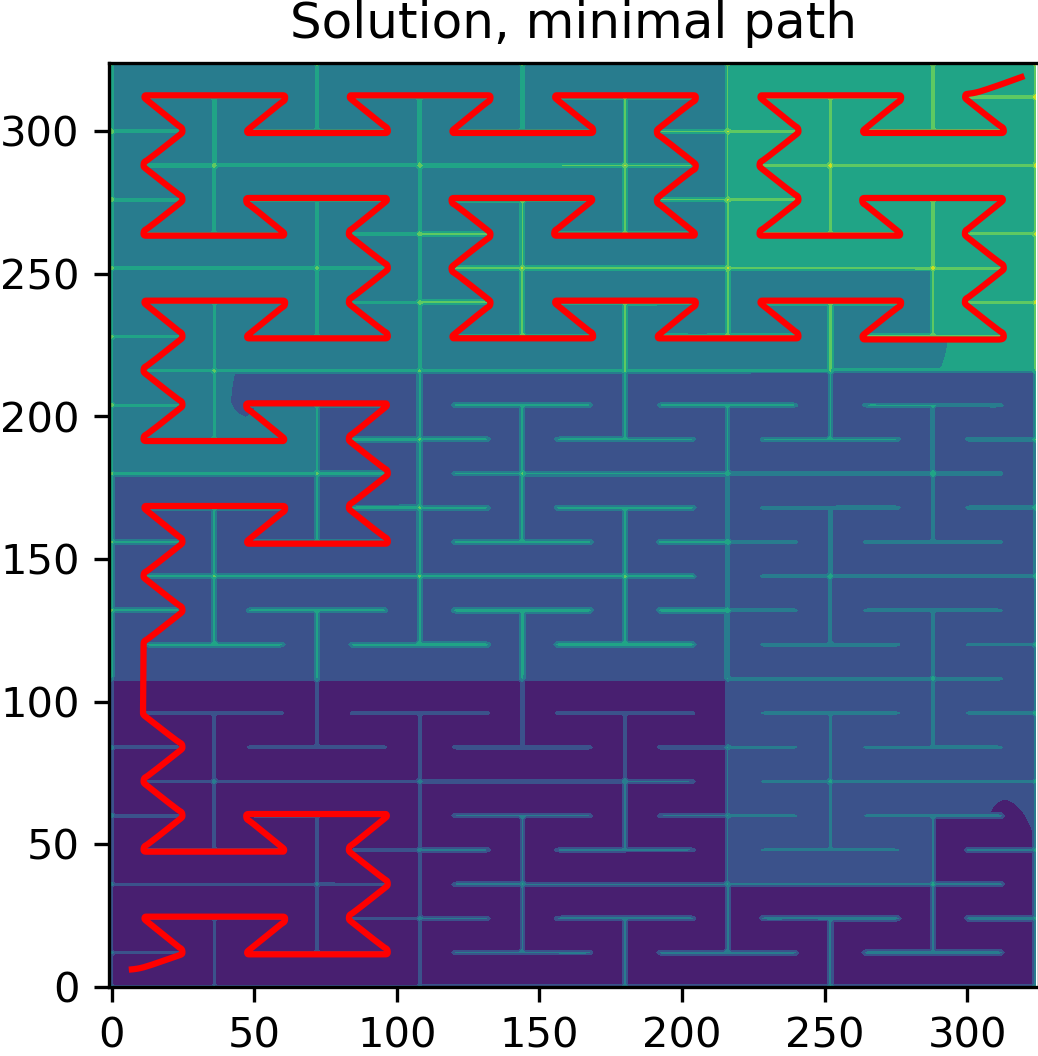}
	\includegraphics[height=\h]{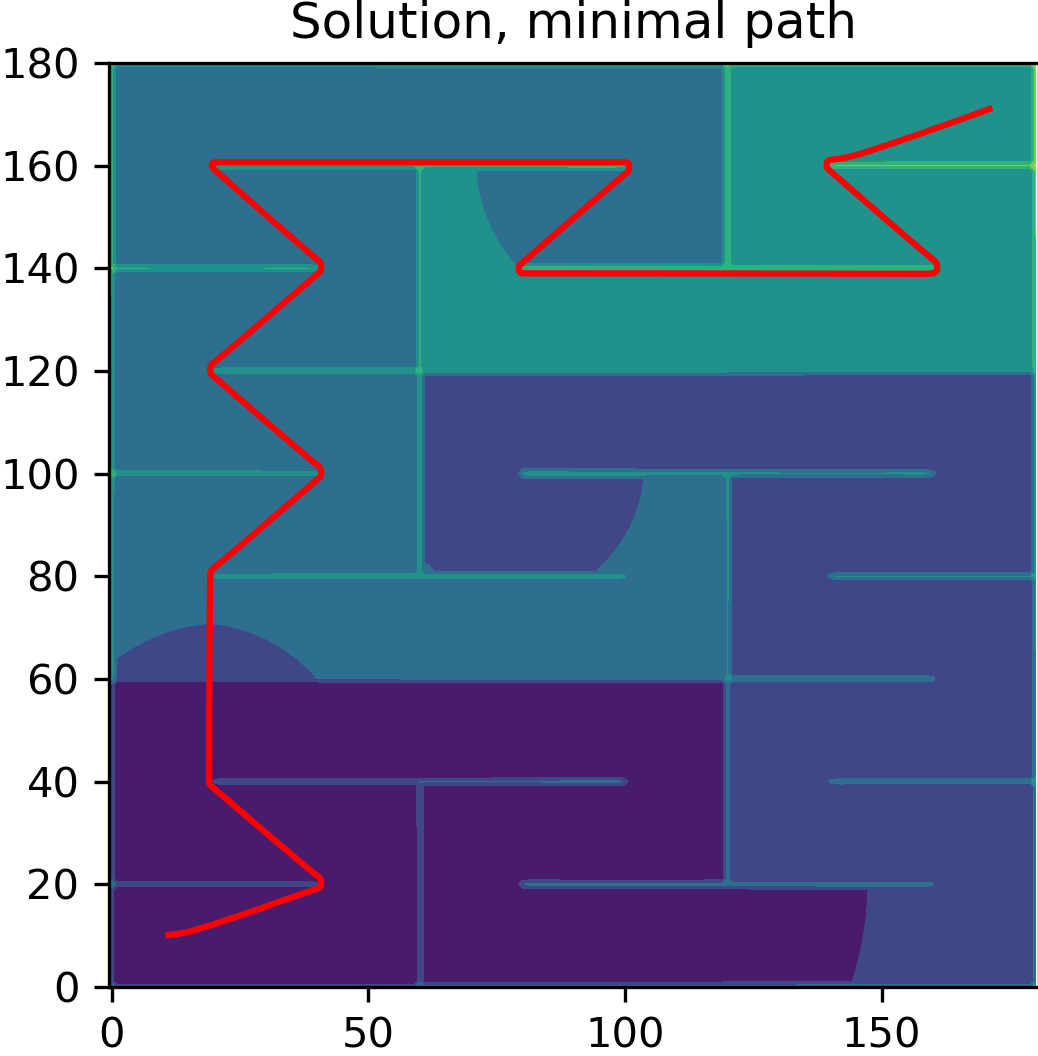}
	\includegraphics[height=\h]{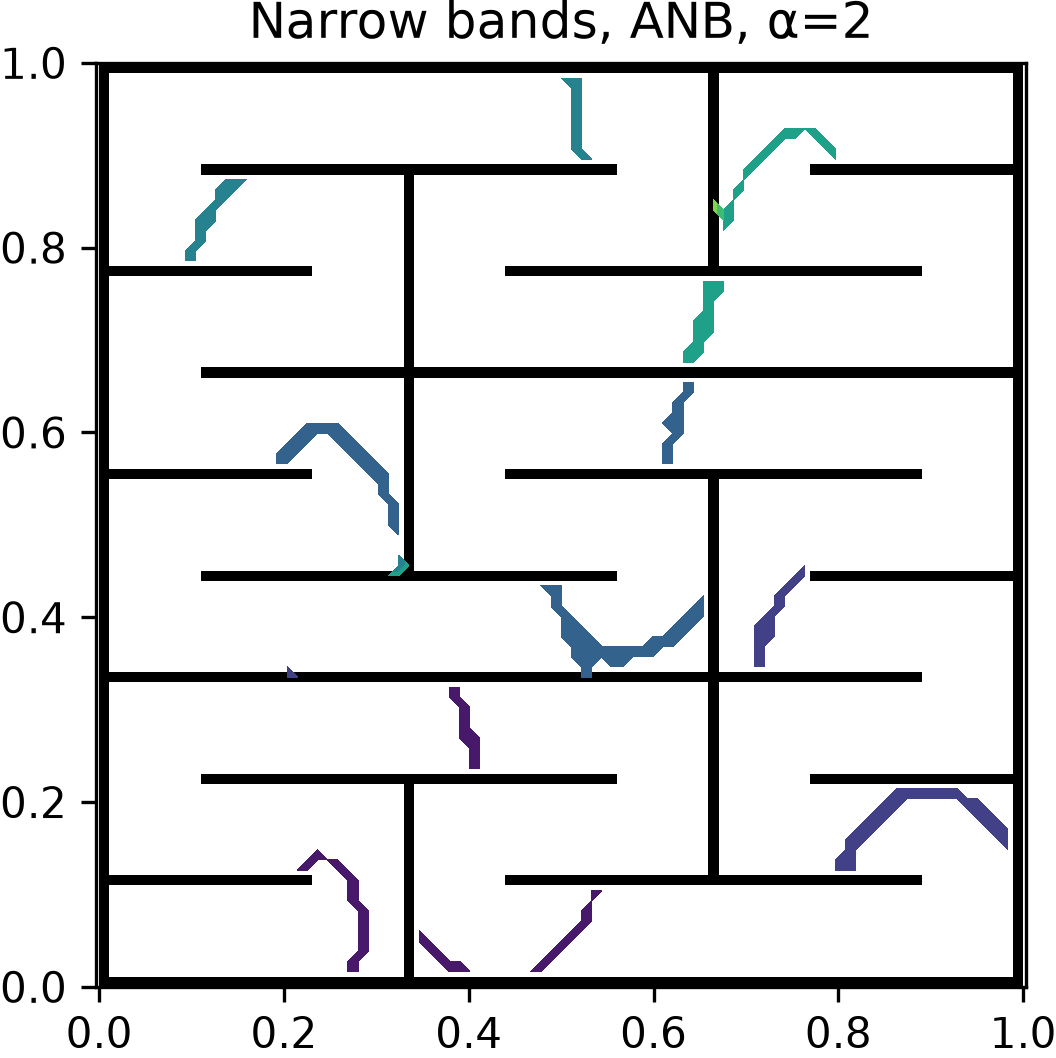}
	\includegraphics[height=\h]{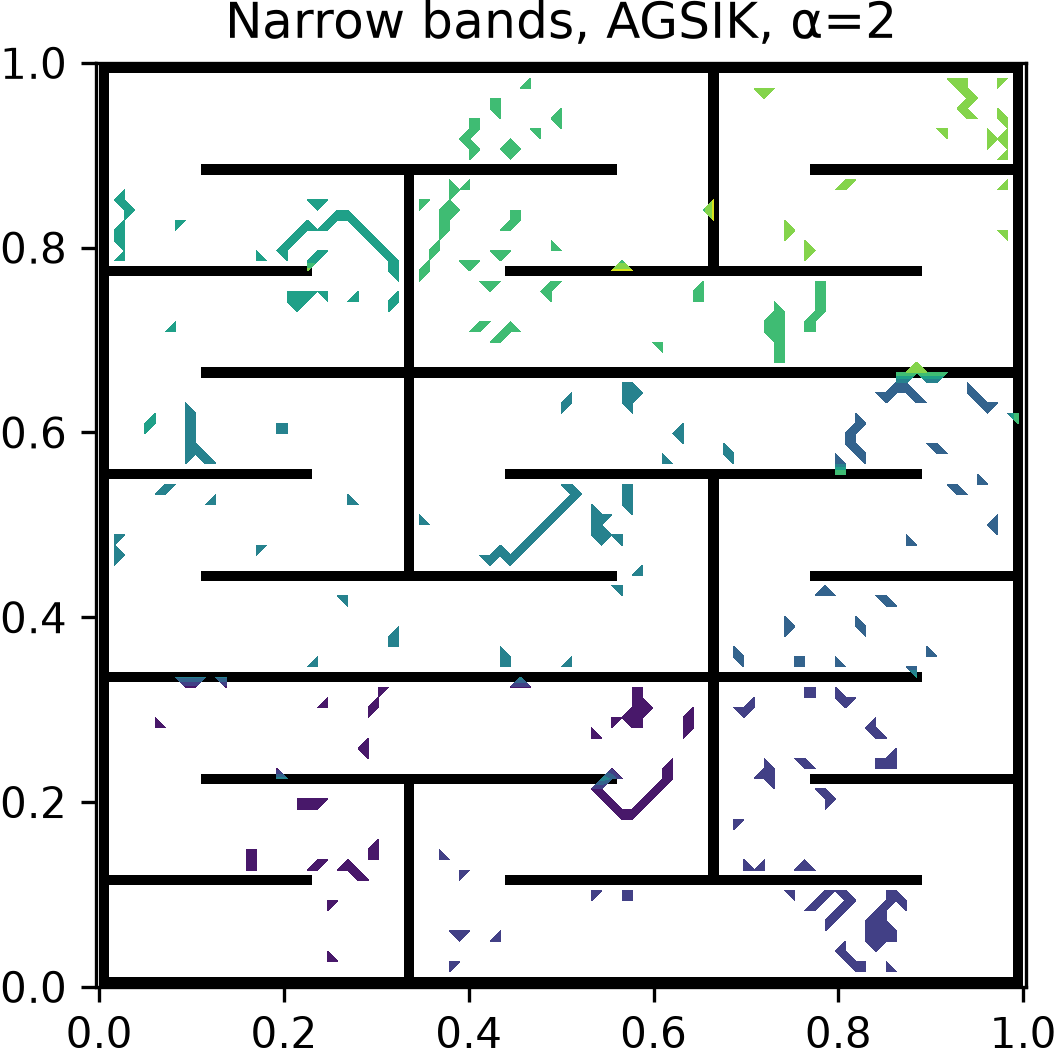}
	\caption{Maze with slightly permeable walls, upwind differences scheme (causal), see \cref{subsec:maze}.
	(i) Level sets of the solution, and a geodesic. (ii) Smaller problem instance, (iii) narrow bands of the ANB, and (iv) of the AGSIK. Problem size $135\times 135$, \Kruzkov\ timescale $\alpha=2h$ in this illustration, and other parameters defaulted \eqref{eq:default_tol}, \eqref{eq:default_rho_tau}.}
	\label{fig:maze}
\end{figure}

\paragraph{Visualization of the narrow bands.}
We illustrate on \cref{fig:pursuit,fig:seismic,fig:maze,fig:snail,fig:swirl} some \emph{level sets} $\{x\in \Omega \mid c_i \leq u(x) \leq c_{i+1}\}$ of numerical solutions to the test cases of interest, where the levels $c_1,\cdots$ are arbitrary and chosen by the visualisation software Matplotlib$^{\text{\textregistered}}$. 
We also display some narrow bands $B_n := \{n\tau - 2 \Delta_n \leq \hat u_n < n\tau\}$ produced at arbitrary steps $n_1,\cdots$ of the ANB method \cref{algo:anb}. 
Our intent is to illustrate that each $B_{n_i}$ visually ``looks like'' a level set corresponding to a thin range of values of the solution $u$, even though it is based on the rough approximation $\hat u_n$ - obtained by extending the sub-solution $\ol u_n$ via Disjktra's algorithm - known at this stage of the numerical resolution.

The narrow bands associated to the ANB$^*$ variant, see \cref{sec:AGSI}, are essentially indistinguishable from those of the ANB. We also illustrate some effective narrow bands $\tilde B_n := \{x \in X \mid \exp(\tfrac{-\Lambda u(x)}\alpha) \geq \exp(\frac{\ve-u(x)} \alpha) + 2^{-n}\}$, at arbitrary iterations $\tilde n_1,\cdots$ of the AGSIK method, see again \cref{sec:AGSI}. The AGSIK bands $\tilde B_n$ are more ``fuzzy'' than the ANB bands $B_n$, since they are based on the residuals of the numerical scheme, but they nevertheless  resemble level sets of $u$.
In contrast, the computational bands of the FS and FSd methods are fixed in advance, as vertically, horizontally or diagonally aligned sets of points; this is effective with the SL4 or SL8 schemes, but not with the LF or SL3 scheme, or in a complex maze structure. The queues underlying the AGSI and FIM may also be regarded as computational bands, built in view of simple front tracking heuristics; while they work well on simple test cases (not shown), they are less relevant in inhomogeneous or anisotropic settings, and they typically become domain-wide hence useless if the scheme is strongly diffusive. 

\begin{table}[!t]

\newcommand\tr[1]{\textcolor{red}{#1}}
\hspace{-0.5cm}
\begin{tabular}{cccccccccc}
&  ANB &  ANB$^*$ &  AGSIK &  NB &  AGSI &  FIM &  FS &  FSd &  GI\\
1D, LF & \tr{148+143} & \tr{97+90} & \tr{112} & \tr{113} & 1242 & 2481 & 3182 & NA & 6868 \\
Pursuit, LF & 163+327 & \tr{104+272} & \tr{120} & \tr{116} & 633 & 1264 & 3484 & 2164 & 4726 \\
Seismic(subs), LF & \tr{41+39} & \tr{41+49} & \tr{61} & \tr{60} & 1509 & 3017 & * & 4864 & * \\
Seismic, LF & \tr{49+50} & \tr{41+54} & \tr{61} & \tr{60} & * & * & * & * & * \\
Maze, UD & \tr{5.1+6.6} & \tr{4.8+10} & \tr{15} & \tr{15} & \tr{48} & \tr{80} & 1156 & 200 & 2824 \\ %
Snail, 2SL3 & 431+2269 & \tr{29+82} & \tr{30} & \tr{60} & 1412 & 2776 & * & 2352 & * \\
Snail, 2SL4 & 264+303 & \tr{18+44} & \tr{30} & \tr{63} & 1544 & 3047 & * & \tr{140} & * \\
Snail, 2SL8 & 174+174 & \tr{20+44} & \tr{34} & \tr{31} & 704 & 930 & \tr{36} & \tr{44} & * \\
Swirl, 2SL3 & 1087+1577 & \tr{25+58} & \tr{28} & \tr{42} & 767 & 1507 & 4016 & 2216 & 7316 \\
Swirl, 2SL4 & 500+907 & \tr{29+53} & \tr{34} & \tr{34} & 810 & 1601 & 2644 & \tr{16} & 4959 \\
Swirl, 2SL8 & 124+263 & \tr{49+98} & \tr{61} & \tr{63} & 836 & 468 & \tr{20} & \tr{16} & * 
\end{tabular}
\caption{Average number of evaluations of $\Lambda$ per point of the various methods (lower is better), in the various test cases, see \cref{subsec:1D,subsec:pursuit,subsec:seismic,subsec:maze,subsec:snail,subsec:zermelo}. 
For the ANB and ANB$^*$, the number of evaluations of $\Lambda$ is reported as $a+b$ where $a$ (resp.\ $b$) counts evaluations on the super-solution (resp.\ sub-solution). 
NA : not applicable. $*$ : computations exceeding one hour were interrupted.
Some ``reasonably low'' numbers highlighted in red by hand. The AGSI and FIM work well in the maze test case, the FS and FSd work well with SL4 and SL8 schemes, and the narrow band type methods ANB$^*$, AGSIK, NB provide consistent performance across the board including with the diffusive SL3 and LF schemes. When semi-Lagrangian schemes are used, the ANB should clearly be discarded in favor of the ANB$^*$, which we recall is a basic modification avoiding many redundant updates, see \cref{sec:AGSI}.
\Kruzkov\ timescale $\alpha = 3h$ (1D, seismic), $\alpha=5h$ (pursuit), $\alpha=h$ (maze), $\alpha=2h$ (snail, swirl, except for the NB method see \cref{subsec:snail,subsec:zermelo}). Other method parameters defaulted \eqref{eq:default_tol}, \eqref{eq:default_rho_tau}. Problem size 2000 (1d), $2001 \times 2001$ (pursuit, snail, swirl), $6298 \times 900$ (seismic subsampled), $12596\times 1801$ (seismic), $135 \times 135$ (maze). See \cref{table:elapsed} for the corresponding computation times.
}
\label{table:niter}
\end{table}

\subsection{One dimensional toy experiment, Lax-Friedrichs scheme}
\label{subsec:1D}

HJB PDEs involving complex anelliptic Hamiltonians, or HJI PDEs featuring non-convex Hamiltonians, are often discretized using the Lax-Friedrichs scheme, see \cref{subsec:pursuit,subsec:seismic} respectively. 
This toy experiment shows that the diffusivity of the LF scheme is the reason why these problems are costly to solve,
rather than the anellipticity or non-convexity of the Hamiltonian.

The standard eikonal equation $|\nabla u(x)|=1$ is numerically solved on the one-dimensional segment $x\in [0,1]$, with the Dirichlet boundary condition $u(0)=0$ and the outflow boundary condition $u(1)=\infty$, using the LF scheme from \cref{sec:effective} with parameters $c_1 = 1/3$ and $C_0=1$. 

\begin{figure}[!t]
\centering
	\def\h{4cm}
	\includegraphics[height=\h]{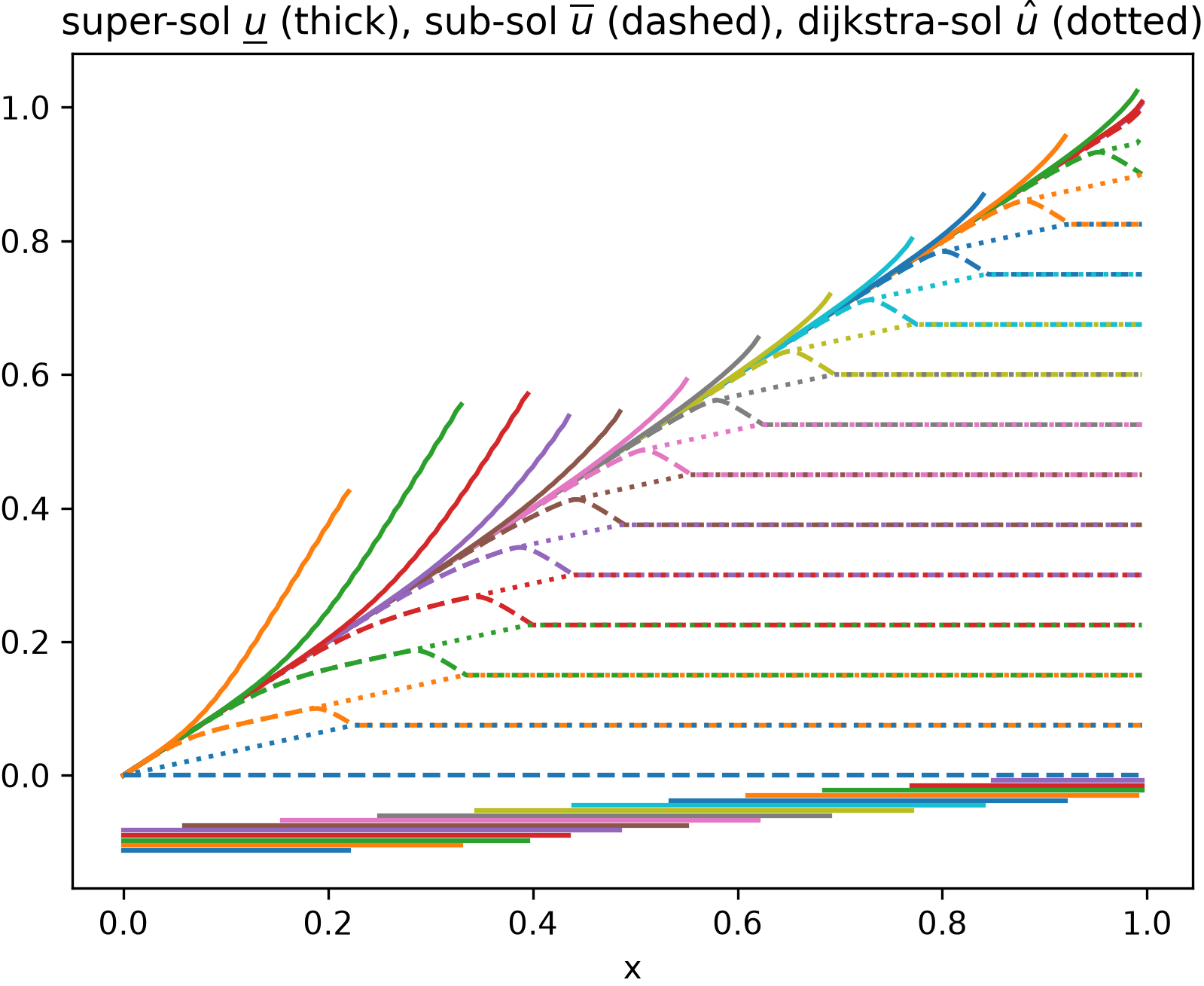}
	\includegraphics[height=\h]{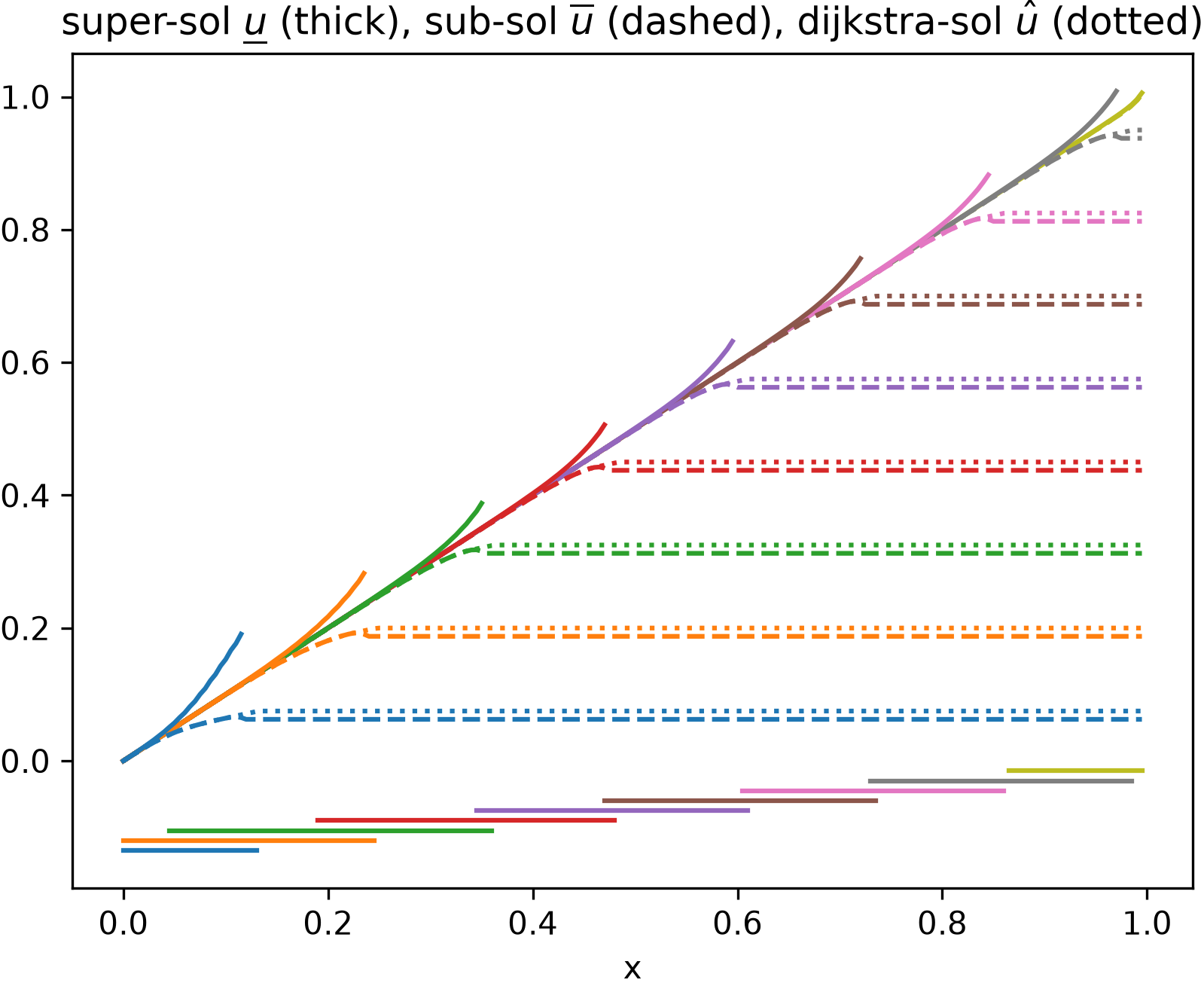}
	\includegraphics[height=\h]{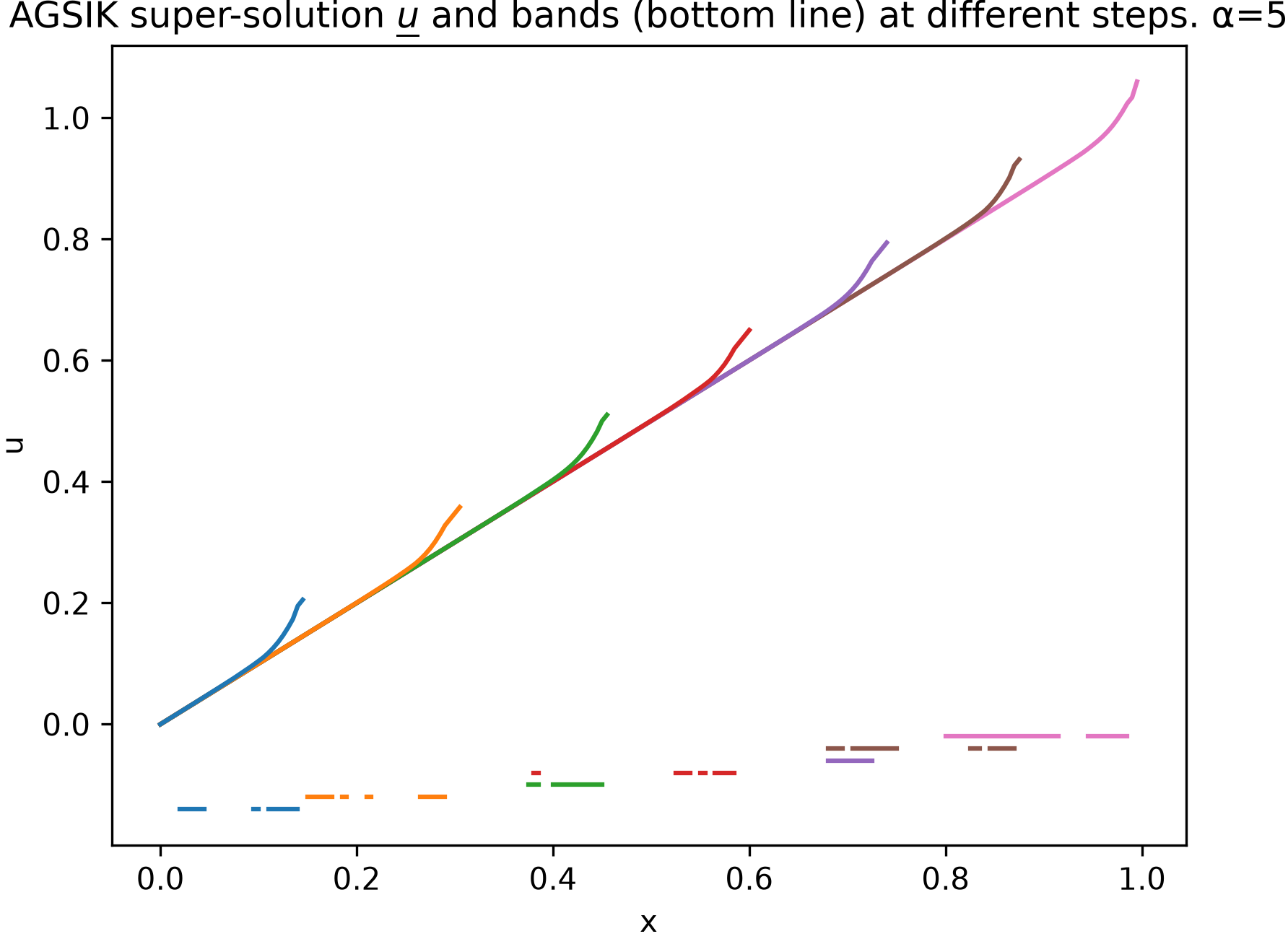}
	\caption{
	One dimensional problem, Lax-Friedrichs scheme, see \cref{subsec:1D}. (Left) Super-solution (solid), sub-solution (dashed), Dijkstra extension of the sub-solution (dotted), and active band (solid at the bottom), in \emph{all} iterations of the ANB with parameter $\tau = 3\alpha$. (Center) Likewise in \emph{some} iterations of the ANB with parameter $\tau = \alpha/2$, following our defaults \eqref{eq:default_rho_tau}.
	(Right) Super-solution (solid), and active band (at the bottom), in several iterations of the AGSIK.
	Problem size 200, \Kruzkov\ scale parameter $\alpha=5h$ in this illustration, $\rho=0.95$, tolerance \eqref{eq:default_tol}.
	}
	\label{fig:1d}
\end{figure}

The intermediate numerical solutions, and corresponding narrow bands, of the ANB and AGSIK are illustrated on \cref{fig:1d}. For simplicity, we use the timescale $\alpha = 3h$ ($\alpha=5h$ in \cref{fig:1d}), and default the other parameters \eqref{eq:default_tol}, \eqref{eq:default_rho_tau}; 
one can show \cite[\S4.1]{mansour2024solving} that the contraction ratio near linear maps \eqref{eq:rho_lin} reaches the minimal value $\rho^{\lin}(\alpha) = \frac{2\sqrt 2} 3 \approx 0.94$ when $\alpha = \frac 2 {\ln 2} h \approx 2.89 h$.
As illustrated on \cref{table:niter,table:elapsed}, even such a simple toy problem puts the popular AGSI, FIM, and FS methods in difficulty: they require hundreds or thousands of updates per point, because the underlying front tracking heuristics cannot cope with the scheme strong diffusivity. 

Of course, it would be trivial to implement a better discretization scheme for such a one-dimensional problem, e.g.\ the causal scheme $\Lambda u(x) = h + \min \{u(x-h),u(x+h)\}$. We are thus admittedly turning a simple problem into a complicated one, for illustration purposes.

\subsection{Pursuit evasion, Lax-Friedrichs scheme}
\label{subsec:pursuit}

We numerically solve an HJI PDE encoding a pursuit evasion game where agent $A$ with speed $v_A(\cos \theta,\sin \theta)$, $\theta \in [-\pi,\pi] \sm ]-\pi/4,\pi/4[$, tries to catch agent $B$ with speed $v_B(\cos \phi,\sin \phi)$, $\phi \in [-\pi,\pi]$, where $v_A=2$ and $v_B = 1$ following \cite{cacace2011local}.
The PDE reads
$\cF^*(\nabla u(x)) = 1$ where
\begin{equation*}
	\cF^*(\eta) := \sup_{\theta \in [-\pi,-\pi/4] \cup [\pi/4,\pi]} \<\eta,v_A (\cos \theta,\sin \theta)\> - \sup_{\phi \in [-\pi,\pi]} \<\eta, v_B (\cos \phi,\sin \phi)\>,
\end{equation*}
for any co-vector $\eta \in \bR^2$, is a \emph{non-convex} Hamiltonian. 
The PDE is posed on the square domain $\Omega = ]0,1[^2$, a point source is introduced in the center,  and outflow boundary conditions are used on $\partial \Omega$. 
The Lax-Friedrichs scheme parameters $C_0$ and $c_1$ are adjusted numerically so as to satisfy $C_0^{-1} \|\xi\|_\infty \leq \cF^*(\xi)$ and $\|\nabla \cF^*(\xi)\|_\infty \leq 0.95\, c_1^{-1}$, for all directions $\xi = (\cos (k \pi/20),\sin (k \pi/20))$, $k=1,\cdots,40$. See \cite[\S4.2]{mansour2024solving} for more discussion.

\begin{figure}
\centering
	\def\h{3.9cm}
	\includegraphics[height=\h]{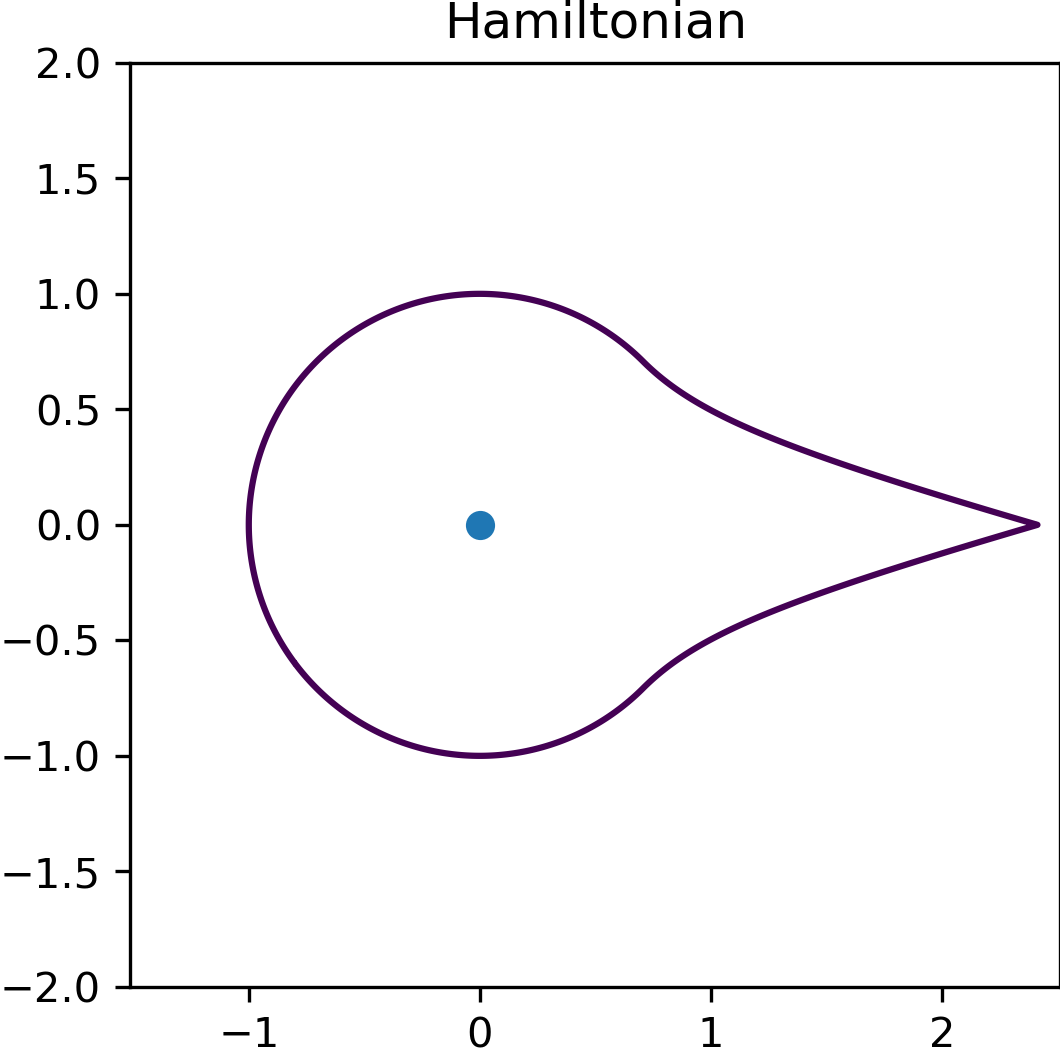}
	\includegraphics[height=\h]{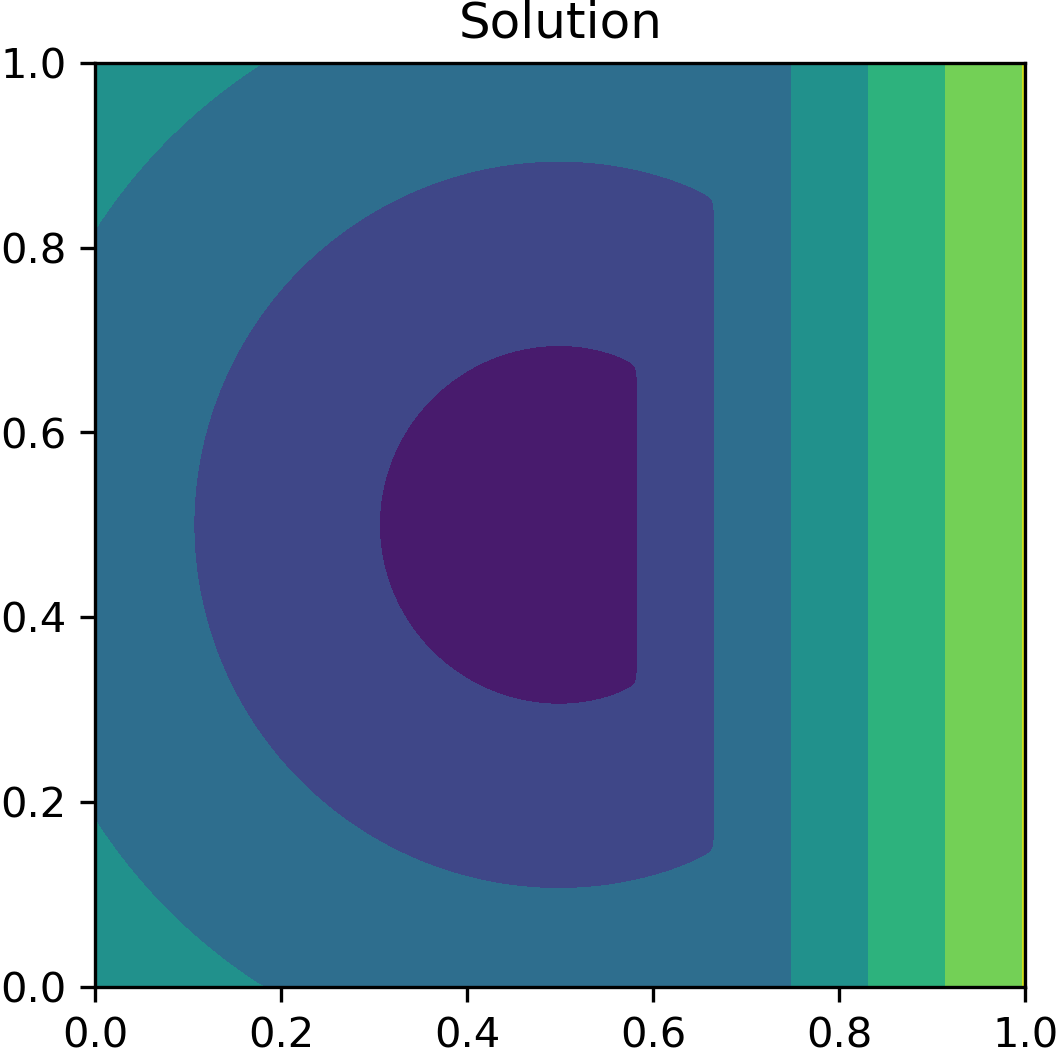}
	\includegraphics[height=\h]{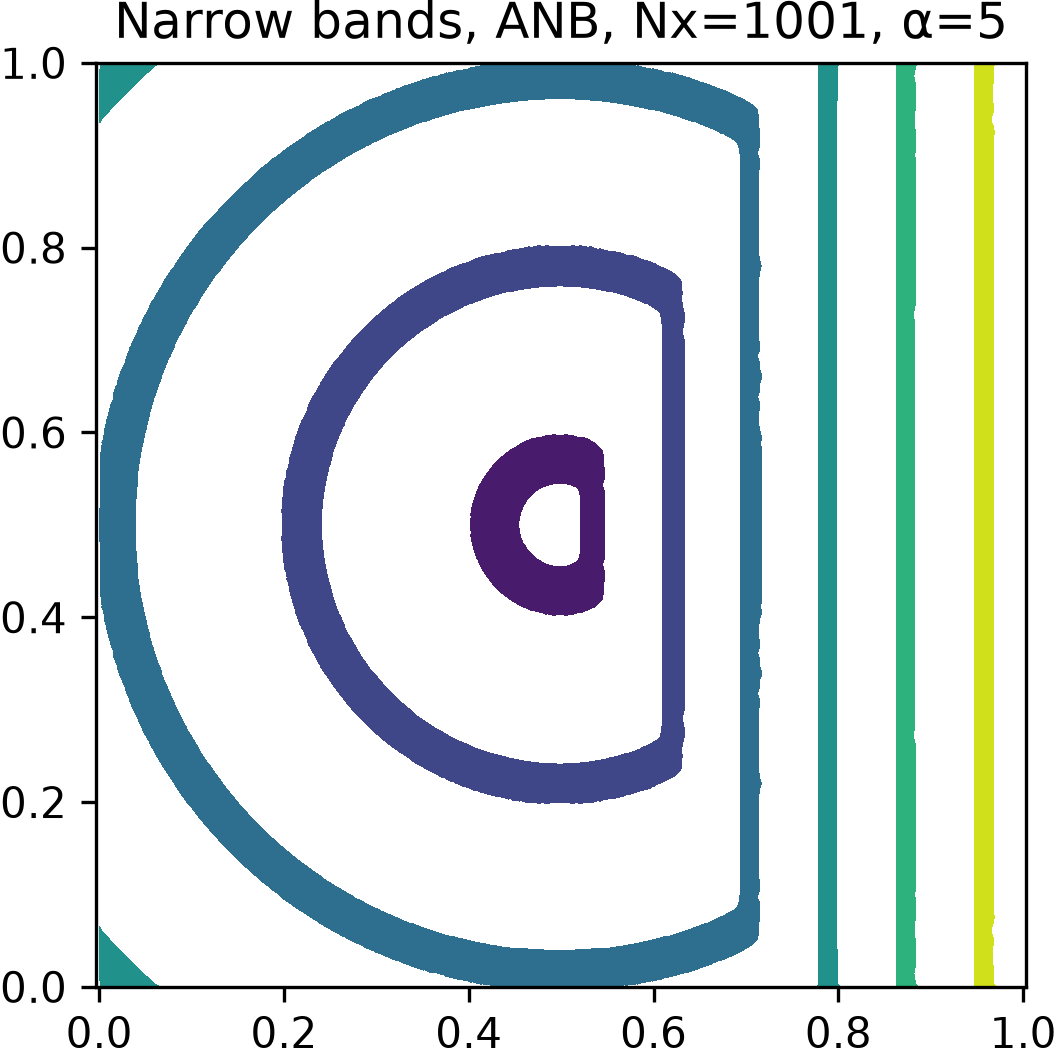}
	\includegraphics[height=\h]{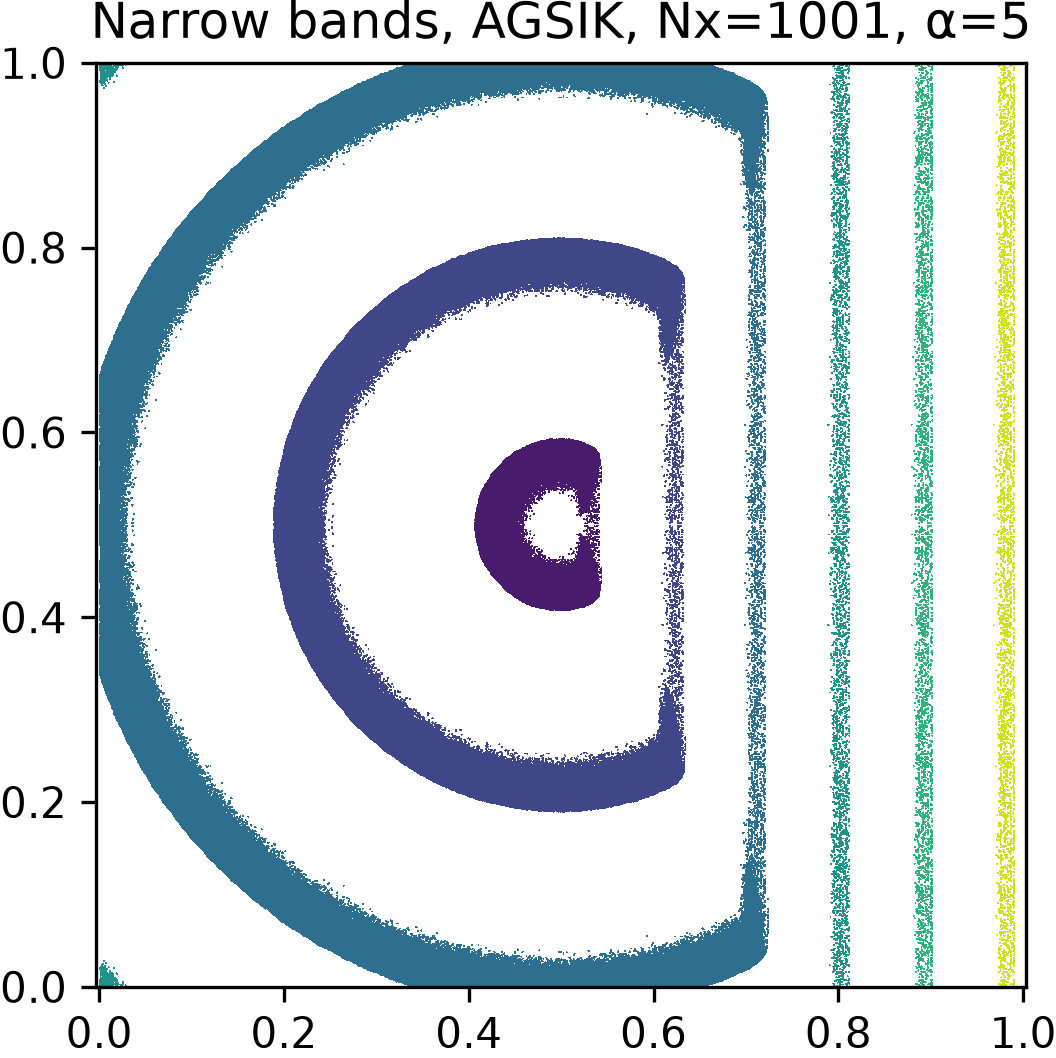}
	\caption{Pursuit evasion test case, Lax-Friedrichs scheme, see \cref{subsec:pursuit}. (i) Level line $\{\eta\in \bR^2 \mid \cF^*(\eta) = 1\}$ of the non-convex Hamiltonian, (ii) Level sets of the solution. (iii) Some narrow bands produced by the ANB, and (iv) by the AGSIK.
	 Problem size $1001 \times 1001$, \Kruzkov\ timescale $\alpha=5h$ and other parameters defaulted \eqref{eq:default_tol}, \eqref{eq:default_rho_tau}.
	}
	\label{fig:pursuit}
\end{figure}

For simplicity, we set $\alpha=5h$ and default the other parameters \eqref{eq:default_tol}, \eqref{eq:default_rho_tau}. 
(The contraction ratio near linear maps takes the minimal value $\rho^\lin(\alpha) \approx 0.98$ for $\alpha \approx 8.37 h$.)
The ANB, ANB$^*$, and AGSIK produce nice looking narrow bands, see \cref{fig:pursuit}. The classical AGSI, FIM, FS and FSd methods use many more iterations per point, since they do not track the propagation front effectively due to the scheme diffusivity, see \cref{table:niter}. 

\subsection{Anelliptic anisotropy, Lax-Friedrichs scheme}
\label{subsec:seismic}
We numerically solve the HJB PDE characterizing the first arrival time of the seismic pressure wave in the complex geological medium described by the BP 2007 TTI two-dimensional test case\footnote{Created by Hemang Shah and provided courtesy of BP Exploration Operation Company Limited.}, as in \cite[\S 4.3]{mansour2024solving}. The PDE reads $\cF^*_x(\nabla u(x)) = 1$ and involves the convex Hamiltonian
\begin{align}
\label{eq:HookeMetric}
	\cF_x^*(\eta) &:= \sqrt{\lambda_{\max}(c(x; \eta))}, &
	\text{ where } c(x; \eta)_{jk} &:= \sum_{i,l} c_{ijkl}(x) \eta_i \eta_l.
\end{align}
defined as the square root of the largest eigenvalue $\lambda_{\max}$ of the symmetric positive semi-definite matrix $c(x; \eta) \in \cS_d^{+}$ obtained by contracting the co-vector $\eta = \nabla u(x)$ with the reduced Hooke tensor $(c_{ijkl}(x))_{ijkl=1\cdots d}$. 
A point source is introduced in the domain interior, and outflow boundary conditions are used on the domain boundary.
The Lax-Friedrichs scheme parameters $C_0(x)$ and $c_1(x)$ are adjusted numerically at each discretization point $x$ so as to satisfy $C_0(x)^{-1} \|\xi\|_\infty \leq \cF^*_x(\xi)$ and $\|\nabla \cF^*(\xi)\|_\infty \leq 0.95\, c_1(x)^{-1}$, for all directions $\xi = (\cos (k \pi/20),\allowbreak \sin (k \pi/20))$, $k=1,\cdots,40$. See \cite[\S4.2]{mansour2024solving} and \cref{sec:effective} for more discussion.

For simplicity, we set $\alpha = 3h$ and default the other parameters \eqref{eq:default_tol}, \eqref{eq:default_rho_tau}. 
(The contraction ratio near linear maps takes the minimal value $\rho^\lin (\alpha) \approx 0.82$ for $\alpha \approx 1.58 h$, see \cref{fig:rho_lin}.)
The classical AGSI, FIM, FS, FSd methods suffer from the scheme diffusivity and the large problem dimension $12596\times 1801$, whereas the ANB, ANB$^*$, AGSIK and NB methods shine and appear to benefit from the rather favorable contraction ratio near linear maps. 

The first arrival time of the \emph{shear} wave propagation mode could in principle be numerically computed  similarly by solving the HJI PDE associated to the \emph{non-convex} Hamiltonian obtained by selecting the second eigenvalue $\lambda_{\min}$ in \eqref{eq:HookeMetric}.
Unfortunately, the chosen dataset describes the Hooke tensor via Thomsen's coefficients $(\delta,\ve,\theta,V_P)$, but omits $V_S$, in such way that $\lambda_{\min}$ cannot be estimated reliably, see \cite[Remark 4.7]{mansour2024solving}.

\subsection{Maze, upwind-differences (causal scheme)}
\label{subsec:maze}

We numerically solve an isotropic eikonal equation on a two-dimensional maze structure, featuring \emph{slightly permeable} walls, following \cite{jeong2008fim} and \cite[\S4.4]{mansour2024solving}. 
The wall's slight permeability is introduced purposely so as to increase the difficulty of front tracking \cite{jeong2008fim}. 
The PDE is posed on the domain $x\in ]0,1[^2$, and reads
\begin{align*}
	|\nabla u(x)| &= c(x), &
	c(x) &:= 
	\begin{cases}
		c^h_n &\text{if } d_\infty(x,\Gamma_n) \leq h/2,\\
		1 &\text{else}
	\end{cases} &
	d_\infty(x,\Gamma) &:= \min_{y \in \Gamma} \|x-y\|_\infty.
\end{align*}
We denoted 
by $\Gamma_n\subset [0,1]^2$ the $n$-th iteration of a Peano curve joining the points $(i+\frac 1 2,j+\frac 1 2)/3^n$, $i,j\in \{0,\cdots, 3^n-1\}$. The domain is discretized on a Cartesian grid with $N^2$ points, where $N = w 3^n$, thus with grid scale $h:=1/N$. The thickness $2(h/2) = h$ of the slightly permeable walls is exactly one pixel, and the cost $c_n^h := \frac w 4 9^n = \frac{3^n}{4 h}$ is chosen in such way that the optimal path joining the two ends of the Peano curve ``cheats'' exactly once by going through them. 

We choose $n=3$ and $w=5$, resulting in a small $135\times 135$ test case, with extremely strong speed contrast $c_n^h \approx 911$. A third, vertical dimension is added in \cite{jeong2008fim,mansour2024solving}, but it does not change the qualitative behavior of the model, and it makes visualization more complicated.

For simplicity, we set $\alpha=h$ and default the other parameters \eqref{eq:default_tol}, \eqref{eq:default_rho_tau}. The contraction ratio near linear maps $\rho^{\lin}(\alpha) \to 0$ as $\alpha\to 0$, because the scheme is causal, but there would be no point in choosing $\alpha \ll h$, since it increases the number $n_{\max}$ of times the narrow band is rebuilt, and the associated costs. 
The ANB and ANB$^*$ work well in this test case, even though the computational time is dominated by the narrow band reconstructions rather than the iteration of $\Lambda$. 
In contrast, the NB method needs a modification - an extension of the narrow band by one-pixel in all directions \cite[\S4.4]{mansour2024solving} - otherwise it cannot find the shortcut through the permeable walls.
We use a straightforward implementation of the AGSIK, based on a priority queue, rather than the bin-based optimized implementation suggested in \cref{sec:AGSI}; indeed, the latter limits the number of simultaneously active bins to $8$ for memory usage considerations, whereas $\approx c_n^h$ bins would be needed here.
This complex maze structure is designed to defeat the FS and FSd methods, which indeed need hundreds of iterations to converge, see \cref{table:niter}.

\subsection{Strongly anisotropic Riemannian test case, semi-Lagrangian scheme}
\label{subsec:snail}

\begin{figure}
	\def\h{3.8cm}
	\includegraphics[height=\h]{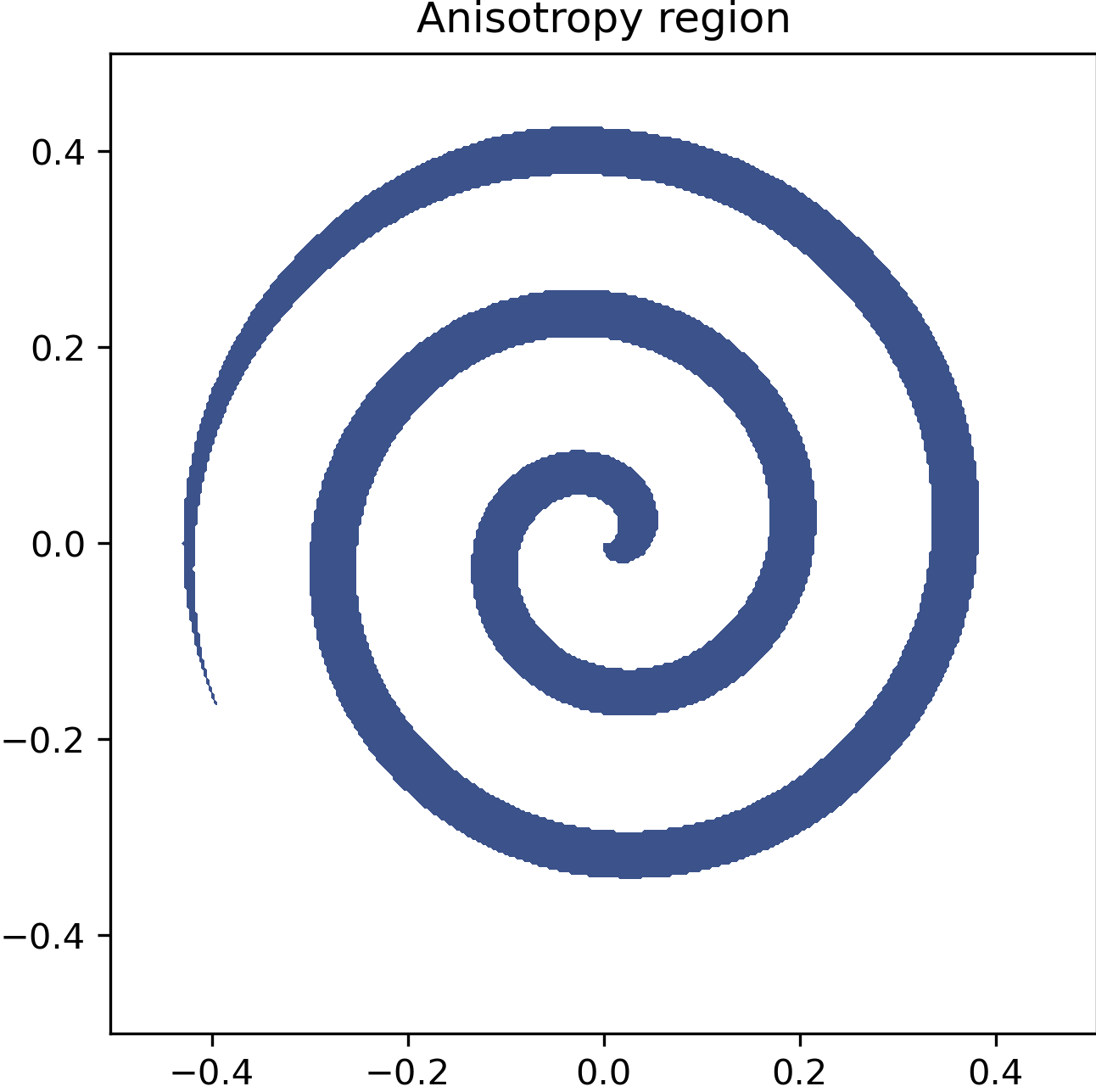}	
	\includegraphics[height=\h]{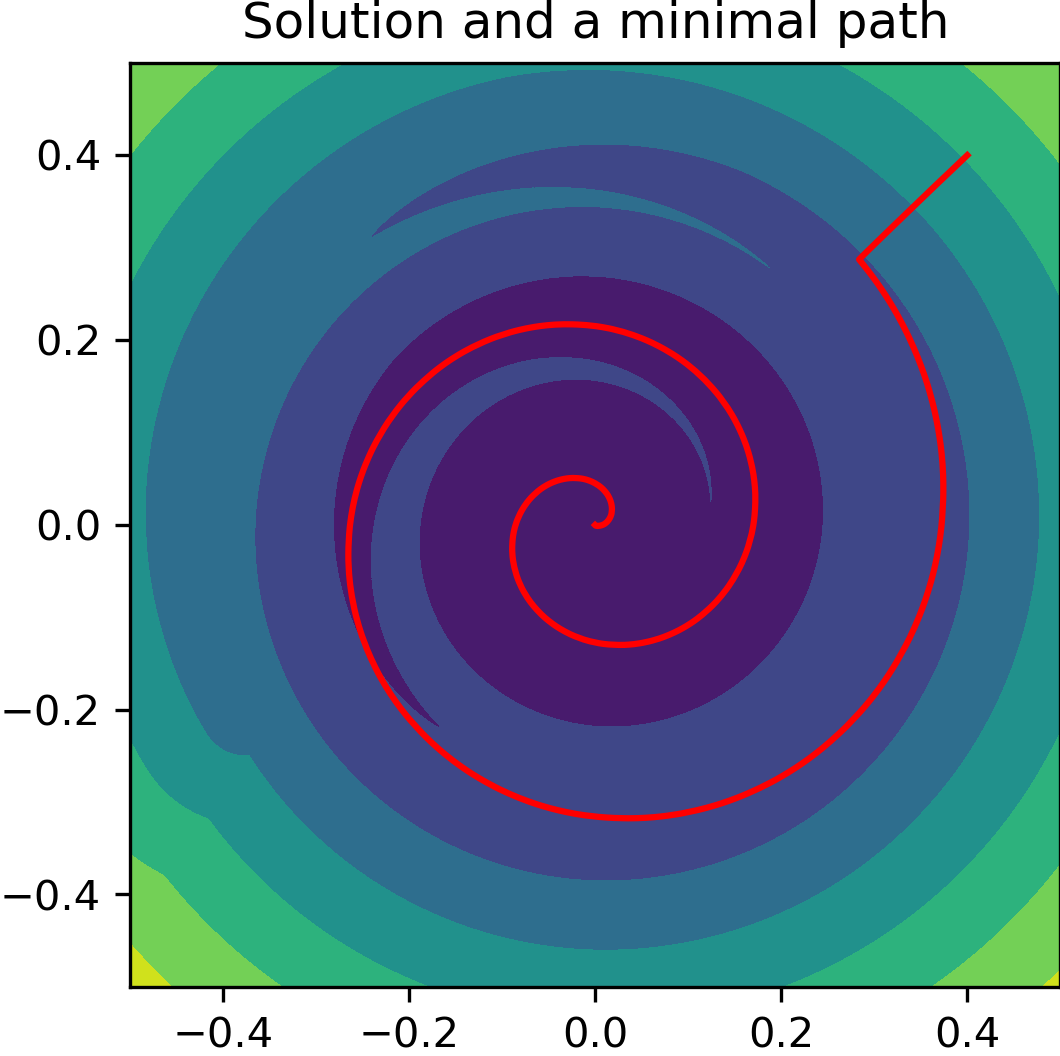}
	\includegraphics[height=\h]{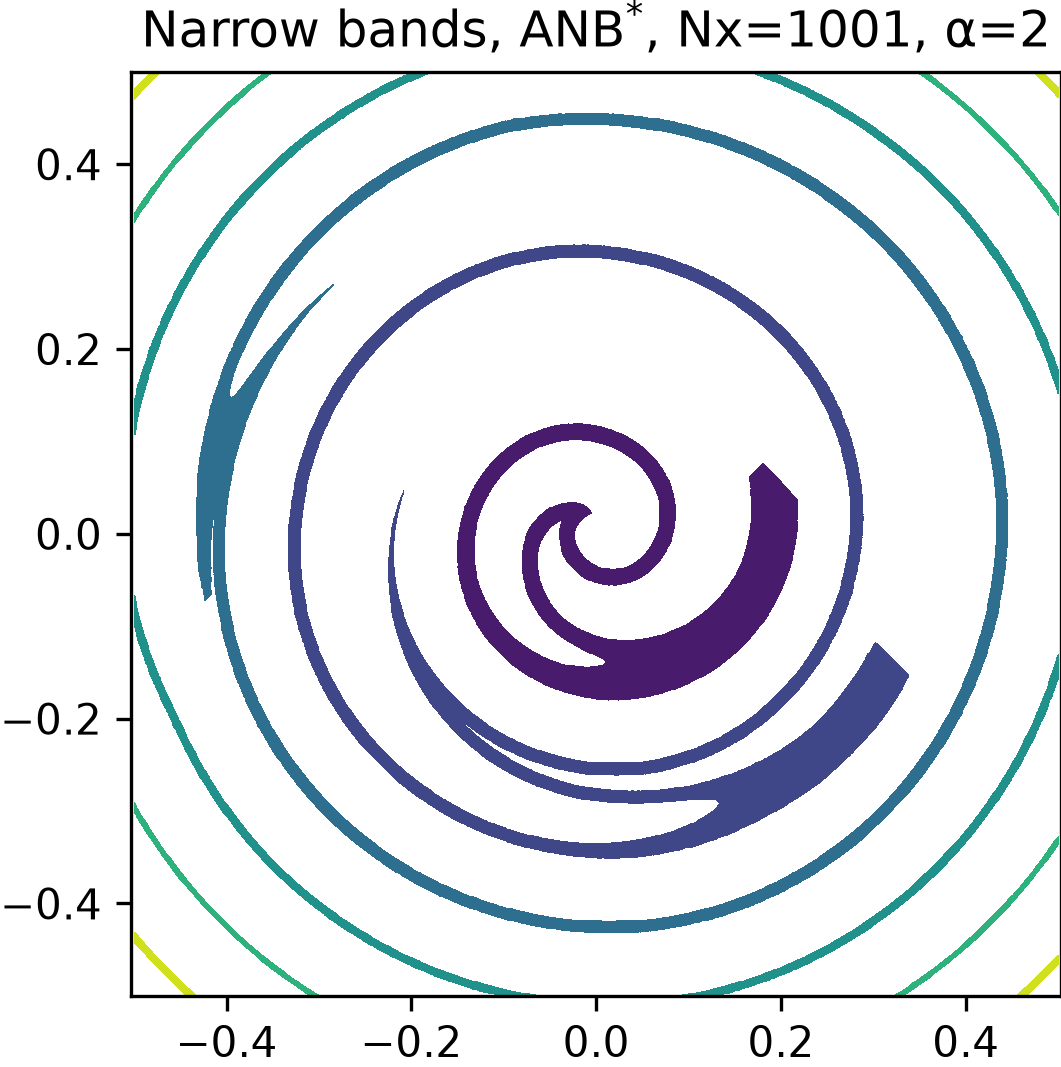}
	\includegraphics[height=\h]{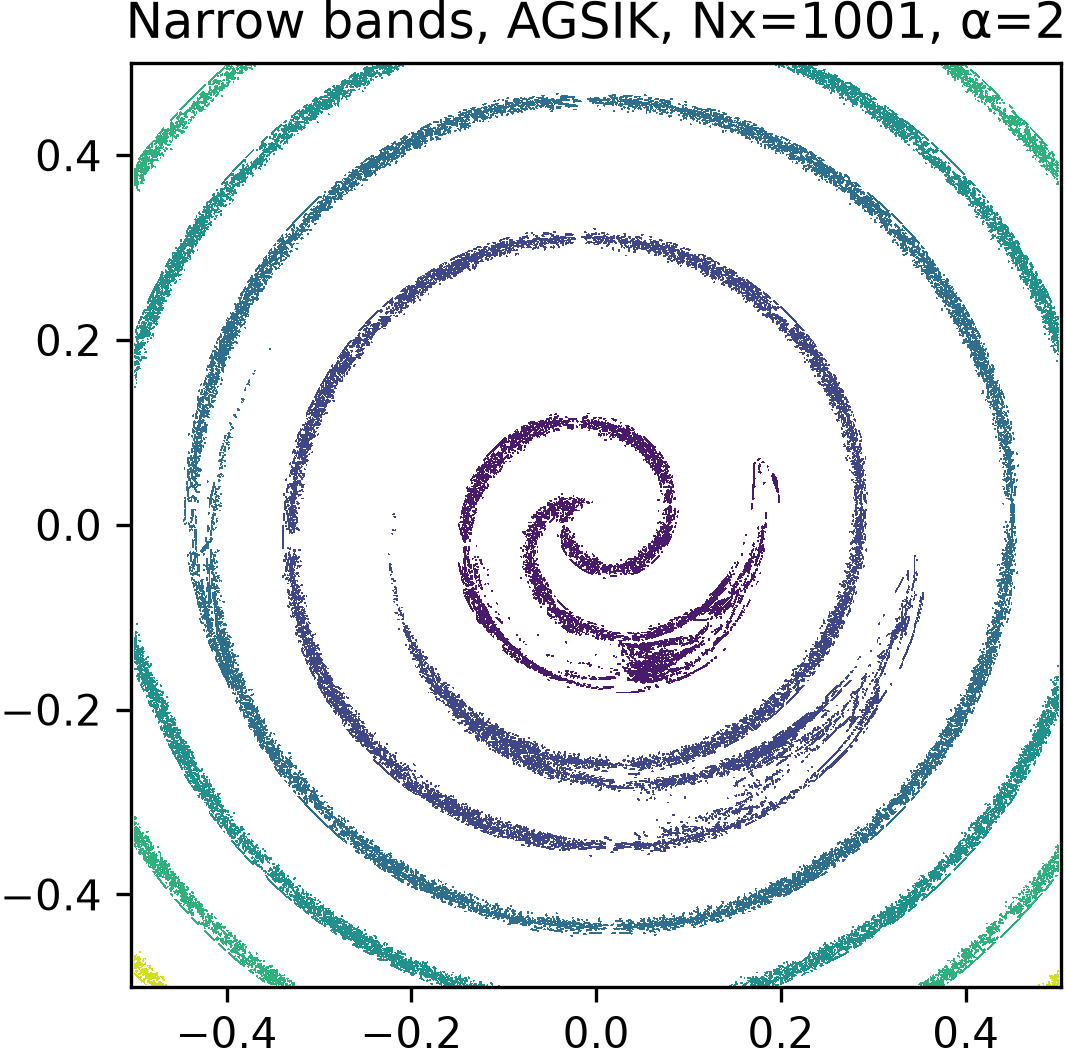}
\caption{Snail test case, SL4 scheme, see \cref{subsec:snail}. (i) Region, in blue, where the metric is anisotropic. (ii) Level sets of the solution, and a minimal path. (iii) Narrow bands of the ANB$^*$, and (iv) of the AGSIK. Problem size $1001 \times 1001$, \Kruzkov\ timescale $\alpha=2h$ and other parameters defaulted \eqref{eq:default_tol}, \eqref{eq:default_rho_tau}.}
\label{fig:snail}
\end{figure}

We numerically solve a strongly anisotropic Riemannian eikonal PDE, inspired by applications to tubular structure segmentation \cite{benmansour2011Tubular}. 
The Riemannian metric is Euclidean except within a small neighborhood of width $w_{\max} := 0.05$ of an Archimedean spiral $\Gamma$ given below, where the cost of motion in the tangential direction $\Gamma'$ is divided by a factor $\kappa:=20$. Because of the shape of the fast region, see \cref{fig:snail} (i), this is often referred to as the snail test case.
More precisely, we numerically solve $\cF^*_x(\nabla u(x)) = 1$ on the domain $x\in \Omega := ]-1/2,1/2[^2\sm\{0\}$ with
\begin{align*}
	\Gamma(r) &:= r (\cos \omega r, \sin \omega r), &
	\cF^*_x(\eta)^2 &:= |\eta|^2 + s(x)(\kappa^2-1) \<\eta,t(x)\>^2.
\end{align*}
The Archimedean spiral $\Gamma$ uses the pulsation $\omega := 12 \pi$, and argument $0 \leq r \leq r_{\max} := 0.43$.
The unit vector $t(x) \propto (-\omega |x| \sin \theta(x) + \cos \theta(x), \omega |x| \cos \theta(x) + \sin \theta(x) )$ approximates the tangent vector to $\Gamma$, where $\theta(x) \in [0,2\pi[$ is the angular coordinate of  $x$.
The  indicator function $s(x)\in \{0,1\}$ characterizes points close to the curve $\Gamma$, and equals one iff 
$0 \leq \mathrm{Mod}(|x|-\theta(x)/\omega,\ 2\pi/\omega) \leq w_{\max}$ 
and $|x| \leq r_{\max}$, where $\mathrm{Mod}(a,b)\in [0,b[$ denotes the remainder of $a$ in the Euclidean division by $b$.
A point source is introduced at the origin, and outflow boundary conditions are used on $\partial \Omega$.
See \cite{mirebeau2014marchingLattice,mansour2024solving} for more discussion. 

For simplicity, we set $\alpha=2h$, and default the other parameters \eqref{eq:default_tol}, \eqref{eq:default_rho_tau}, in the numerical experiment. 
The contraction ratio near linear maps takes the minimal value $\rho^\lin(\alpha)\approx 0.97$ (resp.\ $0.995$, resp.\ $0.99917$) for the timescale $\alpha \approx 1.32 h$ (resp.\ $4.94 h$, resp $19.2 h$), with the SL8 (resp.\ SL4, resp. SL3) scheme. 
The ANB$^*$ and AGSIK work well in all cases. The NB method requires an increased timescale $\alpha = 15h$ with the SL3 and SL4 schemes, otherwise it produces an incorrect solution.
Among the classical AGSI, FIM, FS and FSd methods, all struggle with the SL3 scheme, but the SL4 (resp.\ SL8) scheme is efficiently solved by the FSd (resp.\ FS and FSd) method thanks to the axis-aligned stencil structure, see \cref{table:niter}.

\subsection{Strongly asymmetric Randers test case}
\label{subsec:zermelo}

\begin{figure}
	\def\h{3.75cm}
	\includegraphics[height=\h]{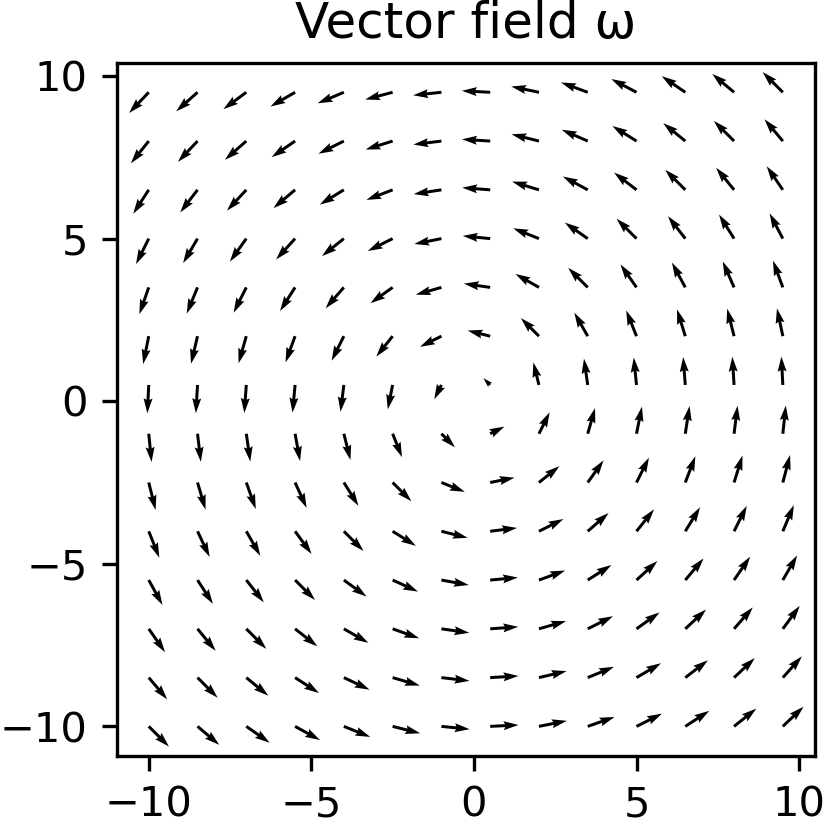}	
	\includegraphics[height=\h]{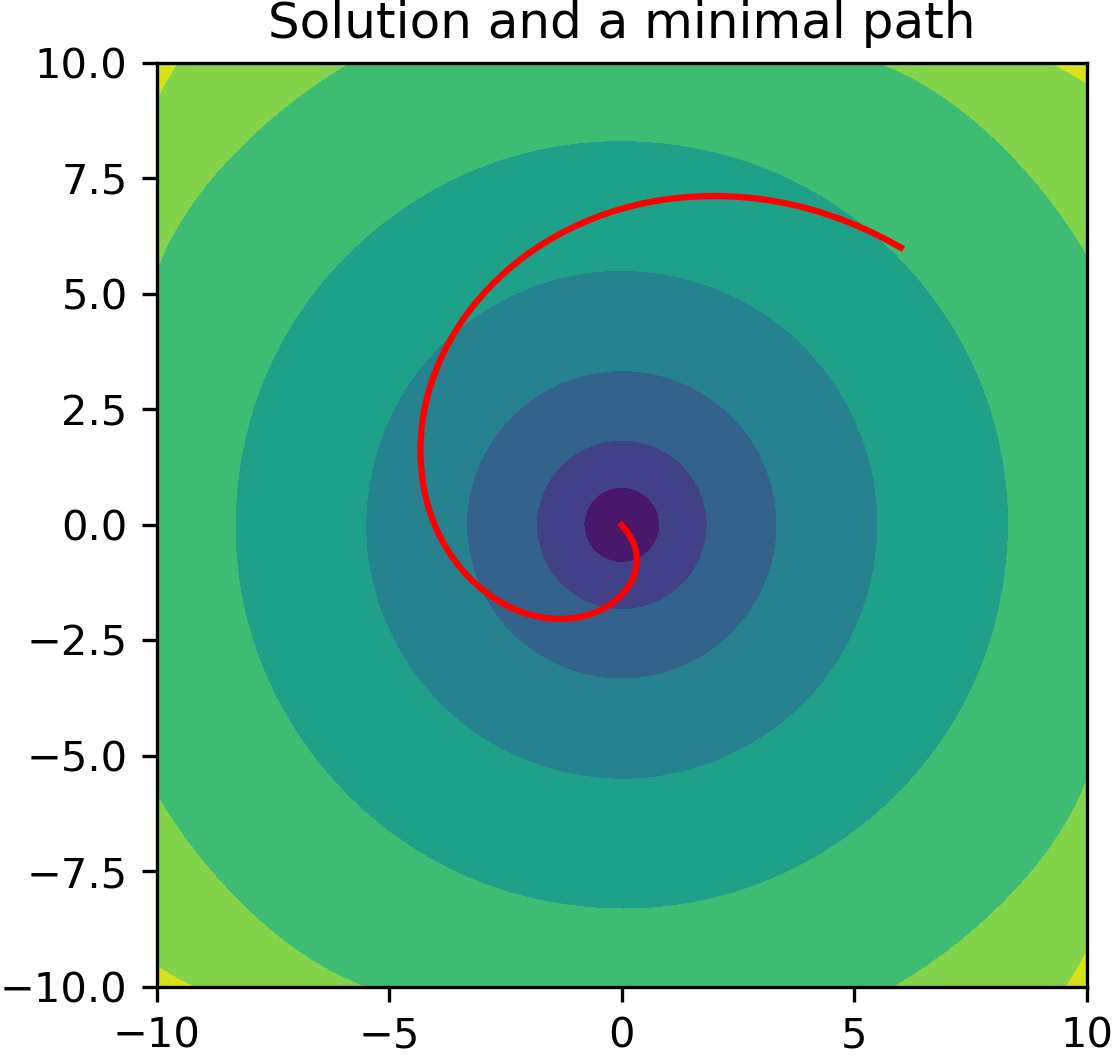}	
	\includegraphics[height=\h]{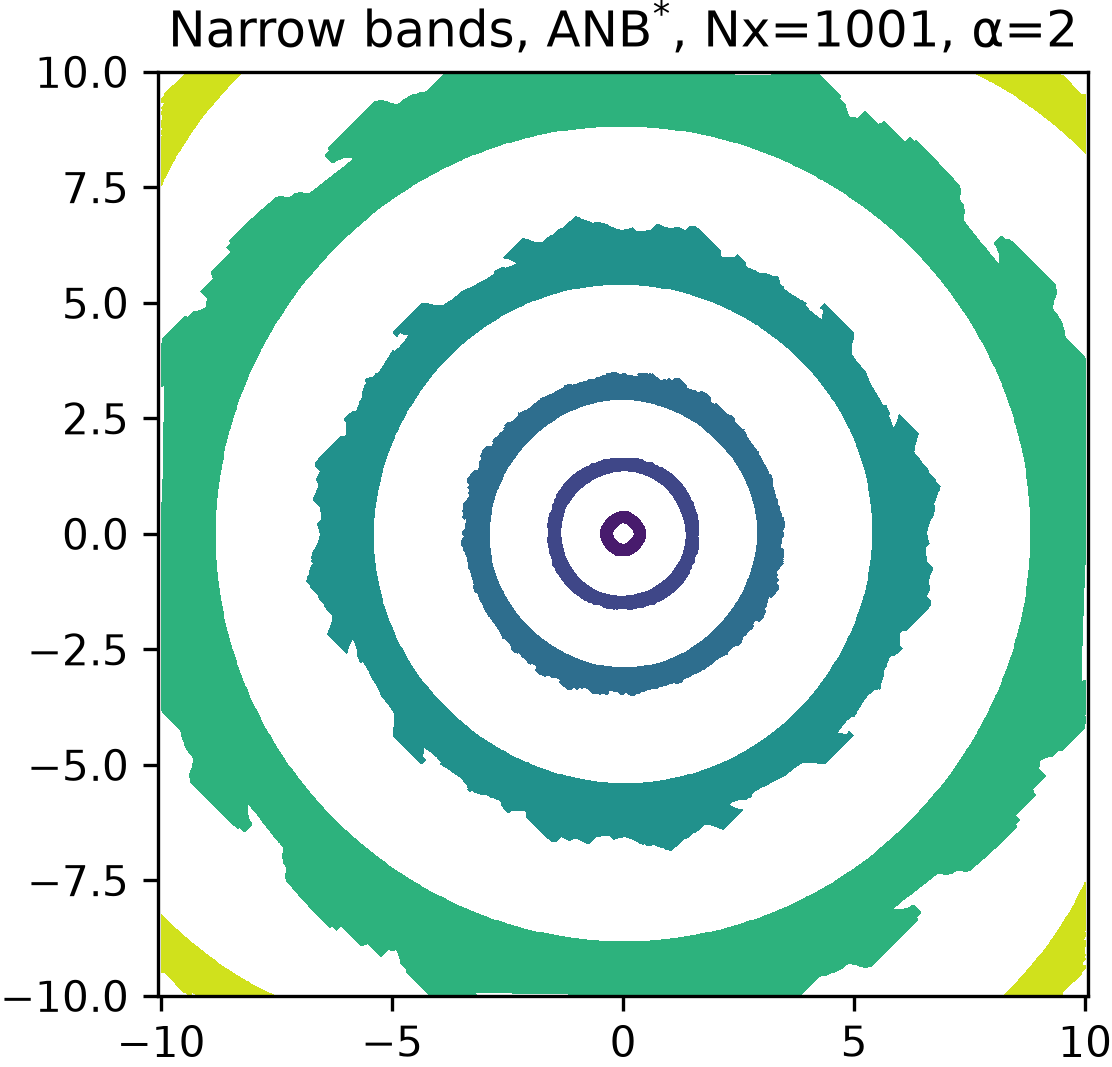}
	\includegraphics[height=\h]{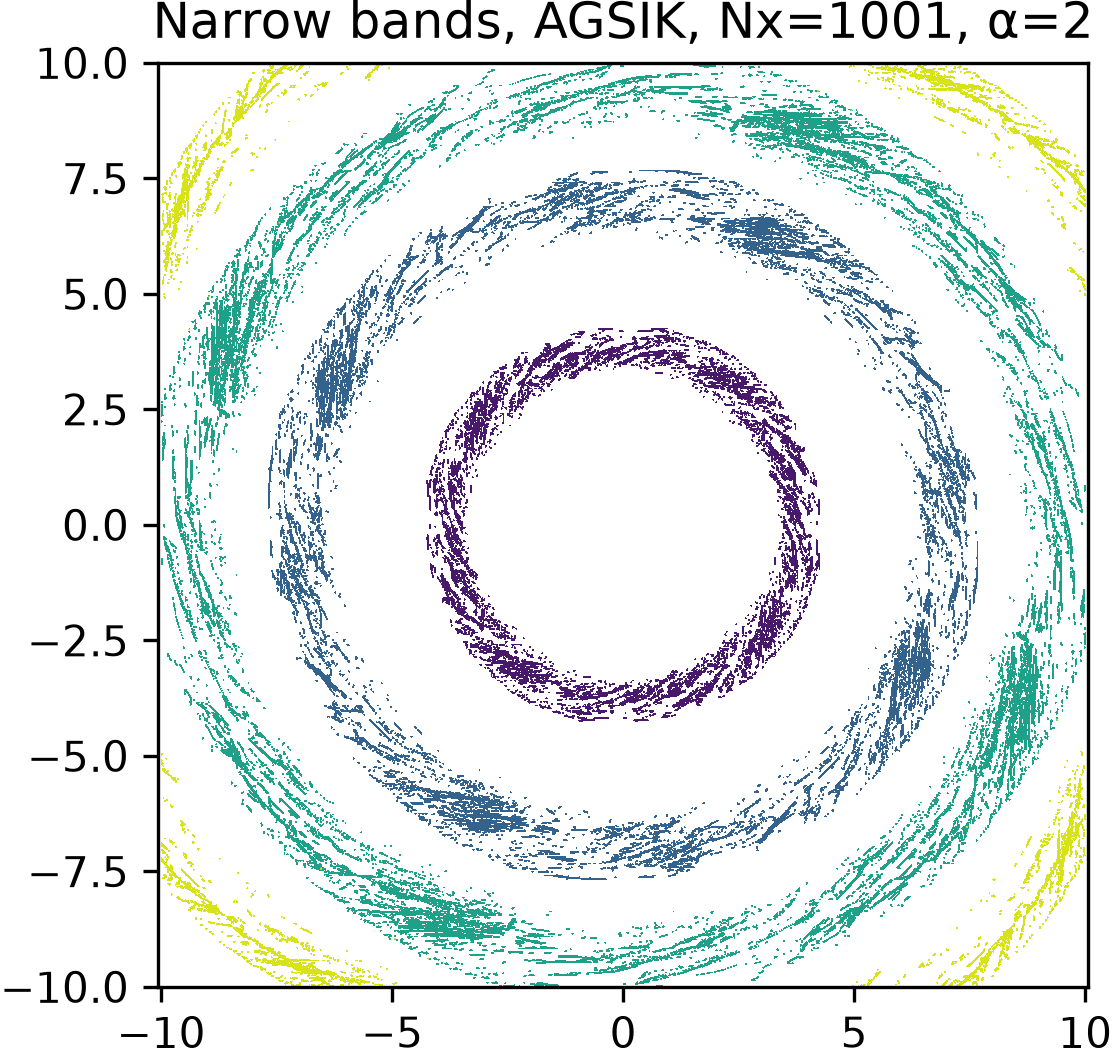}
	\caption{Swirl test case, SL4 scheme, see \cref{subsec:zermelo}. (i) Vector field $\omega$. (ii) Level sets of the solution, and a geodesic. (iii) Narrow bands of the ANB$^*$, and (iv) of the AGSIK. Problem size $1001 \times 1001$, \Kruzkov\ timescale $\alpha=2h$ and other parameters defaulted \eqref{eq:default_tol}, \eqref{eq:default_rho_tau}.}
	\label{fig:swirl}
\end{figure}

We solve an instance of Zermelo's navigation problem, where the vehicle has unit speed and the drift vector field $\omega$ barely satisfies the local controllability condition: $\|\omega(x)\|\leq \omega_{\max} := 0.98 < 1$, following \cite{mirebeau2014EfficientFinsler,mansour2024solving}. 
We often refer to it as the swirl test case, in view of the chosen rotating water current, see (\ref{eqdef:swirl_hamiltonian}, right) and \cref{fig:swirl} (i).
More precisely, we numerically solve the Randers eikonal PDE $\cF^*_x(\nabla u(x)) = 1$ where $x\in\Omega := ]-10,10[^2$ and
\begin{align}
\label{eqdef:swirl_hamiltonian}
	\cF^*_x(\eta) &:= \| \eta - \omega(x)\|, &
	\omega(x) &:= \omega_{\max} \frac{ \|x\|^2}{1+\|x\|^2} \frac{x^\perp}{\|x\|},
\end{align}
denoting $x^\perp := (-x_1,x_0)$ if $x=(x_0,x_1)$.
The Hamiltonian \eqref{eqdef:swirl_hamiltonian} is not positively homogeneous, but a standard reformulation yields a $1$-positively homogeneous Hamiltonian with the same unit level set, see e.g.\ \cite{mirebeau2014EfficientFinsler}.
A point source is introduced at the origin, and outflow boundary conditions are used on $\partial \Omega$.
Let us mention that the SL4 scheme for the Hamiltonian \eqref{eqdef:swirl_hamiltonian} is equivalent to the upwind differences (UD) scheme with a drift, see \cite[Appendix B]{mansour2024solving}.

For simplicity we set $\alpha = 2h$, and default the other parameters \eqref{eq:default_tol}, \eqref{eq:default_rho_tau}, in the numerical experiment. 
The contraction ratio near linear maps takes the minimal value $\rho^\lin(\alpha)\approx 0.85$ (resp.\ $0.97$, resp.\ $0.9956$) for the timescale $\alpha \approx 0.18 h$ (resp.\ $0.67 h$, resp.\ $2.65 h$), with the SL8 (resp.\ SL4, resp. SL3) scheme. The ANB$^*$ and AGSIK work well in all cases. The NB method requires an increased timescale $\alpha = 8h$ with the SL3 scheme, otherwise it produces an incorrect solution. 
Again, among the classical AGSI, FIM, FS and FSd methods, all struggle with the SL3 scheme, but the SL4 (resp.\ SL8) scheme is efficiently solved by the FSd (resp.\ FS and FSd) method thanks to the axis-aligned stencil structure, see \cref{table:niter}.

\section*{Conclusion}

In this paper, we present, to our knowledge, the first computationally efficient method for the progressive construction of \emph{sub-solutions} to (static first order) eikonal, HJB or HJI PDEs. 
The main difficulty, in comparison with the construction of super-solutions, is to find a substitute for the outflow boundary conditions $\ul u = +\infty$ in the intermediate sub-problems; for that purpose, our intermediate boundary conditions are produced by Dijkstra's method applied to a graph-based underestimation of the operator of interest in an obstacle type auxiliary problem.
With both a sub-solution and a super-solution at hand, we have access to a-posteriori error estimates, which are at the foundation of the proposed ANB method with adaptive narrow band width and guaranteed quasi-linear complexity, building and improving on \cite{cristiani2009fast,mansour2024solving}.
Another contribution is a new axiomatization of two mirror properties, the causality and the memory of a discretization scheme, used in the analysis of narrow band methods. 

Numerical experiments, on challenging test cases featuring strong anisotropy or diffusive schemes, illustrate that the proposed ANB method and its variant the ANB$^*$ maintain a quasi-linear complexity, as established. In contrast, the classical FS, FSd, AGSI and FIM iterative methods \cite{zhao2005FastSweeping,bornemann2006FiniteElement,jeong2008fim} become super-linear, while the NB method \cite{mansour2024solving}, which uses a fixed-width narrow band, requires careful tuning of its parameters (in particular the \Kruzkov\ timescale, see \cref{subsec:snail,subsec:zermelo}), otherwise it may yield an incorrect numerical solution.
We also introduce the AGSIK, a progressive variant of adaptive Gauss-Seidel iteration \cite{bornemann2006FiniteElement}, featuring intermediate sub-problems whose stopping criteria are expressed in \Kruzkov\ coordinates. It produces narrow bands similar to those of the ANB, and is often computationally cheaper in practice (since it avoids the sub-solution construction, is simpler and has less overhead in general), but lacks any complexity guarantee.

Several directions are open for future research, such as investigating applications to differential games and HJI PDEs featuring non-convex Hamiltonians (since those are often discretized using diffusive Lax-Friedrichs schemes), applications in which an $\ell^\infty$ error estimate on the numerical solution is crucial (since this can only be certified by providing a sub-solution), parallelization, and extending the proposed methods to problems where the eikonal equation is coupled with e.g.\ a conservation equation.

\paragraph{Acknowledgement.} The authors would like to thank Olaf Gainville (CEA, DAM, DIF) for fruitful discussions.

\appendix

\section{Adaptive Gauss-Seidel Iterations} 
\label{sec:AGSI}
We introduce a variant of the ANB - referred to as the ANB$^*$ - which uses Adaptive Gauss-Seidel Iterations (AGSI) to solve the sub-problems posed on the successive  narrow bands, resulting in greatly improved efficiency in practice, even though the complexity guarantees are slightly weaker, see \cref{th:complexity_with_agsi}. For comparison purposes, we also introduce the AGSIK, see \cref{algo:agsik}, a variant of the AGSI which uses intermediate stopping criteria expressed in \Kruzkov\ coordinates. 
Contrary to the ANB and ANB$^*$, the AGSIK does not benefit from any complexity guarantee, and only produces a super-solution.
However, it is also much simpler to describe and implement, produces similar looking narrow bands, and is quite efficient in practice. 

\paragraph{AGSI.}
Our first step is a complexity estimate for the parallel implementation of the AGSI method, described in \cref{algo:agsi_parallel}, applied to a $\rho$-sub-additive operator, such as $\Phi_\alpha := \rE_\alpha \Lambda_\infty \rE_\alpha^{-1}$ in our context. Note that the AGSI was originally \cite{bornemann2006FiniteElement} applied directly to the operator $\Lambda$, which lacks this property. %
No complexity estimate is known in the latter case to our knowledge, even though a heuristic argument \cite{bornemann2006FiniteElement} suggests the super-linear complexity $\cO(N^{1+1/d})$.

\begin{algorithm}
\caption{Parallel AGSI, for computing a super-solution. \leancheck{algo:agsi_parallel}\newline
Variant for computing a sub-solution: replace all instances of $v - \Phi v$ below with $\Phi v-v$.
}
\label{algo:agsi_parallel}
\begin{algorithmic}[1]
\State 
\textbf{Input:} a target tolerance $\eta>0$, a mapping $v\in \bR^X$ such that $0 \leq v-\Phi v \leq C$. 
\State
\textbf{Ensures at termination:} $0 \leq v-\Phi v\leq \eta$.
\State \textbf{While} $\max(v-\Phi v) > \eta$ \textbf{do}:
\State \label{algoline:agsi_parallel_for}
\phantom{bla} \textbf{For all} $x\in X$ \textbf{such that} $v(x)-\Phi v(x) > \eta$ \textbf{do in parallel}: 
$v(x) \gets \Phi v(x)$.
\end{algorithmic}
\end{algorithm}
\noindent

\emph{Implementation note for \cref{algo:agsi_parallel}:} In line \ref{algoline:agsi_parallel_for}, it is sufficient to test the condition $v(x)-\Phi v(x) > \eta$ at the points $x\in X$ which are neighbors to a point updated in the previous iteration (except of course in the first iteration).
The resulting program is mathematically equivalent, except that the \emph{elementary operations complexity} now matches the \emph{update complexity} estimated in \cref{prop:cplx_agsi_parallel}.

In contrast with the global iteration \cref{algo:global_iteration}, the AGSI may update only a few points in some iterations, hence the \emph{number of updates at any given point}, estimated in the next proposition, can be much smaller than the \emph{number of iterations}, which we do not estimate.

\begin{proposition}[\leancheck{prop:cplx_agsi_parallel}]
\label{prop:cplx_agsi_parallel}
Assume that $\Phi : \bR^X \to \bR^X$ is monotone and $\rho$-sub-additive, with $0<\rho<1$. Then
\cref{algo:agsi_parallel} (and likewise the variant adapted to sub-solutions) terminates using at most 
$N''_\Phi$ elementary updates $v(x) \gets \Phi v(x)$ at any given point $x\in X$. We denoted
\begin{equation*}
	N_\Phi'' := \Big\lceil \frac{\ln ((1-\rho) C/\eta)}{-\ln \rho} \Big\rceil_+ + \frac{2}{(1-\rho)^2}. 
\end{equation*}
\end{proposition}

\begin{proof}
We only address the super-solution case, whose proof applies to the sub-solution variant of \cref{algo:agsi_parallel} with trivial modifications. Denote by $v_n$ the state of map $v$ at the beginning of the $n$-th iteration of \cref{algo:agsi_parallel}, where $n=0,1,\cdots$, and let $B_n := \{v_n - \Phi v_n > \eta\}$ and $C_n := \max (v_n - \Phi v_n)$. 
One has $\Phi v_n \leq v_n$ for all $n\geq 0$, by \cref{prop:partialMonotone}.
By construction, 
\begin{align}	
\label{eq:vnp1_eq_Phi_Bb_vn}
v_n &\leq \Phi v_n + C_n,  &
\Phi v_n \leq v_{n+1} := \Phi_{B_n} v_n & \leq \Phi v_n + \eta,
\end{align}
using (i) the definition of $C_n$, (ii) $\Phi v_n \leq v_n$, and (iii) the definition of $B_n$. Therefore
\begin{equation*}
	v_{n+1} \leq \Phi v_n + \eta 
	\leq \Phi (\Phi v_n + C_n) + \eta 
	\leq \Phi (v_{n+1} + C_n) + \eta 
	\leq \Phi v_{n+1} + \rho C_n + \eta,
\end{equation*}
using successively (i)  (\ref{eq:vnp1_eq_Phi_Bb_vn}, right), (ii) the monotony of $\Phi$ and (\ref{eq:vnp1_eq_Phi_Bb_vn}, left), (iii) (\ref{eq:vnp1_eq_Phi_Bb_vn}, right) again, (iv) the $\rho$-sub-additivity of $\Phi$. It follows that $C_{n+1} \leq \rho C_n + \eta$, hence $C_n \leq \rho^n C_0 + \eta/(1-\rho)$ by an immediate induction argument, and therefore 
\begin{align*}
	v_N &\leq \Phi v_N +  \frac {2\eta}{1-\rho}, \text{ with } N := \Big\lceil \frac{\ln ((1-\rho) C/\eta)}{-\ln \rho} \Big \rceil_+, &
	\text{ thus } v_N &\leq v_* + \xi,
	\text{ with }\xi := \frac{2 \eta}{(1-\rho)^2},
\end{align*}
where $v_*$ denotes the fixed point of $\Phi$.
The second inequality is obtained using the comparison principle, see \cref{prop:comparison_Phi}, and observing that $v_N - \xi$ is a sub-solution, %
see \cref{prop:pertSubSuper}.

Recalling that $\Phi v_n \leq v_n$ in all iterations, by \cref{prop:partialMonotone}, we obtain $v_* \leq v_n$ by the comparison principle. If the update $v(x) \gets \Phi v(x)$ is performed $r$ times overall at a given point $x$ in the iterations $N, \cdots, M-1$, where $M> N$ is an arbitrary integer, then $v_*(x) \leq v_M (x) \leq v_N(x) - r \eta \leq v_*(x) + \xi - r \eta$. Thus $r \leq 2/(1-\rho)^2$ and the result follows.
\end{proof}

\paragraph{ANB$^*$.}
The proposed ANB method, as described in \cref{algo:anb}, solves in each iteration $n$ a discretized eikonal equation restricted to the subdomain $B_n$ defined by the narrow band, up to an adequate convergence tolerance $\eta_n$. Lines \ref{anb:while_sub} and \ref{anb:while_sup} of \cref{algo:anb} show that these sub-problems are addressed using a global iteration (GI) of the operator $\Lambda$ restricted to the narrow band $B_n$. This approach is sufficient to establish the quasi-linear theoretical complexity bound of \cref{th:complexity,th:complexity_withUnknown}, thanks to the contraction property of the update operator $\Phi_\alpha := \rE_\alpha \Lambda_\infty \rE_\alpha^{-1}$ in \Kruzkov\ coordinates, see \cref{prop:contractivity}.

However, GI does not shine by its efficiency in practice, and it is tempting to replace it with the AGSI, which often avoids many useless evaluations of the update operator in regions where the convergence tolerance is already satisfied. This can be especially useful in test cases where the narrow band is relatively wide and inhomogeneous. 
The next result extends the quasi-linear complexity guarantee of \cref{th:complexity} to this setting. %
Unfortunately, we could \emph{not} similarly extend \cref{th:complexity_withUnknown}, which makes the weaker assumption that $\Phi_\beta$ is $\sigma$-sub-additive for some \emph{unknown} constants $\beta>0$ and $\sigma<1$, unrelated to the method parameters $\alpha>0$ and $\rho<1$. %

\begin{theorem}[\leancheck{th:complexity_with_agsi}]
\label{th:complexity_with_agsi}
The statement of \cref{th:complexity} still holds if lines \ref{anb:while_sub} and \ref{anb:while_sup} of \cref{algo:anb} are replaced with the three lines below, with the following modified constant which now bounds the number of updates at any given point in the inner loops (rather than the number of inner loops)
\begin{equation*}
N''_\Lambda :=  \Big\lceil \frac{2 \Delta_{\max}}{-\alpha \ln \rho}\Big\rceil_+ + \frac 2{(1-\rho)^2}.	
\end{equation*}
\begin{algorithmic}[1]
\setcounter{ALG@line}{10}
\State
\textbf{Until} $\exp(-\ol u / \alpha) \leq \exp(-\Lambda \ol u / \alpha) + \eta_n$ and $\exp(-\Lambda \ul u / \alpha) \leq \exp(-\ul u / \alpha) + \eta_n$ on $B_n$ \textbf{do:}
\State \phantom{bl} 
\textbf{For all} $x \in B_n$ \text{such that} $\exp(-\ol u / \alpha) > \exp(-\Lambda \ol u / \alpha) + \eta_n$ \textbf{do in parallel}: $\ol u(x) \gets \Lambda \ol u(x)$ %
\State \phantom{bl}
\textbf{For all} $x \in B_n$ \text{such that} $\exp(-\Lambda \ul u / \alpha) > \exp(-\ul u / \alpha) + \eta_n$ \textbf{do in parallel}: $\ul u(x) \gets \Lambda \ul u(x)$
\end{algorithmic}
\end{theorem}

\begin{proof}
In the proof of \cref{th:complexity}, see the last paragraph of \cref{subsec:cplx_subadd}, replace \cref{prop:niter_while} with: \cref{prop:cplx_agsi_parallel} applied to the operators $\Phi:=[\Phi_\alpha]_{B_n}^{\hat v_n}$ (resp.\ $[\Phi_\alpha]_{B_n}^{\ul v_n}$) restricted to $B_n$ and using frozen values on $X \sm B_n$, similar to  \cref{prop:niter_while,prop:niter_while_withUnknown}, and to the parameters $\eta := (1-\rho) \exp(-n\tau/\alpha)$ and $C := \exp(-(n\tau - 2 \Delta_{\max})/\alpha)$.
\end{proof}

\paragraph{AGSIK method.}
As announced, we introduce a progressive variant of the AGSI \cite{bornemann2006FiniteElement}, featuring intermediate objectives whose stopping criterion is expressed in \Kruzkov\ coordinates, see line \ref{agsik:until} of \cref{algo:agsik}. They can be reformulated as follows, in the $n$-th main iteration:
\begin{equation*}
	u < \Lambda u + \ve - \alpha \ln (1- \exp(\frac{\Lambda u - n \alpha \ln 2} \alpha)).
\end{equation*}
By convention $\ln t = -\infty$ for all $t \leq 0$, and this criterion is thus automatically satisfied at all $x \in X$ such that $\Lambda u(x) > n\alpha \ln 2$.
In contrast, if $\Lambda u (x) \ll n \alpha \ln 2$, then it approximates the standard super-solution convergence criterion $u(x) \leq \Lambda u(x) + \ve$, based on the scheme residual.

\begin{algorithm}
\caption{AGSI with \Kruzkov\ based intermediate stopping criteria (AGSIK) \leancheck{algo:agsik}}
\label{algo:agsik}
\begin{algorithmic}[1]
\State
\textbf{Input:} a target tolerance $\ve>0$, a timescale $\alpha>0$, a monotone operator $\Lambda : \bR^X \to \bR^X$.
\State 
\textbf{Initialization}: Set the mutable array $u \in ]-\infty,\infty]^X$ to $u = +\infty$  
\State
\textbf{Ensures at termination:} $\Lambda u\leq u\leq \Lambda u+\ve$ 
\State
\textbf{For} $n=0,1,\cdots$ \textbf{ until} $u \leq \Lambda u + \ve$ on $X$ \textbf{do}:
\State \phantom{bla} \label{agsik:until}
\textbf{Until} $\exp(-\frac{\Lambda u} \alpha) < \exp(\frac{\ve-u} \alpha) + 2^{-n}$ on $X$ \textbf{do}:
\State \phantom{blabla} \label{agsik:parallel_for}
\textbf{For all} $x\in X$ \textbf{s.t.} $\exp(-\frac{\Lambda u(x)}\alpha) \geq \exp(\frac{\ve-u(x)}\alpha) + 2^{-n}$ \textbf{do in parallel}:
$u(x) \gets \Lambda u(x)$.
\end{algorithmic}
\end{algorithm}

\noindent\emph{Implementation notes.}
In order to efficiently select the points in line \ref{agsik:parallel_for} of \cref{algo:agsik}, one can maintain a sorted list (implemented using a heap) containing the following elements sorted w.r.t.\ the second variable $\eta$
\begin{align}
\label{eqdef:agsik_L}
	L &= \big\{(x,\eta(u,x)) \mid x\in X,\ \eta(u,x)>0\big\}, &
	\eta(u,x) &:= \exp(-\frac{\Lambda u(x)} \alpha) - \exp(\frac{\ve-u(x)} \alpha)
\end{align}
Alternatively, one can maintain the collection of bins $L_k := \{x  \mid 2^{-k} \leq \eta(u,x) < 2^{-(k-1)}\}$, where $k \in \bN$, similarly to Dial's algorithm \cite{tsitsiklis1995EfficientTrajectories}.
In our experience the latter is usually more computationally efficient, and $8$ simultaneously active bins are enough, \emph{except} in the maze example with slightly permeable walls of \cref{subsec:maze} due to the high speed contrast. If there are too few bins, points are updated early rather than lost.

To avoid floating point underflow, the expression (\ref{eqdef:agsik_L}, right) should not be numerically evaluated directly, but one should rather compute its base two logarithm as follows: 
\begin{align}
\label{eq:log_approx}
-\ln_2 \eta(u,x) &= \frac{\Lambda u(x)}{\alpha\ln 2} -
\mathrm{L}_2(\frac{u(x)-\Lambda u(x)-\ve} \alpha), &
\mathrm{L}_2 (x) &:= \ln_2(1-e^{-x}) \approx
\begin{cases}
0 \hfill \text{if } x\geq 1,\\
\lceil\ln_2 x\rceil \text{ else}.
\end{cases}\!\!\!\!\!\!\!
\end{align}
Finally, the approximation (\ref{eq:log_approx}, right) directly extracts the exponent from the floating point representation of $x$, which eliminates the computational cost of evaluating the transcendental $\log$ and $\exp$ functions, without compromising the method's efficiency in practice.

\section{Discretizations of the anisotropic eikonal equation}
\label{sec:effective}
We describe in this appendix three classical discretizations of the anisotropic eikonal, HJB and HJI PDEs of interest, namely the semi-Lagrangian, Lax-Friedrichs, and Upwind-Differences schemes, and we show that they obey the assumptions of our main result \cref{th:complexity}. 
For the most part, these scheme properties are established in \cite{mansour2024solving}, hence we simply recall them.
The new results from this section are the proof of \cref{prop:constr_Gamma}, which was announced in the introduction and is presented below, and the explicit expressions \eqref{eq:rho_loc_LF}, \eqref{eq:rho_loc_SL} and \eqref{eq:rho_loc_UD} of the local contraction ratio \eqref{eq:subadd_diff} of the considered schemes near a linear mapping.

\paragraph{Preliminary: construction of a causal under-estimation $\Gamma \leq \Lambda$.}
Both the presented ANB method and the NB method \cite{mansour2024solving} use a causal approximation of the discretisation scheme $\Lambda$, in order to propagate information into regions where the numerical solution is not yet stabilized, and generate a useful narrow band. However, the NB method needs an \emph{over-estimation} $\widetilde \Gamma \geq \Lambda$.
Conveniently, the causal \emph{under-estimation} $\Gamma \leq \Lambda$ required by the ANB can be built in a straightforward and systematic manner via \cref{prop:constr_Gamma}, whose proof follows. 

\begin{proof}[Proof of \cref{prop:constr_Gamma}]
	\emph{Proof of $C_\past$-memory \eqref{eqdef:Cpast_memory}.} 
	Let $u,v\in \bR^X$, $t \in \bR$, and let us assume that $\max \{u,t\} = \max\{v,t\}$ following the premise of \eqref{eqdef:Cpast_memory}.
	Fix $x \in \Int X$, and assume that $\Lambda u(x) \geq \Lambda v(x)$ w.l.o.g.\ up to exchanging the roles of $u$ and $v$. 
	We distinguish two cases: If $t + C_\past \geq \Lambda u(x)$, then (\ref{eqdef:Cpast_memory}, right) holds by construction. On the other hand, if $t + C_\past < \Lambda u(x)$, then $t < \min_{\cV(x)} u$  by (\ref{eqdef:Cpast_dependency}, right).
	Recalling that $\max \{u,t\} = \max \{v,t\}$, we obtain $u_{|\cV(x)} = v_{|\cV(x)}$, hence $\Lambda u(x) = \Lambda v(x)$ by (\ref{eqdef:Cpast_dependency}, left), which establishes (\ref{eqdef:Cpast_memory}, right) and concludes.
	
	\emph{Construction of a super-solution to $\Lambda$.} 
	Consider the mapping $\ul u(x) := C_\past d_\cV(x,\partial X) + \max(u_\bd)$, where $d_\cV(x,\partial X)$ denotes the number of edges from $x$ to $\partial X$ on the oriented graph defined by the stencils $\cV$, and $u_\bd \in \bR^{\partial X}$ is the boundary data. Then $\ul u(x) = C_\past + \min_{\cV(x)} \ul u \geq \Lambda \ul u(x)$ for all $x\in \Int X$ by \eqref{eqdef:Cpast_dependency}, and $\ul u = \max (u_\bd) \geq u_\bd = \Lambda \ul u$ on $\partial X$, thus $\ul u\geq \Lambda \ul u$ as required.%
	
	\emph{Properties of the operator $\Gamma$.}
	The monotonicity, additive invariance and $\delta$-strict-causality of $\Gamma$ are clear. 
	\emph{Proof of the underestimation property $\Gamma \leq \Lambda$}: given $u \in \bR^X$, one has $\Lambda u = u_\bd = \Gamma u$ on $\partial X$, and
	\begin{equation*}
		\Gamma u(x) = \Lambda[0](x) + \min_{\cV(x)} u = \Lambda\big[\min_{\cV(x)} u\big] (x) \leq \Lambda u(x),
	\end{equation*}
	for all $x\in \Int X$,
	using successively (i) the definition \eqref{eqdef:constr_Gamma} of $\Gamma$, (ii) the additive invariance of $\Lambda$, and (iii) the stencil property \eqref{eqdef:Cpast_dependency} along with monotony. \emph{Proof of $(C_\Gamma,C'_\Gamma)$-comparison}: $\forall x \in \Int X$
	\begin{equation*}
		C_\Gamma \Lambda u(x) \leq C_\Gamma (C_\past + \min_{\cV(x)} u) = C_\Gamma' + \delta + C_\Gamma \min_{\cV(x)} u \leq C_\Gamma' + \Gamma [C_\Gamma u](x) \leq C_\Gamma' + \Lambda [C_\Gamma u](x)
	\end{equation*}
	using successively (i) assumption \eqref{eqdef:Cpast_dependency}, (ii) the expression of $C'_\Gamma$ (\ref{eqdef:constr_Gamma}, iv), (iii) $\delta \leq c(x)$, (iv) $\Gamma \leq \Lambda$ as established above. We conclude the proof noting that
	\begin{equation*}
		\Lambda[C_\Gamma u](x) \leq C_\past + C_\Gamma \min_{\cV(x)} u = C_\Gamma (\delta + \min_{\cV(x)} u) \leq C_\Gamma \Gamma u(x)
	\end{equation*}
	using successively (i) assumption \eqref{eqdef:Cpast_dependency}, 
	(ii) the expression \eqref{eqdef:constr_Gamma} of $C_\Gamma$, %
	and (iii) $\delta \leq c(x)$. 
	Note also that $\delta\leq \Lambda[0]\leq C_\past$ on $\Int X$ by (\ref{eqdef:Cpast_dependency}, right) with $u=0$, hence $C_\Gamma\geq 1$ and $C'_\Gamma=C_\past^2/\delta-\delta\geq 0$, consistently with \cref{assum:all}.
\end{proof}

Before turning to the different discretization schemes $\Lambda$, let us mention that the constants associated with their properties always obey by construction $\alpha,\nu = \Theta(h)$, $C_\past = \cO(h)$, $\rho\leq \rho_{\max}<1$, where $h$ denotes the grid scale, as required for our final complexity analysis \eqref{eq:homog_scheme_csts}.
Thus $\delta := \min_{\Int X} \Lambda[0] \geq \nu = \Theta(h)$ by \cref{prop:cst_subsol} (alternatively $\delta \geq -\alpha \ln \rho = \Theta(h)$, again by \cref{prop:cst_subsol}).
It follows that $C_\Gamma = C_\past/\delta = \cO(1)$ and $C_\Gamma' = C_\Gamma C_\past-\delta = \cO(h)$ by \cref{prop:constr_Gamma}, which concludes the proof of the homogeneity properties \eqref{eq:homog_scheme_csts}.

\paragraph{Semi-Lagrangian scheme.}  
We refer to \cite{falcone2013semi} for an overview of this class of methods, and follow \cite{tsitsiklis1995EfficientTrajectories,kimmel1998Manifold,sethian2003OUM,bornemann2006FiniteElement,alton2012MAOUM,mirebeau2014EfficientFinsler,mirebeau2014marchingLattice,dahiya2018ordered}, generalizing \cref{rem:FE}. Let
\begin{equation}
\label{eqdef:SLupdate}
  \Lambda u(x) := \min_{(v_1,\cdots,v_d) \in \bV(x)}
  \ \min_{(\xi_1,\cdots,\xi_d) \in \Xi_d} \ 
  \ol \cF_x\Big(\sum_{1 \leq i \leq d} \xi_i v_i\Big) + \sum_{1 \leq i \leq d} \xi_i u(x+v_i), 
\end{equation}
for all $u \in \bR^X$ and $x \in \Int X$, and by convention $\Lambda u = u_\bd$ on $\partial X$. 
We denoted $\Xi_d := \{\xi \in [0,1]^d\mid \<\xi,\bbone\> = 1\}$ with $\bbone := (1,\cdots,1) \in \bR^d$, 
see \cref{assum:stencil} for the other terms.
The semi-Lagrangian scheme \eqref{eqdef:SLupdate} discretizes the PDE $\cF_x^*(\nabla u(x)) = 1$,  where the Hamiltonian $\cF_x^*(\eta) := \sup\{\<\xi,\eta\> \mid \cF_x(\xi) \leq 1\}$ is dual to the metric $\cF$, hence necessarily convex. The optimization problem \eqref{eqdef:SLupdate} has a closed form solution for Isotropic, Riemannian, and Randers metrics, but may require an inner e.g.\ Newton solver for more complex metrics \cite{desquilbet2021single}.

\begin{assumption}
\label{assum:stencil}
$X = \Int X \sqcup \partial X \subset \bR^d$ is a finite set. $\cF_x$ is a (possibly asymmetric) given norm on $\bR^d$ and $\ol \cF_x(v) := \cF_x(-v)$, for each $x \in \Int X$, $v \in \bR^d$. 
$\bV(x) \subset \bR^{d \times d}$ is a given finite collection of square matrices, and the corresponding stencil $\cV(x) := \{x+v_i \mid (v_1,\cdots,v_d)\in \bV(x),\ 1\leq i\leq d\}$ satisfies $\cV(x) \subset X$ for all $x \in \Int X$.
We further assume:
\begin{itemize}
	\item (Connectedness to the boundary) 
For all $x_0 \in \Int X$, 
there exists $n\geq 1$, 
$x_1,\cdots,x_{n-1} \in \Int X$, $x_n \in \partial X$, such that $x_{i+1} \in \cV(x_i)$ for all $0 \leq i <n$.
	\item  (Non-degeneracy of simplices) For all $x \in \Int X$, all $(v_1,\cdots,v_d) \in \bV(x)$, 
	$\det(v_1,\cdots,v_d) \neq 0$.
\end{itemize}
\end{assumption}
Monotonicity and additive invariance hold by construction.
\cite[Corollary 3.11]{mansour2024solving} 
establishes $\nu$-sub-multiplicativity with the parameter given below.
The operator  $\Phi_\alpha := \rE_\alpha \Lambda_\infty \rE_\alpha^{-1}$, in \Kruzkov\ coordinates, is $\rho$-sub-additive \cite[Proposition 2.27]{mansour2024solving}, with the parameters:
\begin{align*}
	\alpha &:= \max_{(x,V) \in \bV} \max_{|\xi|=|\zeta|=1} \frac{2\,\ol\cF_x(V\xi)^2}{\ol\cF_x(V\zeta)}, &
	\rho &:= 1-\frac{\nu} {2\alpha}, &
	\nu &:= \min_{(x,V) \in \bV} \min_{\xi \in \Xi_d} \ol\cF_x(V \xi),
\end{align*}
where $(x,V) \in \bV$ is a shorthand for: $x \in \Int X$ and $V \in \bV(x)$.
Define as follows the intrinsic grid scale $h>0$, and the anisotropy ratio $\kappa \geq 1$ which measures the geometric degeneracy of the eikonal PDE and of its discretization: denoting by $|\cdot|$ the Euclidean norm
\begin{align*}
	h &:= \max_{(x,V) \in \bV,\, |\xi|=1} \ol\cF_x(V\xi), &
	\kappa &:= \max_{(x,V) \in \bV,\, |\xi|=1} \ol\cF_x(V\xi) \Big/ \min_{(x,V) \in \bV,\, |\xi|=1} \ol\cF_x(V\xi).
\end{align*}
Assumption \eqref{eqdef:Cpast_dependency} of \cref{prop:constr_Gamma} holds with $C_\past := h$, which yields a causal underestimation $\Gamma$ and the remaining required properties of $\Lambda$.

Finally, the homogeneity properties $2 h \leq \alpha \leq 2 \kappa h$, $\rho \leq 1-\tfrac 1 {4 \kappa^2\sqrt d}$, $\tfrac h {\kappa\sqrt d}\leq \nu \leq h$, established in \cite[Proposition 2.27]{mansour2024solving}, are as required \eqref{eq:homog_scheme_csts}.

\paragraph{Lax-Friedrichs scheme.}
Define the following scheme $\Lambda$, inspired by \cite{kao2004lax}:
\begin{align}
\nonumber
	\Lambda_0 u(x) &
	:= h C_0 +\!\! \min_{\substack{1 \leq i \leq d\\ \chi\in\{-1,1\}}} u(x+\chi h b_i), &
	\Lambda_1 u(x) 
	&:= \frac {c_1} d\big(h - \cF_x^*(h\nabla_h u(x))\big) 
	+ \!\!\! \sum_{\substack{1 \leq i \leq d\\ \chi\in\{-1,1\}}} \frac{u(x+\chi hb_i)}{2d}, 
\end{align}
\vspace{-7mm}
\begin{align}
	\label{eqdef:LFupdate}
	\nabla_h u(x) &:= \sum_{1 \leq i \leq d} \frac{u(x+h b_i)-u(x-h b_i)}{2h} b_i,&
	\Lambda u(x) &:= \min \{\Lambda_0 u(x), \Lambda_1 u(x)\},
\end{align}
for all $u \in \bR^X$ and all $x \in \Int X$, where $(b_i)_{1 \leq i \leq d}$ denotes the canonical basis of $\bR^d$. By convention $\Lambda_0 u = \Lambda_1 u = \Lambda u = u_\bd$ on $\partial X$. 
The scheme is the combination of a slope limiter $\Lambda_0$, and of a discretization $\Lambda_1$ of the HJB or HJI (for an non-convex Hamiltonian $\cF^*$) PDE $\cF^*_x(\nabla u(x))=1$  using centered differences and a Lax-Friedrichs stabilization term, see \cite[Appendix A]{mansour2024solving} for more discussion.  %
\begin{assumption}
	\label{assum:LF}
	The domain 
	$\Omega \subset \bR^d$ is open and bounded, and the Hamiltonian $\cF^* \in C^0(\overline \Omega\times \bR^d,[0,\infty[)$ obeys 
\begin{align}
\label{eq:c0C1}
	\cF^*_x(\lambda \xi) &= \lambda \cF_x^*(\xi), &
	C_0^{-1} \|\xi\|_\infty &\leq \cF_x^*(\xi), &
	|\cF_x^*(\xi)-\cF_x^*(\zeta)| \leq C_1 \|\xi-\zeta\|_1,
\end{align}
for all $\lambda\geq 0$ and $\xi,\zeta\in \bR^d$, where $C_0,C_1$ are positive constants.
	The grid scale $h>0$ is fixed, and the discrete domain $X := \Int X \sqcup \partial X$ is defined by 
\begin{align}
\label{eqdef:XhLF}
	\Int X &:= \Omega \cap h \bZ^d, &
	\partial X := \{x+\chi h b_i\mid x\in \Int X,\ 1 \leq i \leq d,\ \chi = \pm 1\}\sm \Int X.
\end{align}
We also fix $0 < c_1 < C_1^{-1}$ and denote $\delta_\LF := (1-C_1c_1)/(2d)$. 
\end{assumption}
Note that the  Hamiltonian $\cF^*$ needs \emph{not} be convex w.r.t.\ the second variable.
The monotonicity, additive invariance, %
and $\nu$-sub-multiplicativity of $\Lambda$, and the $\rho$-sub-additivity of $\Phi_\alpha:= \rE_\alpha \Lambda_\infty \rE_\alpha^{-1}$ in \Kruzkov\ coordinates, are established in \cite[Theorem A.9]{mansour2024solving} with the parameters
\begin{align}
\label{eq:LF_alpha_rho_nu}
	\alpha &:= 2h \frac {d C_0^2}{c_1 \delta_\LF^2}, & 
	\rho &:= 1- \frac 1 2 \Big(\frac {c_1 \delta_\LF}{C_0 d}\Big)^2, & 
	\nu &:= h \frac {c_1} d.
\end{align}
Assumption \eqref{eqdef:Cpast_dependency} of \cref{prop:constr_Gamma} holds with $C_\past := C_0 h$ and $\cV(x) = \{x + \chi h b_i\}_{\chi=\pm 1}^{1 \leq i \leq d}$, which yields a causal underestimation $\Gamma$ and the remaining required properties of $\Lambda$. 

\emph{Pointwise scheme stabilization parameters:} one may consider $C_0(x)$, $C_1(x)$ and $c_1(x)$ (thus also $\delta_\LF(x)$) depending on the current point $x$, see \cref{subsec:seismic}. In that case $\alpha = \max_{x \in \Int X} \alpha(x)$, $\rho = \max_{x \in \Int X} \rho(x)$, and $\nu = \min_{x \in \Int X} \nu(x)$, where $\alpha(x)$, $\rho(x)$ and $\nu(x)$ are defined by \eqref{eq:LF_alpha_rho_nu}.

\paragraph{Upwind-differences scheme.}
Following \cite{mirebeau2019riemannian}, and generalizing the seminal work \cite{sethian1996fast}, define $\Lambda u(x)$ as the unique solution $\lambda$ to the non-linear finite differences equation
\begin{align}
\label{eqdef:UDupdate}
	\sum_{e \in E} \mu(x,e) \max \big\{0,\, \lambda-u(x-h e)+h\<\omega(x),e\>,\, \lambda - u(x+ h e)-h\<\omega(x),e\>\big\}^2 = h^2, 
\end{align}
for all $u \in \bR^X$ and $x \in \Int X$, and by convention $\Lambda u = u_\bd$ on $\partial X$. 
The scheme discretizes Randers eikonal equation $\|\nabla u(x) + \omega(x)\|^2_{\cD(x)}=1$, for all $x\in \Omega$. The symmetric positive definite matrix $\cD(x)$ should be regarded as the inverse of a Riemannian metric, and is related to the scheme coefficients $\mu(x,e)$ and offsets $e\in E$ by \eqref{eqdef:selling}. The vector field $\omega(x)$ models a drift and is subject to a smallness condition ensuring local controllability, see below. Randers metrics arise, for instance, in Zermelo's navigation problem, see \cite{bardi2008viscosity} and \cref{subsec:zermelo}. 
Solving \eqref{eqdef:UDupdate} w.r.t.\ $\lambda$ amounts to sorting $\#(E)$ values and solving $\#(E)$ univariate quadratic equations.

\begin{assumption}
\label{assum:Euler}
The domain $\Omega\subset \bR^d$ is open and bounded. The Randers eikonal equation parameters $\cD \in C^0(\overline \Omega,\cS_d^{++})$ and $\omega \in C^0(\overline\Omega, \bR^d)$ obey the compatibility condition $\|\omega(x)\|_{\cD(x)} \leq \omega_{\max} < 1$. The 
stencil $E \subset \bZ^d \sm \{0\}$ is finite and the coefficients $\mu \in C^0(\overline \Omega \times E,]0,\infty[)$ obey 
\begin{align}
\label{eqdef:selling}
	\cD(x) &:= \sum_{e \in E} \mu(x,e) e e^\top, &
	0 < \mu_{\min} \leq \mu(x,e) \leq \mu_{\max}.
\end{align}
The discrete domain $X = \Int X \sqcup \partial X$ is defined in terms of the grid scale $h>0$ by:
\begin{align*}
	\Int X &:= \Omega \cap h \bZ^d, &
	\partial X &:= \{x+\chi h e\mid x \in \Int X,\ e \in E,\ \chi = \pm 1\} \sm \Int X.
\end{align*}
\end{assumption}
The upwind differences (UD) scheme $\Lambda$ defined by \eqref{eqdef:UDupdate} is monotone, additively invariant, $\nu$-sub-multiplicative, and 
$\Phi_\alpha := \rE_\alpha \Lambda_\infty \rE_\alpha^{-1}$ is $\rho$-sub-additive with  \cite[Corollary B.6]{mansour2024solving}
\begin{align}
	\alpha &:= \frac{8 h \kappa}{\sqrt{\mu_{\min}}}, &
	\rho &:= 1-\frac{\nu}{2\alpha}, &
	\nu &:= \frac{h}{\kappa \sqrt{\mu_{\min} \#(E)}}, & 	
	\text{where } \kappa := \frac{\sqrt{\mu_{\max}/\mu_{\min}}}{1-\omega_{\max}}.
\end{align}	
These properties are established in \cite[Corollary B.6]{mansour2024solving} by reformulating the upwind-differences scheme as a semi-Lagrangian scheme associated with a diagonal metric in higher dimension $d' := \#(E)$. 
For all $x \in \Int X$ one has %
\begin{align*}
	\Lambda u(x) & \leq C_\past + \min_{\cV(x)} u, &
	\text{where } \cV(x) &:= \{x+\chi h e\mid e \in E, \chi = \pm 1\}, &
	\text{and } C_\past  &:= \frac {2 h}{\sqrt{\mu_{\min}}},
\end{align*}
using the defining expression \eqref{eqdef:UDupdate}, along with $|\<\omega(x),e\>| \leq \|\omega(x)\|_{\cD(x)} \|e\|_{\cD(x)^{-1}} \leq 1/\sqrt{\mu_{\min}}$ using (i) Cauchy-Schwarz's inequality and (ii) the smallness assumption on $\omega$ and (\ref{eqdef:selling}, left). This establishes assumption \eqref{eqdef:Cpast_dependency} of  \cref{prop:constr_Gamma}, which yields a causal underestimation $\Gamma$ and the remaining properties of $\Lambda$.

\paragraph{Local contraction ratio near linear maps.}
Consider a Hamiltonian $\cF^* : \Omega \times \bR^d \to [0,\infty[$, a discretization scheme $\Lambda : \bR^X \to \bR^X$, an interior point $x \in \Int X$, and a co-vector $\eta \in \bR^d$ such that  $\cF_x^*(\eta) = 1$, i.e.\ the linear mapping $\<\eta,\cdot\>$ satisfies the PDE of interest at the point $x$. 
We evaluate the local contraction ratio \eqref{eq:subadd_diff}, near the linear mapping $\<\eta,\cdot\>$, 
for the three schemes of interest: denoting by $(b_1,\cdots,b_d)$ the canonical basis of $\bR^d$
\begin{align}
\label{eq:rho_loc_LF}
\rho^\loc_\LF(\alpha; \<\eta,\cdot\>, x) &= 
\frac 1 {2d} \sum_{\substack{1 \leq i \leq d \\ \chi=\pm}} \big(1-c_1 \chi \< \nabla \cF_x^*(\eta), b_i\> \big) \exp\big(\frac{\chi \<\eta, b_i\>}{\alpha/h}\big). \\
\label{eq:rho_loc_SL}
\rho^\loc_\SL(\alpha; \<\eta,\cdot\>, x) &= \!
\sum_{1 \leq i \leq d} \frac{\<V^{-1} \nabla \cF_x^*(\eta), b_i\>}{\<V^{-1} \nabla \cF_x^*(\eta), \bbone\>} 
\exp\big(\frac{\<V b_i, \eta\>} {\alpha}\big), %
\qquad
V \in \bV(x),\, V^{-1} \nabla \cF_x^*(\eta) < 0.
\\
\label{eq:rho_loc_UD}
\rho^\loc_\UD (\alpha; \<\eta,\cdot\>, x)\! &=  \!\!
\sum_{e \in E} 
\frac{\mu(x,e) \zeta(e)}{\<\mu(x),\zeta\>} 
\exp\!\big(\frac{- \<\eta,e\>\, \sign(\<\eta+\omega,e\>)}{\alpha/h}\big),
\ \zeta(e) := |\<\eta+\omega,e\>|,\, e \in E.
\end{align}
We only sketch the proof, and denote $u := \<\eta,\cdot\>$. First observe that $\Lambda u(x) = u(x)$, since the considered schemes are consistent on affine maps, 
and since $\cF_x^*(\nabla u(x)) = 1$ by assumption; 
thus $u(y) - \Lambda u(x) = u(y)-u(x) = \<\eta, y-x\>$, and we recognize in \eqref{eq:rho_loc_LF}, \eqref{eq:rho_loc_SL} and \eqref{eq:rho_loc_UD} the exponential term $\exp(\frac{u(y) - \Lambda u(x)}\alpha)$ from \eqref{eq:subadd_diff}.
We next evaluate the partial derivatives $\tfrac{\partial \Lambda u(x)}{\partial u(y)}$, $y\in X$.

For the Lax-Friedrichs scheme, consider a stencil element $y = x+\chi h b_i$,
where $\chi\in \{-1,1\}$ and $1 \leq i \leq d$; then $\tfrac{\partial \Lambda u(x)}{\partial u(y)} = \tfrac{\partial \Lambda_1 u(x)}{\partial u(y)} = - \frac {c_1} d \frac \chi 2 \<\nabla \cF^*(\eta), b_i\> + \frac 1 {2d}$, using successively (\ref{eq:c0C1}, center) and \eqref{eqdef:LFupdate}, thus establishing \eqref{eq:rho_loc_LF}.
For the Semi-Lagrangian scheme, we make the assumption that the angular sectors cover all directions: $\bR^d = \cup_{V \in \bV(x)} V \bR_+^d$, in a \emph{disjoint} manner up to their shared faces, and that the opposite of the group velocity $v_g := \nabla \cF^*_x(\eta)$ is contained in the interior of a unique angular sector, attached to some $V =  (v_1,\cdots,v_d) \in \cV(x)$ which is thus characterized by $V^{-1} v_g < 0$ componentwise; then this sector realizes the minimum in \eqref{eqdef:SLupdate}, and the active stencil elements are therefore $y_i = x + V b_i = x+v_i$, $1\leq i \leq d$, which satisfy $\tfrac {\partial \Lambda u(x)}{\partial u(y_i)} = \xi_i \propto -\<V^{-1} v_g , b_i\>$ using (i) the envelope theorem and (ii) the Lagrange optimality conditions, establishing \eqref{eq:rho_loc_SL}. 
Finally, for the Upwind-Differences scheme, 
substitute $\lambda = u(x)$ in  \eqref{eqdef:UDupdate}, observe that each maximum in \eqref{eqdef:UDupdate} equals $h |\<\eta+\omega,e\>| = h\zeta(e)$ and is attained for the neighbor $y_e = x-\chi h e$ where $\chi := \sign(\<\eta+\omega,e\>)$. In addition, $\tfrac{\partial \Lambda u(x)}{\partial u(y_e)} \propto \mu(x,e) |\<\eta+\omega,e\>|$ by the implicit function theorem, which yields \eqref{eq:rho_loc_UD}.

\section{\Kruzkov-contractivity and causality}
\label{sec:subadd}

We establish in this appendix several results announced in \cref{rem:causality}, which relate the \emph{causality} property of an operator $\Lambda$ with the \emph{$\rho$-sub-additivity} property of $\Phi_\alpha := \rE_\alpha \Lambda_\infty \rE_\alpha^{-1}$ in \Kruzkov\ coordinates. 
Our first step is a derivative-free version of the differential characterization \eqref{eq:subadd_diff} of the latter property.

\begin{proposition}[\leancheck{prop:subadd_local}]
	\label{prop:subadd_local}
Let $\Lambda : \bR^X \to \bR^X$ be monotone and continuous. Denote by $\Phi_\alpha := \rE_\alpha \Lambda_\infty \rE_\alpha^{-1}$ this operator in \Kruzkov\ coordinates, where $\alpha>0$. 
Given $\rho > 0$, the following properties are equivalent: $\Phi_\alpha$ is $\rho$-sub-additive
\begin{align}
	\label{eq:subadd_local}
	 &\Leftrightarrow 
	\forall u \in \bR^X,\ \Lambda\big[u + t \exp\big(\frac u \alpha\big)\big] \leq \Lambda u + \rho t \exp\big(\frac{\Lambda u} \alpha\big) + o(t) \text{ as } t \to 0^+ \\
	\label{eq:subadd_local_x}
	& \Leftrightarrow \forall u \in \bR^X,\ \forall x \in X,\ 
	\Lambda\big[u + t \exp\big(\frac{u-\Lambda u(x)} \alpha\big)\big](x) \leq \Lambda u(x) + \rho t + o(t) \text{ as } t \to 0^+
\end{align}
The equivalence also holds if the remainder $o(t)$ is replaced with $\cO(t^2)$. Note that \eqref{eq:subadd_local_x} involves the mapping $u - \Lambda u(x) := [u(y) - \Lambda u(x)]_{y \in X} \in \bR^X$, where $x \in X$ is fixed.
\end{proposition}

\begin{proof}
We focus on \eqref{eq:subadd_local}, since \eqref{eq:subadd_local_x} is equivalent up to the change of variables $t' = te^{\Lambda u(x)/\alpha}$.
We assume $\alpha = 1$, w.l.o.g.\ up to considering the rescaled operator $\alpha^{-1} \Lambda (\alpha \cdot)$. Then 
\begin{align}
\label{ineq:E1}
	\rE_1 u + e^{-\ve} t \leq \rE_1(u+t e^u) &\leq \rE_1 u +t, &
	 \rE_1^{-1} v - t/v &\leq 
	\rE_1^{-1} (v+t) \leq \rE_1^{-1} v - e^\ve t/v
\end{align}
where $u\in \bR$ and $v<0$ are fixed, $\ve>0$ is arbitrary, and the estimates involving $\ve$ hold for $t\in [0,T_\ve[$ for some $T_\ve>0$ by a Taylor expansion. The other estimates hold unconditionally by concavity of $\rE_1 = -\exp(-\cdot)$ and convexity of $\rE_1^{-1} = -\ln(-\cdot)$.
\emph{Proof of the direct implication}: one has 
\begin{equation}
\label{ineq:E1_Lambda}
	\rE_1 \Lambda (u+t e^u) 
	= \Phi_1(\rE_1(u+t e^u)) 
	\leq \Phi_1(\rE_1 u+t) 
	\leq \Phi_1 \rE_1 u + \rho t
	= \rE_1 \Lambda u + \rho t,
\end{equation}
using successively (i) the defining expression $\Phi_1 = \rE_1 \Lambda_\infty \rE_1^{-1}$, (ii) the monotony of $\Phi_1$ and (\ref{ineq:E1}, left, r.h.s.), (iii) the $\rho$-sub-additivity of $\Phi_1$, (iv) the defining expression of $\Phi_1$, again. Applying (\ref{ineq:E1}, right, r.h.s.) yields \eqref{eq:subadd_local} since $\ve>0$ is arbitrary. \emph{Proof of the reverse implication}: we obtain
\begin{equation*}
	\Phi_1(\rE_1 u+t) \leq \Phi_1(\rE_1 (u+t e^{\ve+u})) = \rE_1 \Lambda(u+t e^{\ve+u}) \leq \rE_1(\Lambda u + \rho t e^{\ve+ \ve'+ \Lambda u}) \leq \rE_1 \Lambda u + \rho t e^{\ve+\ve'},
\end{equation*}
where $\ve,\ve'>0$ are arbitrary, using successively (i) (\ref{ineq:E1}, left) for $t$ small enough, (ii) the defining expression of $\Phi_1$, (iii) \eqref{eq:subadd_local}, for $t\geq 0$ small enough, and (iv) (\ref{ineq:E1}, left, r.h.s). 
This shows that $\Phi_1(v+t) \leq \Phi_1 v + \rho t + o(t)$ as $t \to 0^+$, for all $v \in ]-\infty,0[^X$. Noting that $\Phi_1$ is continuous, since $\Lambda$ is continuous, and that the Lipschitz property is local (specifically, a continuous function with non-positive upper right Dini derivative is non-increasing), we obtain $\Phi_1(v+t) \leq \Phi_1 v + \rho t$ for all $v<0$ and $t \geq 0$, such that $v+t < 0$. The restrictions $v<0$ and $v+t<0$ are then lifted by noting that $\Phi_1 v = \sup\{\Phi_1 w\mid  w\leq v \text{ and } w<0\}$, $v\in \bR^X$, and the result follows.
\end{proof}

The continuity assumption in \cref{prop:subadd_local}, used for the reverse implication, follows in \cref{prop:subadd_alpha_beta,prop:strict_causal_implies_subadd,prop:causal_implies_subadd} from the scheme monotonicity and $1$-sub-additivity, since these imply the $1$-Lipschitz property by \cref{prop:contractivity}, and from the scheme concavity in \cref{corol:subadd_concave}.

We can slightly simplify the characterization \eqref{eq:subadd_local} if the scheme $\Lambda$ is concave w.r.t.\ $u\in \bR^X$.
The semi-Lagrangian scheme \eqref{eqdef:SLupdate} is concave in this sense, since it is defined as a minimum of affine functions. The finite-element scheme \eqref{eqdef:FE}, and the upwind-differences scheme \eqref{eqdef:UDupdate}, are therefore also concave, since they are special cases of the former \cite[Appendix B]{mansour2024solving}. The Lax-Friedrichs scheme \eqref{eqdef:LFupdate} is concave if the Hamiltonian $\cF^*_x$ is convex for all $x \in X$, but this scheme is not concave in the general case of a non-convex Hamiltonian.

\begin{corollary}[\leancheck{corol:subadd_concave}]
\label{corol:subadd_concave}
Let $\Lambda : \bR^X \to \bR^X$ be monotone and such that $u \in \bR^X \mapsto \Lambda u(x) \in \bR$ is concave, for all $x \in X$. Then one has the equivalence, where $\alpha,\rho> 0$ and $\Phi_\alpha := \rE_\alpha \Lambda_\infty \rE_\alpha^{-1}$
\begin{equation*}
	\Phi_\alpha  \text{ is $\rho$-sub-additive}
	\ \Leftrightarrow \
	\forall u \in \bR^X,\ \forall t \geq 0,\ \Lambda\big[u + t \exp(\frac u \alpha)\big] \leq \Lambda u + \rho t \exp(\frac{\Lambda u} \alpha).
\end{equation*}
\end{corollary}
\begin{proof}
	Consider a concave $f : \bR \to \bR$ satisfying $f(a+t) \leq f(a)+\rho t + o(t)$ as $t \to 0^+$, where $a,\rho \in \bR$ are fixed.
	Given $0 < t \leq T$, one has $(f(a+T) - f(a))/T \leq (f(a+t)-f(a))/t \leq \rho + o(1)$ as $t \to 0^+$, by concavity, hence $f(a+T) - f(a) \leq \rho T$ as desired. The result follows by choosing $f(t):=\Lambda[u+t\exp(\frac{u-\Lambda u(x)} \alpha)](x)$ from \eqref{eq:subadd_local_x}.
\end{proof}

Our next result shows that the $\alpha$-\Kruzkov-contractivity property, defined as: ($\exists \rho\in [0,1[$, $\Phi_\alpha$ is $\rho$-sub-additive), weakens as $\alpha$ increases, as announced in \cref{rem:causality}. 

\begin{proposition}[\leancheck{prop:subadd_alpha_beta}]
\label{prop:subadd_alpha_beta}
Let $\Lambda : \bR^X \to \bR^X$ be monotone and $1$-sub-additive. If $\Phi_\alpha$ is $\rho$-sub-additive for some $\alpha,\rho>0$, then $\Phi_\beta$ is $\sigma$-sub-additive for any $\beta \geq \alpha$ and $\sigma \geq (1-\frac \alpha \beta) + \rho \frac \alpha \beta$. 
(Note in addition that $(1-\frac \alpha \beta) + \rho \frac \alpha \beta<1$ if $\rho<1$.)
\end{proposition}

\begin{proof}
	Denoting $a:=1/\alpha$, $b:=1/\beta$ and $s:=\alpha/\beta \in [0,1]$, we obtain	
	\begin{equation*}
		\Lambda\big[u + t \exp(b (u - \Lambda u(x)) )\big](x) 
		\leq \Lambda\big[u + t (s \exp(a (u - \Lambda u(x)) ) + (1-s))\big](x) 
		\leq \Lambda u(x) + \rho t s + t(1-s) + o(t),
	\end{equation*}
	as $t\to 0^+$. We used successively (i) the convexity of the exponential applied to $b = s a + (1-s) \times 0$, and the monotonicity of $\Lambda$, and (ii) the 1-sub-additivity of $\Lambda$ and \eqref{eq:subadd_local_x}. Using again the equivalence \eqref{eq:subadd_local_x}, we obtain that $\Phi_\beta$ is $(\rho s + (1-s))$-sub-additive, and the result follows.
\end{proof}

\paragraph{Relation with the causality property.}
Our final result \cref{prop:causal_implies_subadd} shows that the causality property implies the $\rho$-sub-additivity of $\Phi_\alpha$, for any $\alpha>0$ and with an explicit constant $\rho = \rho(\alpha)$, under mild additional assumptions, as announced in \cref{rem:causality}. This generalizes \cite[Proposition 2.25]{mansour2024solving}, which establishes a similar result in the special case of the semi-Lagrangian scheme. As a preliminary, we consider in \cref{prop:strict_causal_implies_subadd} the stronger assumption of strict-causality.

\begin{proposition}[\leancheck{prop:strict_causal_implies_subadd}]
\label{prop:strict_causal_implies_subadd}
	Let $\Lambda : \bR^X \to \bR^X$ be $\delta$-strictly-causal, monotone, and $1$-sub-additive. 
	Then $\Phi_\alpha$ is $\rho := \exp(-\delta/\alpha)$-sub-additive, for any $\alpha>0$. 
\end{proposition}

\begin{proof}
We assume $\alpha=1$, up to considering the operator $\alpha^{-1} \Lambda (\alpha\cdot)$, which is $(\delta/\alpha)$-strictly-causal, monotonic and $1$-sub-additive.
Observe that $\Lambda [u + t e^u] \leq \Lambda u + C t$ with $C := \exp(\max u)$, by monotonicity and $1$-sub-additivity. 
	Hence, for any $\eta < \delta$ and any $0 \leq t \leq (\delta-\eta)/C$,
	\begin{align*}
		\Lambda[u+te^u](x) &= \Lambda \min \{u + t e^u,\, \Lambda u(x) -\eta\} (x)\\
		&\leq \Lambda [u+t e^{\Lambda u(x)-\eta}](x) \\
		&\leq \Lambda u(x) + t e^{-\eta} e^{\Lambda u(x)}.
	\end{align*}
	We used successively (i) the $\delta$-causality property, with parameter $t' := \Lambda u(x)-\eta$, noting that $t'+\delta \geq \Lambda u(x) + Ct \geq  \Lambda [u+te^u](x)$ by construction, (ii) monotony and 
	$\min \{u + t e^u,\, \Lambda u(x) -\eta\} \leq u+t e^{\Lambda u(x)-\eta}$ by distinguishing the cases $u \leq \Lambda u(x)-\eta$ and $u> \Lambda u(x)-\eta$, (iii) the $1$-sub-additivity of $\Lambda$. 
	Thus $\Phi_1$ is $\exp(-\eta)$-sub-additive, by \cref{prop:subadd_local} \eqref{eq:subadd_local}, and letting $\eta\to \delta$ we conclude the proof.
\end{proof}

\begin{proposition}[\leancheck{prop:causal_implies_subadd}]
\label{prop:causal_implies_subadd}
Let $\Lambda : \bR^X \to \bR^X$ be monotone, causal, additively invariant, 
$\nu$-sub-multiplicative, and have $C_\past$-memory, for some $C_\past,\nu>0$.
Then $\Phi_\alpha := \rE_\alpha \Lambda_\infty \rE_\alpha^{-1}$ is $\rho$-sub-additive for any $\alpha>0$, where $\rho \in ]0,1[$ is defined as 
\begin{equation}
\label{eqdef:rho_causal}
\rho := 1- \frac {\nu}{C_\past} \Big(1-\exp(-\frac{C_\past} \alpha)\Big).
\end{equation} 
\end{proposition}

\begin{proof}
Clearly $\rho < 1$, and noting that $\nu \leq \Lambda[0] \leq C_\past$, see \cref{prop:cst_subsol} and the remark below it, we find that $\rho>0$.
We assume that $\alpha = 1$, w.l.o.g.\ up to considering the operator $\alpha^{-1} \Lambda (\alpha\cdot)$ which is monotone, causal, additively invariant, $(\nu/\alpha)$-sub-multiplicative, and has $(C_\past/\alpha)$-memory. 
Our objective is to establish \eqref{eq:subadd_local_x}, and for that purpose we fix $u \in \bR^X$ and $x\in X$.
We assume that $x \in \Int X$, since the case $x \in \partial X$ is trivial in view of the Dirichlet boundary conditions, which are always implicitly assumed.
We further assume that $\Lambda u(x) = 0$, w.l.o.g.\ by additive invariance up to adding a constant to $u$.
	Let us introduce 
	\begin{align*}
		c &:= 1-\exp(-C_\past),& 
		C & := \tfrac {c}{C_\past} \exp(-C_\past) 
	\end{align*}\vspace{-5mm}
	\begin{align*}
		u_t &:= u+t \exp(u), &
		r_t &:= \frac{c u}{C_\past} + 1 + C t &
		v_t &:= u + t\, r_t = \Big(1+\frac{c t}{C_\past}\Big) u + t + C t^2.
	\end{align*}
	Note that $0 = \Lambda u(x) \leq \Lambda u_t(x) \leq \Lambda u(x) + t \exp(\max(u))$, for all $t \geq 0$. We consider $t > 0$ small enough such that $\Lambda u_t(x) \leq u(y)$ for all $y\in X$ such that $u(y)>0$. We claim that
	\begin{align}
	\label{ineq:vt_ut}
		v_t &\geq u_t \text{ on } \{-C_\past \leq u_t \leq \Lambda u_t(x)\}, &
		v_t \geq \Lambda u_t(x) \text{ on } \{ u_t \geq \Lambda u_t(x)\}.
	\end{align}
	Indeed, consider an arbitrary $y\in X$, denote $s := -u(y)$, and observe that for any $t>0$
	\begin{align}
	\label{ineq:uv_st}
		u_t(y) &\leq v_t(y) \Leftrightarrow e^{-s} + \tfrac {c}{C_\past} s \leq 1+C t,&
		-C_\past &\leq u_t(y) \Leftrightarrow s - C_\past \leq t e^{-s}.
	\end{align}
	We distinguish three cases:
	\begin{itemize}
		\item If $s < 0$, then $u(y) > 0$, thus $u(y) \geq \Lambda u_t(x)$ and therefore both $u_t(y) \geq \Lambda u_t(x)$ and $v_t(y) \geq \Lambda u_t(x)$, fulfilling (\ref{ineq:vt_ut}, right).
		\item If $0 \leq s \leq C_\past$, then $e^{-s} + \frac {c} {C_\past} s \leq 1$ by convexity of the l.h.s., since this inequality holds for $s\in \{0,C_\past\}$. This implies (\ref{ineq:uv_st}, left), which concludes in both cases of \eqref{ineq:vt_ut}.
		\item If $s \geq C_\past$, we assume that $u_t(y) \geq -C_\past$, w.l.o.g.\  since otherwise $y$ belongs to neither set in \eqref{ineq:vt_ut}, and obtain that  
		\begin{equation*}
			e^{-s} + \tfrac {c}{C_\past} s \leq e^{-C_\past} + c + \tfrac {c} {C_\past} (s-C_\past) \leq 1+C t,
		\end{equation*}
		using successively (i) $e^{-s} \leq e^{-C_\past}$ and a rearrangement of terms, and (ii) the  expression of $c$ and  (\ref{ineq:uv_st}, right).
		This yields (\ref{ineq:uv_st}, left), which concludes in both cases of \eqref{ineq:vt_ut}.
	\end{itemize}
	Having now established \eqref{ineq:vt_ut}, we obtain 	
	\begin{equation}
	\label{ineq:Lambda_ut_vt}
		\Lambda u_t(x) 
		\leq \Lambda v_t(x) \leq (1+\frac{c t}{C_\past}) \Lambda u(x) - \frac {\nu c t}{C_\past} + t + C t^2  
		= \rho t + C t^2,
	\end{equation}
	using successively (i) \cref{corol:mon_min_max} applied to $u_t$ and $v_t$ with $t_0 = -C_\past$ and $t_1 = \Lambda u_t(x)$, noting that $t_0 + C_\past = 0 = \Lambda u(x) \leq \Lambda u_t(x) = t_1$ as required,
	 and (ii) the additive invariance and $\nu$-sub-multiplicativity of $\Lambda$, as defined in \cref{rem:submult_strictcaus}, and (iii) the assumption $\Lambda u(x) = 0$ and the defining expression \eqref{eqdef:rho_causal} of $\rho$.
	The $\rho$-sub-additivity of $\Phi_\alpha$ then follows by \eqref{eq:subadd_local_x}, \cref{prop:subadd_local}, which concludes.
\end{proof}

\section{Sensitivity of the numerical cost to the parameters $(\alpha,\rho,\tau,\ve)$}
\label{sec:num_secondary}
\begin{figure}
\def\w{3.7cm}\def\W{5cm}
\includegraphics[width=\w]{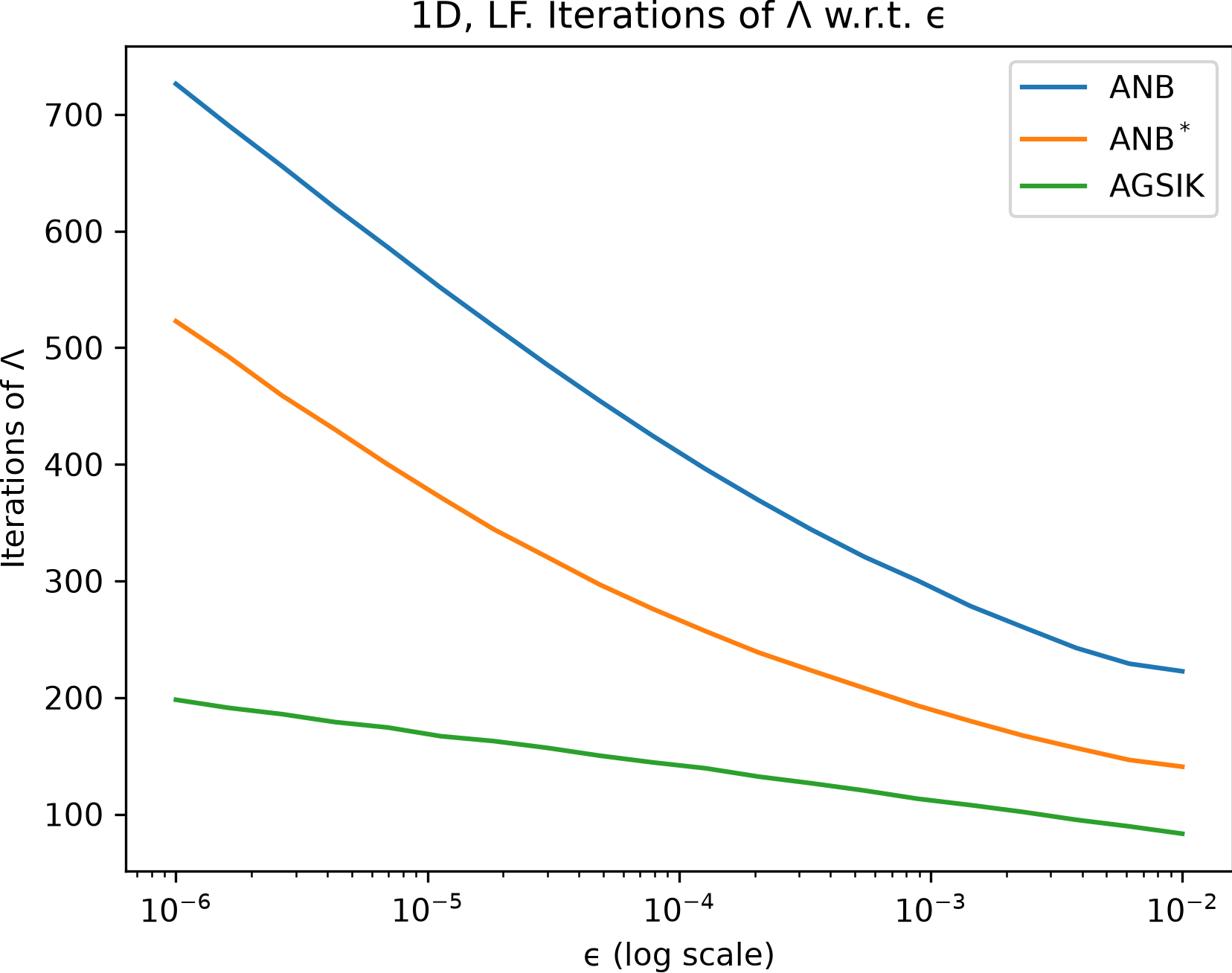}
\includegraphics[width=\w]{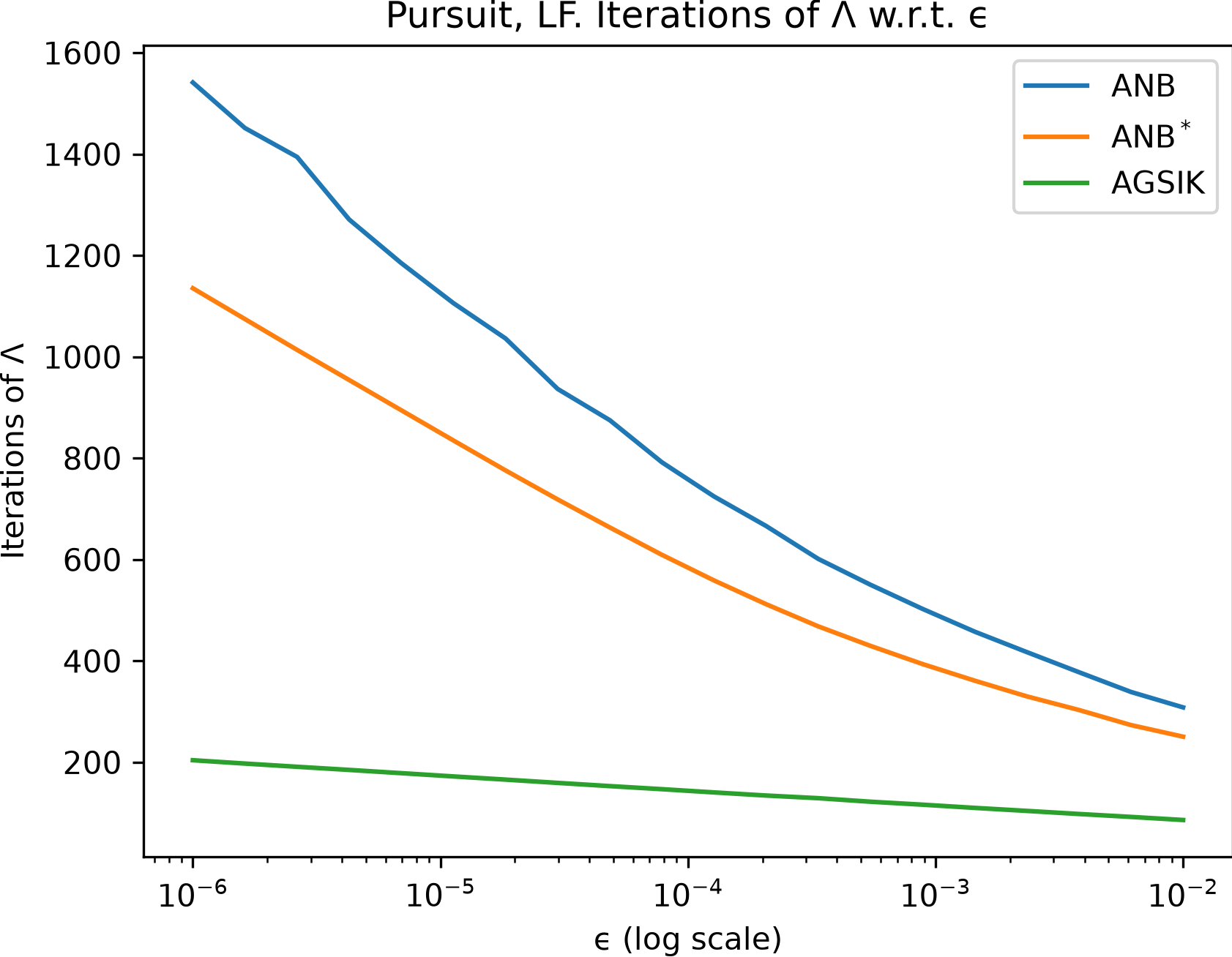}
\includegraphics[width=\w]{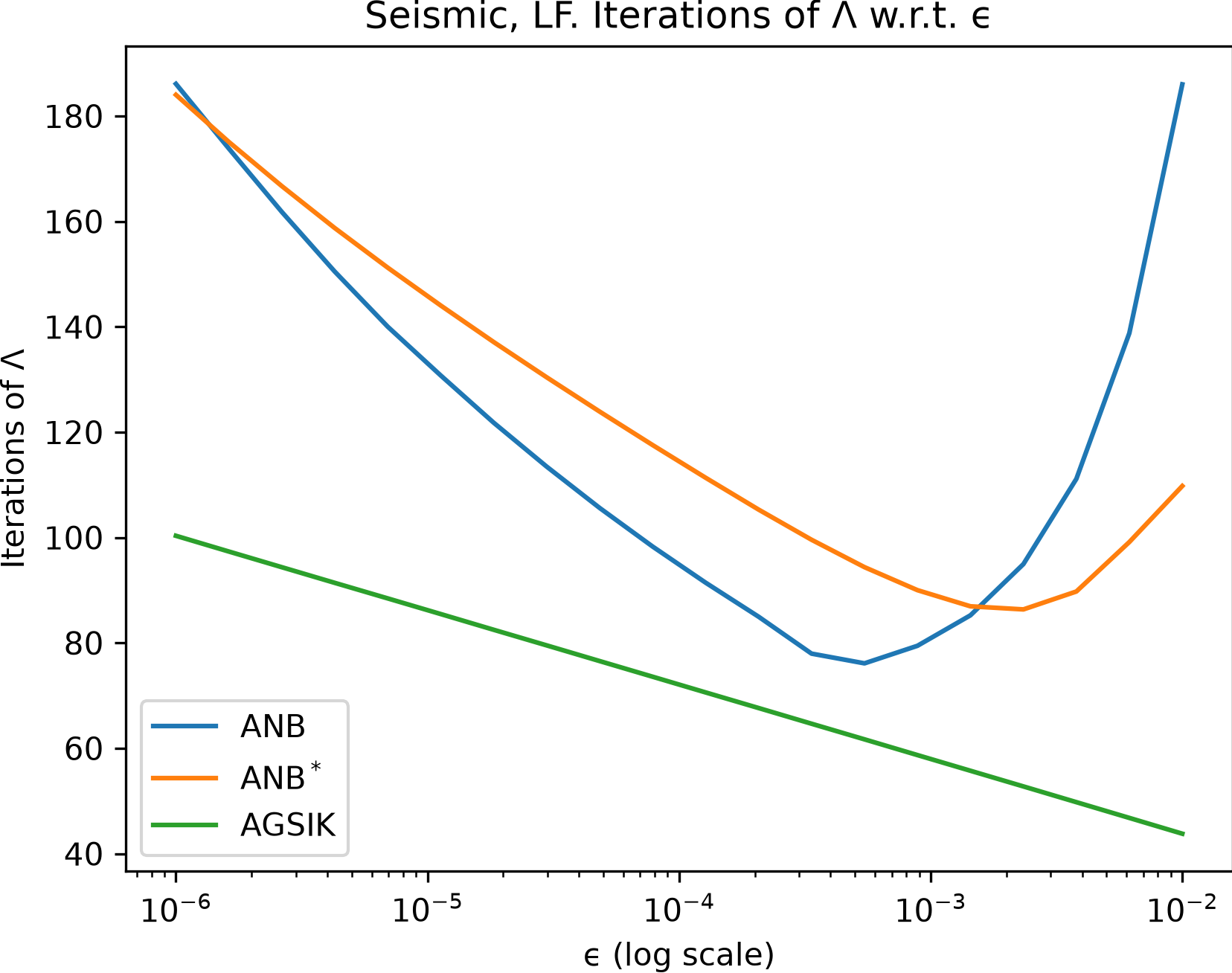}
\includegraphics[width=\w]{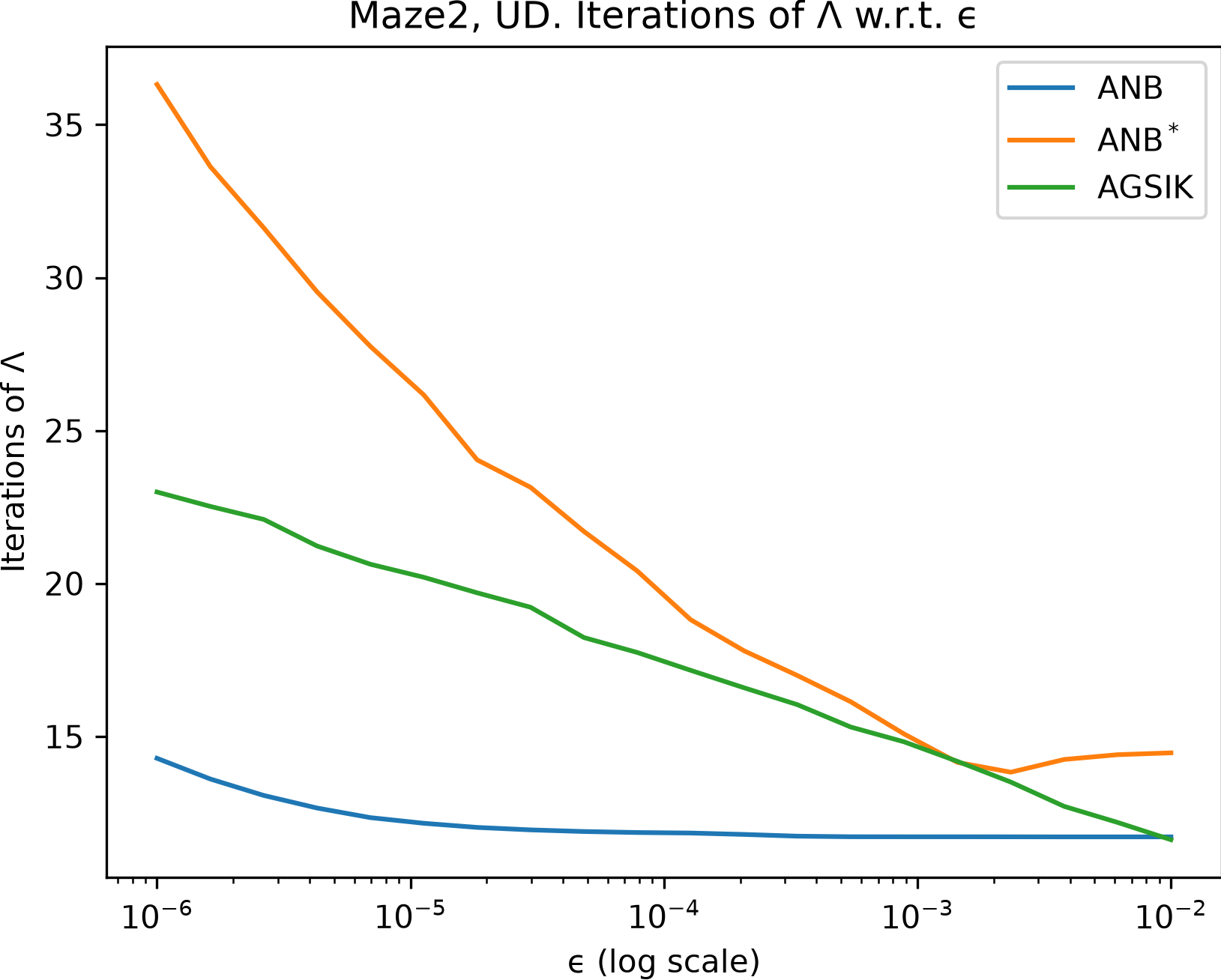}\\
\includegraphics[width=\W]{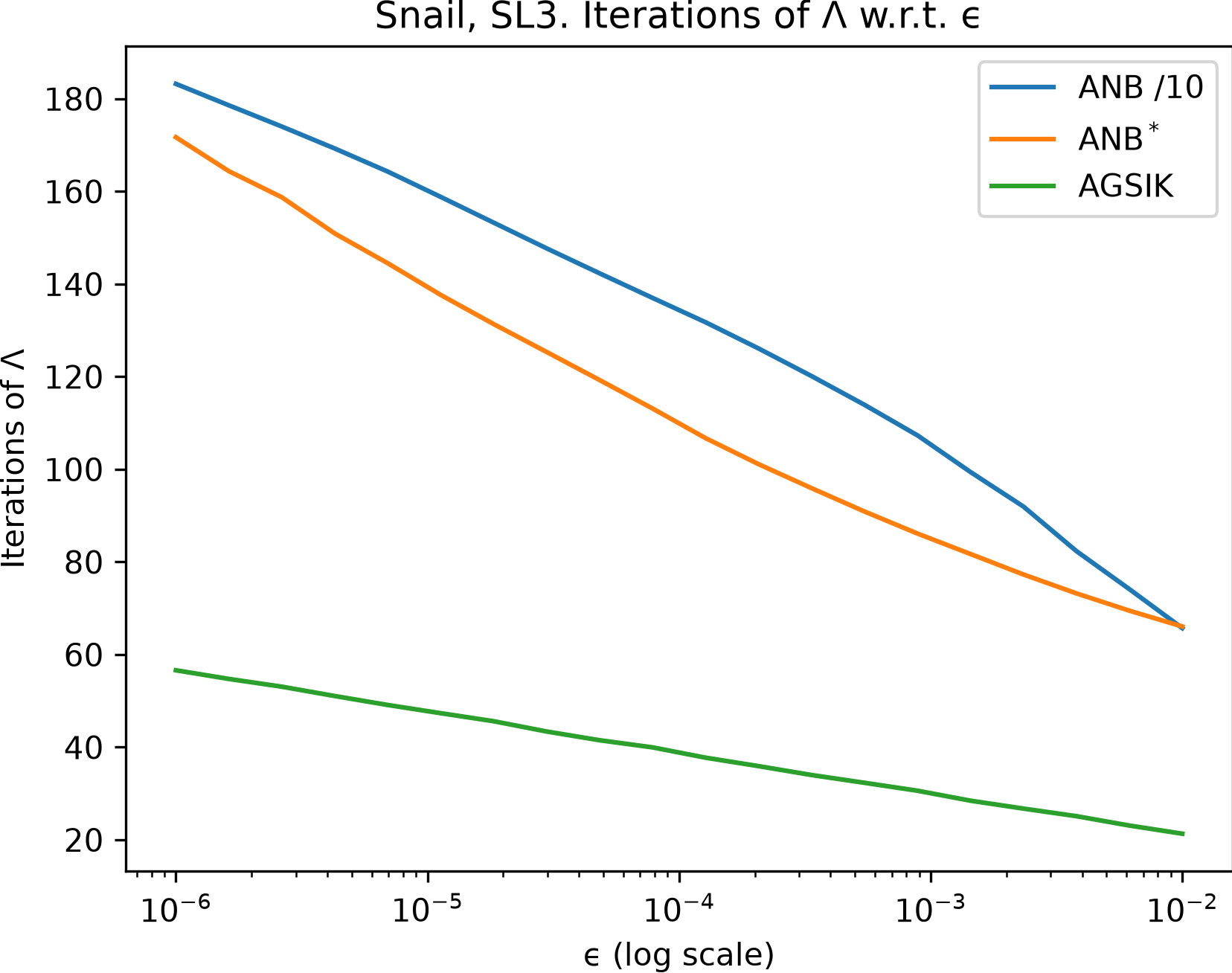}
\includegraphics[width=\W]{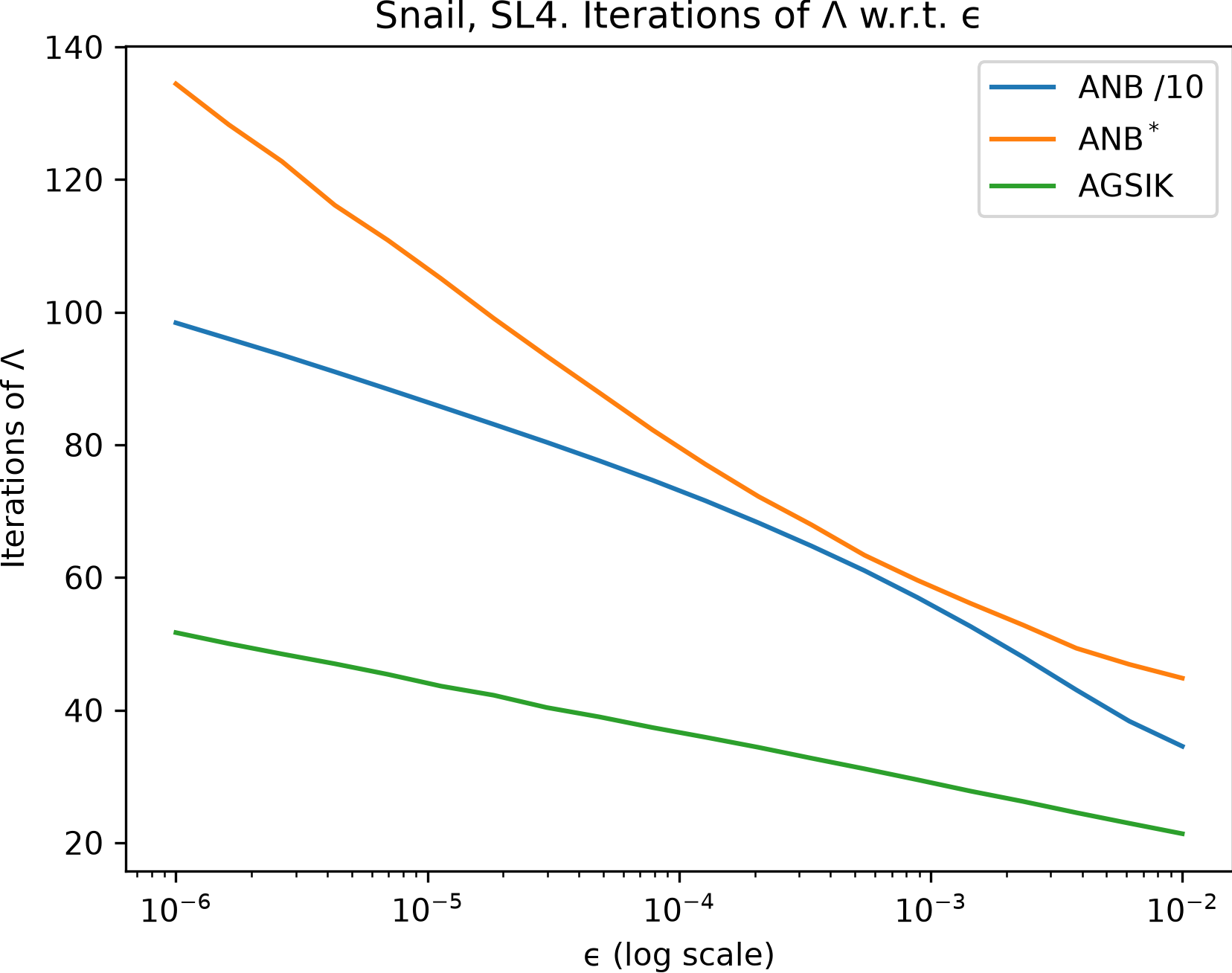}
\includegraphics[width=\W]{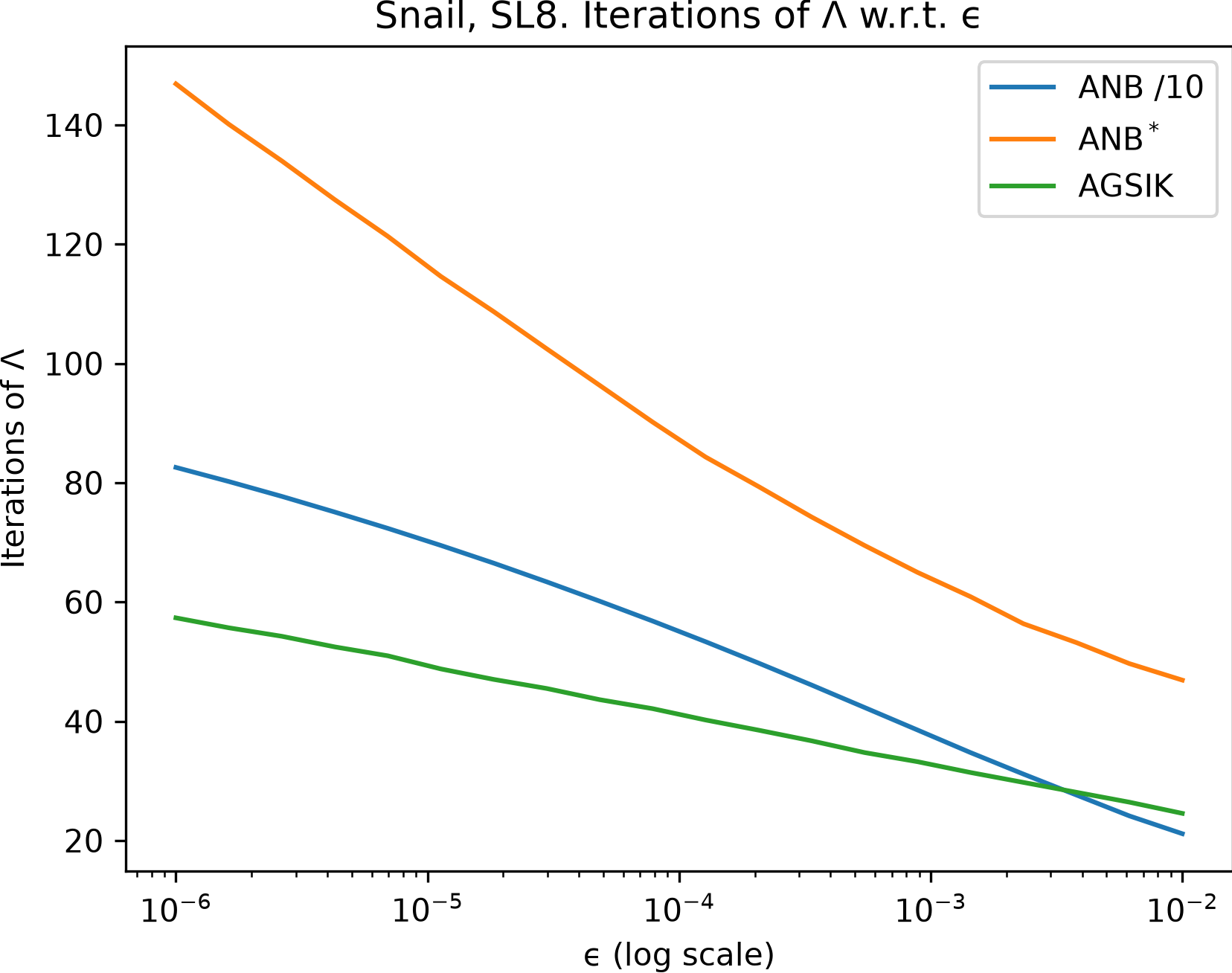}\\
\includegraphics[width=\W]{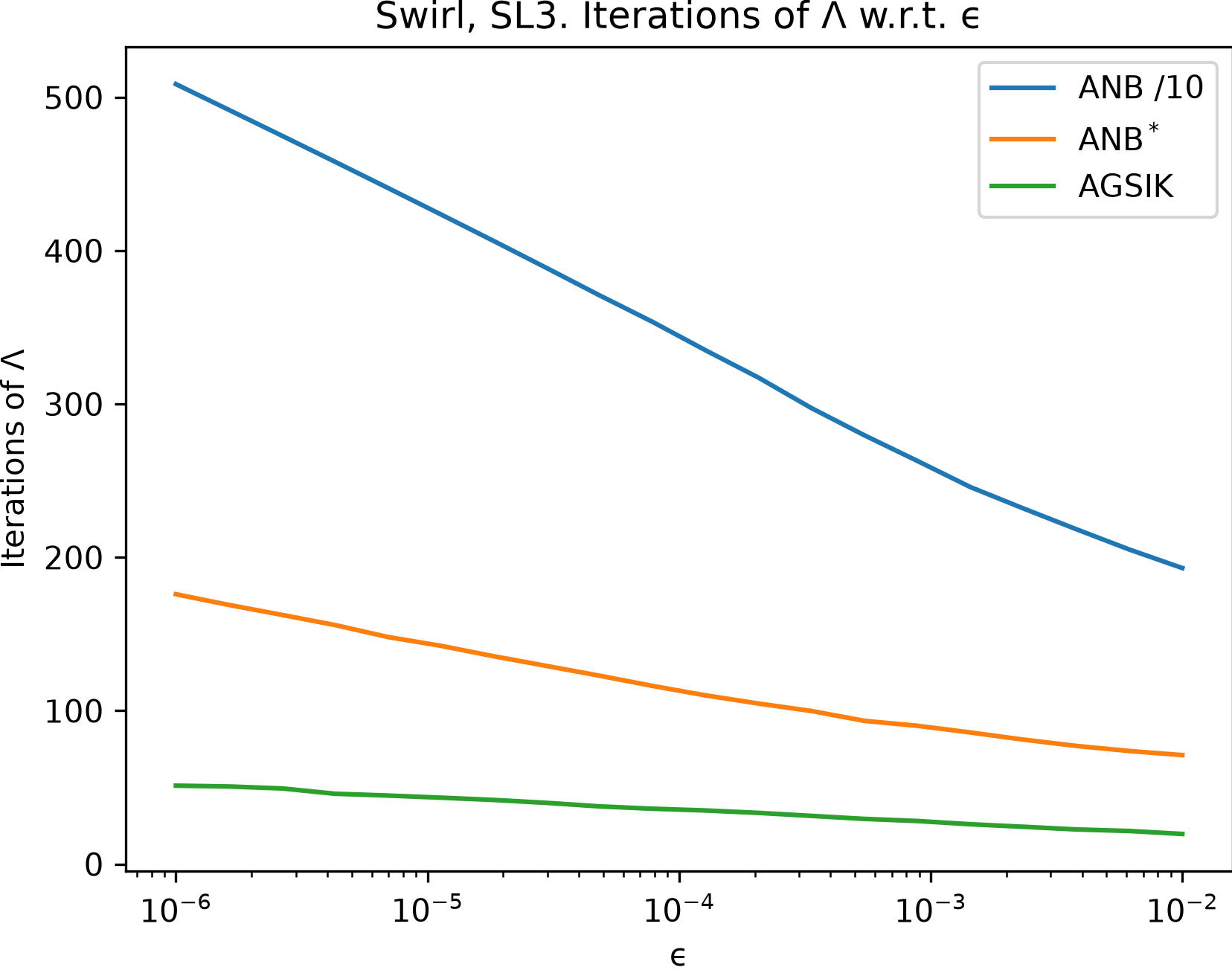}
\includegraphics[width=\W]{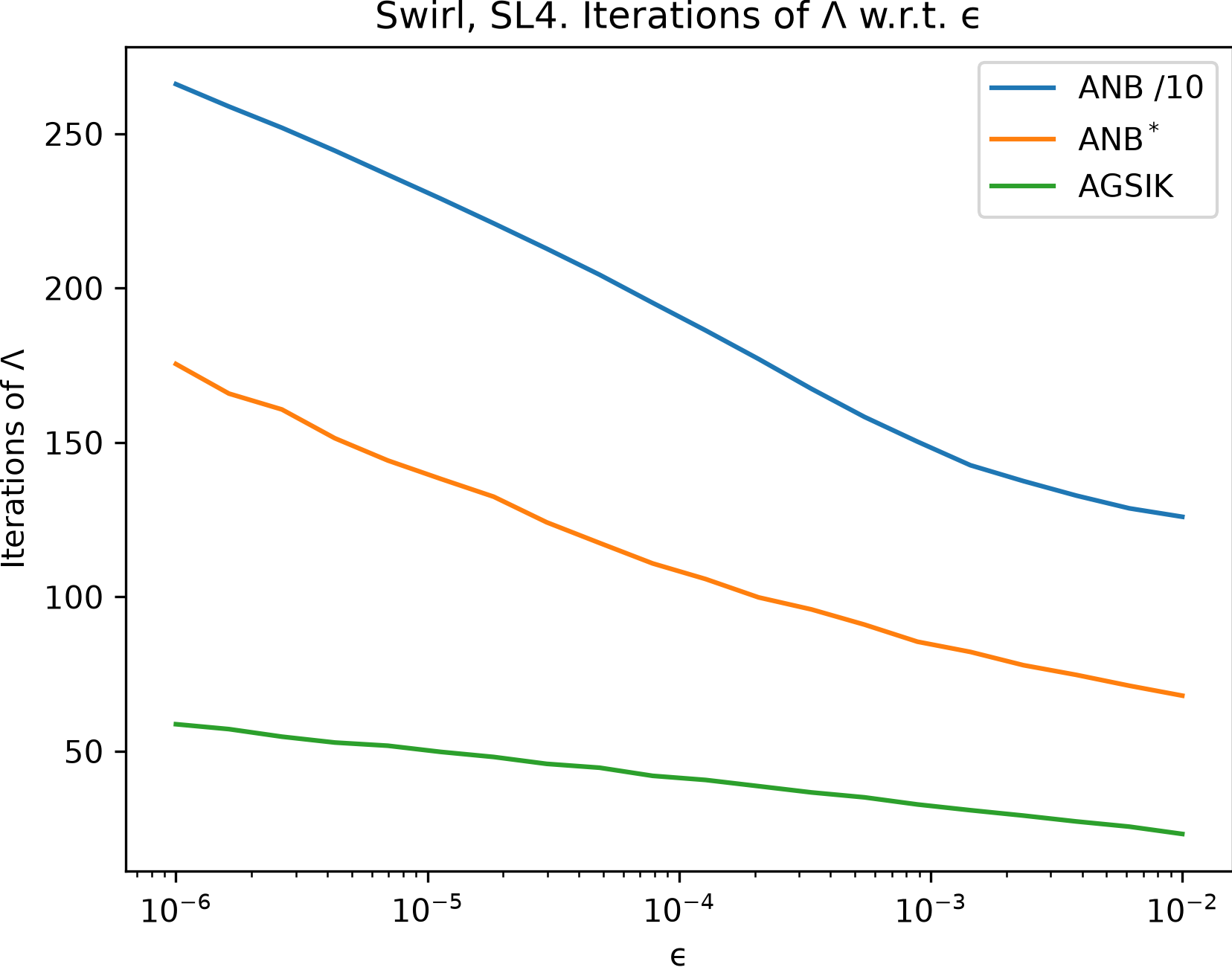}
\includegraphics[width=\W]{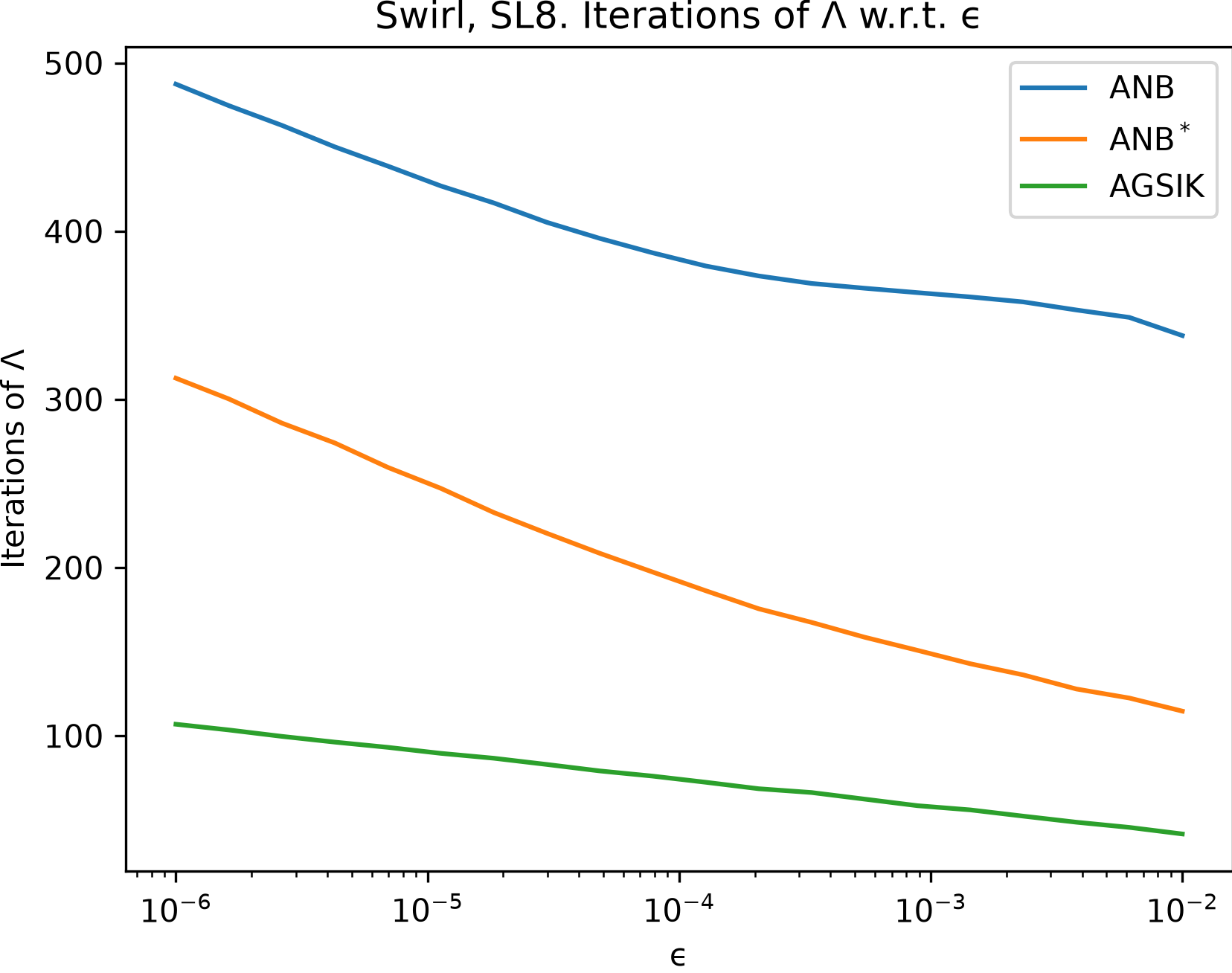}

\caption{Average number of iterations $N_\Lambda(\ve)$ of the update operator $\Lambda$ at each point of the domain, w.r.t.\ the tolerance $\ve$ (log scale).
Our complexity estimates \cref{th:complexity,th:complexity_withUnknown} of the ANB show that $N_\Lambda(\ve) = \cO(\ln^2 \ve)$ as $\ve \to 0^+$, yet in most cases we empirically observe $N_\Lambda(\ve) = \cO(\ln \ve)$.
In the Seismic test case, $N_\Lambda(\ve)$ also increases when the tolerance $\ve$ is \emph{excessively large}. 
Indeed, our complexity estimate \cref{th:complexity} assumes an upper bound $\ve < \ve_{\max} = \Theta(h)$, see \eqref{eq:bigO_h_ln_eps}; if the tolerance is too large, then the front tracking is not effective, see \cref{rem:tol_range}.
}
\label{fig:iter_wrt_eps}
\end{figure}

\begin{figure}
\def\w{3.7cm}\def\W{5cm}
\includegraphics[width=\w]{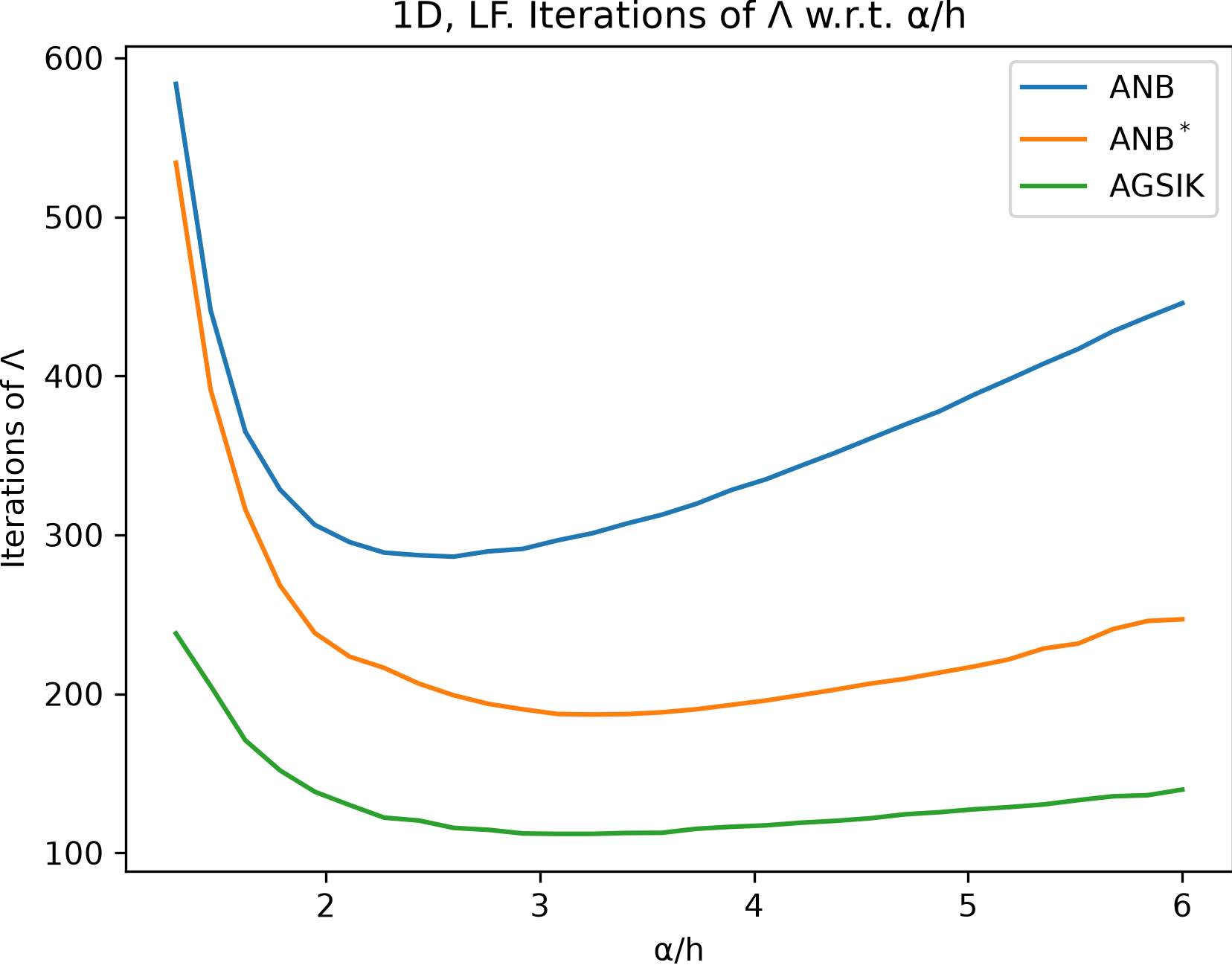}
\includegraphics[width=\w]{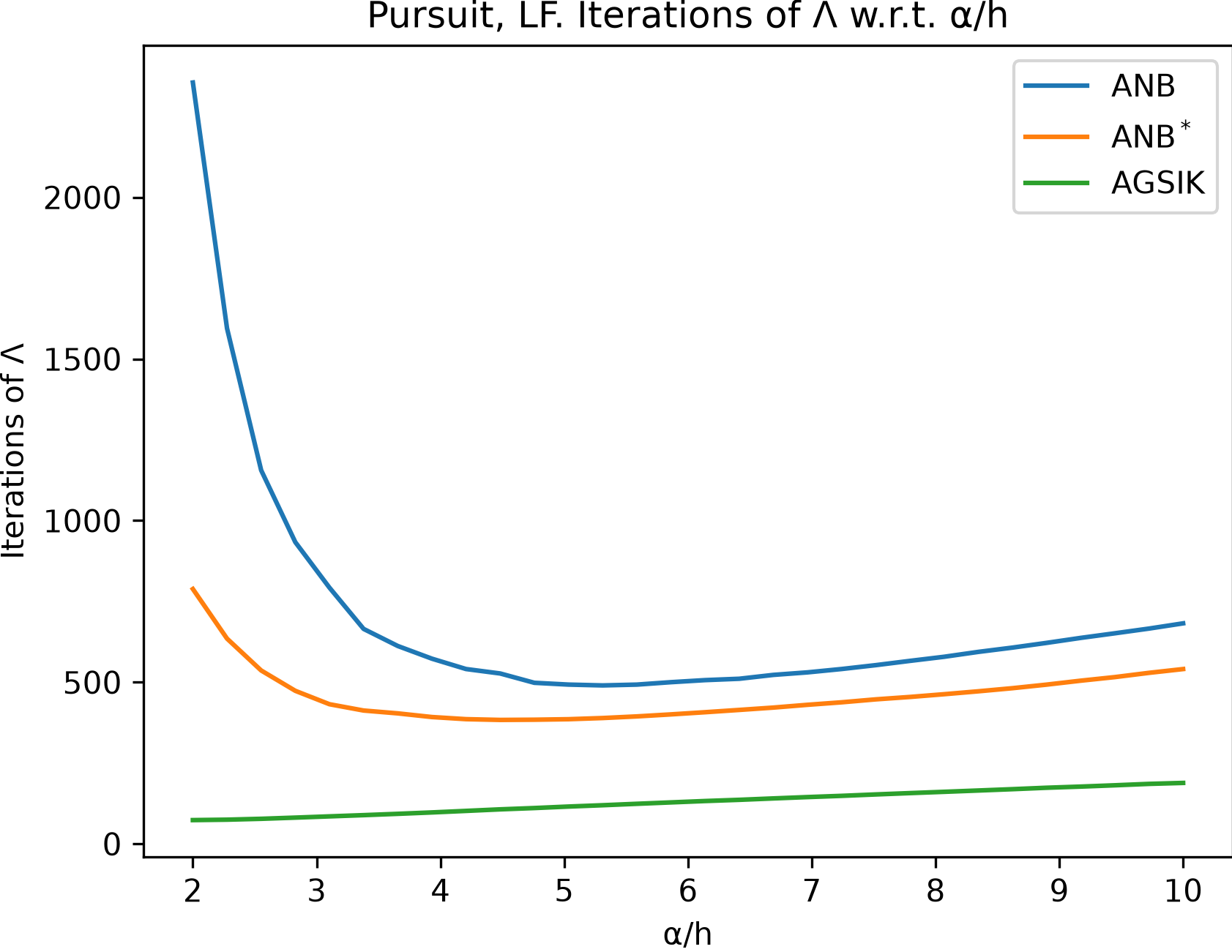}
\includegraphics[width=\w]{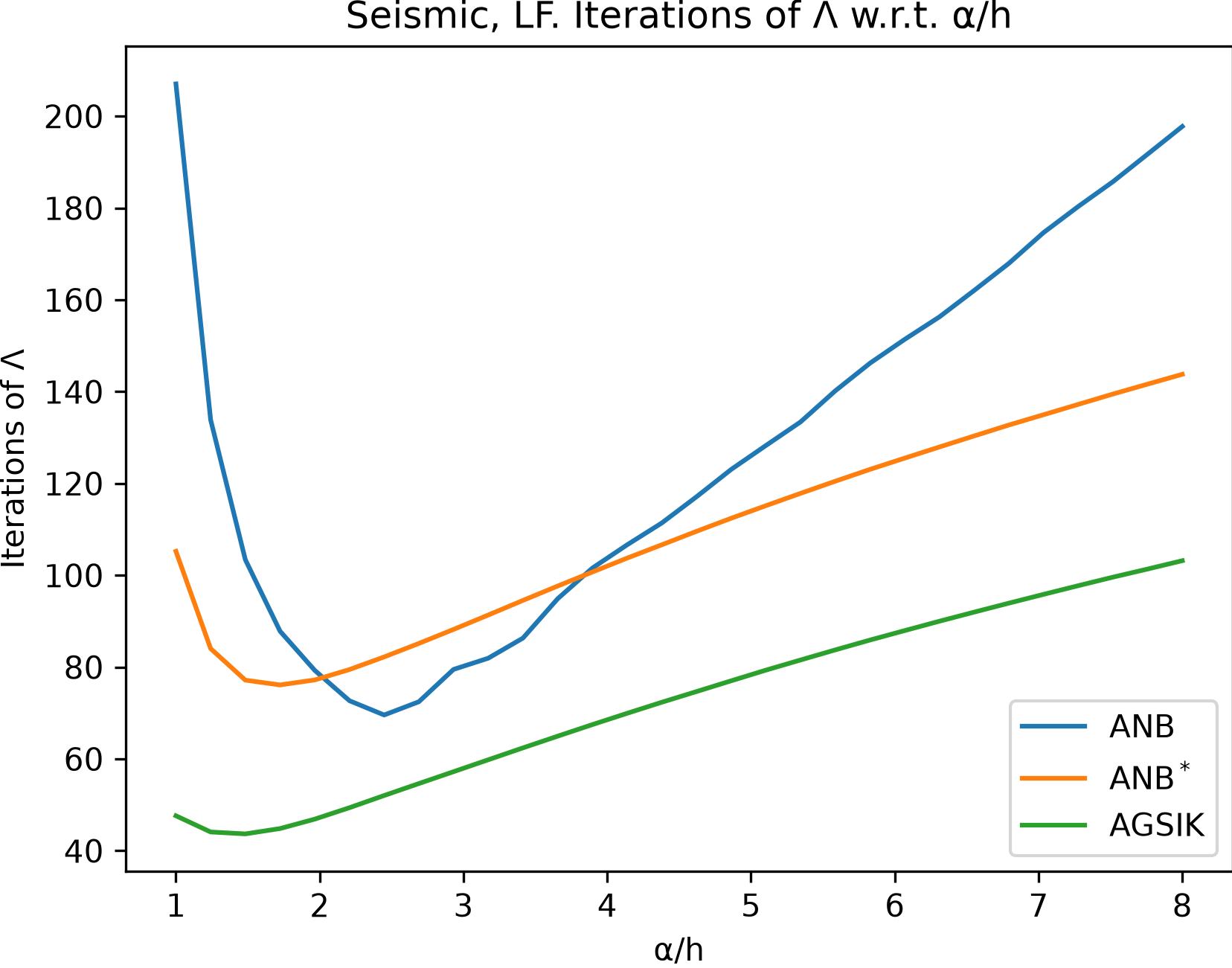}
\includegraphics[width=\w]{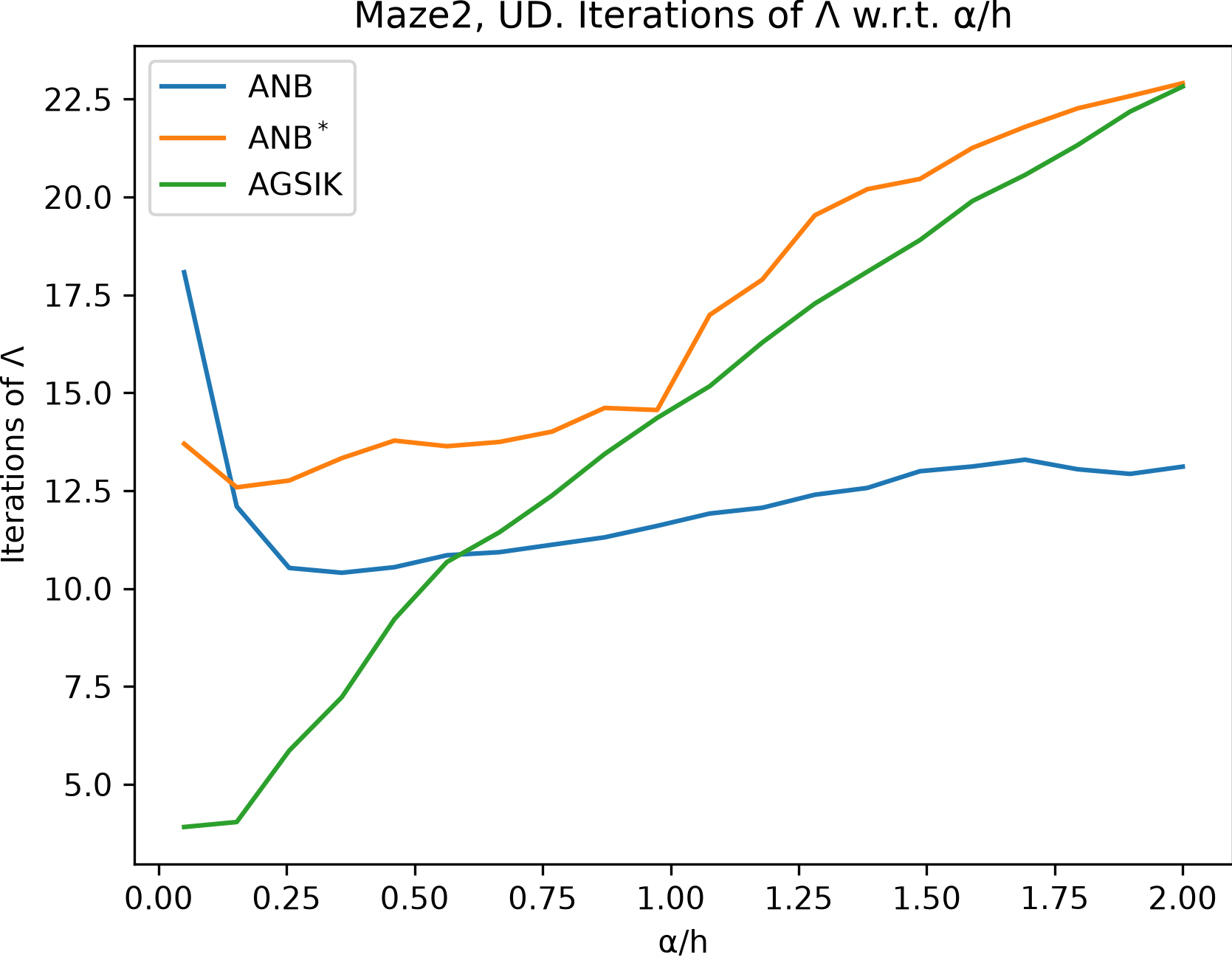}\\
\includegraphics[width=\W]{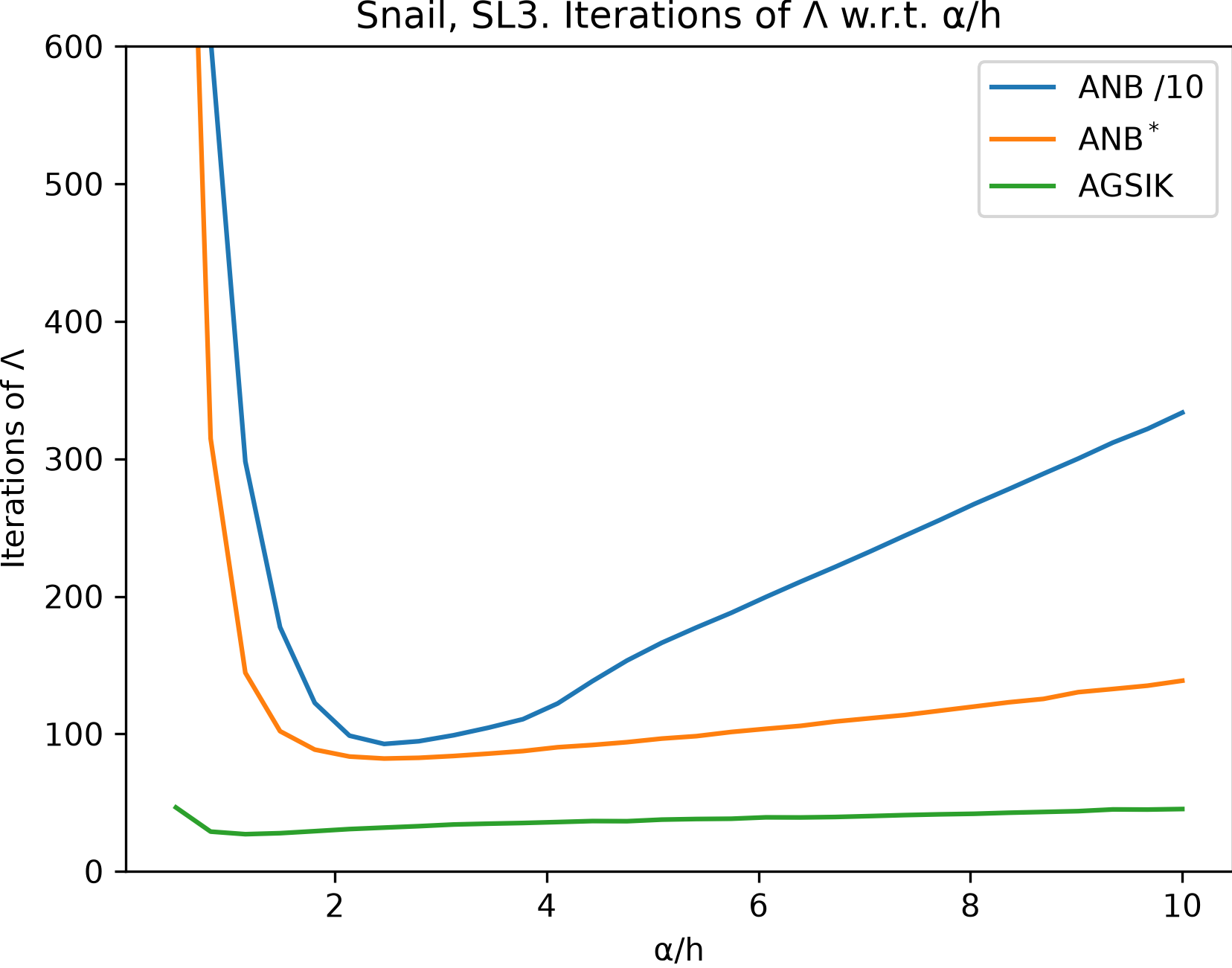}
\includegraphics[width=\W]{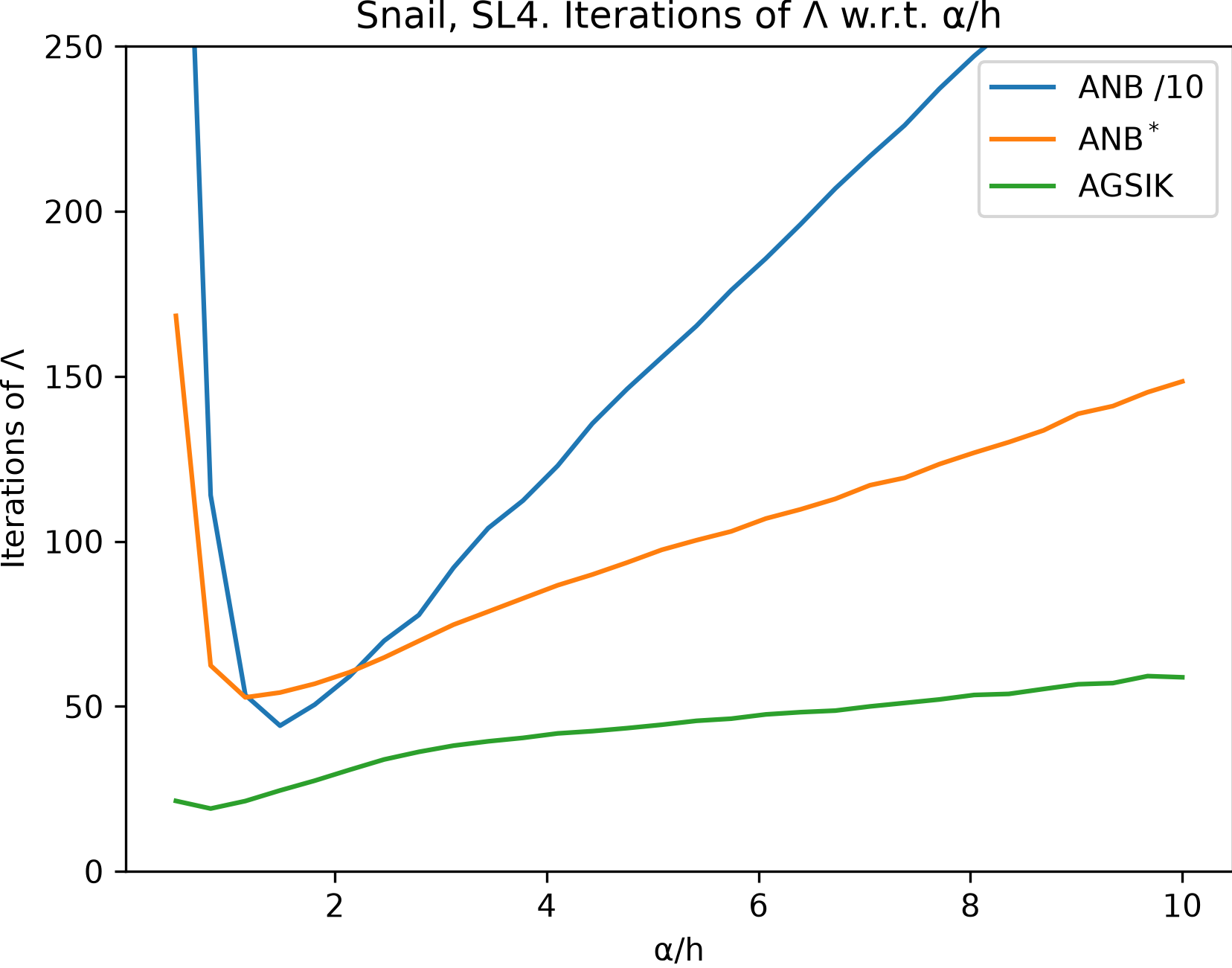}
\includegraphics[width=\W]{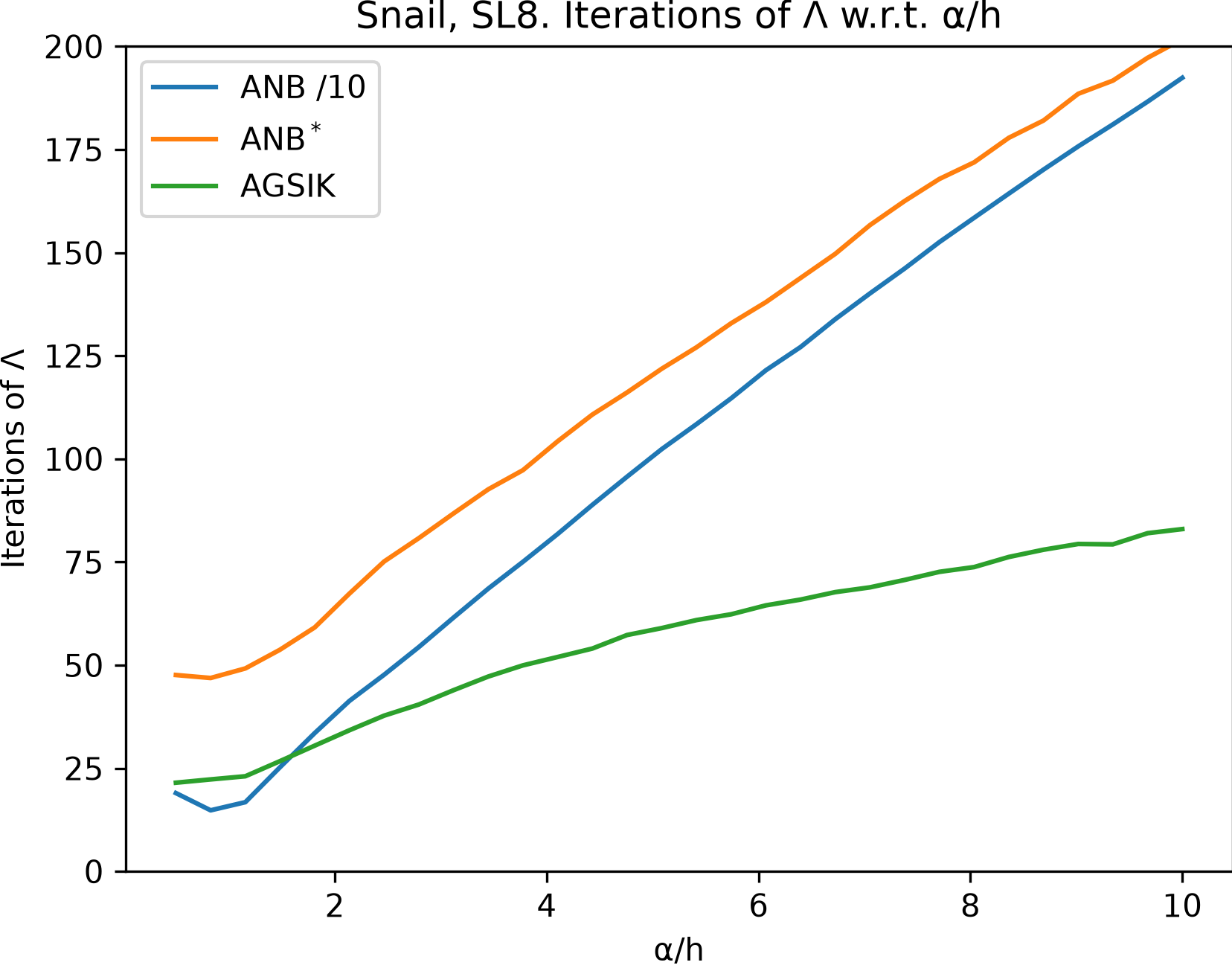}\\
\includegraphics[width=\W]{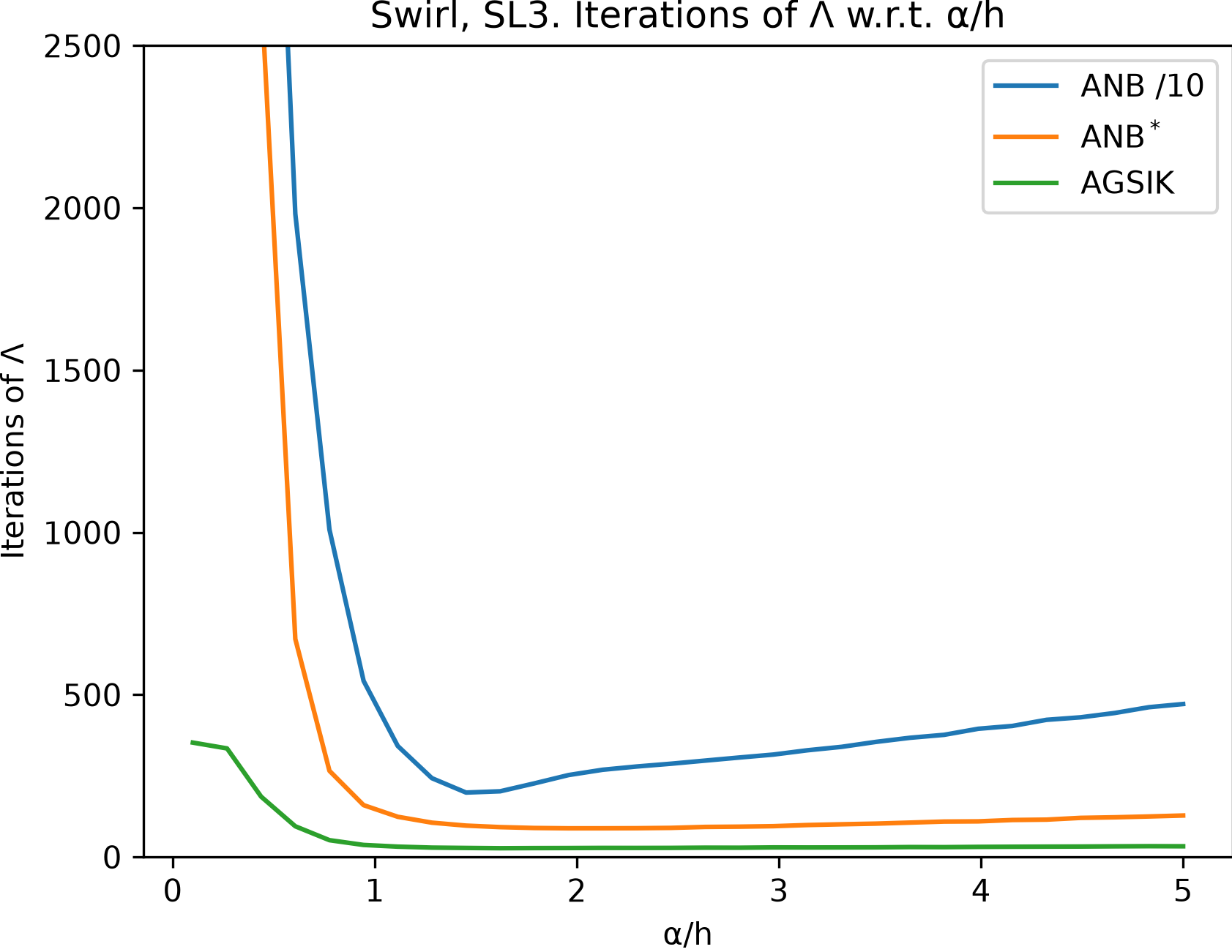}
\includegraphics[width=\W]{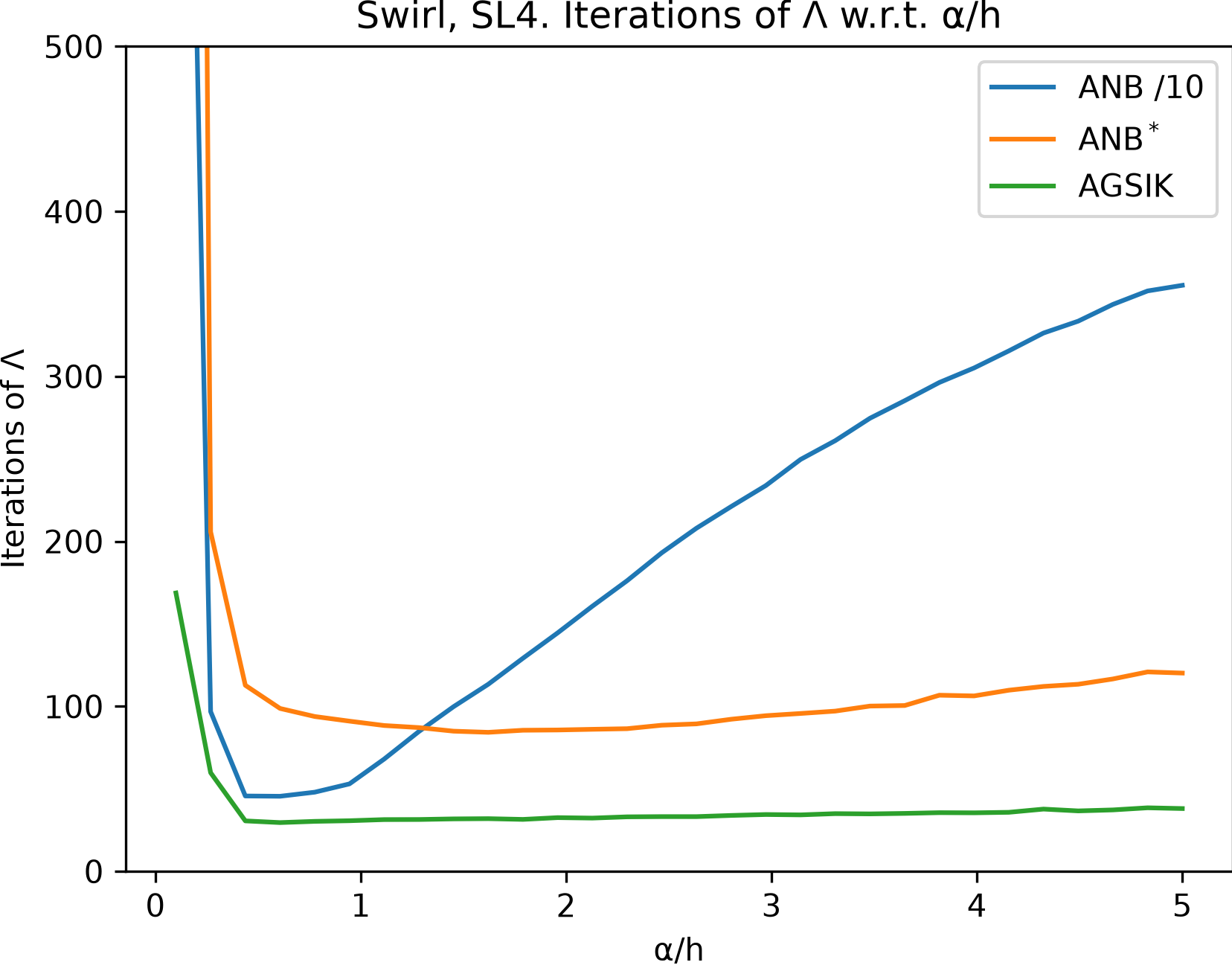}
\includegraphics[width=\W]{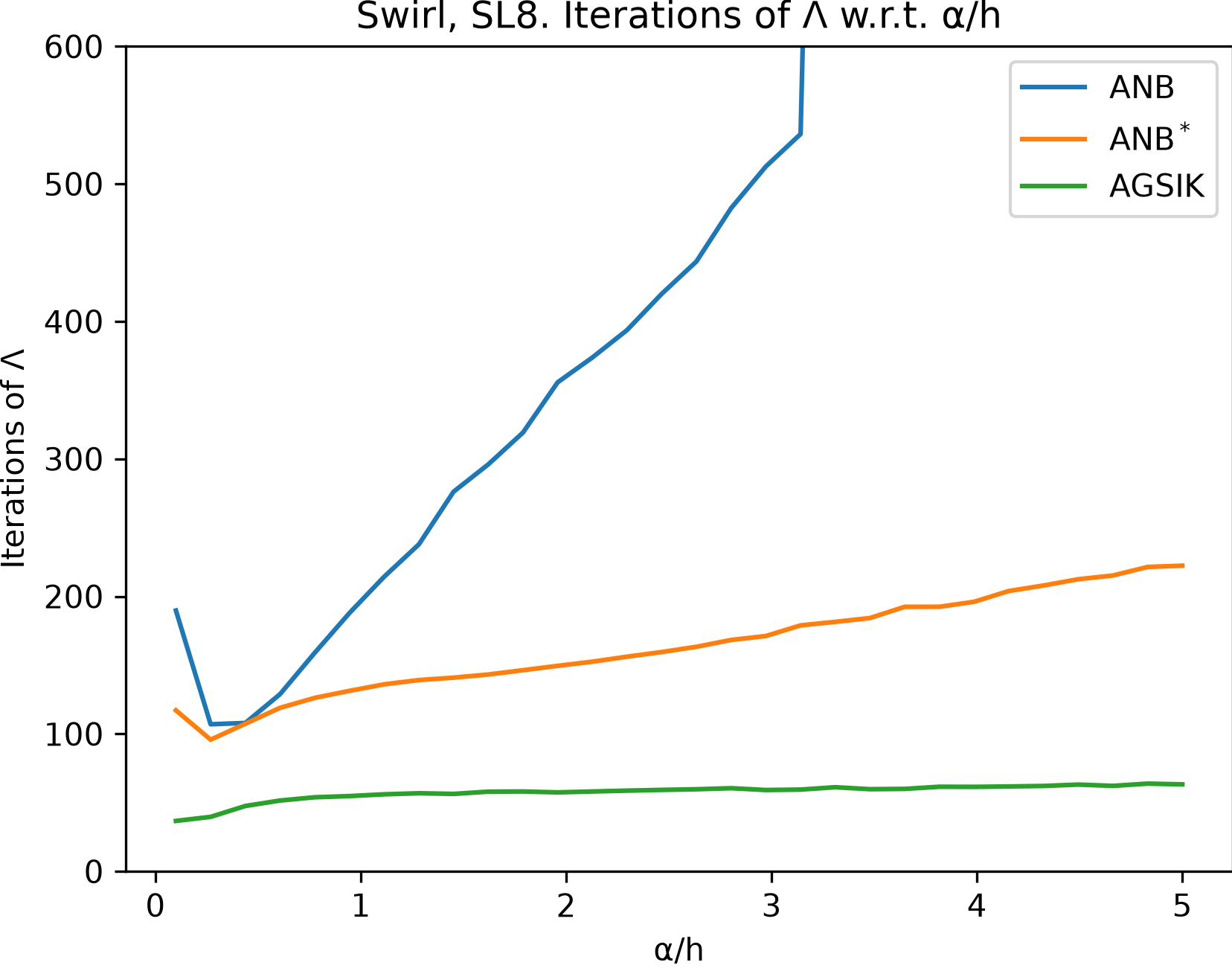}

\caption{Average number of iterations $N_\Lambda(\alpha)$ of the update operator $\Lambda$ at each point of the domain, w.r.t.\ the \Kruzkov\ timescale $\alpha$ (specifically: $r_\alpha = \alpha/h$), in the various test cases. 
In order to show all curves on the same graph, we often divide the ANB mean iteration number by 10, and write ANB/10 in the legend.
The choice of the parameter $\alpha$ is unfortunately quite sensitive for the ANB and ANB$^*$, whose numerical cost $N_\Lambda(\alpha)$ may increase significantly when $\alpha$ is too small or too large. 
Fortunately, it can be guided by the linearized contraction ratio, which is easy to compute a priori: $N_\Lambda(\alpha)$ is often approximately minimal when $\rho^\lin(\alpha)$ is minimal, compare with \cref{fig:rho_lin}. 
The numerical cost of the AGSIK appears to be less sensitive to $\alpha$.
}
\label{fig:iter_wrt_alpha}
\end{figure}

\begin{figure}
\def\w{3.7cm}\def\W{5cm}
\includegraphics[width=\w]{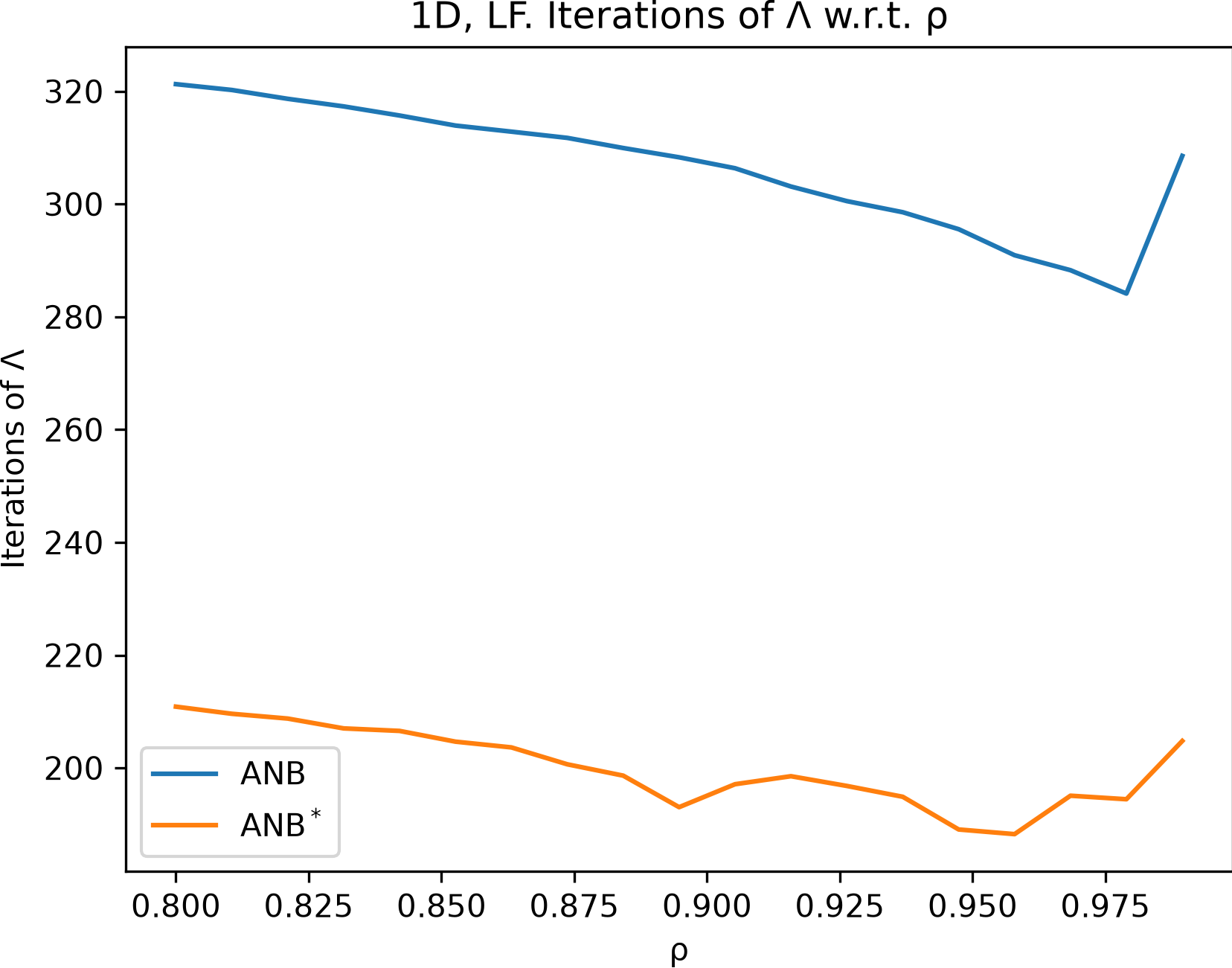}
\includegraphics[width=\w]{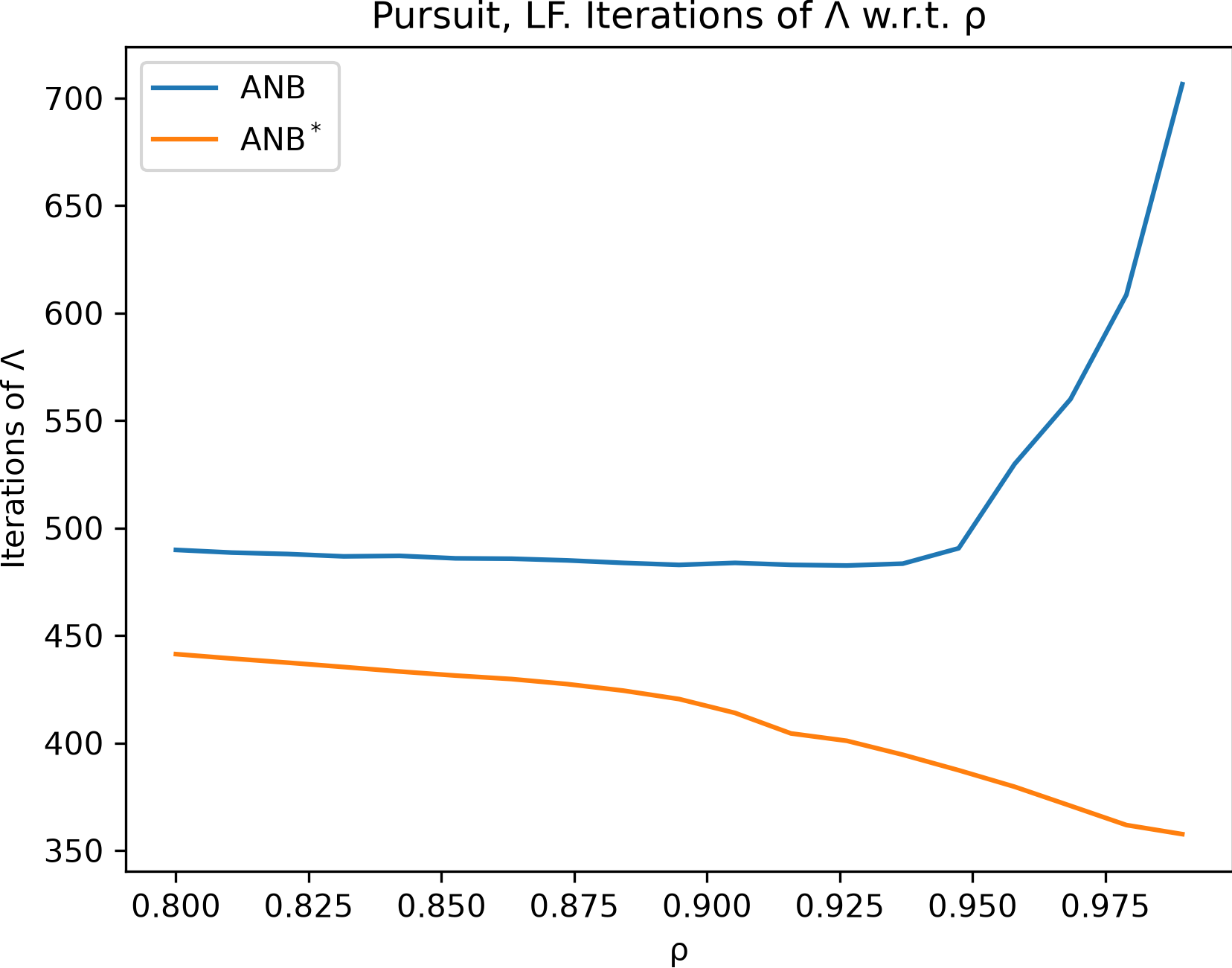}
\includegraphics[width=\w]{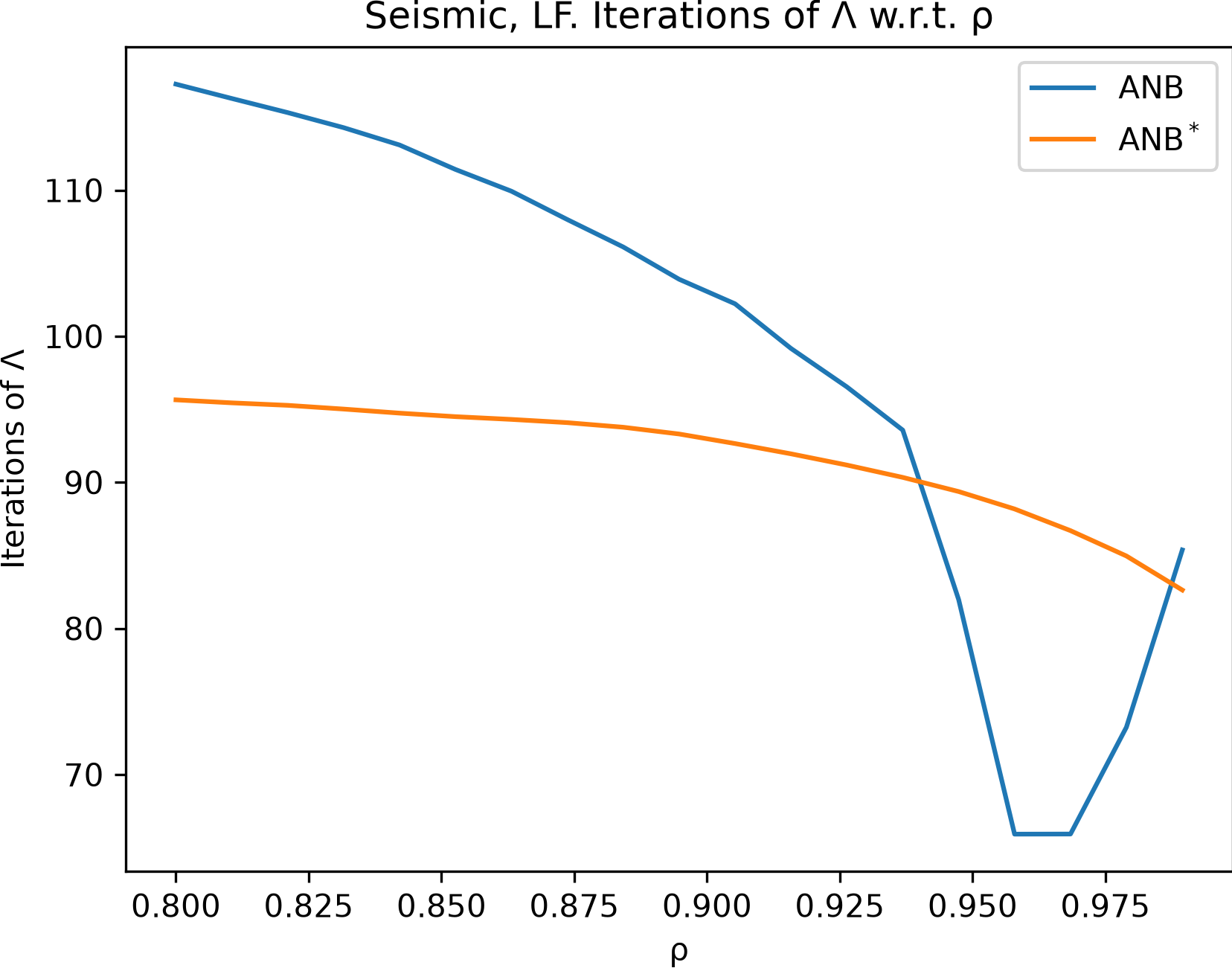}
\includegraphics[width=\w]{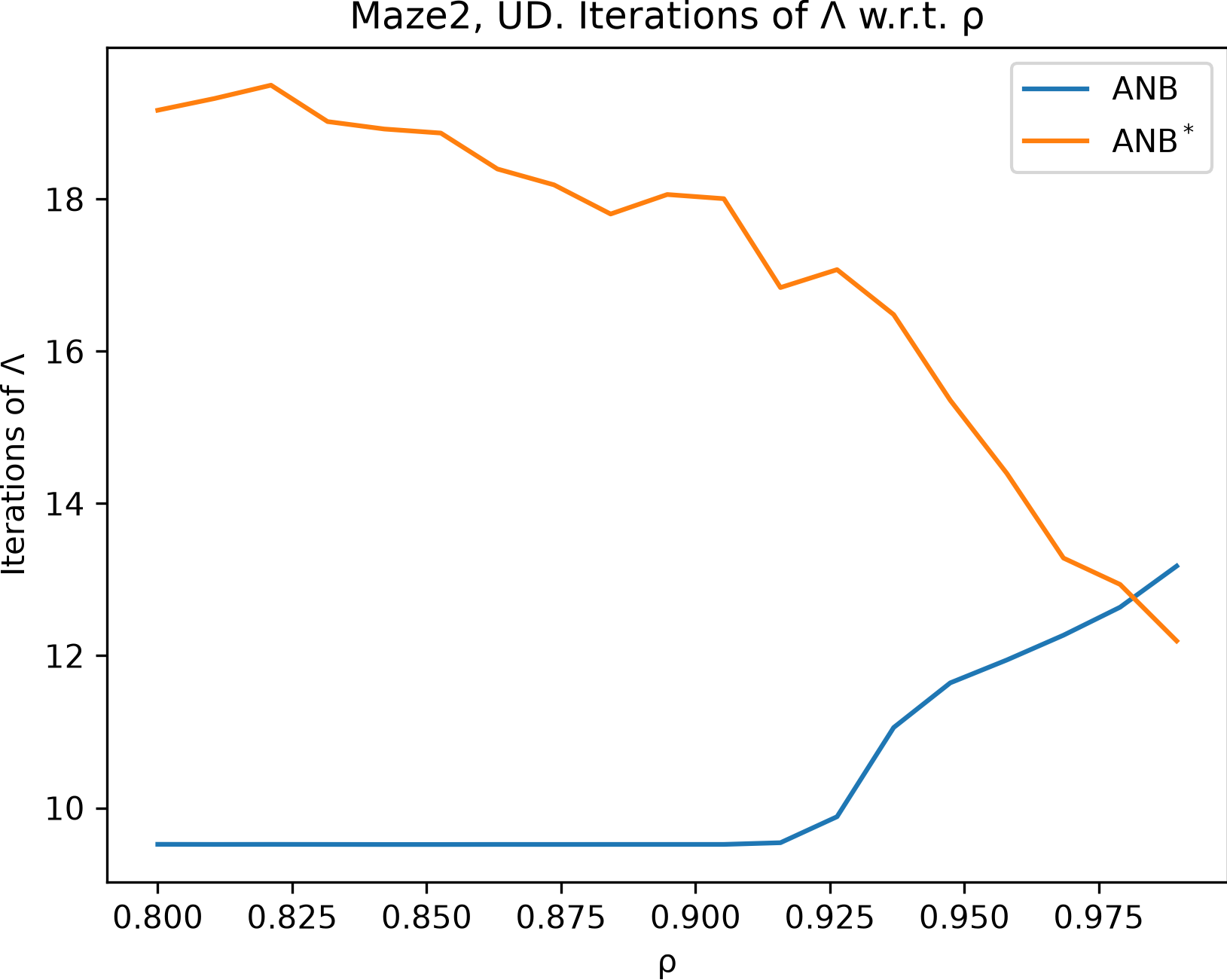}\\
\includegraphics[width=\W]{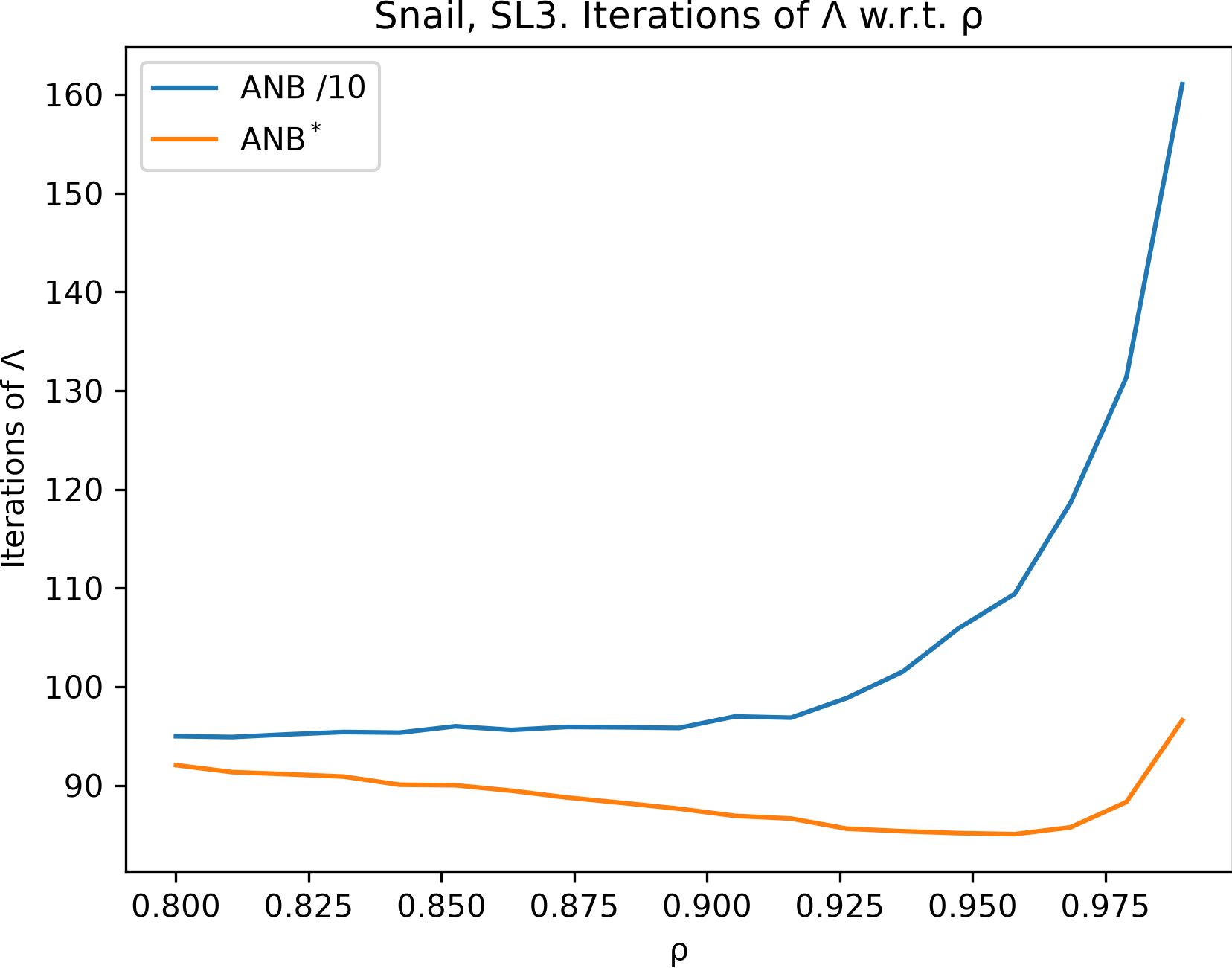}
\includegraphics[width=\W]{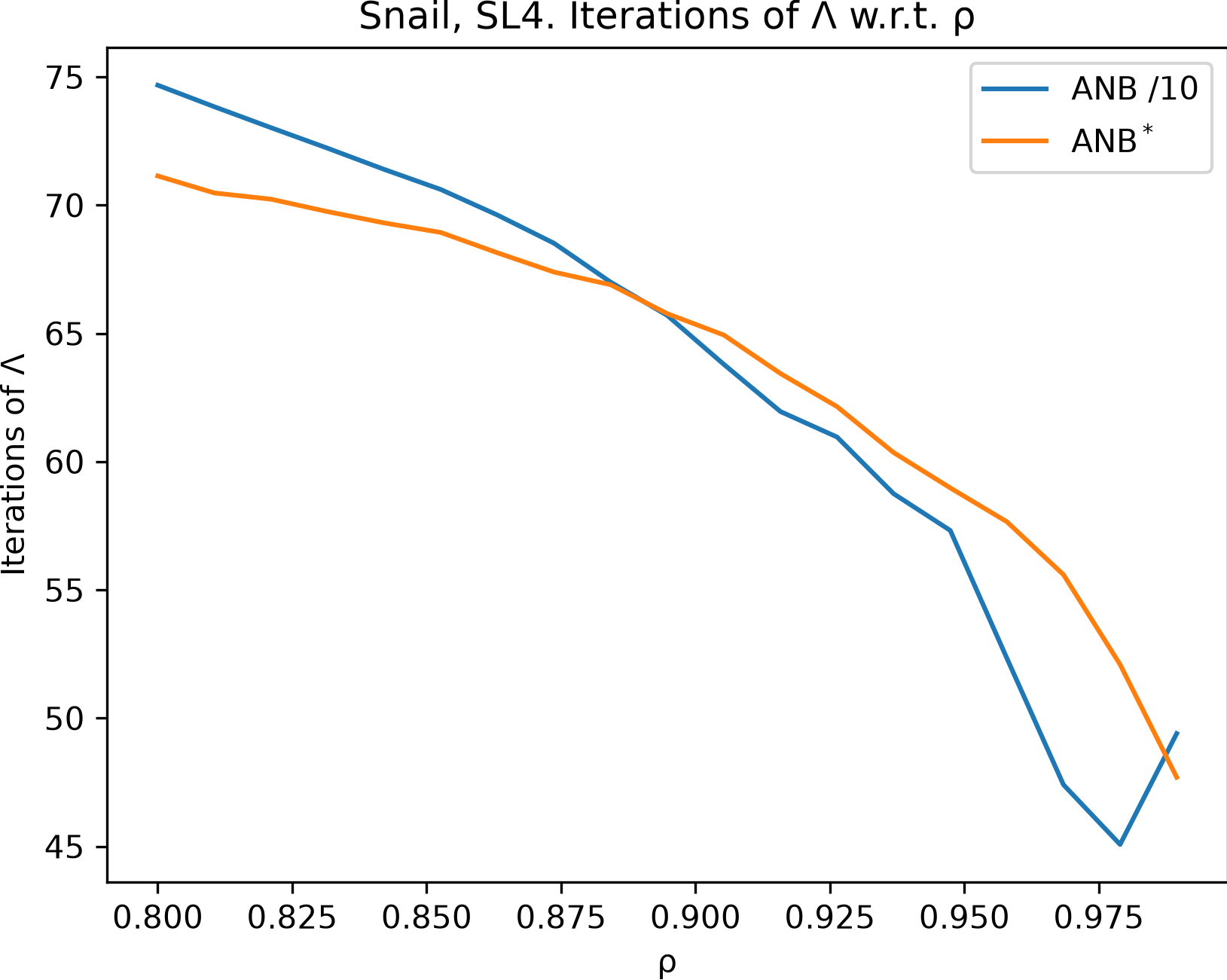}
\includegraphics[width=\W]{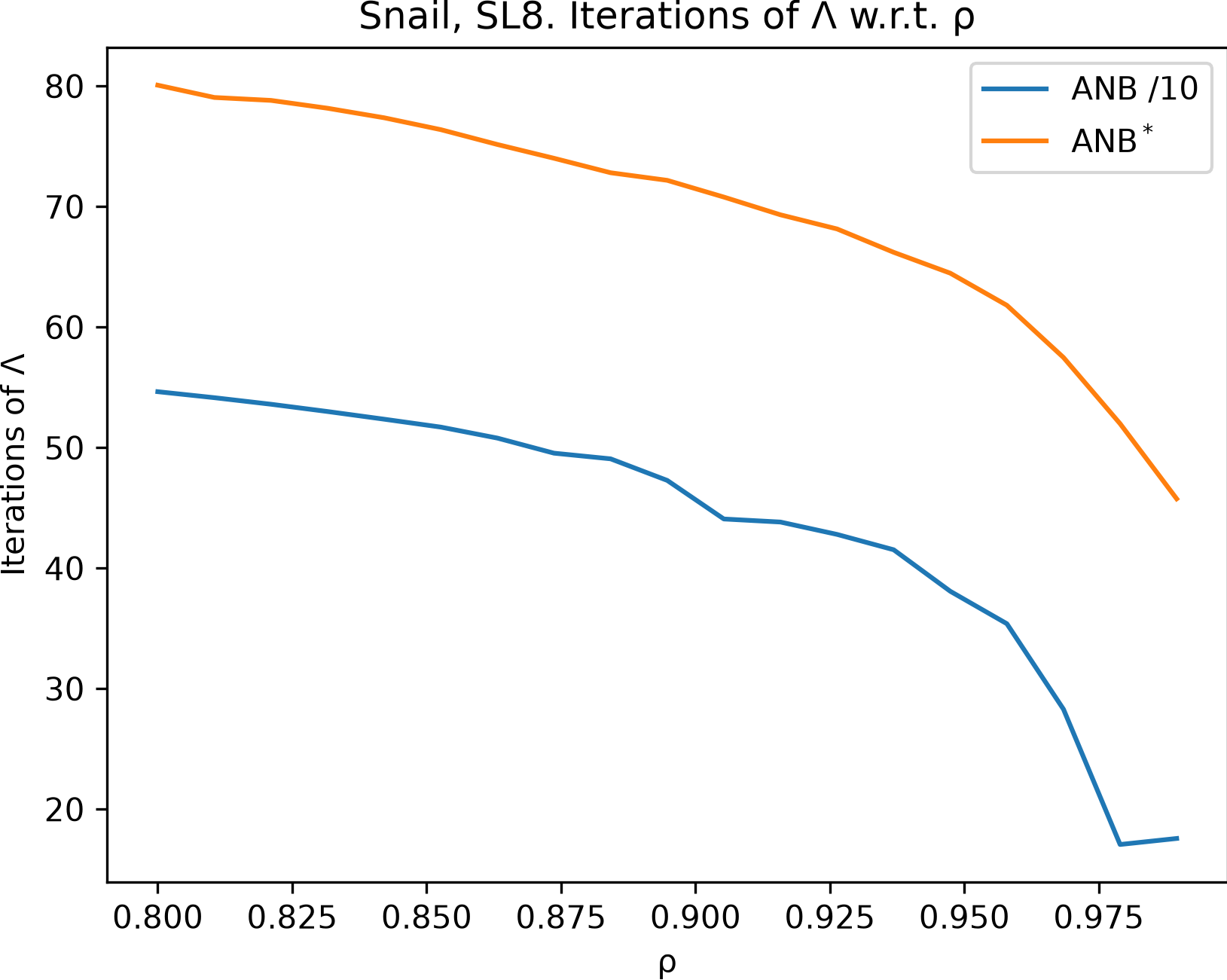}\\
\includegraphics[width=\W]{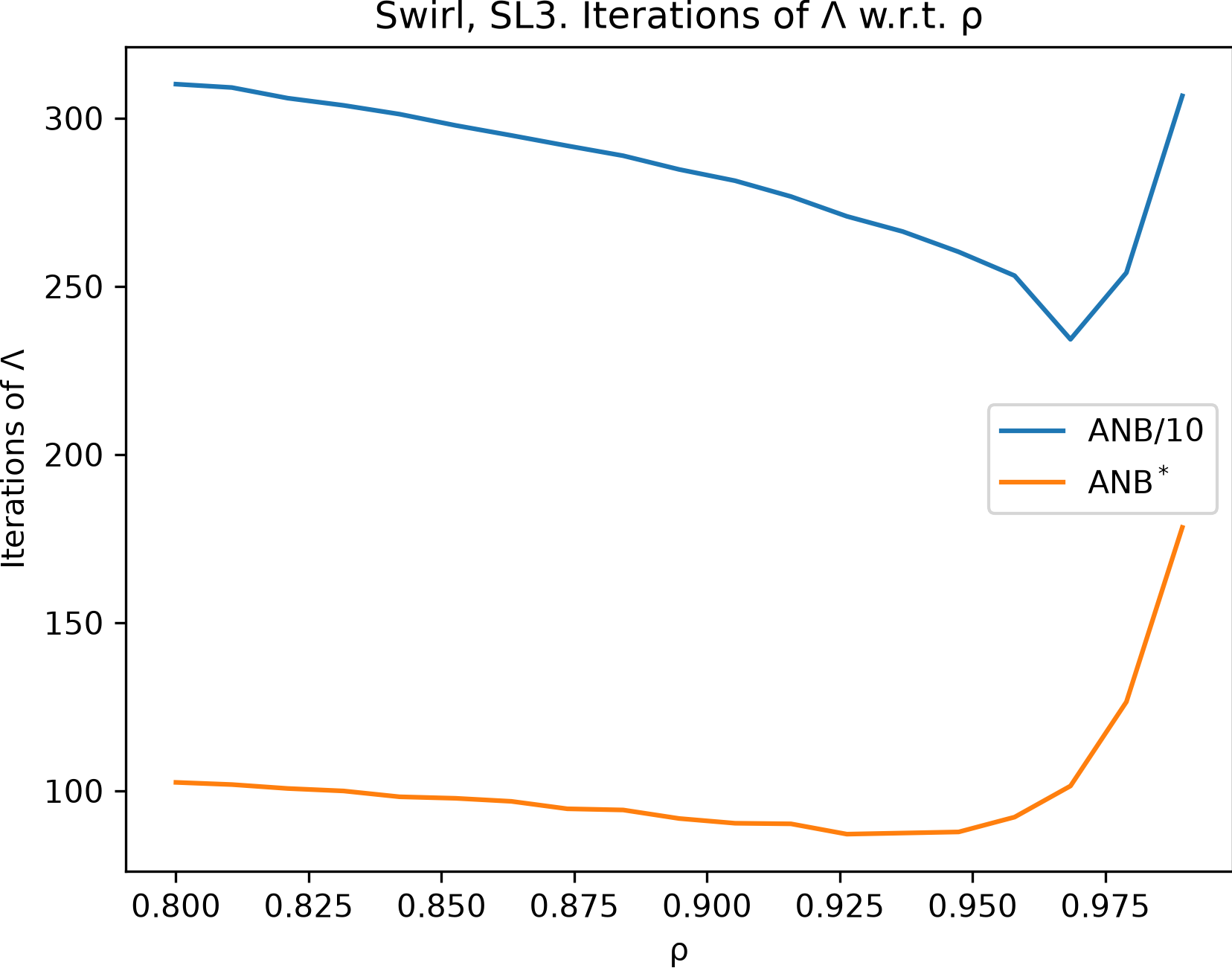}
\includegraphics[width=\W]{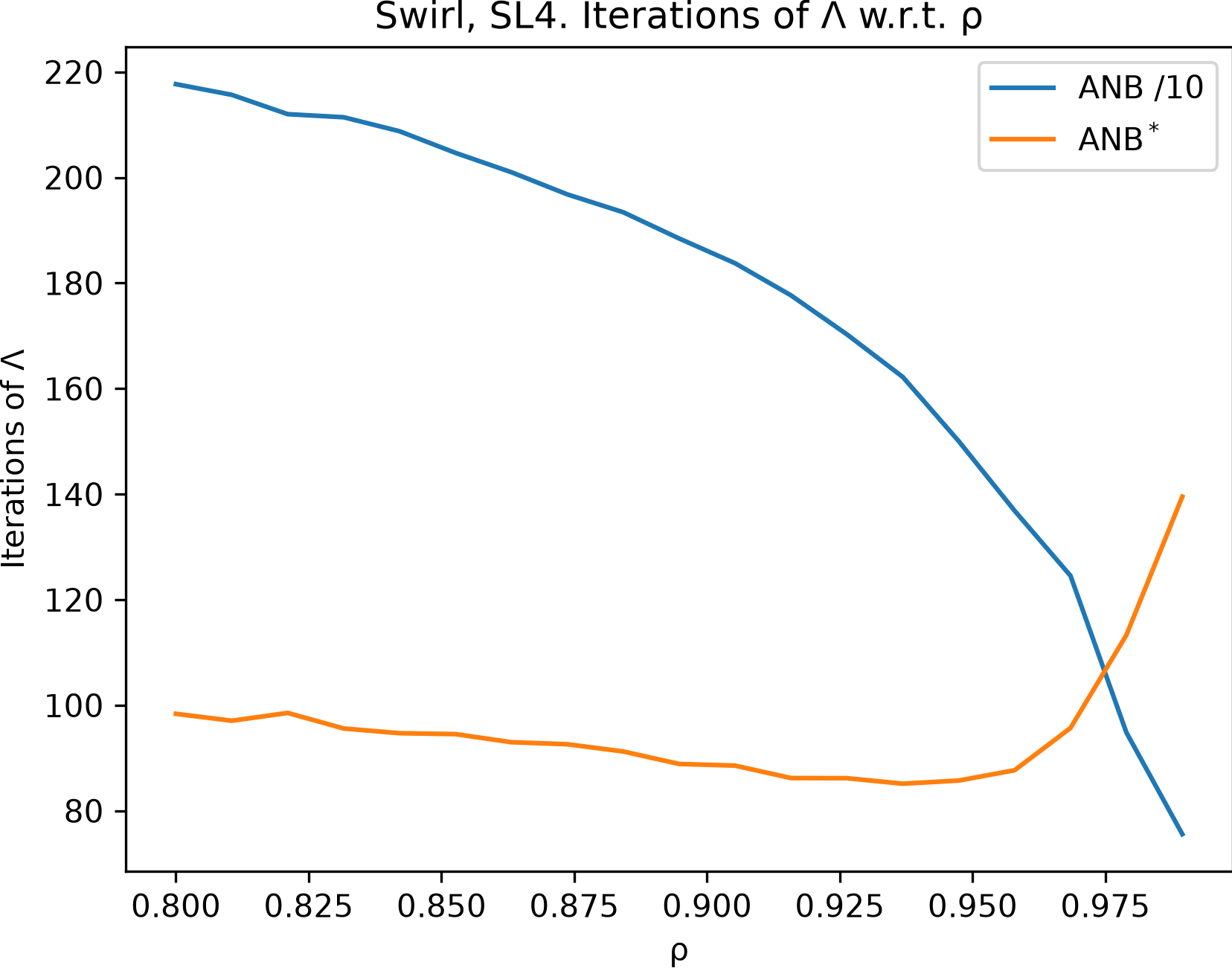}
\includegraphics[width=\W]{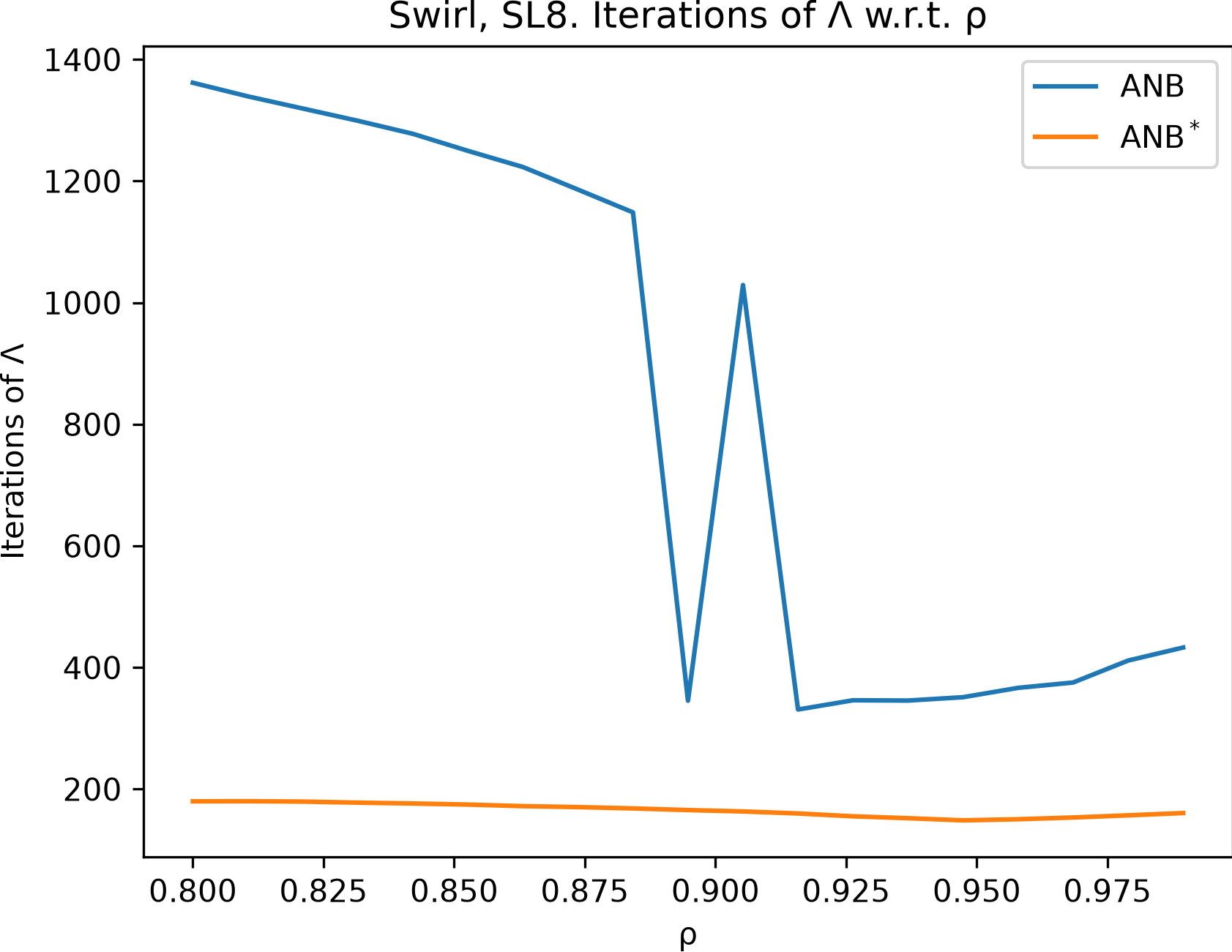}

\caption{Average number of iterations $N_\Lambda(\rho)$ of the update operator $\Lambda$ at each point of the domain, w.r.t.\ the parameter $\rho\in [0.8,0.99]$, in the considered test cases, other parameters defaulted.
The value of $\rho$ minimizing $N_\Lambda(\rho)$ depends on the discretization scheme (LF, UD, SL3-4-8), the fixed point solver (ANB, ANB$^*$), and the test case, in a rather unpredictable manner. 
Nevertheless, $N_\Lambda(\rho)$ rarely varies by a factor larger than $2$ as $\rho\in [0.8,0.99]$ in a given configuration.  This means that the method is not too sensitive to this parameter, which is arbitrarily fixed to $\rho = 0.95$ in all other experiments, see \eqref{eq:default_rho_tau}.
}
\label{fig:iter_wrt_rho}
\end{figure}

\begin{figure}
\def\w{3.7cm}\def\W{5cm}
\includegraphics[width=\w]{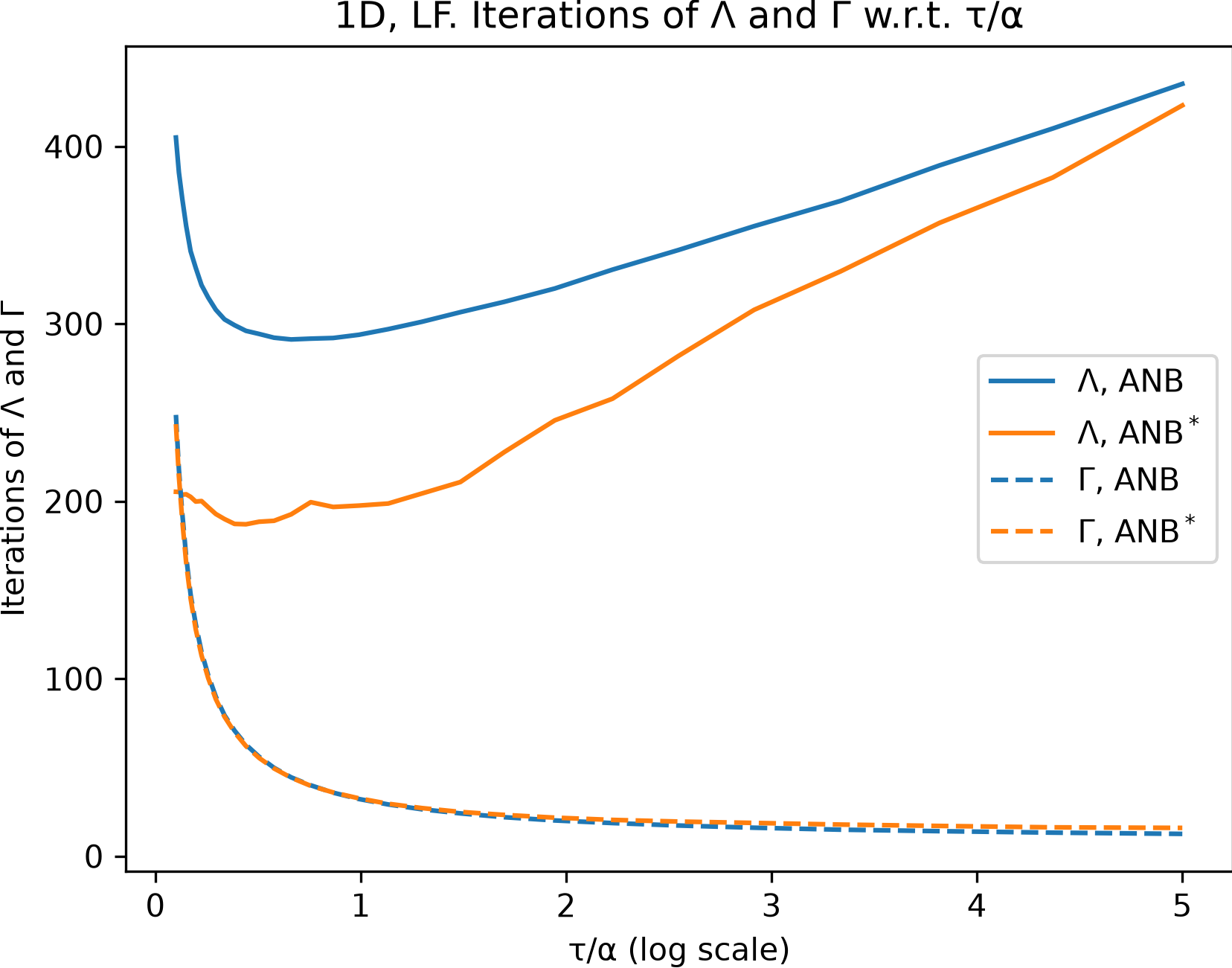}
\includegraphics[width=\w]{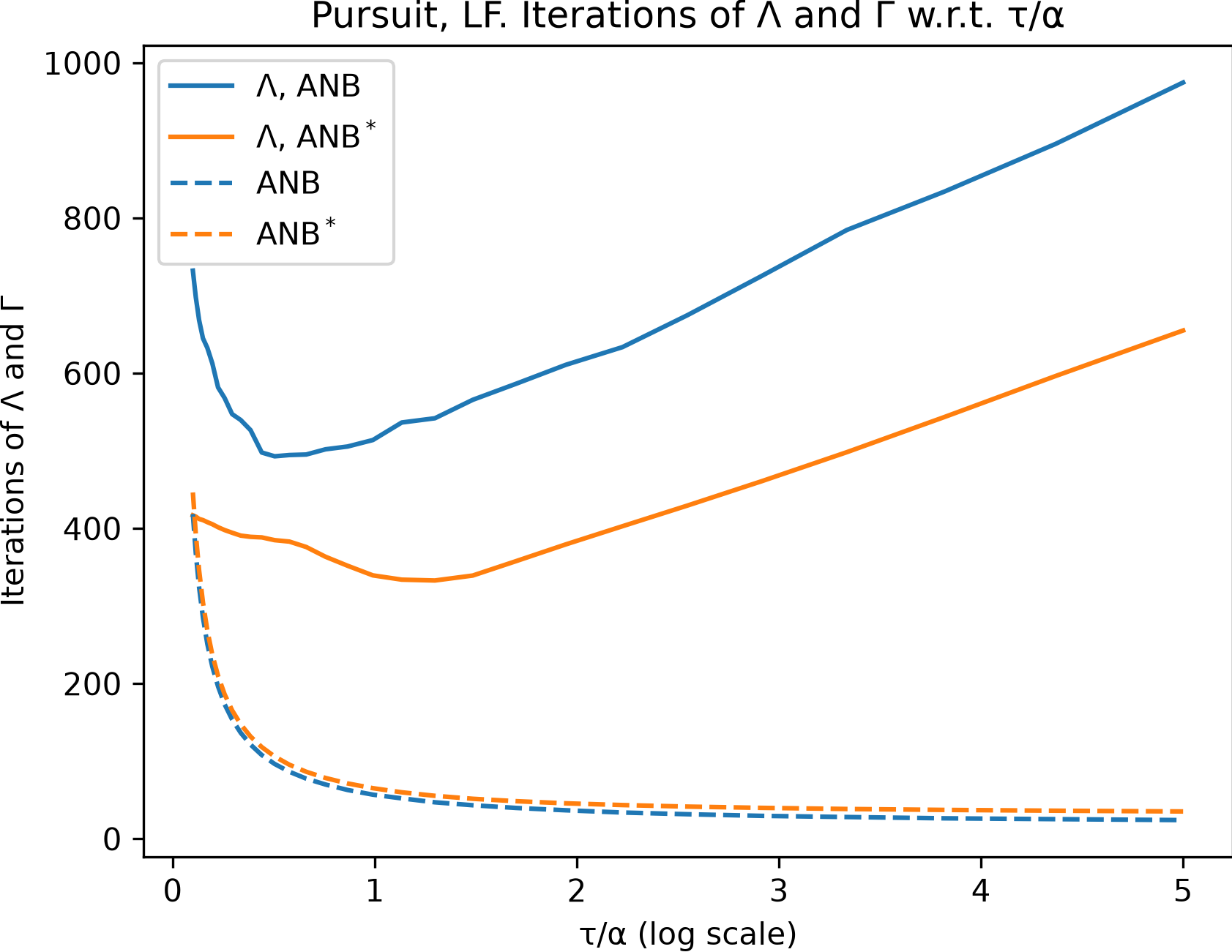}
\includegraphics[width=\w]{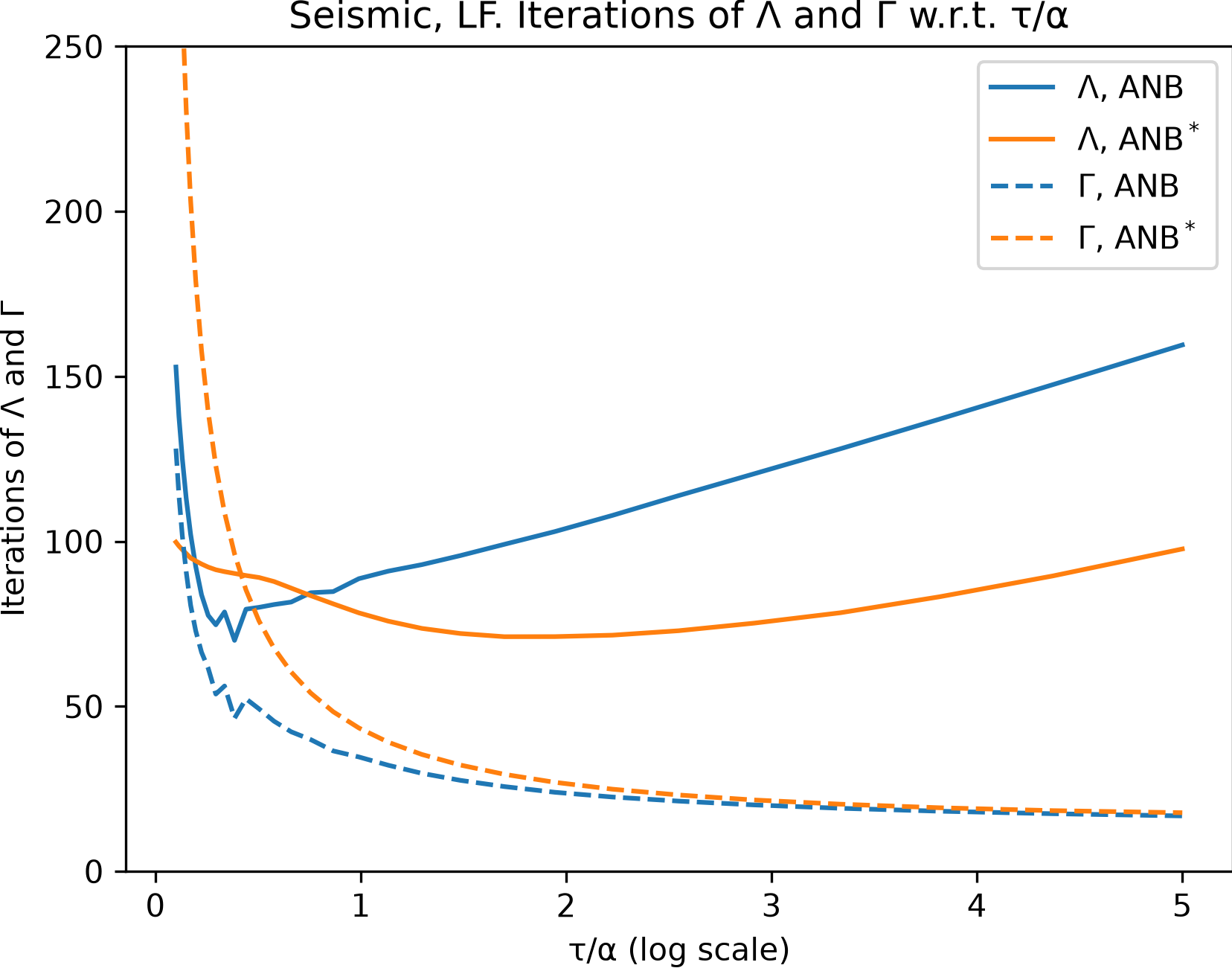}
\includegraphics[width=\w]{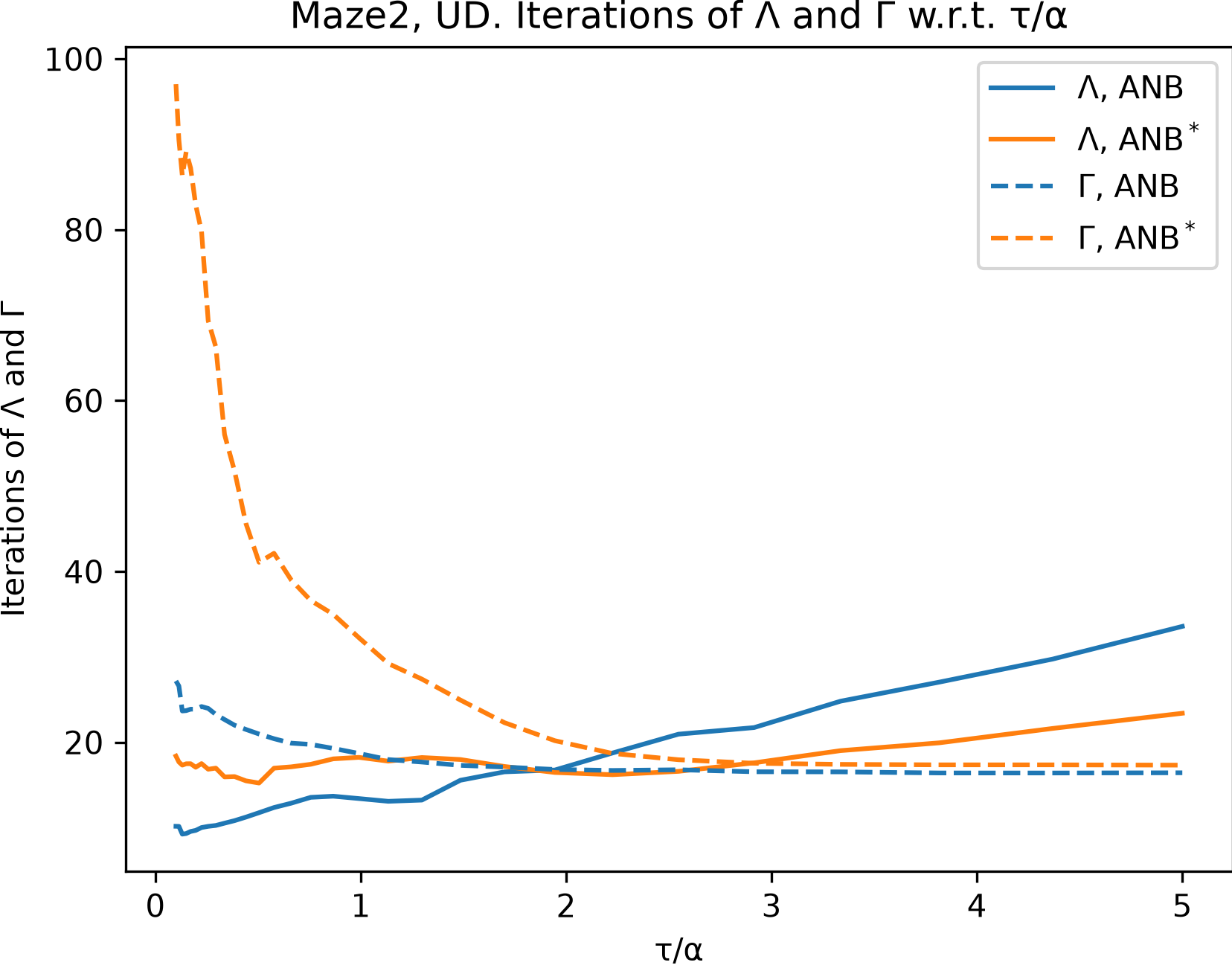}\\
\includegraphics[width=\W]{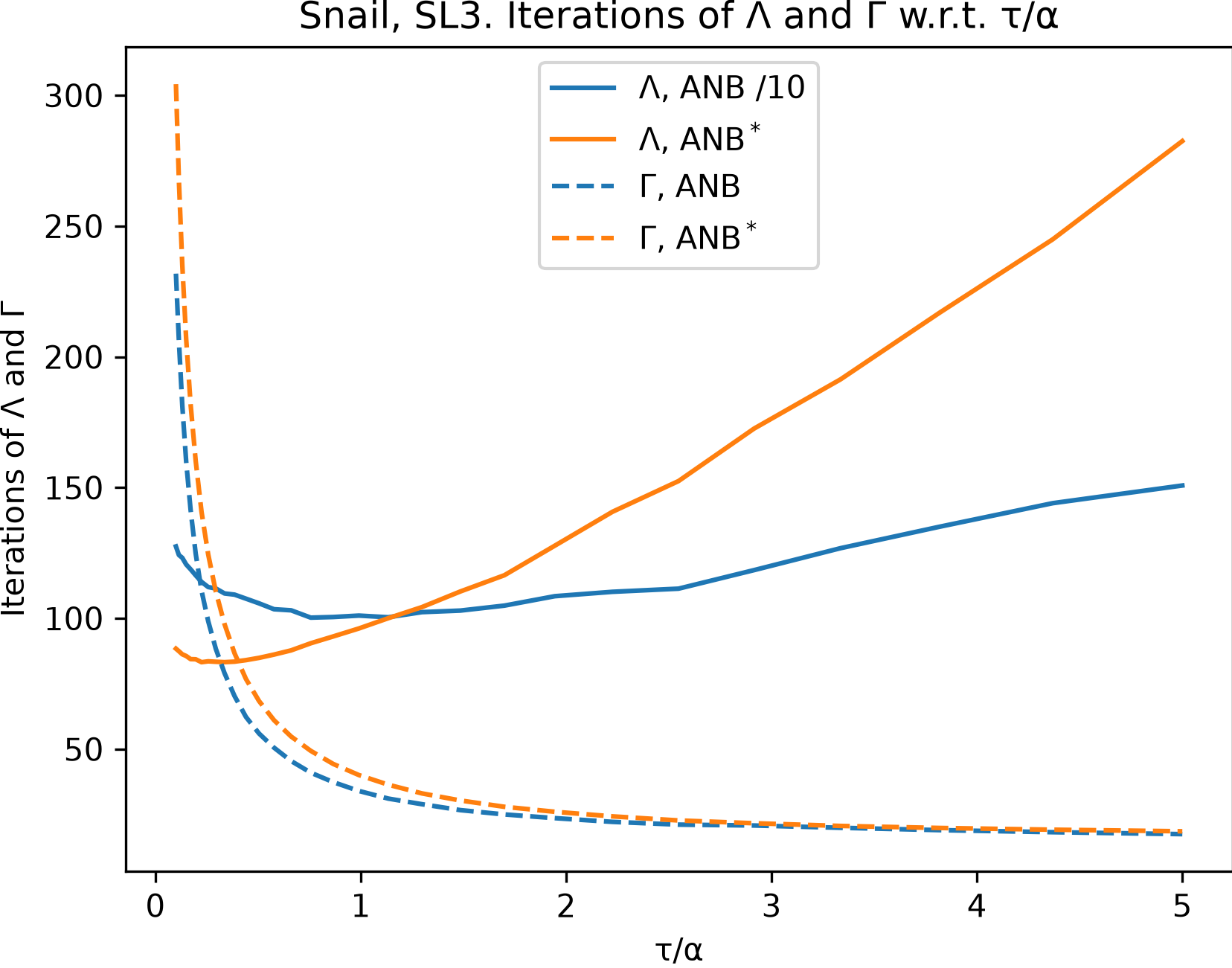}
\includegraphics[width=\W]{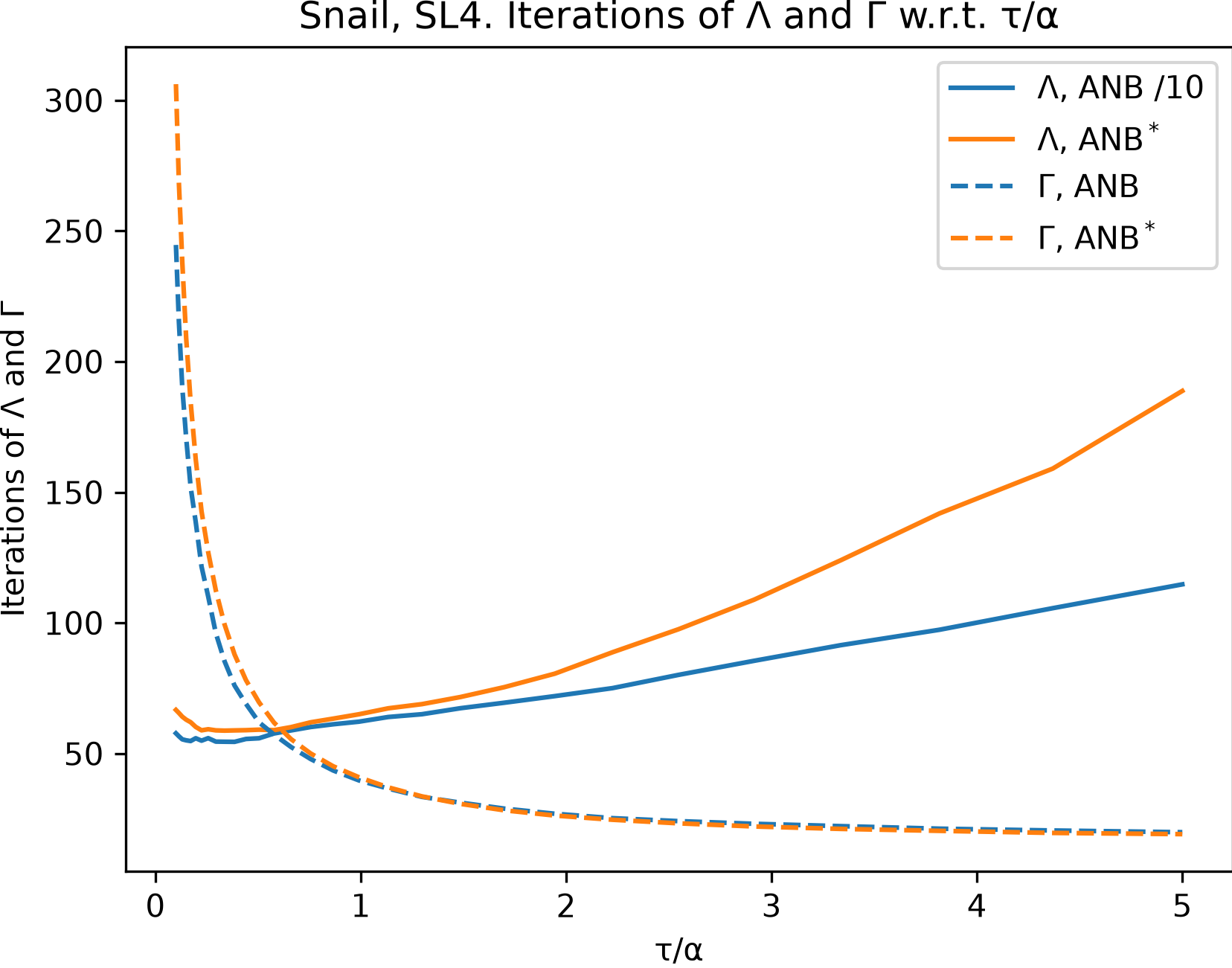}
\includegraphics[width=\W]{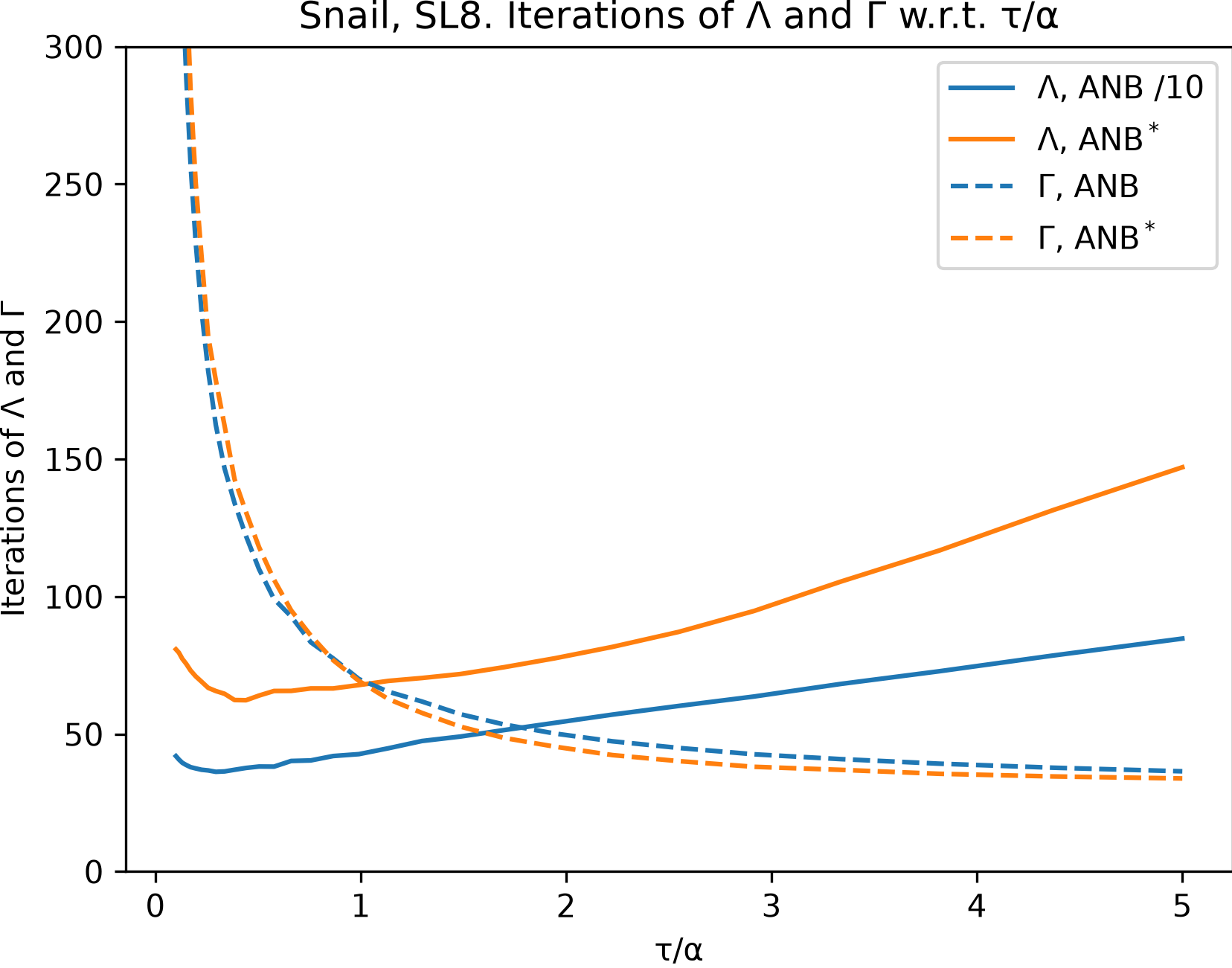}\\
\includegraphics[width=\W]{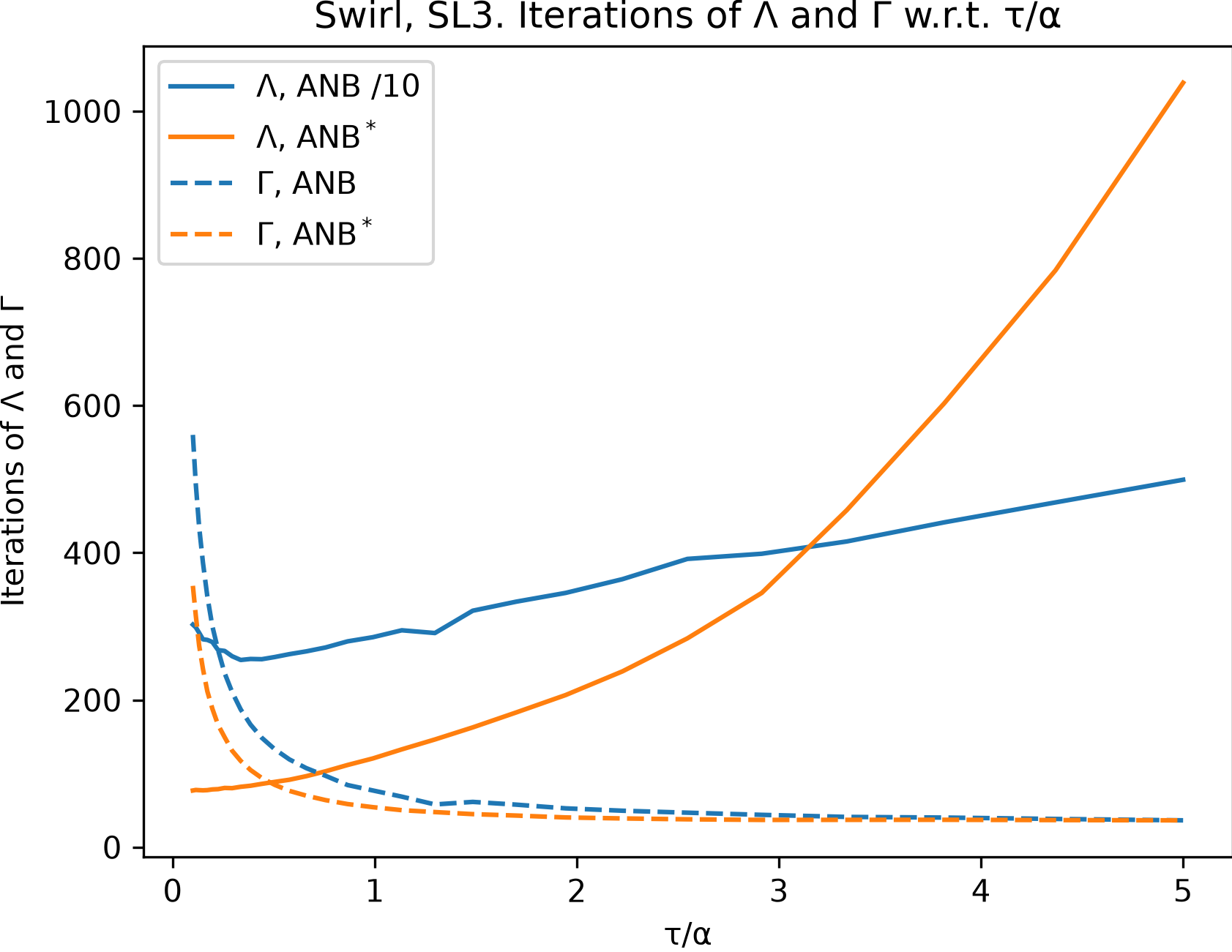}
\includegraphics[width=\W]{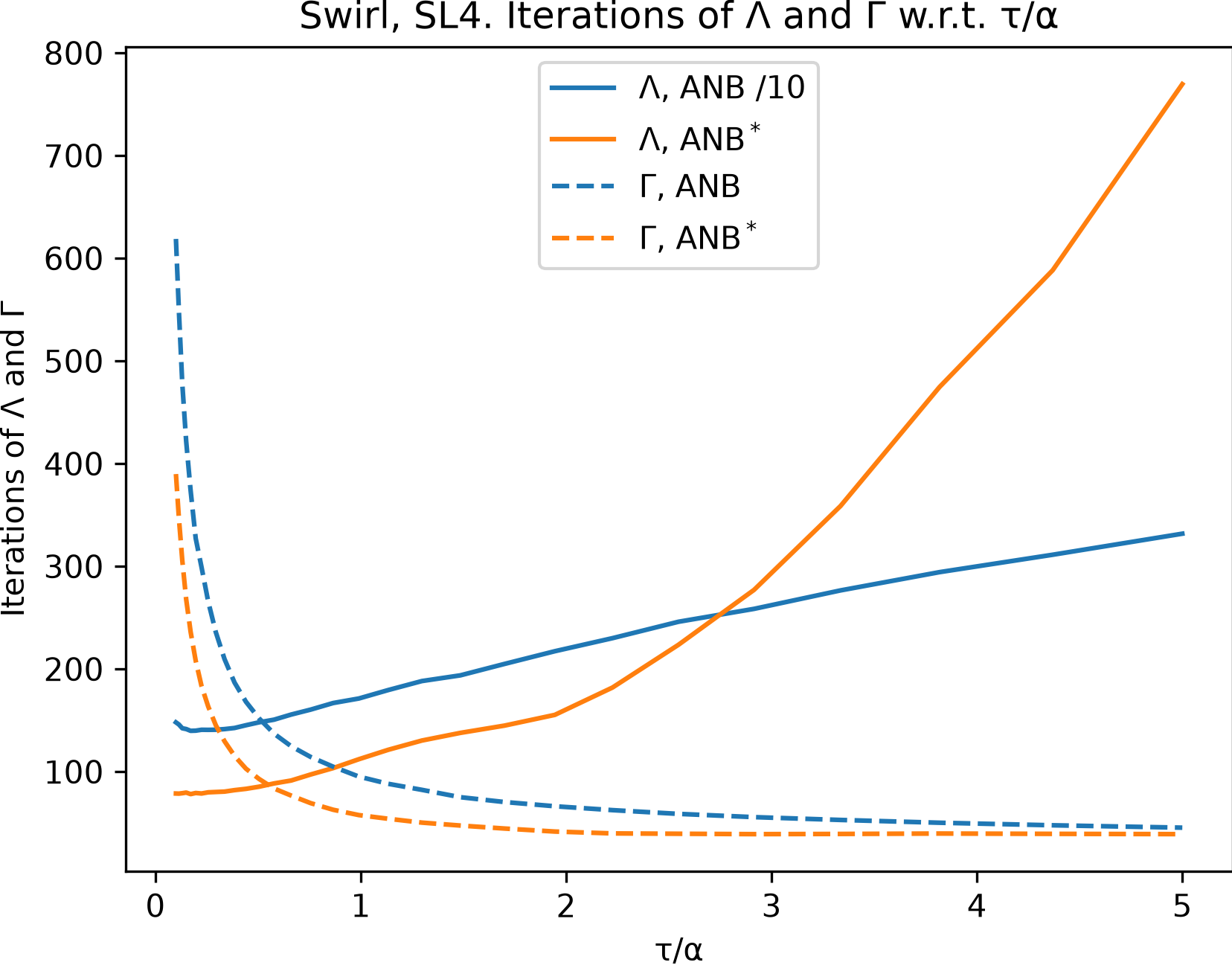}
\includegraphics[width=\W]{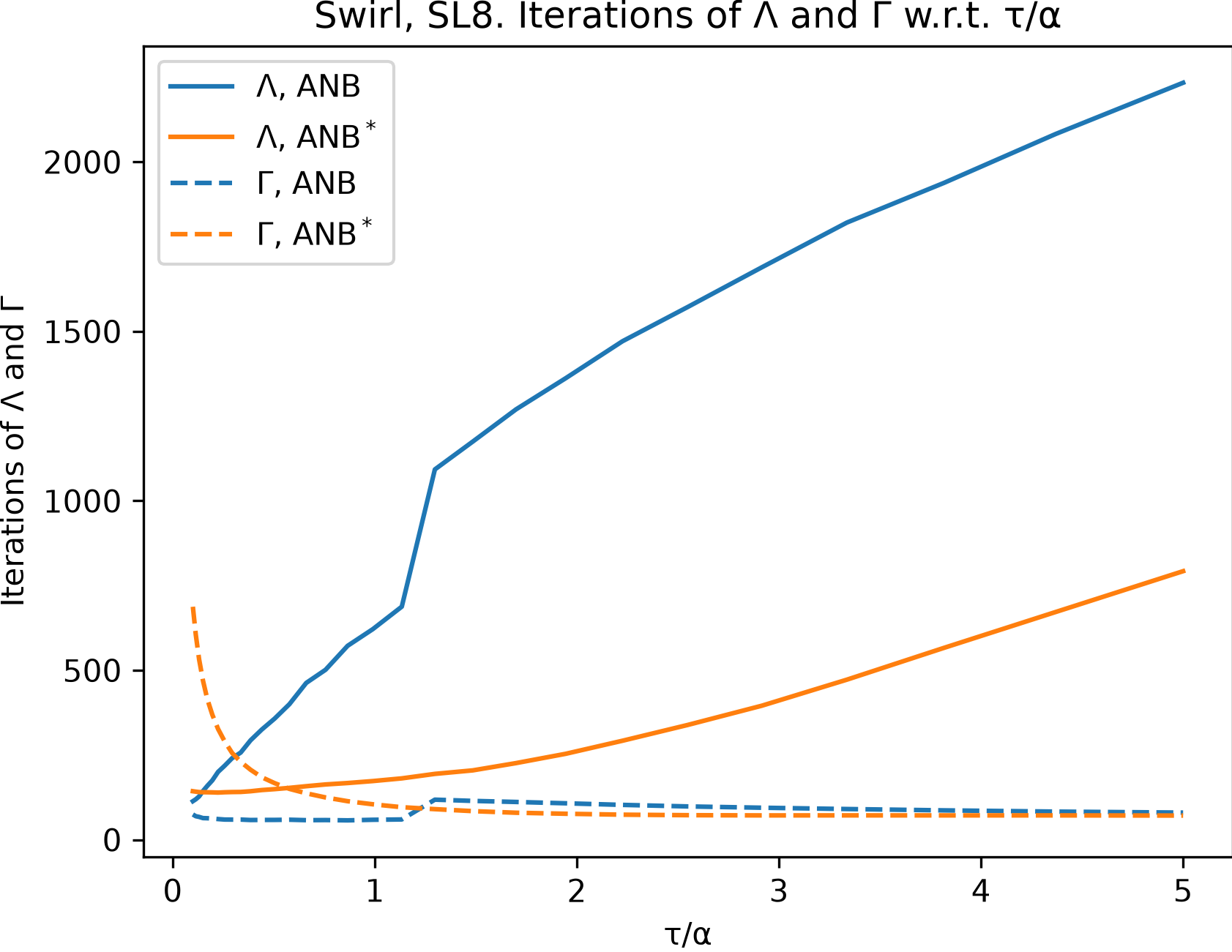}

\caption{Average number of iterations of the update operator $\Lambda$ at each point of the domain, w.r.t.\ the timestep $\tau$ (specifically: $r_\tau = \tau/\alpha$), for the ANB and ANB$^*$ methods, in the various test cases. Number of iterations of the scheme $\Lambda$ (solid lines), and of the causal approximation $\Gamma$ used for the narrow band construction (dashed lines).
The parameter $\tau$ determines how much the narrow band advances in each iteration.
The band is rebuilt too often if $\tau/\alpha \ll 1$, thus requiring more evaluations of $\Gamma$, and insufficiently often if $\tau/\alpha \gg 1$, thus becoming too large for effective front tracking and requiring more evaluations of $\Lambda$. 
Choosing $\tau/\alpha \in [0.5,2]$ is usually adequate; we set $\tau = \alpha/2$ in all other experiments \eqref{eq:default_rho_tau}.
}
\label{fig:iter_wrt_tau}
\end{figure}

We gather in this appendix a series of figures illustrating the computational cost of the proposed ANB, ANB$^*$ and AGSIK methods as their parameters and the gridscale and convergence tolerance vary. 
Cost is measured in iterations of the numerical scheme $\Lambda$ in \cref{fig:iter_wrt_alpha,fig:iter_wrt_eps,fig:iter_wrt_rho,fig:iter_wrt_tau}, but we also report some computation times in \cref{table:elapsed}. 
The first two parameters, the gridscale $h \approx N_x^{-1} \approx N^{-\frac 1 d}$ and the tolerance $\ve>0$ appear in the quasi-linear cost estimate $\cO(N|\ln \ve| \ln (N/\ve))$ established in \cref{th:complexity,th:complexity_withUnknown,th:complexity_with_agsi,prop:complexity_lazy} for the ANB and ANB$^*$, see \eqref{eq:ln_N_eps}.

\begin{figure}
\def\h{4cm}
\includegraphics[height=\h]{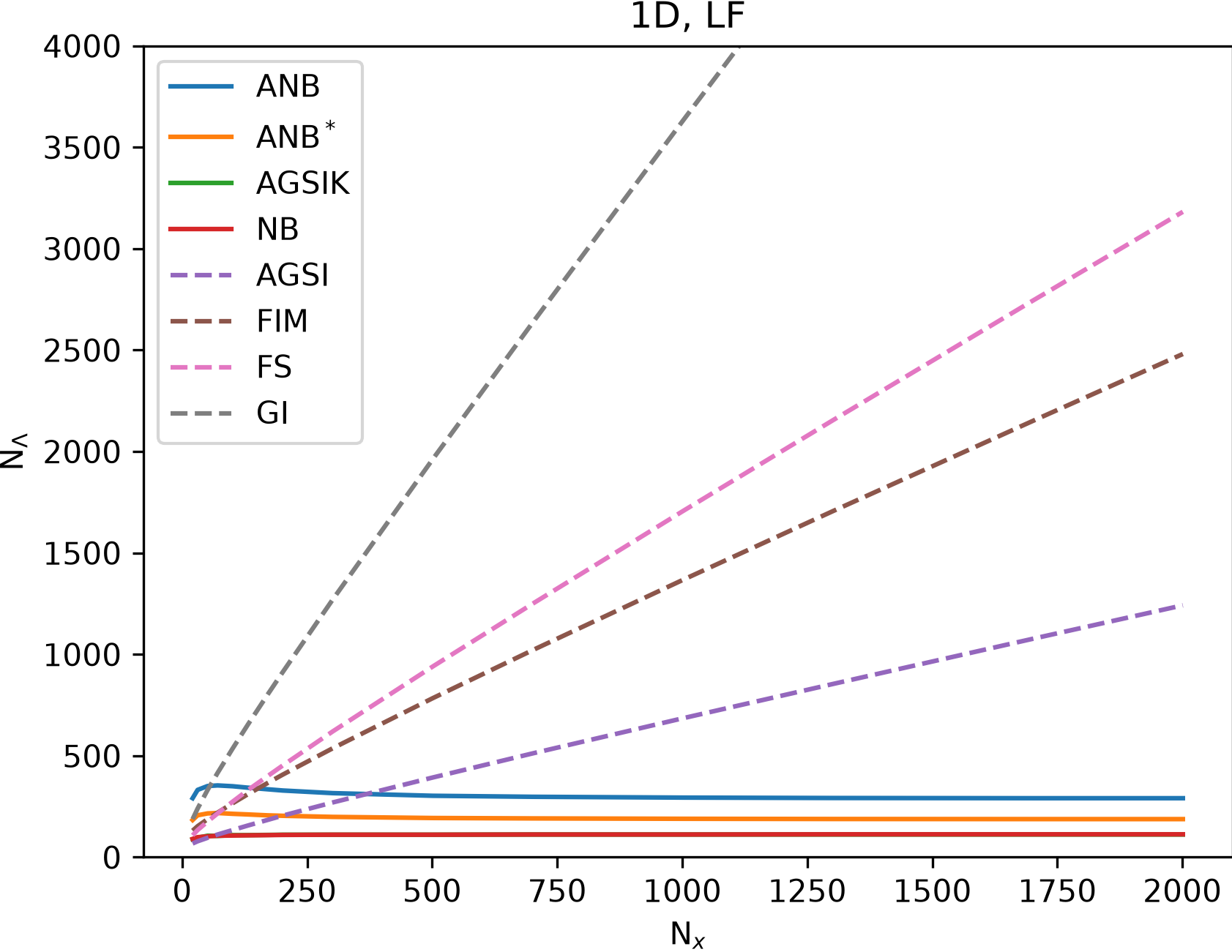}
\includegraphics[height=\h]{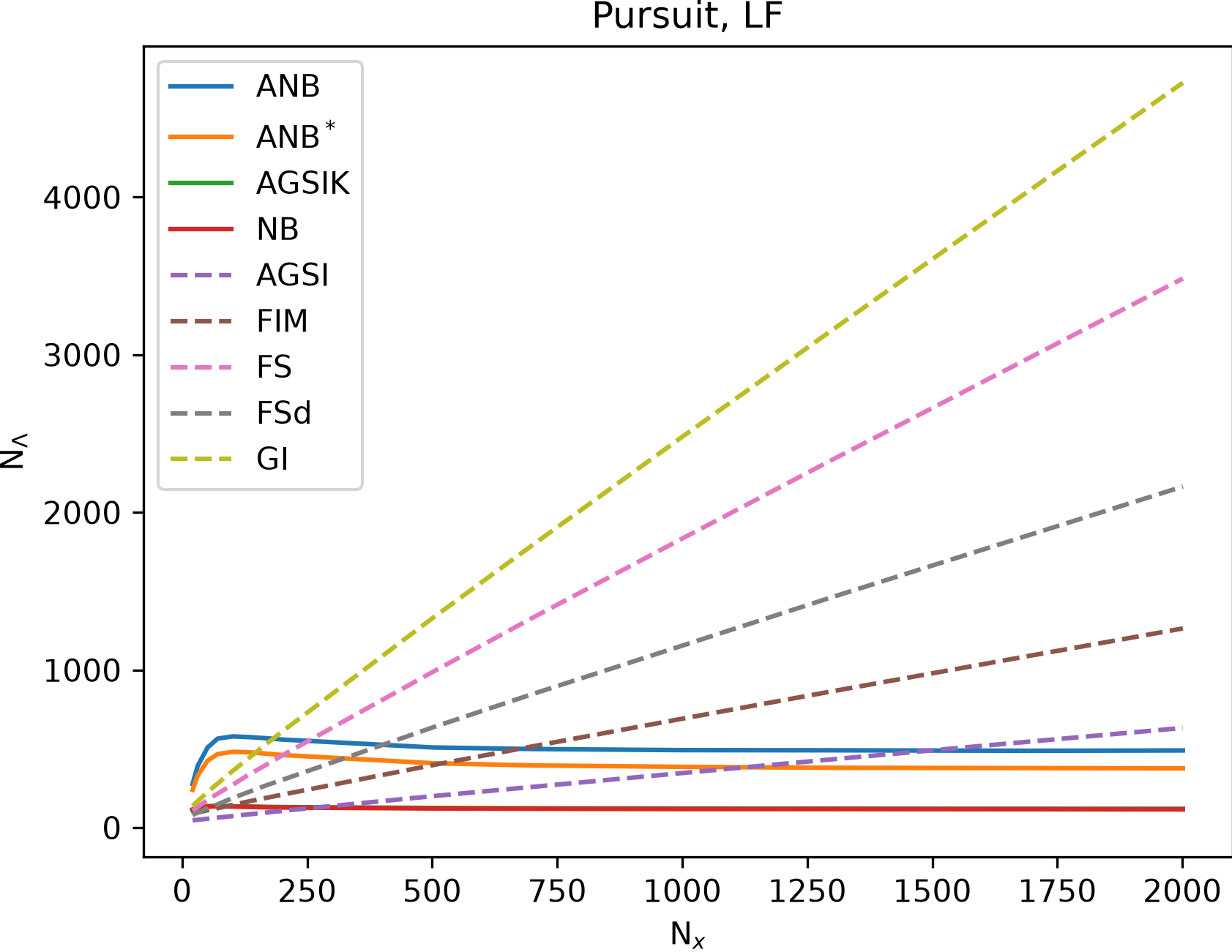}
\includegraphics[height=\h]{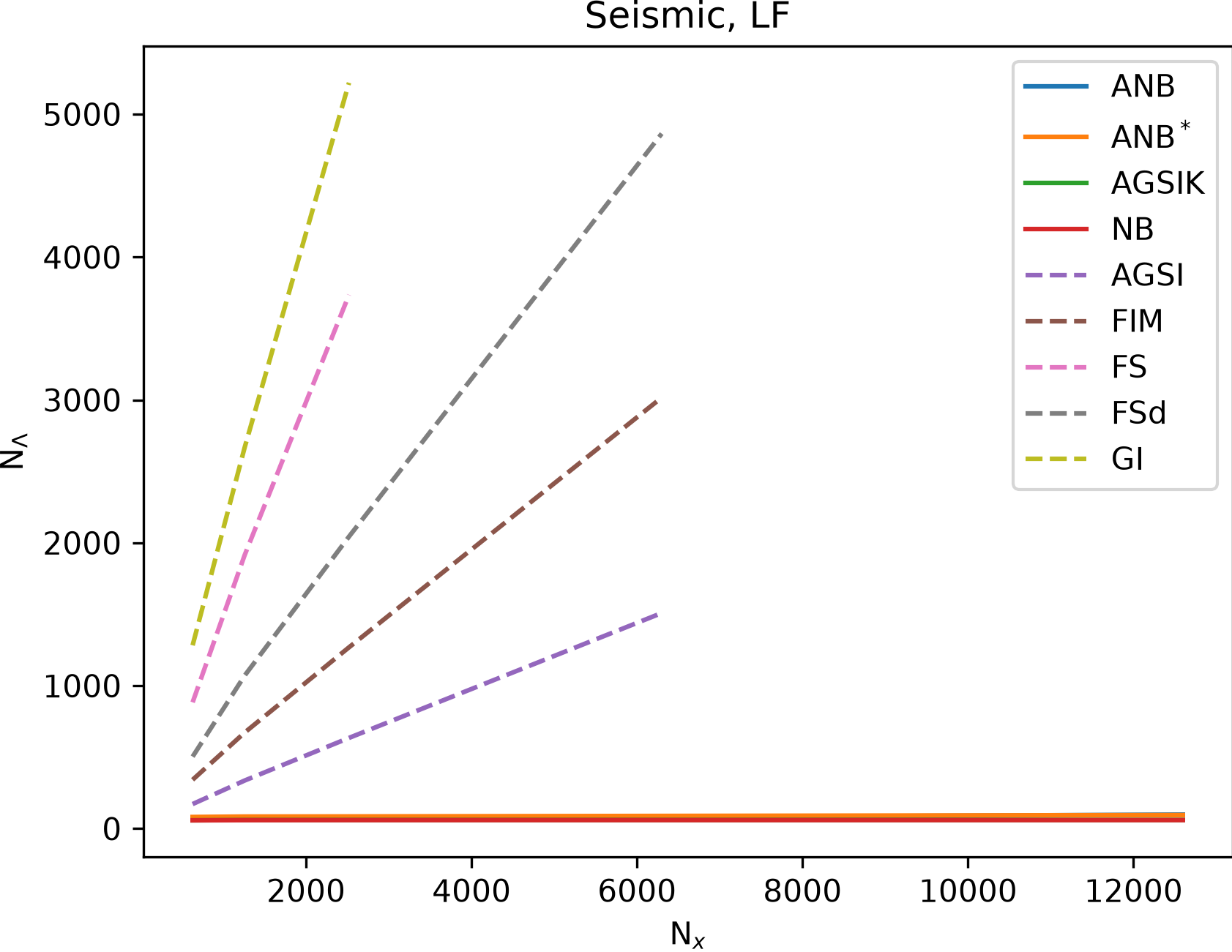}\\
\includegraphics[height=\h]{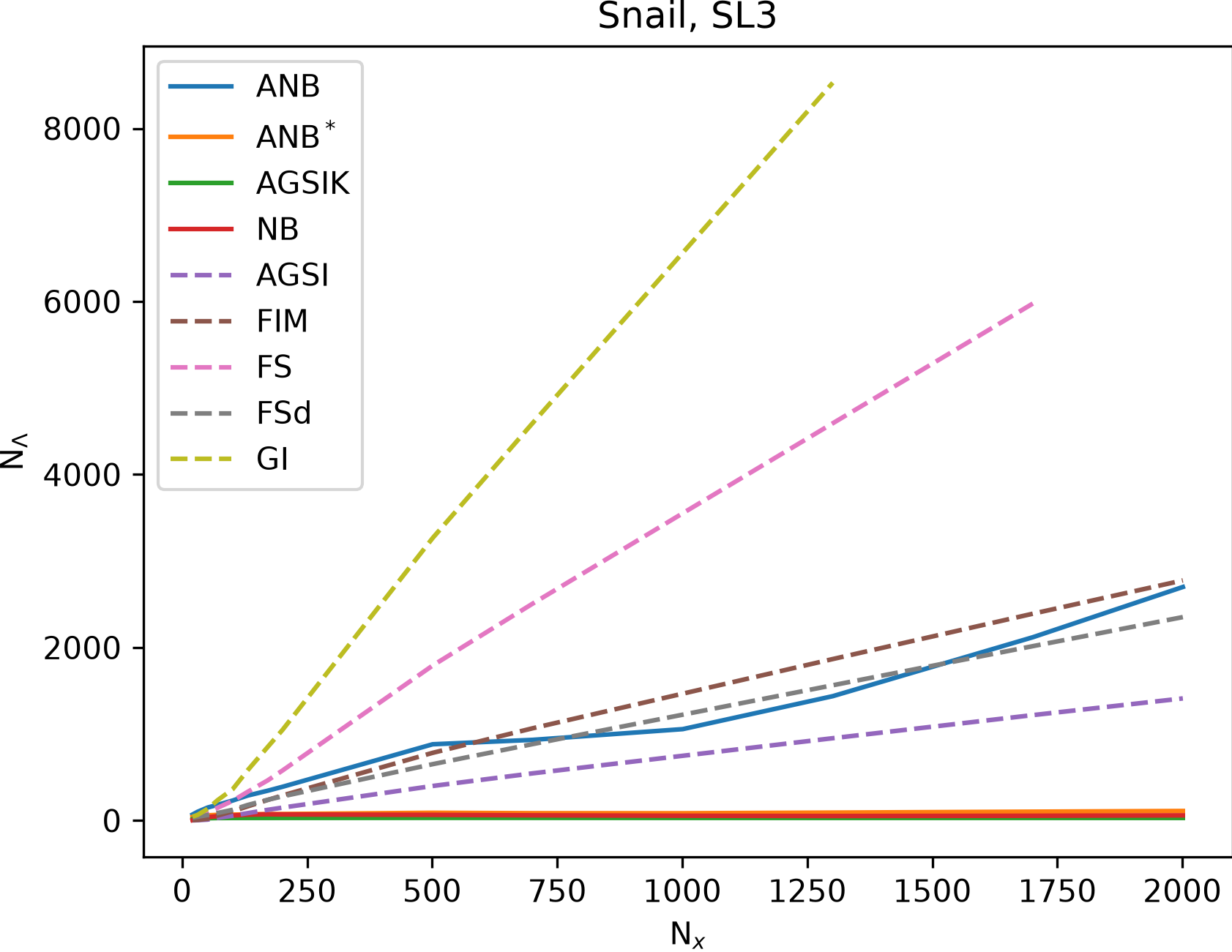}
\includegraphics[height=\h]{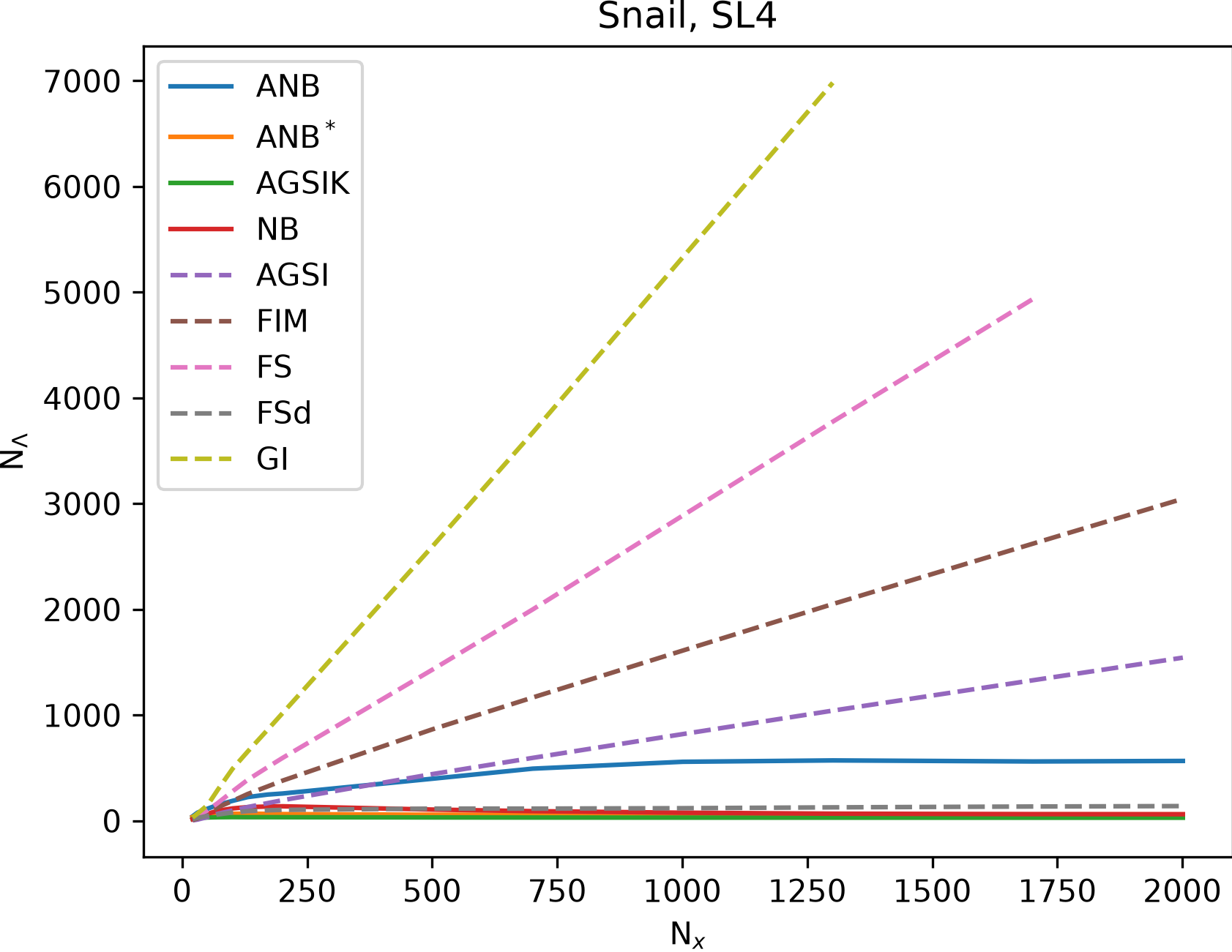}
\includegraphics[height=\h]{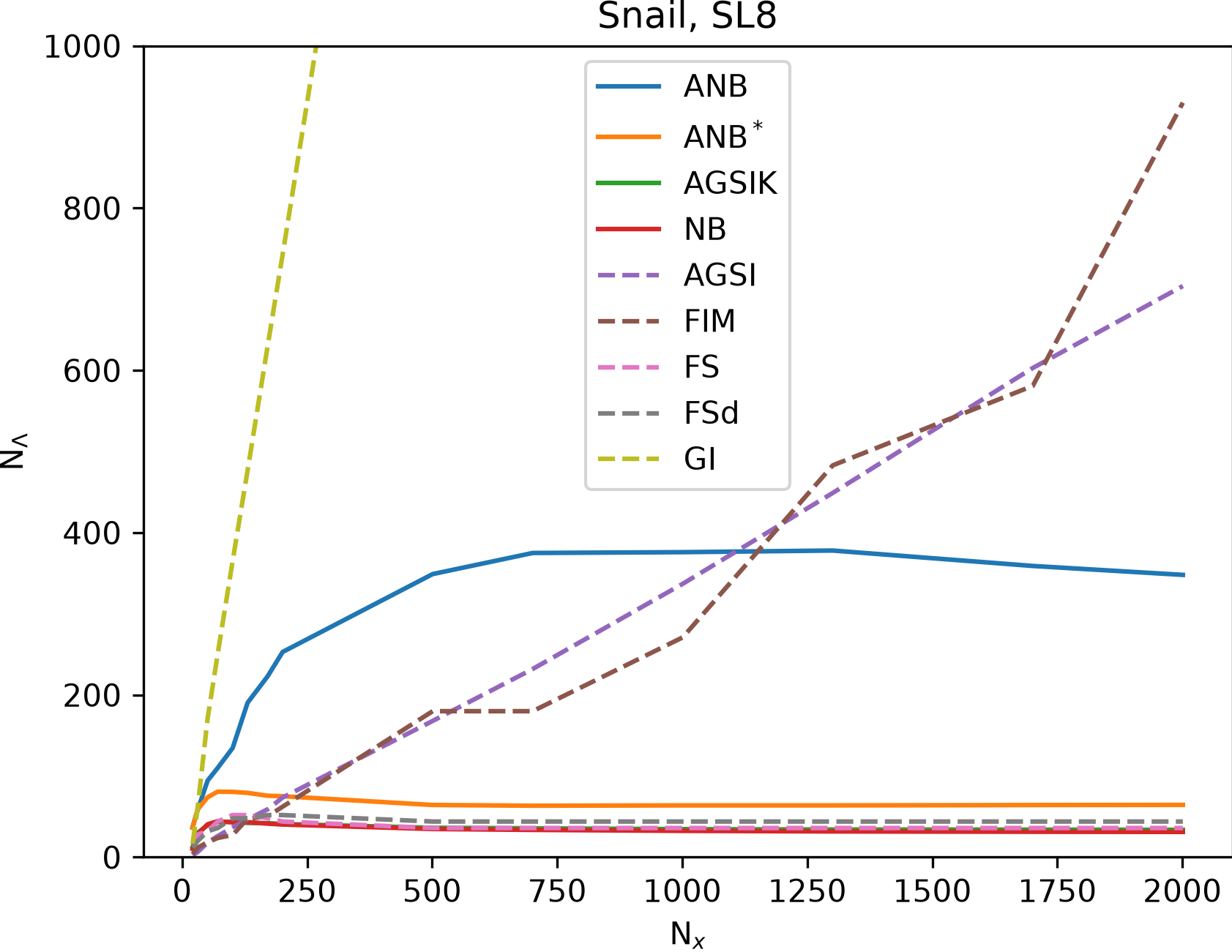}\\
\includegraphics[height=\h]{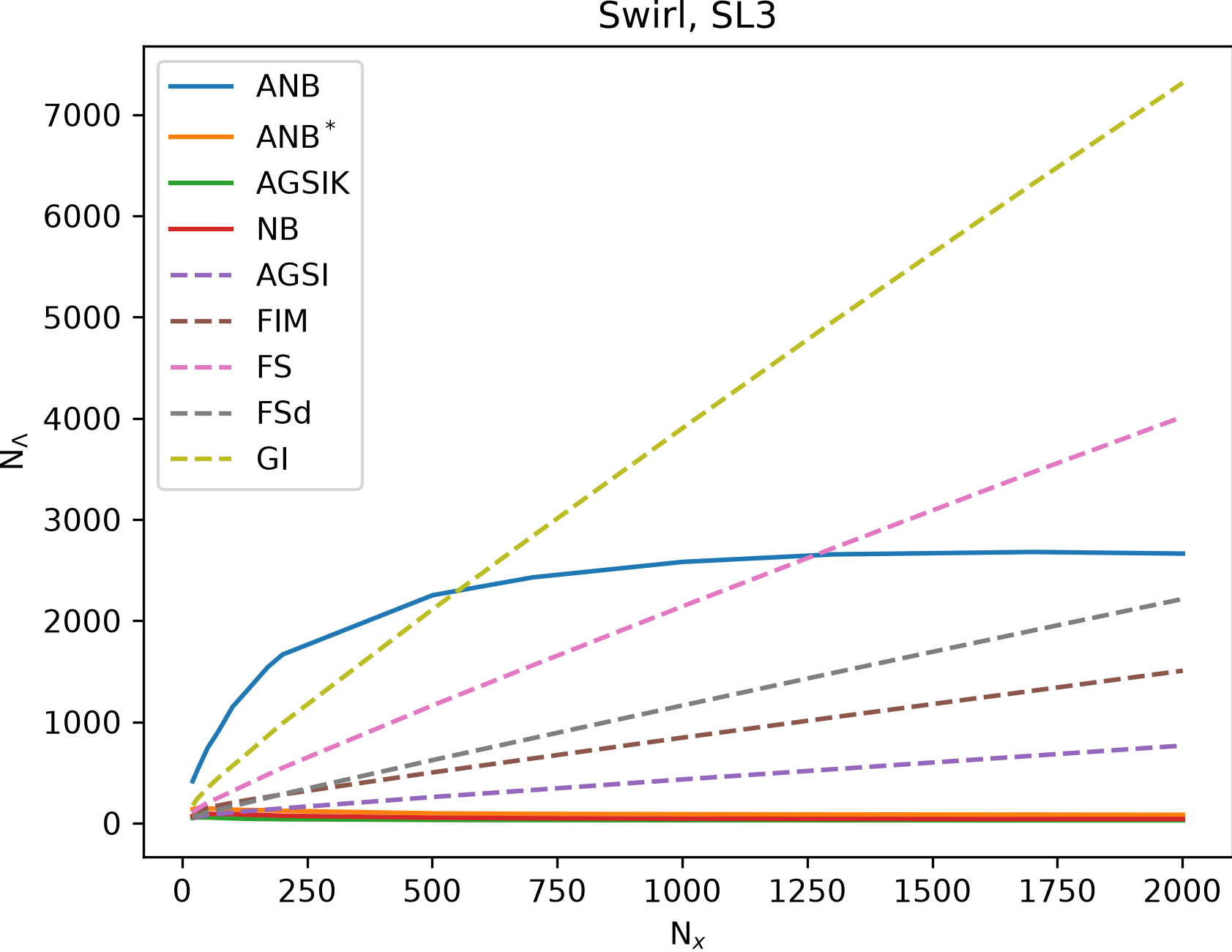}
\includegraphics[height=\h]{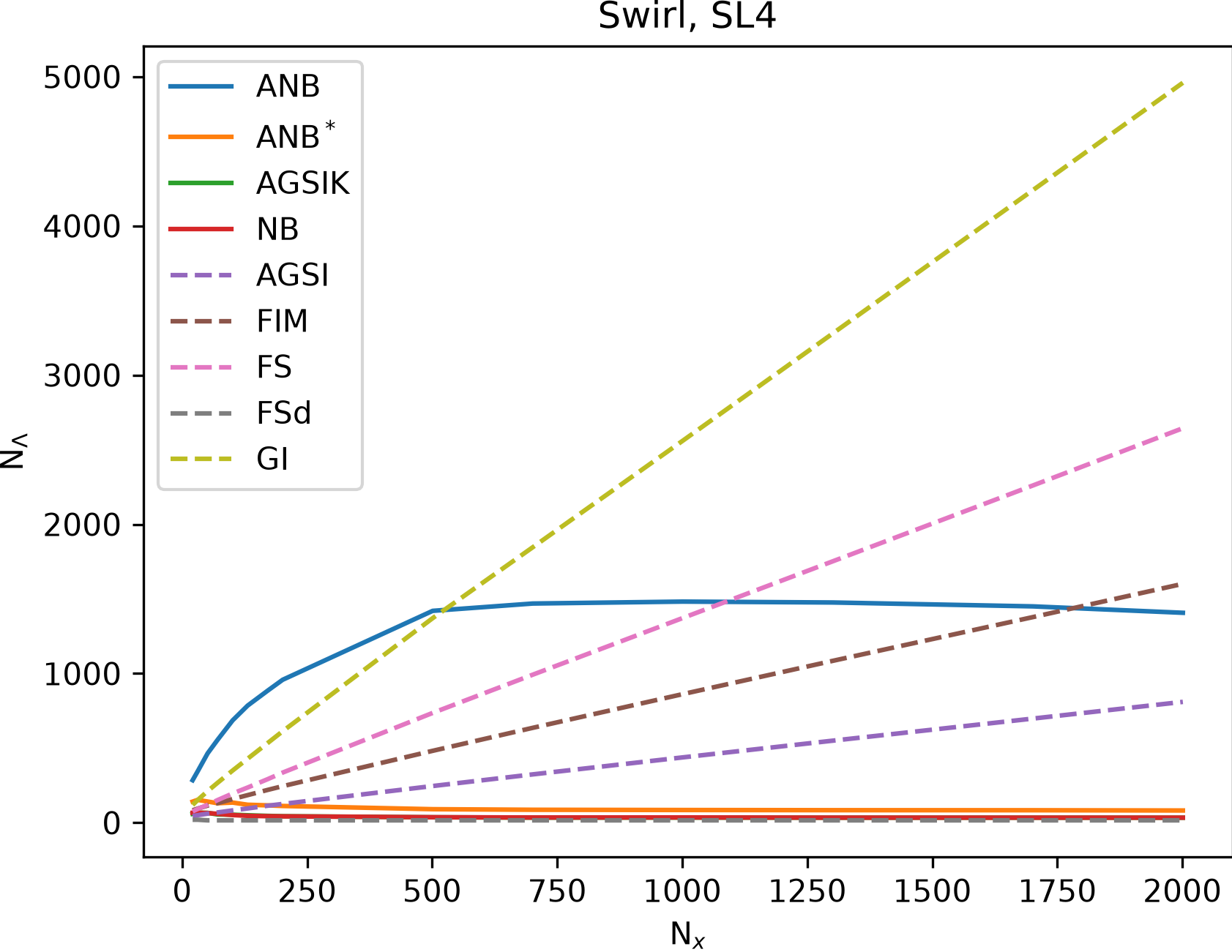}
\includegraphics[height=\h]{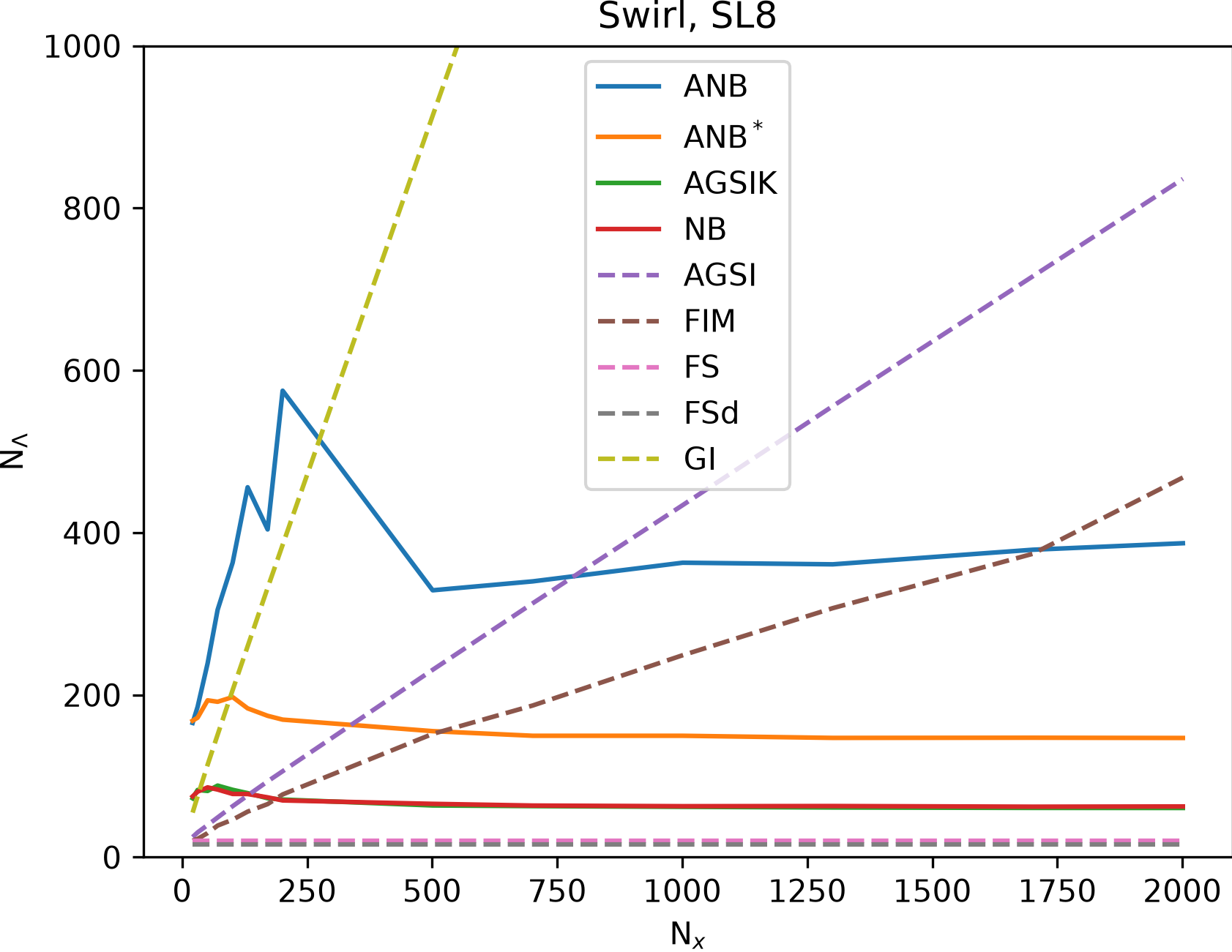}

\caption{Average number of iterations per point $N_\Lambda$, as a function of the number $N_x$ of discretization points in the horizontal dimension (inversely proportional to the gridscale), in the various test cases, for the various methods.
$N_\Lambda$ is independent of $N_x$ for the ANB, ANB$^*$, AGSIK, and NB method. Likewise for the FSd method with the SL4 scheme, and the FS and FSd methods with the SL8 scheme. In the other cases, $N_\Lambda$ grows linearly with $N_x$, in such way that the method cost is super-linear.
}
\label{fig:iter_wrt_h}
\end{figure}

\paragraph{Gridscale $h$.}
The average number of pointwise updates of the ANB, ANB$^*$, AGSIK, and NB methods is on \cref{fig:iter_wrt_h} empirically independent of the number of discretization points $N_x$ along the horizontal dimension, i.e.\ the inverse gridscale. (Up to one exception: the standard ANB in the Snail test case with the SL3 scheme. This is an extremely degenerate test case, where we expect stabilisation to occur for larger $N_x$.) This is consistent with the quasi-linear cost estimate, and with the uniform pointwise bound on evaluations of $\Lambda$ established in \cref{th:complexity,th:complexity_withUnknown,th:complexity_with_agsi}. Let us acknowledge nevertheless that the number of iterations of the ANB is often quite high in the semi-Lagrangian test cases, which is why we recommend the ANB$^*$ variant in practice.
The classical methods see their iteration count increasing linearly with $N_x$, reflecting their super-linear cost $\cO(N^{1+1/d})$, except for the FS (resp.\ FSd) method with the SL4 (resp.\ SL4 and SL8) scheme as already noted.

\paragraph{Convergence tolerance $\ve$.} %
	The computational cost of the ANB, ANB$^*$ and AGSIK empirically appears to depend logarithmically on the tolerance parameter $\ve>0$, see \cref{fig:iter_wrt_eps}, which is compatible with the $\cO(N|\ln \ve| \ln (N/\ve))$ estimate as $\ve \to 0$.
	On the other hand, front tracking appears to fail when $\ve$ is \emph{too large} in the seismic test case, as anticipated in \cref{rem:tol_range}.\\

The other parameters $(\alpha,\rho,\tau)$ are internal to the ANB, ANB$^*$ and AGSIK. 
We vary the following quantities $(r_\alpha,\rho,r_\tau)$, one by one and leaving the others constant to the value:
\begin{align}
\label{eq:default_ratios}
	r_\alpha & \in \{1,2,3,5\}, &
	\rho &= 0.95, &
	r_\tau &= 0.5, 
\end{align}
where the specific choice of $r_\alpha$ depends on the test case and is summarized in \cref{table:niter}.
We let
\begin{align*}
	\alpha &= r_\alpha h, &
	\tau &= r_\tau \alpha, 
\end{align*}
and run the ANB and ANB$^*$ with parameters $(\alpha,\rho,\tau,\ve)$, and the AGSIK with parameters $(\alpha,\tilde \ve)$, following \eqref{eq:default_tol} and \eqref{eq:default_rho_tau}.
The overall picture is that one must pay some attention to the choice of the timescale $\alpha$ (controlled via $r_\alpha$), but there is not much point in changing $\rho$ or $\tau$ (controlled via $r_\tau$) from their defaults.

\paragraph{Timescale $\alpha$ (controlled via $r_\alpha = \alpha/h$).} 	This parameter has the strongest impact on the computational cost of the ANB and ANB$^*$, see \cref{fig:iter_wrt_alpha}. 
We obtained good results with $h \leq \alpha \leq 5 h$ in our experiments, see \cref{table:niter} for the specific values used. 
In practice, choosing $\alpha$ which minimizes the contraction ratio near linear maps $\rho^\lin(\alpha)$ \eqref{eq:rho_lin} is often almost optimal, compare \cref{fig:rho_lin,fig:iter_wrt_alpha}. In contrast, in our test cases, the  theoretical estimates of $\rho(\alpha)$ obtained in \cref{sec:effective} appear to be too pessimistic to be useful.
The complexity of the AGSIK also depends on the parameter $\alpha$, but less strongly than the ANB and ANB$^*$, see again \cref{fig:iter_wrt_alpha}.

\paragraph{Parameter $\rho$.}
The dependency of the computational cost w.r.t.\ this parameter is generally unpredictable: non-monotone, non-continuous, and different for the ANB and ANB$^*$, see \cref{fig:iter_wrt_rho}. 
Fortunately, this dependency is not very strong overall: the number of iterations rarely varies by a factor larger than $2$ as $\rho \in [0.8, 0.99]$. The default value $\rho = 0.95$ \eqref{eq:default_rho_tau} is thus good enough for our experiments, and we make no attempt to use the theoretically estimated contraction ratio $\rho(\alpha)$, see \cref{sec:effective}, or its numerically estimated variant $\rho^\lin(\alpha)$ near linear maps, see \cref{fig:rho_lin}.

\paragraph{Timestep $\tau$ (controlled via $r_\tau = \tau/\alpha$)} This parameter determines how often the narrow band is rebuilt.
If $\tau \ll \alpha$, then the band is rebuilt too frequently, which increases the number of iterations of the causal operator $\Gamma$ used for that purpose. 
If $\tau \gg \alpha$, then the width of the narrow band increases, and the method becomes global.
We choose the default $\tau = \alpha/2$, but note that slightly less frequent narrow band reconstructions $\tau/\alpha \in [1/2,2]$ work as well, see \cref{fig:iter_wrt_tau}.

\paragraph{Computational time vs number of iterations.}
Finally, we present in \cref{table:elapsed} the computational time of the various methods, in the various test cases, obtained using a single thread on an Apple$^{\small\textregistered}$ macbook pro equipped with an M1 pro$^{\small\textregistered}$ processor. 
Runtimes are not perfectly reproducible, and can vary substantially between runs, which is why some numbers in \cref{table:elapsed} differ from those reported in \cite{mansour2024solving}.
Finally, let us acknowledge that much faster computational times can be achieved using fast marching methods with modified schemes and adaptive stencils \cite{mirebeau2014EfficientFinsler,mirebeau2019riemannian,mirebeau2014marchingLattice} on a sequential architecture, or by massively parallel methods on GPU architectures \cite{kimmel2008ParallelDistance,jeong2008fim,mirebeau2023massively}, but they are outside the scope of this paper.

\begin{table}
\hspace{-0.5cm}
\begin{tabular}{cccccccccc}
 &  ANB &  ANB$^*$ &  AGSIK &  NB &  AGSI &  FIM &  FS &  FSd &  GI\\
1D, LF & 0.01 & 0.01 + 0.01 & 0.01 & 0.01 & 0.04 & 0.05 & 0.07 & NA & 0.09\\
Pursuit, LF & 53 + 41 & 62 + 68 & 36 & 30 & 152 & 204 & 550 & 305 & 299\\
Seismic(subs), LF & 37 + 23 & 57 + 33 & 45 & 44 & 1100 & 1598 & * & 2695 & *\\
Seismic, LF & 198 + 126 & 268 + 155 & 191 & 193 & * & * & * & * & *\\
Maze, UD & 0.49 + 0.03 & 0.59 + 0.04 & 0.15 & 0.16 & 0.15 & 0.19 & 2.6 & 0.37 & 5.4\\
Snail, 2SL3 & 324 + 1540 & 56 + 65 & 25 & 54 & 1060 & 1707 & * & 1533 & *\\
Snail, 2SL4 & 246 + 252 & 44 + 41 & 29 & 69 & 1363 & 2306 & * & 115 & *\\
Snail, 2SL8 & 309 + 272 & 72 + 69 & 55 & 70 & 1143 & 1489 & 53 & 67 & *\\
Swirl, 2SL3 & 883 + 1096 & 69 + 46 & 22 & 36 & 550 & 913 & 2716 & 1400 & 3578\\
Swirl, 2SL4 & 510 + 769 & 70 + 49 & 32 & 40 & 710 & 1210 & 2209 & 14 & 3085\\
Swirl, 2SL8 & 215 + 416 & 139 + 156 & 97 & 119 & 1394 & 776 & 29 & 24 & *\\

\end{tabular}
\caption{Elapsed time, in seconds. For the ANB and ANB$^*$, we report the total time as $a+b$ where $b$ is devoted to the sub-solution updates, and $a$ includes everything else (super-solution updates, narrow band construction via Dijkstra's method). 
}
\label{table:elapsed}
\end{table}

\newcommand{\etalchar}[1]{$^{#1}$}

\end{document}